\documentclass[aap]{amsart}

\usepackage{url}
\usepackage[colorlinks=true, linkcolor=blue, citecolor=blue, urlcolor=blue]{hyperref}

\usepackage{mdframed}
\usepackage{enumitem}

\usepackage{amssymb,amsmath,amsthm}
\usepackage{mathabx} %
\usepackage{empheq}
\usepackage{scalerel}
\usepackage{mathrsfs} %
\usepackage{bm}
\usepackage{bbm}
\usepackage[retainorgcmds]{IEEEtrantools}
\usepackage{textgreek}
\usepackage{stmaryrd}

\usepackage{algorithm}
\usepackage{algorithmic}
\usepackage[algo2e,linesnumbered,ruled,vlined]{algorithm2e}

\usepackage{graphicx}
\usepackage{caption}
\usepackage{subcaption}
\usepackage[all,cmtip]{xy} %
\usepackage{tikz} %
\usepackage{adjustbox}

\usepackage{epstopdf}

\usepackage[margin=2.5cm]{geometry}

\usepackage{hyperref}
\usepackage[square,numbers]{natbib}
\usepackage[titletoc]{appendix}

\usepackage{color}
\usepackage{xcolor}
\usepackage{todonotes}

\newcommand\bI{ {\mathbf I} }

\newcommand\bM{ {\mathbf M} }

\newcommand\bU{ {\mathbf U} }

\newcommand\bX{ {\mathbf X} }
\newcommand\bY{ {\mathbf Y} }

\newcommand\bp{{\mathbf p}}

\newcommand\bu{{\mathbf u}}

\newcommand\by{{\mathbf y}}
\newcommand\bz{{\mathbf z}}

\newcommand{\cB}{{\mathcal{B}}}

\newcommand{\cF}{{\mathcal{F}}}

\newcommand{\cL}{{\mathcal{L}}}
\newcommand{\cM}{{\mathcal{M}}}

\newcommand{\cS}{{\mathcal{S}}}
\newcommand{\cT}{{\mathcal{T}}}
\newcommand{\cU}{{\mathcal{U}}}

\newcommand{\cX}{{\mathcal{X}}}

\newcommand{\bbeta}{  {\bm{\beta}} }

\newcommand{\bvarepsilon}{  {\bm{\varepsilon}} }

\newcommand\EE{ {\mathbb E} }

\newcommand\PP{ {\mathbb P} }

\newcommand\RR{ {\mathbb R} }

\def\SS{\mathbb{S}}

\newcommand\rA{ {\mathrm A} }
\newcommand\rB{ {\mathrm B} }

\newtheorem{theorem}{Theorem}
\newtheorem{definition}[theorem]{Definition}

\newtheorem{lemma}[theorem]{Lemma}
\newtheorem{proposition}[theorem]{Proposition}
\newtheorem{remark}[theorem]{Remark}
\newtheorem{corollary}[theorem]{Corollary}

\newtheorem{assumption}{Assumption}

\newcommand{\vertiii}[1]{{\left\vert\kern-0.25ex\left\vert\kern-0.25ex\left\vert #1 
		\right\vert\kern-0.25ex\right\vert\kern-0.25ex\right\vert}}

\newcommand\bone{\mathbf 1}
\newcommand\bbone{\mathbbm{1}}

\usepackage{amsopn}
\DeclareMathOperator{\diago}{diag}

\usepackage{empheq}
\usepackage{physics}

\newcommand{\idynoise}[1]{ \varepsilon_{#1} }
\newcommand{\comnoise}[1]{ \varepsilon_{#1}^0 }
\newcommand{\Xnt}{ {X_t^{(n)}} }
\newcommand{\unt}{ {u_t^{(n)}} }

\newlength\myindent
\title[Policy Gradient Convergence for LQ Mean-Field RL]{Linear-Quadratic Mean-Field Reinforcement Learning: Convergence of Policy Gradient Methods}
\author{Ren\'e Carmona$^{1}$}
\thanks{$^{1}$Department of Operations Research and Financial Engineering \& Program in Applied and Computational Mathematics, Princeton NJ 08544, USA, (\href{mailto:rcarmona@princeton.edu}{rcarmona@princeton.edu})}
\author{Mathieu Lauri\`ere$^{2}$}
\thanks{$^{2}$Shanghai Frontiers Science Center of Artificial Intelligence and Deep Learning; NYU-ECNU Institute of Mathematical Sciences at NYU Shanghai; NYU Shanghai, 567 West Yangsi Road, Shanghai, 200126, People’s Republic of China, (\href{mailto:mathieu.lauriere@nyu.edu}{mathieu.lauriere@nyu.edu})}
\author{Zongjun Tan$^{3}$}
\thanks{$^{3}$Department of Operations Research and Financial Engineering \& Program in Applied and Computational Mathematics, Princeton NJ 08544, USA; Walmart Global Tech, USA (\href{mailto:zongjun.tan@walmart.com}{zongjun.tan@walmart.com})}
\thanks{{\bf Acknowledgements.}  The authors were partially supported through funding from AFOSR under the grant FA9550-23-1-0324. M.L. and Z.T. thank Zhuoran Yang for fruitful discussions in the early stage of this project. }
\date{}

\begin{document}

\begin{abstract}
We investigate reinforcement learning in the setting of Markov decision processes for a large number of exchangeable agents interacting in a mean field manner. Applications include, for example, the control of a large number of robots communicating through a central unit dispatching the optimal policy computed by maximizing an aggregate reward. An approximate solution is obtained by learning the optimal policy of a generic agent interacting with the statistical distribution of the states and actions of the other agents. We first provide a full analysis this discrete-time mean field control problem. We then rigorously prove the convergence of exact and model-free policy gradient methods in a mean-field linear-quadratic setting and establish bounds on the rates of convergence. We also provide graphical evidence of the convergence based on implementations of our algorithms. 
\end{abstract}

\maketitle



\section{\textbf{Introduction}}
\label{se:introduction}

Reinforcement Learning (RL) is one of the most powerful frameworks in learning theory, and it is particularly useful when feedback information like penalties or rewards is available to a learning agent. A \emph{learning agent} is generally understood as a computational mechanism with a high degree of autonomy that collects feedback information by interacting with an environment.
A typical issue in RL is  to connect observations of the environment to actions. 
In the standard literature on RL algorithms are intended to identify policies for a single agent. When multiple agents interact within a shared environment,  \emph{Multi-Agent Reinforcement Learning} (MARL) algorithms are proposed. In this paper, we consider MARL models with a large number of agents, and in the limit when this number of agents tends to infinity, reducing the analysis to the optimization of a single agent interacting with an aggregated state generated by all the other individual agents. We extend the single-agent RL paradigm to the framework of Mean Field Control (MFC) theory, and call it \emph{Mean Field Reinforcement Learning} (MFRL). In this learning framework, we need to acquire feedback signals both from a representative agent and the aggregated population.

MFRL leverages the power of RL and MFC to efficiently find optimal policies in very large populations of agents. 
Learning with mean field interactions has recently gained momentum, in particular with model-free techniques, see e.g.~\cite{SubramanianMahajan-2018-RLstatioMFG,guo2019learning,fu2019actorcriticMFG,elie2020convergence,gu2023dynamic,gu2021mean,anahtarci2023q,perrin2020fictitious,motte2022mean,carmona2020policyCDC}.
The theory behind RL can be viewed as a development of optimal control in which the model is (partially or fully) unknown. On the other hand, MFC appears as a generalization of optimal control for cooperative agents in mean field interaction when the number of agents grows to infinity. MFRL lies at the intersection of these two extensions. It aims at describing how a large number of agents can collectively learn a reasonable approximation of the optimal solution of a control problem.  It is worth noticing that when the number of agents is finite, we are dealing with an $N$-agent optimal control problem, or we are learning the optimal policies using a MARL framework. Note that the convergence from the $N$-agent problem to its mean-field approximation depends on the equilibrium concepts. See for example~\cite{carmona2023nash, lacker2018convergence, cardaliaguet2019master} and~\cite{CarmonaDelarue_book_I} for a discussion of the connection between MFG and MFC and the convergence from $N$-agent problems to mean-field problems.

MFRL for the cooperative setting is rooted in the theory of mean field Markov decision processes (MFMDP), which have been studied in~\cite{gast2012mean}. \cite{motte2022mean,carmona2023model} studied MFMDPs with common noise and proved dynamic programming principles under generic assumptions. \cite{gu2023dynamic,gu2021mean} focused on dynamic programming principles and RL for the state-action value function for finite-state MFC without common noise. \cite{carmona2023model} proposed tabular and neural-network based RL algorithms for MFMDP with common noise and common randomization, on more general spaces. Continuous-time RL for MFC has been studied in~\cite{frikha2023actor,frikha2024full,wei2025continuous}. 
Another line of research focuses on the non-cooperative setting with mean-field games~\cite{fu2019actorcriticMFG, wang2021global, lauriere2022scalable, muller2022learning}, when the agents in the system play against each other and optimize their own rewards. 
Other works related to MFRL include~\cite{yang2018mean} in the framework of stochastic games, ~\cite{arabneydi2016new,shao2024reinforcement} for mean-field teams and mean field type games, ~\cite{angiuli2022reinforcement,angiuli2023convergence} for mean field control games, or~\cite{gu2024mean} for a network structure. We refer to~\cite{lauriere2022learning} for a recent survey on this topic.

In this paper, we restrict ourselves to fully cooperative models whose mean-field nature can be understood as a limiting case of both the team learning and the concurrent learning models in MARL. 
The interactions between a representative agent and the agent statistical distribution can be viewed as an asymptotic approximation of the interactions between homogeneous agents and the centralized planner in team learning. 
Each individual agent is embedded in a large system
and needs to react to the changes caused by the field of all the agents, but on their own, they have essentially no direct impact on other agents' policies like in concurrent learning.
It is important to highlight that the mean-field distribution is not posed apriori. The agents are learning together with the evolution of the mean-field distributions generated from these agents. 
As a consequence, a generic agent should also learn a macroscopic level policy in order to evaluate the impact of  overall field of agents.
Our model includes one more special feature: in addition to the idiosyncratic noise felt by  each agent, a common noise is added to help the model be applicable to more realistic applications. Despite the increased complexity die to its inclusion, we show that this extra source of randomness helps the exploration in the space of distributions.

The present version of the paper contains significant differences with the first version, which was posted on arXiv in 2019. First, we added Section~\ref{sec:assm-defi}, which establishes the well-posedness of the problem and clarifies the connection between admissible controls and admissible parameters. We also added  Section~\ref{sec:analysis}, which provides a complete analysis of the theoretical framework of the problem.  Moreover, the assumptions have been reorganized, and several proofs are presented in more details. In particular, for each of the three algorithms, we show that the learned parameters remain admissible after each update. 
Since the first arXiv version of this paper, several follow-ups appeared online. To cite just a few closely related works, let us mention that \cite{angiuli2022unified} studied a unified RL approach for MFC and MFG and the convergence was proved in~\cite{angiuli2023convergence}, 
\cite{carmona2021linear,uz2024independent,zaman2024robust} studied convergence of policy gradient methods for LQ mean field type games (i.e., a Nash equilibrium between MFC-type players),  \cite{li2023policy} and \cite{wang2024online} respectively studied a policy iteration and a value iteration algorithm for LQ MFC, and \cite{xu2024model,xu2025mean} studied LQ models with multiplicative noise. Moreover, \cite{cui2024partially,cui2024learning} studied RL for MFC with partial observation. Notice that none of these works supersedes the results presented in this paper.

In this paper, we investigate  Linear Quadratic (LQ) models with Mean Field interactions in the presence of a common noise.
In general, the mean field limit provides a streamlined description of the interactions between a large group of collaborative agents. However these problems remain extremely difficult to solve~\cite{MR3343705,MR3501391,MR3631380}. Fortunately, in mean-field LQ problems, simpler optimality conditions can be formulated in terms of algebraic Riccati equations. Because of its tractability, the family of LQ models proved to be of great importance in optimal control, and they have also been studied from a RL viewpoint (see for example \cite{recht2018tour,fazel2018global, yang2019provably}). 

The main thrust of the present paper is to generalize the global convergence of the Policy Gradient (PG) algorithm given in~\cite{fazel2018global} for standard deterministic LQ problem  to a LQ MFC problem with common noise. Incidentally, we take advantage of the impact of the common noise to explore the unknown environment, even if its presence significantly complicates the convergence analysis. In particular, we present two model-free algorithms based on a derivative-free gradient estimation technique called the \emph{zeroth order optimization} method. This type of derivative-free method is related to another research line called \emph{Evolution Strategy} (ES), which has been studied for a long time \cite{rechenberg1973evolutionsstrategie} and recently attracted more attention in the RL community~\cite{salimans2017evolution}. The main idea is to evaluate the function to minimize on a set of perturbed parameters so that one can move along a weighted average direction which improves the value of the target function~\cite{MR3627456}. It has been shown that a careful implementation of ES can provide performance improvements when compared with classical RL algorithms based on back-propagation. Moreover, the method can be easily scaled up, and be adapted to distributed computer systems \cite{salimans2017evolution,choromanski2018structured}.

\vskip 6pt
The main contributions of the paper can be organized in two separate buckets.
\begin{itemize}
\item First, we prove existence and uniqueness of an optimal control for a mean-field linear-quadratic Markov decision process with a common noise. Beyond existence and uniqueness, we prove that the optimal admissible control is linear in the state and its conditional mean process. This necessary form of the optimal control allows us to reduce the search for optimality in the second part of the paper to a specific parametric class of controls. It is derived by an application of necessary and sufficient conditions in a discrete-time Pontryagin's maximum principle tailored to our LQMFC set-up. This provides a crucial link between the MFC problem and learning by a finite population of agents as it gives an approximately optimal control for the problem with a finite number of learners.

\item Second, we investigate the convergence properties of the PG algorithm, both for model-based and for model-free problems. We prove three global convergence results: 1) with an exact computation of the policy gradient terms, 2) with a gradient estimator based on an idealized MKV simulator for the MFC cost, and finally, 3) a more realistic gradient estimator with a population simulator for a system with a finite number of agents. 
In particular, we show how a finite number of agents can collaborate to learn a control which is approximately optimal for the MFC problem. 
\end{itemize}
Finally, we stress that our numerical experiments show that the method is very robust and that  it can be used for finitely many agent models, even if the agents are not exchangeable and have noisy dynamics.

\vskip 12pt\noindent 
\textbf{Frequently used notations. }

The set of $d \times d$ real symmetric matrices is denoted by $\cS\cM_d$, or $\cS\cM$ for short when the dimension is clear from the context. For $U \in \cS\cM$, $\lambda_{min}(U)$ denotes the minimum eigenvalue of $U$. $U \succeq 0$ means that $U$ is positive semi-definite (PSD), namely $\lambda_{min}(U) \geq 0$, and $U \succ 0$ means that $U$ is positive definite (i.e. $\lambda_{min}(U) > 0$). The symbol $\otimes$ denotes the Kronecker product between matrices.  The transpose of a vector or a matrix $U$ is denoted by $U^\top$. The inner product between two vectors $x$ and $y$ is written as $x \cdot y$ or $x^\top y$. 

For matrices $U, V \in \RR^{\ell \times d}$, we denote by $\langle U, V \rangle_{tr} = Tr(U^\top V)$ the trace of the product $U^\top V$, and  we denote by $\| U \|$, $\| U \|_F$ and $\| U \|_{tr}$ the spectral norm, the Frobenius norm and the trace norm respectively. In other words, $\| U \| = \max \{ \sigma_i(U) \}$, $\| U \|_{F} = \sqrt{\langle U^\top, U \rangle_{tr} } = \sqrt{ \sum \sigma_i(U)^2 }$ and $\| U \|_{tr} = \sum \sigma_i(U)$, where $\big( \sigma_i(U) \big)_{i=1}^{\min\{\ell, d\}}$ are the singular values of $U$.  
Let $\tilde U \in \RR^{(\ell + d) \times (\ell + d)}$ be a symmetric matrix with $U$ and $U^\top$ on its off-diagonal blocks, called the dilatation matrix of $U$ (see~\cite{MR2946459}). We have $\| \tilde U \| = \| U \|$.
We also recall that 
$$\| U \| \leq \| U \|_{F} \leq \sqrt{\min\{ \ell, d\}} \| U \|
\qquad\text{and}\qquad
| \langle U, V \rangle_{tr} | \leq \| U\|. \| V \|_{tr}.
$$
Moreover, if $U, V \in \cS\cM$ and $U \succeq 0$, then 
$$
\lambda_{min}(U) Tr(V) \leq | \langle U, V \rangle_{tr} |. %
$$
For a vector $\xi \in \RR^d$, we denote by $\| \xi \| = \sqrt{ \xi^\top \xi }$ its Euclidean norm, and $\diago(\xi) \in \RR^{d\times d}$ denotes the diagonal matrix with entries of $\xi$ on its diagonal. When $\xi$ is a random vector, $\| \xi \|_{L^p(\RR^d)} = \EE[ \| \xi \|^p]^{1/p}$ denotes the $L^p$ norm of $\xi$. Let . 
The sub-gaussian norms (see~\cite{MR3837109}) of a sub-gaussian random variable $\eta \in \RR$ or of a random vector $ \xi \in \RR^d$, denoted by $\| \eta \|_{\psi_2}$ and $\| \xi \|_{\psi_2}$, are defined by
$$
\| \eta \|_{\psi_2} = \inf \left\{ s > 0 : \mathbb{E}[ \exp( \eta^2 / s^2) ] \leq 2 \right\}
\qquad\text{and}\qquad
\| \xi \|_{\psi_2} = \sup_{ v \in \RR^d: \| v \| = 1} \| \xi^\top v \|_{\psi_2}. 
$$ 
Similarly, the sub-exponential norm of a sub-exponential random variable $\zeta \in \RR$ is defined by 
$$
\| \xi \|_{\psi_1} = \inf \{ s > 0: \EE[ \exp( \xi / s) ] \leq 2 \}.
$$
The following inequalities connect the $L^p$ norms with the sub-exponential norm and the sub-gaussian norm (see~\cite[Proposition 2.7.1]{vershynin2018high}): 
$$
    \| \xi \|_{L^2(\RR^d)} \leq \sqrt{2 d} \| \xi \|_{\psi_2}
\qquad
\text{and}
\qquad
    \| \zeta \|_{L^2(\RR)} \leq 2 \| \zeta \|_{\psi_1}.
$$
Moreover, for any $p \geq 2$, we have 
$
    \| \zeta \|_{L^p(\RR)} \leq (2 p !)^{1/p} \| \zeta \|_{\psi_1} \leq p \| \zeta \|_{\psi_1}.
$
We use bold letters for sequences,  so that $\bX = (X_t)_{t \geq t}$ stands for a sequence of random vectors $X_t \in \RR^d$. 
For $N$ random vectors $(X^1, \ldots, X^N) \in (\RR^d)^N$, we denote the concatenated random vector in $\RR^{dN}$ by $\mathrm{Vect}(X^1, \ldots, X^N) = [(X^1)^\top, \ldots, (X^N)^\top]^\top$. We write $I_d$ for the $d \times d$ identity matrix, $\bone_d$ and $\bbone_d$ for the $d \times d$ matrix and the $d-$vector whose entries are all ones. We shall drop the subscripts when the dimension is clear from the context. 

\section{\textbf{Models}}
\label{sec:def-model}

In this section, we introduce an LQ model for a finite number of agents and an LQ MFC model with common noise. We work with a probability space $(\Omega, \cF, \PP)$ that supports an initial random variable, an i.i.d. idiosyncratic noise sequence $\bvarepsilon = (\varepsilon_{t+1})_{t \geq 0}$, and an i.i.d. common noise sequence $\bvarepsilon^0 = (\varepsilon_{t+1}^0)_{t \geq 0}$. 
The two noise sequences $\bvarepsilon$ and $\bvarepsilon^0$ considered here are independent of each other in the sense that $\varepsilon_{t+1}$ is independent of $\varepsilon_{s+1}^0$ for every $s, t \geq 0$. These noise terms are of mean zeros.
We also assume that the probability space support two \emph{initial} random variables, denoted by $\varepsilon_0^0$ and $\varepsilon_0$, that are independent of each other and independent of the noise sequences $\bvarepsilon^0$ and $\bvarepsilon$, $\varepsilon_0^0$ representing the initial random common shock at initial time $0$. 
The initial perturbations $\varepsilon_0$ and $\varepsilon_0^0$ are not assumed to have mean-zero distributions like the noise terms $\varepsilon^0_{t+1}$ and $\varepsilon_{t+1}$. 
\begin{assumption}
    \label{ass:subgaussian}
We shall assume that the initial perturbations $\varepsilon_0$ and $\varepsilon_0^0$, and all noise terms $(\varepsilon_{t+1})_{t \geq 0}, (\varepsilon_{t+1}^0)_{t \geq 0}$ are sub-gaussian random vectors. 
\end{assumption}

We shall use the filtrations $(\cF_t)_{t \geq 0}$ and $(\cF^{0}_t)_{ t \geq 0}$ defined by:
\begin{align}
    \cF^{0}_t & = \sigma \big( \{ \varepsilon_0^0, \varepsilon_1^0, \ldots, \varepsilon_{t}^0  \} \big),
    \\
    \cF_t & = \sigma \big( \{ \varepsilon_0, \varepsilon_0^0, \varepsilon_1, \varepsilon_1^0, \ldots, \varepsilon_{t}, \varepsilon_t^0 \} \big).
\end{align}
The building blocks of the model are:
\begin{itemize}
    \item a measurable \emph{state space} $(S, \cB_S) = (\RR^d, \cB_{\RR^d})$; 
    \item a measurable \emph{action space} $(A, \cB_A) = (\RR^\ell, \cB_{\RR^\ell})$; 
    \item two measurable spaces $(E, \cB_E) = (E^0, \cB_{E^0}) =(\RR^d, \cB_{\RR^d})$ in which the idiosyncratic noise and common noise sequences take their values;
    \item a discount factor $\gamma \in (0, 1)$;
    \item a linear \emph{system function}, denoted by $F : S \times S \times A \times A  \times E^0 \times E \to S$, with model parameters $(\rA, \bar\rA, \rB, \bar\rB) \in \RR^{d\times d} \times \RR^{d \times d} \times \RR^{d \times \ell} \times \RR^{d \times \ell}$ given by 
    \begin{equation}
        \label{eq:lq_sys_function_F}
        F(x, \bar x, u, \bar u, e^0, e) = \rA x + \bar\rA \bar{x} + \rB u + \bar\rB \bar{u} + e^0 + e;
    \end{equation}
    \item a quadratic one-step cost function, denoted by $c: S \times S \times A \times A \to \RR$, with model parameters $(Q, \bar Q, R, \bar Q) \in \RR^{d \times d} \times \RR^{d \times d} \times  \RR^{\ell \times \ell} \times \RR^{\ell \times \ell}$ given by 
    \begin{equation}
        \label{eq:lq_one_step_cost}
        c(x, \bar{x}, u, \bar{u}) = (x-\bar{x})^\top Q (x-\bar{x}) + \bar{x}^\top (Q + \bar{Q}) \bar{x} + (u-\bar{u})^\top R (u-\bar{u}) + \bar{u}^\top (R + \bar{R}) \bar{u}. 
    \end{equation}
\end{itemize}

\subsection{Model for  a finite number of agents}
\label{subsec:OC-N-Agents}

We consider $N$ agents with symmetric interactions.  
Their states at time $t \geq 0$, denoted by $(X^n_t)_{n=1}^N \in (\RR^d)^N$, satisfy the system of equations: 
\begin{align}
\label{fo:N-multi_state}
    X_{t+1}^{(n)} & = F \big( X_t^{(n)}, \bar{X}_t^{N}, u_t^{(n)}, \bar{u}_t^N, \comnoise{t+1}, \idynoise{t+1}^{(n)} \big) 
    \nonumber \\
    & = \rA \Xnt + \bar\rA \bar{X}_t^{N}  + \rB \unt + \bar{\rB} \bar{u}_t^N + \comnoise{t+1} + \idynoise{t+1}^{(n)}
\end{align}
with initial conditions $X_0^{(n)} = \varepsilon_0^0 + \varepsilon_0^{(n)}$, for $n = 1, \ldots, N$, where 
$
\bar{X}^N_t = \frac1N \sum_{n=1}^N X_t^{(n)}
$ 
and 
$
\bar{u}^N_t = \frac1N\sum_{n=1}^N u_t^{(n)}
$
denote the sample averages of the individual states and actions at time $t$, and $( \varepsilon_0^{(n)}, (\idynoise{t+1}^{(n)})_{t \geq 0})_{n=1,\ldots, N}$ are $N$ identical and independent copies of the initial perturbation and the idiosyncratic noise sequence $( \varepsilon_0, (\varepsilon_{t+1})_{t \geq 0} )$.

The time evolution of the sample average of the state reads:
\begin{equation}
\label{fo:N-mean_state_dynamics}
     \bar{X}^N_{t+1} = (\rA + \bar\rA) \bar{X}^N_t  +  (\rB + \bar \rB) \bar{u}^N_t + \varepsilon^0_{t+1} + \frac{1}{N}\sum_{n=1}^N \varepsilon^{(n)}_{t+1},
\end{equation}
with initial average $\bar{X}^N_0 = \comnoise{0} + \frac{1}{N} \sum_{n=1}^N \idynoise{0}^{(n)}$.
The last term in~\eqref{fo:N-mean_state_dynamics} is expected to be small when $N$ is large due to the law of large numbers, but $(\varepsilon^0_{t+1})_{t \geq 0}$ affects all the agents as it represents aggregate shocks. See for example  \cite{MR3325083,alasseur2020extended,graber2016linear} for applications to systemic risk and energy management in the mean field game setting.   

Letting $\underline{X}_t= \mathrm{Vect}\big(X_t^{(1)}, \ldots, X_t^{(N)} \big)$ and $\underline{U}_t= \mathrm{Vect} \big( u_t^{(1)}, \ldots, u_t^{(N)} \big)$ for $t\ge 0$, we obtain a vector dynamics as a representation of the $N$ agents' dynamics~\eqref{fo:N-multi_state}: 
\begin{equation}
\label{fo:N-vector_state}
    \underline{X}_{t+1} = \rA^N \underline{X}_t  + \rB^N \underline{U}_t + \underline{\varepsilon}^0_{t+1} + \underline{\varepsilon}_{t+1}
\end{equation}
with 
$\underline{\varepsilon}_{t+1} = \mathrm{Vect}( \idynoise{t+1}^{(1)}, \ldots ,\idynoise{t+1}^{(N)} )$, $\underline{\varepsilon}^0_{t+1} = \mathrm{Vect}( \comnoise{t+1},\ldots,\comnoise{t+1} )$, $ \rA^{N} = I_N \otimes \rA + \frac{1}{N} \bone_N \otimes \bar{\rA}$ and $\rB^{N} = I_N \otimes \rB + \frac{1}{N} \bone_N \otimes \bar{\rB} $. In other words, the matrix $\rA^{N}$ is the sum of the block-diagonal matrix with $N$ blocks equal to $\rA$, and the matrix with $N \times N$ blocks $\frac{1}{N}\bar{\rA}$, and similarly for the matrix $\rB^{N}$. 

The goal is to minimize the so-called \textit{social cost} of the population, defined as
\begin{equation}
    \label{eq:social_cost_of_population}
    J^N( \underline{\bU} ) = \EE \Big[ \sum_{t \geq 0} \gamma^{t} \bar c^N(\underline{X}_t, \underline{U}_t) \Big]
\end{equation}
over the set of admissible control process $\underline{\bU} = (\underline{U}_t)_{t \geq 0}$ consisting of $\RR^{\ell N}$-value process adapted to the filtration $(\cF_t)_{t \geq 0}$ that makes $J^N(\underline{\bU}) < \infty$.
Here, the instantaneous \textit{social cost} function $\bar c^N : \RR^{dN} \times \RR^{\ell N} \to \RR$ is
\begin{equation}
    \label{fo:N-cost}
    \bar c^{N}( \underline{X}_t, \underline{U}_t) := \frac{1}{N} \sum_{n=1}^N c \big( X_t^{(n)}, \bar{X}_t^N, u_t^{(n)}, \bar{u}_t^N \big) = \underline{X}_t^\top Q^N \underline{X}_t + \underline{U}_t^\top R^N \underline{U}_t,
\end{equation}
where $Q^N = \frac{1}{N} I_N \otimes Q + \frac{1}{N^2} \boldsymbol{1}_{N} \otimes \bar{Q} \in \RR^{dN \times dN}$ and $R^N =  \frac{1}{N} I_N \otimes R + \frac{1}{N^2} \boldsymbol{1}_{N} \otimes \bar{R} \in \RR^{\ell N \times \ell N}$. 

This is now a classical infinite-horizon discounted-cost linear-quadratic control problem with $(d N)-$dimensional state process $\underline{\bX} = (\underline{X}_t)_{t \geq 0}$ following dynamics~\eqref{fo:N-vector_state}. 
On the basis of the standard theory of finite dimensional LQ control problems, one should expect that, under Assumptions~\ref{as:finit_cost} and~\ref{as:positivity-qr} stated below, the optimal control $\underline{\bU}^{*,N}$ exists and is of the form
$\underline{U}^{*,N}_t = \Phi^{*,N} \underline{X}_t$ 
for a deterministic constant matrix $\Phi^{*,N} \in \RR^{(\ell N) \times (d N)}$.

\subsection{The limit mean-field control model}
\label{subsec:OC-MKV}

The model introduced in this section should be understood as the limit when $N\to\infty$ of the $N$-agent model introduced above. The time evolution of the state $X_t\in\RR^d$ of a representative agent who uses action $u_t \in \RR^{\ell}$ at time $t$, should be given by the equation:
\begin{align}
    X_{t+1} &= F( X_t, u_t, \bar{X}_t, \bar{u}_t, \comnoise{t+1}, \idynoise{t+1})
    \nonumber \\
    & = \rA X_t + \bar{\rA}\bar{X}_t  + \rB u_t + \bar{\rB}\bar{u}_t + \varepsilon^0_{t+1} + \varepsilon_{t+1},
\label{fo:MKV-state}
\end{align}
for $t \geq 0$, with an initial state $X_0 =\varepsilon^0_{0}+\varepsilon_{0}$.
Recall that we do not assume that the initial perturbations $\varepsilon_0^0$ and $\varepsilon_0$ are mean zero. Let $\tilde \mu^0_0 = \cL(\varepsilon_0^0)$ and $\tilde \mu_0 = \cL(\varepsilon_0)$ be their respective distributions.
The terms $\bar{X}_t = \EE[X_t|\cF^{0}_t]$ and $\bar{u}_t = \EE[u_t|\cF^{0}_t]$ appearing in the dynamics~\eqref{fo:MKV-state} of the state are random variables. They are the conditional expectations of $X_t$ and $u_t$ given $\cF^{0}_t$, the past of the common noise up to and including time $t$. 
The mean field (MF) cost function is
\begin{equation}
    \label{fo:MKV-discounted_cost}
    J(\bu)=\EE \Big[ \sum_{t \geq 0} \gamma^t c(X_t, \bar{X}_t, u_t, \bar{u}_t) \Big].
\end{equation}
Notice that the expected instantaneous cost at time $t$ reads
\begin{equation}
\label{eq:expression_expected_instateneous_cost}
\EE[ c(X_t, \bar{X}_t, u_t, \bar{u}_t)] = \EE[ X_t^\top Q X_t + (\bar{X}_t)^\top \bar Q \bar{X}_t + u_t^\top R u_t + (\bar{u}_t)^\top \bar R \bar{u}_t ].
\end{equation}
The set of \textit{admissible control processes} will be  denoted by $\cU_{ad}$ and defined as the collection of sequences of random variables on $(\Omega,\cF,\PP)$ adapted to the filtration $(\cF_t)_{t \geq 0}$ and square integrable (a.k.a. $L^2-$discounted integrable):

\begin{equation}
\label{def:admissible_control}
    \cU_{ad} := \left\{ (u_t)_{t \geq 0} \left|\, u_t \in \RR^{\ell} \text{ is } \cF_t-\text{measurable}, \quad \EE \Big[ \sum_{t \geq 0} \gamma^t \Vert u_t \Vert^2 \Big] <\infty  \right. \right\}.
\end{equation}
We also consider the set $\cX$ of $L^2-$discounted integrable processes in $\RR^{d}$ defined by:
\begin{equation}
\label{def:L2_discoutned_integratble_process}
    \cX := \left\{ (X_t)_{t \geq 0} \left|\, X_t \in \RR^{d} \text{ is } \cF_t-\text{measurable}, \quad \EE \Big[ \sum_{t \geq 0} \gamma^t \Vert X_t \Vert^2 \Big] <\infty  \right. \right\}.
\end{equation}
Note that both $\cU_{ad}$ and $\cX$ are linear spaces.  
The MFC problem is then:
\begin{equation}
    \label{pb:MFC_L2_admissible_control}
    \inf_{\bu \in \cU_{ad}} J(\bu).
\end{equation}

The above stochastic control problem has some unique characteristics: 1) it is set in infinite horizon; 2) it includes a \emph{common noise} $\bvarepsilon^0$, and 3) the interaction is not only through the conditional mean of the state, but also through the conditional mean of the control. 

\subsection{Standing assumptions and well-definiteness}
\label{sec:assm-defi}

We shall also use the following three assumptions.
\begin{assumption}
\label{as:finit_cost}
    The discount coefficient $\gamma \in (0,1)$ and the matrices $\rA$ and $\bar\rA$ satisfy $\gamma \| \rA \|^2 < 1$ and $\gamma \| \rA + \bar \rA \|^2 < 1$.
\end{assumption}

\begin{assumption}
\label{as:positivity-qr}
    The matrices $Q, \bar Q, R, \bar R$ are symmetric and they satisfy $Q \succeq 0, Q + \bar Q \succeq 0, R \succ 0, R + \bar R \succ 0$.
\end{assumption}
 
\begin{assumption}[Non-degeneracy] The following holds
\label{as:non-deg}
    \begin{equation}
    \label{eq:non_degenerate_condition}
        \max \Big\{ \lambda_{min} \big(\Sigma_{y_0} \big), \lambda_{min} \big( \Sigma^1 \big) \Big\} > 0,
        \quad
        \max \Big\{ \lambda_{min} \big( \Sigma_{z_0} \big), \lambda_{min} \big(\Sigma^0 \big) \Big\} > 0,
    \end{equation}
    where $y_0 = \varepsilon_0 - \EE[ \varepsilon_0 ]$, $z_0 = \varepsilon_0^0 + \EE[ \varepsilon_0]$, and 
    $
    \Sigma_{y_0} = \mathbb{E} \left[y_0 (y_0)^\top \right]
    $,
    $ 
      \Sigma_{z_0} = \mathbb{E} \left[z_0 (z_0)^\top \right]
    $,
    and where the variance matrices of $(\varepsilon_{t+1}, \varepsilon_{t+1}^0)$ for any $t \geq 0$ in the i.i.d. noise sequences $(\bvarepsilon, \bvarepsilon^0)$ are given by
    $
      \Sigma^1 = \EE[\varepsilon_{t+1} (\varepsilon_{t+1})^\top]
    $,
    $
      \Sigma^0 = \EE[\varepsilon^0_{t+1} (\varepsilon^0_{t+1})^\top]
    $.
\end{assumption}
As we shall see in the following sections, Assumption~\ref{as:finit_cost} is required for the well-definiteness of the state process and the finiteness of the MF cost~\eqref{fo:MKV-discounted_cost}. We show in Section~\ref{sec:analysis} that together with~\ref{as:positivity-qr}, it implies existence and uniqueness of the optimal control process for the MFC problem. Assumption~\ref{as:non-deg} is critical for the global convergence results of PG algorithms discussed in Section~\ref{ref:PGconv}.

We will sometimes use the notation $\bX^\bu$  to emphasize the fact that the state process is controlled by the admissible control process $\bu \in \cU_{ad}$.
Lemma~\ref{lemma:L2_discounted_integrable_of_state_process} and Proposition~\ref{proposition:admissible_of_X} below justify our choices for the space $\cU_{ad}$ of admissible control processes and the set $\cX$ for controlled state processes.

\begin{lemma}
\label{lemma:L2_discounted_integrable_of_state_process}
    If $\gamma \| \rA \|^2 < 1$, if $ (\beta_t)_{t \geq 0} \in \cU_{ad}$, and if $Y_0$ is $\cF_0$-measurable and square-integrable, a  process $\bY= (Y_t)_{t \geq 0}$ satisfying 
    \begin{equation}
    \label{eq:sys_dyn_Y}
        Y_{t+1} = \rA Y_t +  \beta_{t+1},\qquad t \geq 0,
    \end{equation}
    belongs to  $\cX$.
\end{lemma}

\begin{proof}
$Y_t$ is $\cF_{t}$-measurable by construction. For every $t \geq 1$,
\begin{align*}
    \EE\left[ \gamma^t \| Y_t \|^2 \right] &= \EE \bigg[ \gamma^t \Big\| \rA^t Y_0 + \sum_{s=0}^{t-1} \rA^{t-1-s} \beta_{s+1} \Big\|^2 \bigg]  
    \\
    & =  \EE \bigg[ \Big\| (\gamma^{1/2} \rA)^t Y_0 + \sum_{s=1}^{t} (\gamma^{1/2} \rA)^{t-s} \big(\gamma^{s/2} \beta_{s} \big)  \Big\|^2 \bigg] %
\end{align*}
Because $\gamma \| \rA \|^2 < 1$, there exists a real number $1 > \rho > \gamma$ such that $\xi := \rho \| A \|^2 < 1$. Let $\eta := \gamma \|A \|^2 / \xi $, then $\eta = \gamma / \rho < 1$. 
From the Fubini's theorem, we switch the two summations in the following equality:
$$
    \sum_{t=1}^\infty \Big( \sum_{s=1}^{t} \xi^{t-s} \EE\big[ \gamma^s  \| \beta_s \|^2 \big] \Big) =  \sum_{s=1}^{\infty} \Big( \sum_{t=0}^\infty  \xi^t \Big) \EE\big[ \gamma^{s} \| \beta_s \|^2  \big] = \frac{1}{1 - \xi} \sum_{s \geq 0} \gamma^{s+1} \EE[ \| \beta_{s+1} \|^2 ] < \infty.
$$
Consequently, by summing up $\EE[ \gamma^t \| Y_t\|^2 ]$ from $t=0$ to $\infty$, we obtain that
\begin{align*}
    \sum_{t \geq 0}  \EE\left[ \gamma^t \| Y_t \|^2 \right]  & \leq \EE\big[ \| Y_0\|^2 \big] + \sum_{t \geq 1} \EE \bigg[ \Big( \| \gamma^{1/2} \rA \|^t  \| Y_0 \|  + \sum_{s=1}^{t} \| \gamma^{1/2} \rA \|^{t-s} \gamma^{s/2} \| \beta_s \| \Big)^2 \bigg] 
    \\
    & =  \EE\big[ \| Y_0\|^2 \big] + \sum_{t\geq 1}  \EE \bigg[ \Big( \eta^{t/2} (\xi^{t/2} \| Y_0 \| ) + \sum_{s=1}^{t} \eta^{(t-s)/2} \xi^{(t-s)/2} \gamma^{s/2} \| \beta_s \| \Big)^2 \bigg]
    \\
    & \leq  \EE\big[ \| Y_0\|^2 \big] +  \sum_{t\geq 1}  \EE \bigg[ \Big(\eta^t + \sum_{s=1}^{t} \eta^{t-s} \Big). \Big( \xi^t \|Y_0 \|^2 + \sum_{s=1}^{t} \xi^{t-s} \gamma^{s} \| \beta_s \|^2 \Big)   \bigg]
    \\
    & \leq \EE \big[ \|Y_0 \|^2 \big] + \frac{1}{1- \eta} \bigg( \sum_{t = 1}^\infty \xi^t \EE \big[ \|Y_0 \|^2 \big] +  \sum_{t = 1}^\infty \sum_{s=1}^{t} \xi^{t-s}  \EE \big[ \gamma^{s} \| \beta_s \|^2 \big] \bigg) 
    \\
    & \leq \frac{1}{1 - \eta} \frac{1}{1 - \xi} \bigg( \EE\big[ \| Y_0\|^2 \big] + \sum_{s \geq 0} \gamma^{s+1} \EE \big[ \| \beta_{s+1} \|^2 \big] \bigg)
    \\
    & < \infty,
\end{align*}
where the first inequality is due to the matrix norm inequality $\| \rA^t \| \leq \| \rA \|^t$ for any $t \geq 1$, and the third inequality is justified by Cauchy-Schwarz' inequality.
\end{proof}

\begin{proposition}
\label{proposition:admissible_of_X}
Under Assumption~\ref{as:finit_cost}, for any admissible control $\bu \in \cU_{ad}$, the corresponding state process $\bX^\bu$ is in $\cX$and the MF cost $J(\bu)$ is finite.
\end{proposition}

\begin{proof}
The dynamics of the processes $(X_t^\bu - \bar X^{\bu}_t)_{t \geq 0}$ and $(\bar X_t^\bu)_{t \geq 0}$ are given by:
    \begin{align}
        X_{t+1}^\bu - \bar X_{t+1}^\bu &= \rA \big(X_t^\bu - \bar X^{\bu}_t\big) + \rB (u_t - \bar u_t) + \varepsilon_{t+1},
        \label{eq:Xt_minus_barXt_u}
        \\
       \bar X_{t+1}^\bu &= (\rA + \bar \rA)  \bar X^{\bu}_t + (\rB + \bar \rB) \bar u_t + \varepsilon_{t+1}^0,
       \label{eq:barXt_u}
    \end{align}
    for $t \geq 0$, with the initial conditions $X_0^\bu - \bar X^{\bu}_0 = \varepsilon_0 - \EE[ \varepsilon_0]$ and $\bar X^{\bu}_0 = \varepsilon_0^0 + \EE[ \varepsilon_0]$. Let $\beta_{t+1} = \rB (u_t - \bar u_t) + \varepsilon_{t+1}$ for every $t \geq 0$ and $\beta_0 = 0$. Because $\bu \in \cU_{ad}$ and $(\varepsilon_{t+1})_{t \geq 0}$ are i.i.d., by Cauchy-Schwarz' inequality and Jensen's inequality for conditional expectations we have
    $$
        \EE \big[ \| \beta_{t+1} \|^2 \big] \leq \big( \| \rB \|^2 + \| \rB \|^2 + 1 \big) \big( \EE \big[ \| u_t \|^2 \big] + \EE\big[ \EE[ \| u_t \|^2 \, | \, \cF_t^{0} ] \big] + \EE\big[ \| \varepsilon_{t+1} \|^2 \big] \big) 
    $$
    for every $t \geq 0$, and so
    $$
        \sum_{t \geq 0} \gamma^t \EE\big[ \| \beta_t \|^2 \big] \leq  ( 2 \| \rB \|^2 + 1)  \Big( 2 \sum_{t \geq 0} \gamma^t \EE[ \| u_t \|^2 ]  + \sum_{t \geq 0} \gamma^t \EE[ \| \varepsilon_{t+1} \|^2 ] \Big) < \infty.
    $$
    Thus, $(\beta_t)_{t\geq 0} \in \cU_{ad}$ because $\beta_t$ is $\cF_t$-measurable so we can apply Lemma~\ref{lemma:L2_discounted_integrable_of_state_process}. Under Assumption~\ref{as:finit_cost}, we have $\gamma \| \rA \|^2 < 1$, then we have the process $(X_t^\bu - \bar{X}_t^\bu)_{t \geq 0} \in \cX$.  
    
    Similarly, if we now let $\beta_{t+1} = (\rB + \bar \rB) \bar u_{t} + \varepsilon_{t+1}^0$ and $\beta_0 = 0$, then we have $(\beta_t)_{t \geq 0} \in \cU_{ad}$, and we can use again Lemma~\ref{lemma:L2_discounted_integrable_of_state_process} with $Y_t=\bar X_t^\bu$ for every $t \geq 0$ and $\rA$ replaced by $\rA + \bar \rA$ in Lemma~\ref{lemma:L2_discounted_integrable_of_state_process} and conclude that $\bX^\bu \in \cX$.
   From the definition \eqref{fo:MKV-discounted_cost} of the MF cost $J(\bu)$ and the expression for the expectation of the instantaneous cost given in \eqref{eq:expression_expected_instateneous_cost}, we see  that 
    $$
        J(\bu) \leq (\| Q \| +  \| \bar Q \|)   \sum_{t \geq 0} \gamma^t \EE \big[ \| X_t^\bu \|^2 \big] + ( \| R \| +  \| \bar R \|) \sum_{t \geq 0} \gamma^t \EE \big[ \| u_t \|^2 \big]  < \infty,
    $$
    by Jensen's inequality. This concludes the proof.
\end{proof}

\begin{corollary}
\label{corollary:uniqueness_of_controlled_state_process}
Under Assumption~\ref{as:finit_cost}, for any admissible control process $\bu \in \cU_{ad}$,  there exists a unique state process $\bX \in \cX$ controlled by $\bu$
    starting from $X_0 = \varepsilon_0 + \varepsilon_0^0$ satisfying \eqref{fo:MKV-state} with noise processes $(\bvarepsilon, \bvarepsilon^0)$.
\end{corollary}

\begin{proof}
    The existence of controlled state process $\bX \in \cX$ of $\bu \in \cU_{ad}$ is justified by Proposition~\ref{proposition:admissible_of_X}. For the uniqueness, suppose that there exists another process $\bX' \in \cX$ starting from the same initial state $X_0' = X_0$ and following the dynamics~\eqref{fo:MKV-state} with the same noise processes $\bvarepsilon$ and $\bvarepsilon^0$. Then, the difference between $X_t$ and $X'_t$ satisfies the dynamics
   \begin{equation}
       \label{fo:diff}
        X_{t+1} - X'_{t+1} = \rA (X_t - X_t') + \bar{\rA} \EE\big[ X_t - X_t' \, | \cF_t^{0} \big]
\end{equation}
    for every $t \geq 0$ with $X_0 - X_0' = 0$ almost surely. Taking the conditional expectation $\EE[ \cdot | \cF_{t}^{0}]$ on both sides of  the equation \eqref{fo:diff}, we deduce that $\EE[ X_t - X_t' \, | \, \cF_t^{0} ] = 0$ almost surely for every $t \geq 0$. Plugging this back into equation \eqref{fo:diff}, we obtain that $X_t = X_t'$ almost surely for every $t \geq 0$. 
\end{proof}

\begin{remark}
\label{remark:equivalence_of_welldefineness_and_admissibility}
Under~\ref{as:positivity-qr}, we have $R, R + \bar R \succ 0$, and thus $\EE[ (u_t - \bar {u}_t)^\top R (u_t - \bar u_t) ] \geq  \lambda_{min}(R) \EE[ \| u_t - \bar u_t \|^2 ] $ and $\EE[ (\bar{u}_t)^\top (R + \bar R) \bar{u}_t] \geq \lambda_{min}(R + \bar R) \EE[ \| \bar u_t \|^2 ]$ for every $t \geq 0$. 
Thus, if the MF cost $J(\bu)$ if finite for an $(\cF_t)_{t \geq 0}$-adapted control process $\bu = (u_t)_{t \geq 0}$, we infer that $(u_t - \bar u_t)_{t \geq 0}$ and $(\bar u_t)_{t \geq 0}$ are $L^2-$discounted integrable, and that the control process $\bu$ is admissible in the sense that $\bu \in \cU_{ad}$.  Together with Proposition~\ref{proposition:admissible_of_X}, we see that under Assumptions~\ref{as:finit_cost} and~\ref{as:positivity-qr}, the condition of the $L^2-$discounted integrability in $\cU_{ad}$ is equivalent to the finiteness of the MF cost $J(\bu)$. 
\end{remark}

In Section~\ref{sec:analysis}, we will show existence and uniqueness of an optimal admissible control process for the MFC problem~\eqref{pb:MFC_L2_admissible_control} using a discrete-time version of Pontryagin's maximal principle. Moreover, we will also show that the optimal control process $\bu$ is in a feedback form and is linear in $(\bX^\bu, \bar \bX^\bu)$, with coefficients independent of time. To rigorously establish these results with the aforementioned notion of admissible control, we define the following set of matrices used for parameterizing an admissible control process. 

\begin{definition}
    The set of \textbf{admissible control parameters} is defined as
    \begin{equation}
        \label{eq:admissible_parameters_LQMF}
        \Theta := \left\{ \theta = (K,L) \in \mathbb{R}^{\ell \times d} \times \mathbb{R}^{\ell \times d} \  | \  \gamma \| \rA - \rB K\|^2 < 1, \, \gamma \|  \rA + \bar \rA  - (\rB + \bar \rB) L\|^2 < 1 \right\}.
    \end{equation} 
\end{definition}

It is worth noting that in our discrete-time set-up, contrary to continuous-time analogs, there is no issue of existence of state process controlled by a control process given in feedback form. The following lemma shows that when the feedback form is defined through a policy function $\pi_\theta : (x, \bar x) \mapsto \pi_\theta(x, \bar x)$ and the control parameter $\theta$ is admissible in the sense of the above definition, then the control process $\bu$ is also admissible.

\begin{lemma}
\label{lemma:admissible_of_control_with_feedback_form_on_theta}
    Consider an admissible control parameter $\theta = (K, L) \in \Theta$ and a control process $\bu = (u_t)_{t \geq 0}$ in $\RR^{\ell}$ in a feedback form satisfying 
    \begin{equation}
    \label{eq:feedback_linear_form_of_control}
        u_t = - K (X_t^\bu - \bar X_t^\bu) - L \bar X_t^{\bu}
    \end{equation}
    for every $t \geq 0$. Then $\bu \in \cU_{ad}$, and $\bX^{\bu} \in \cX$.
\end{lemma}

\begin{proof}
The idea of the proof is to replace $u_t - \bar u_t$ and $\bar u_t$ in the dynamics~\eqref{eq:Xt_minus_barXt_u}--\eqref{eq:barXt_u} with expressions $u_t - \bar u_t = - K (X_t^\bu - \bar X_t^\bu)$ and $\bar u_t = - L \bar X_t^\bu$, so that
\begin{align}
    X_{t+1}^\bu - \bar X_{t+1}^\bu & = \big( \rA - \rB K \big) (X_t^\bu - \bar{X}_t^\bu) + \varepsilon_{t+1}
    \\
    \bar X_{t+1}^\bu & = \big( \rA + \bar \rA - (\rB + \bar \rB) L \big) \bar X_t^\bu + \varepsilon_{t+1}^0
\end{align}
for every $t \geq 0$. The initial conditions are $X_0^\bu - \bar X_0^\bu = \varepsilon_0 - \EE[ \varepsilon_0 ]$, $\bar X_0^\bu = \varepsilon_0^0 + \EE[ \varepsilon_0 ] $.  Now, we choose $(Y_t, \beta_t) = (X_t^\bu - \bar X_t^\bu, \varepsilon_t)$ for every $t \geq 1$, $\beta_0 = 0$, and $Y_0 = \varepsilon_0 - \EE[ \varepsilon_0]$ in Lemma~\ref{lemma:L2_discounted_integrable_of_state_process}. Because the process $(\varepsilon_{t+1})_{t \geq 0}$ is i.i.d. with square-integrable terms, the process $(\beta_t)_{t \geq 0} \in \cU_{ad}$ and the process $(Y_t)_{t \geq 0} = (X_t^\bu - \bar X_t^\bu)_{t \geq 0}$ is adapted to the filtration $(\cF_t)_{t \geq 0}$. 
Since the control parameter $K$ satisfies $\gamma \| \rA - \rB K \|^2 < 1$, Lemma~\ref{lemma:L2_discounted_integrable_of_state_process} implies that $(X_t^\bu - \bar X_t^{\bu})_{t \geq 0} \in \cX$. Similar arguments yield that $(\bar X_t^\bu) \in \cX$. Hence, 
$$
    \sum_{t \geq 0} \gamma^t \EE[ \| u_t \|^2 ] \leq 2 \| K \|^2 \sum_{t \geq 0} \gamma^t \EE[ \| X_t^\bu - \bar X_t^\bu \|^2 ] + 2 \| L \|^2 \sum_{t \geq 0} \gamma^t \EE[ \| \bar X_t^\bu \|^2 ] < \infty. 
$$
We conclude that  $\bu = (u_t)_{t \geq 0} \in \cU_{ad}$, and  $\bX^\bu = (X_t^\bu)_{t \geq 0} \in \cX$.
\end{proof}

For an admissible control parameter $\theta = (K,L) \in \Theta$, we define its corresponding policy function by:
$
\pi_{\theta}: (x,\bar{x}) \mapsto -K (x - \bar{x}) - L \bar{x}.
$
We say that an admissible control process $\bu = (u_t)_{t \geq 0} \in \cU_{ad}$ is \textbf{parameterized} by $\theta$ if, for every $t \geq 0$, $u_t = \pi_{\theta}(X_t^\bu, \bar{X}_t^\bu)$. In this case, we may use the notations $\bu^\theta$ and $\bX^\theta$ to stress the dependence on $\theta$.  
The collection of \textbf{(admissible) parametrized control processes} is denoted by $\cU^\Theta_{ad} = \{ \bu^\theta \, | \, \theta \in \Theta \} \subset \cU_{ad}$.

\begin{remark}
$\Theta$ is a convex subset of $\RR^{\ell \times d} \times \RR^{\ell \times d}$ because the spectral norm $\| \cdot \|$ is convex. In Section~\ref{ref:PGconv}, we will work on PG algorithms within the set $\Theta$ to search for the optimal parameterized control process. Along the iterations of the algorithm, we not only focus on the admissibility of $\theta \in \Theta$, but we would like it to be in a level-set $\{ \theta \, | \, J(\bu^\theta) \leq C_0 \}$ for some predefined constant $C_0 \in \RR$ for convergence purposes. Unfortunately, the set of control parameters $\{ \theta \, | \, \theta \in \Theta, J(\bu^\theta) \leq C_0\}$ may not be convex because the mapping $\Theta \ni \theta \mapsto J(\bu^\theta)$ is not convex. The importance of level-set considerations in the global convergence of PG algorithms has been discussed in~\cite{frikha2024full} for a continuous-time LQMFC problem with common noise.
\end{remark}

It is worth pointing out that our definition of the set $\cU_{ad}$ of admissible controls is slightly different from the one considered in~\cite{fazel2018global}. Remark~\ref{remark:equivalence_of_welldefineness_and_admissibility} argues that under our standing assumptions, we have $\bu \in \cU_{ad} \Leftrightarrow J(\bu) < \infty$. The authors of \cite{fazel2018global} focus on a control coefficient $K \in \RR^{\ell \times d}$ such that the \emph{spectral radius} of $\rA - \rB K$ satisfies $\rho(\rA - \rB K) < 1$. Unfortunately, the set $\{ K \,| \, \rho(\rA - \rB K) < 1\}$ is not convex. Given that $\rho(\rA - \rB K) \leq \| \rA - \rB K \|$ for any matrix $K \in \RR^{\ell \times d}$, our admissibility condition for $\Theta$ is more restrictive than the condition considered in~\cite{fazel2018global} with discount factor $\gamma$. As a result, a proper analysis of the admissibility of the optimal control is necessary. We will provide more details in Section~\ref{sec:admissibility_of_opt_control}.

\section{\textbf{Analysis of the MFC Problem}}
\label{sec:analysis}

The following theorem is the main result of this section. Not only does it includes existence and uniqueness of an optimal control, but it also identifies a specific feedback form this optimal control needs to have. This necessary structure of the optimal control will allow us to concentrate the numerical implementations of the policy gradient search for optimality to an explicit set of parameters on which we will be able to prove convergence of the gradient descent algorithm, even in model free circumstances. Not surprisingly, the search for a necessary form of the optimal control will require us to prove a form of the Pontryagin maximum principle. In fact, before proving the final feedback form \eqref{eq:optimal_control_linear_in_x_barx}, we will show that the optimal control is of the form
\begin{equation}
\label{eq:control_ut_expression_2}
    u_t = \Gamma ( p_t - \bar{p}_t) + \Lambda \bar{p}_t, \qquad \bar u_t = \Lambda \bar{p}_t
\end{equation}
for every $t \geq 0$, where $(p_t)_{t\ge 0}$ is the adjoint process corresponding to $(u_t, X_t)_{t \geq 0}$ (see Definition~\ref{def:adjoint_process} below).

\vskip 4pt
In order to simplify the notations, we introduce the quantities $\tilde\rA = \rA + \bar \rA$, $\tilde \rB = \rB + \bar \rB$, $\tilde Q = Q + \bar Q$, and $\tilde R = R + \bar R$.

\begin{theorem}
\label{thm:existence_linear_control}
Under Assumptions \ref{as:finit_cost} and~\ref{as:positivity-qr}, there exists a unique optimal control for the MFC problem~\eqref{pb:MFC_L2_admissible_control}. The optimal control $\bu = (u_t)_{t \geq 0}$ is in feedback form, and is linear in $(\bX^\bu, \bar{\bX}^\bu)$, specifically:
\begin{equation}
\label{eq:optimal_control_linear_in_x_barx}
    u_t = - K^* (X_t^\bu - \bar X_t^\bu) - L^* \bar X_t^\bu
\end{equation}
for every $t \geq 0$, with $K^* =-\Gamma P$ and $L^* = -\Lambda \bar P$ where
\begin{equation}
\label{eq:Gamma_Lambda}
    \Gamma = - \frac{1}{2} R^{-1} \rB^\top, \qquad \Lambda = - \frac{1}{2} \tilde{R}^{-1} \tilde{B}^\top,
\end{equation}
and $P$ and $\bar P$ satisfy the two matrix Riccati equations
\begin{equation}
\label{eq:Riccati}
\left\{
    \begin{array}{rcl}
        P & = & \gamma \big( \rA^\top P + 2Q \big) \big( \rA -  \rB R^{-1} \rB^\top P / 2 \big) 
        \\
        \bar P & = & \gamma \big( \tilde{\rA}^\top \bar P + 2 \tilde Q \big) \big( \tilde\rA -  \tilde\rB \tilde{R}^{-1} \tilde{\rB}^\top \bar P / 2 \big).
    \end{array}
\right.
\end{equation}
\end{theorem}

\subsection{Admissibility of the optimal control}
\label{sec:admissibility_of_opt_control}

Before we show the optimality of the control process $\bu$ defined by equation~\eqref{eq:optimal_control_linear_in_x_barx}, we justify its admissibility with the help of Lemma~\ref{lemma:admissible_of_control_with_feedback_form_on_theta}. The idea is to show that the pair of matrices $\theta^* = (K^*, L^*) \in \RR^{\ell \times d} \times \RR^{\ell \times d}$ defined in Theorem~\ref{thm:existence_linear_control} is such that $\theta^* \in \Theta$. 
More precisely, the first step is to construct a pair of solution matrices $(P, \bar P)$ satisfying the Riccati equations~\eqref{eq:Riccati}, and then we show that the matrices $K^*$ and $L^*$ satisfy the ``strictly smaller than 1 norm'' condition. 
To start, we construct $(P, \bar P)$ from the solutions $(P^{*, y}, P^{*, z})$ of another pair of discrete-time Algebraic Riccati equations (DAREs):
\begin{equation}
    \label{eq:DARE_closeloop_feedback_excat_PG}
        \left\{
        \begin{array}{rcl}
            P^{*, y} &=& Q  + \gamma \rA^\top P^{*, y} \rA - \gamma^2 \rA^\top P^{*, y} \rB ( R + \gamma \rB^\top P^{*, y} \rB )^{-1} \rB^\top P^{*, y} \rA,
            \\
            P^{*, z} &=& \tilde{Q}  + \gamma \tilde{\rA}^\top P^{*, z} \tilde{\rA} - \gamma^2 \tilde{\rA}^\top P^{*, z} \tilde{\rB} ( \tilde{R} + \gamma \tilde{\rB}^\top P^{*, z} \tilde{\rB} )^{-1} \tilde{\rB}^\top P^{*, z} \tilde{\rA}.
        \end{array}
        \right.
\end{equation}

\begin{lemma}
\label{lemma:solution_to_ARE_openloop_from_closefeedback_DAREs}
    Under Assumptions \ref{as:finit_cost} and \ref{as:positivity-qr}, there exist symmetric and positive semi-definite matrices $P^{*, y}$ and $P^{*, z}$ in $\RR^{d \times d}$ solving the (DAREs)~\eqref{eq:DARE_closeloop_feedback_excat_PG}.
    Moreover, if we define two matrices $P, \bar P \in \RR^{d \times d}$ by
    \begin{equation}
    \label{eq:expression_of_solution_openloop_DARE}
        \left\{
            \begin{array}{ll}
            P &:= 2 \gamma P^{*, y} \big( I_d + \gamma \rB R^{-1} \rB^\top P^{*, y} \big)^{-1} \rA  
            \\
            \bar P &:= 2 \gamma  P^{*, z} \big( I_d + \gamma \tilde{\rB} \tilde{R}^{-1}  \tilde{\rB}^\top  P^{*, z} \big)^{-1} \tilde{\rA},
            \end{array}
        \right.
    \end{equation}
    then $P$ and $\bar P$ satisfy the matrix Riccati equations~\eqref{eq:Riccati}.
\end{lemma}

\begin{proof}
    Under Assumption~\ref{as:finit_cost}, we know that the two pairs of matrices $(\sqrt{\gamma} \rA, \sqrt{\gamma} \rB)$ and $(\sqrt{\gamma} \tilde \rA, \sqrt{\gamma} \tilde \rB)$ are stabilizable in the sense that there exist $\tilde K, \tilde L \in \RR^{\ell \times d}$ satisfying that all eigenvalues of $\sqrt{\gamma} (\rA - \rB \tilde K)$ and $\sqrt{\gamma} (\tilde \rA - \tilde \rB \tilde L)$ fall strictly inside the unit circle, i.e. $\Theta$ is not empty (take for example  $\tilde K = \tilde L = 0$). Under Assumption~\ref{as:positivity-qr}, we have $Q, \tilde Q \succeq 0$ and $R, \tilde R \succ 0$, so by~\cite[Theorem 3.2 (iii)]{ran1988existence} we know that the (DAREs)~\eqref{eq:DARE_closeloop_feedback_excat_PG} admit solutions, and the solutions $P^{*, y}$ and $P^{*, z}$ are also symmetric and positive semi-definite.

    Now, to show that $P, \bar P$ given by~\eqref{eq:expression_of_solution_openloop_DARE} are solutions to~\eqref{eq:Riccati}, we first recall a standard matrix identity: for any matrix $W \in \RR^{d \times d}$ satisfying that $I_d + W$ is invertible,
    $$
        (I_d + W)^{-1} = I_d - W (I_d + W)^{-1} = I_d - W + W (I_d + W)^{-1} W.
    $$
    Let $W = \gamma \rB^\top R^{-1} \rB P^{*, y}$, $U = \rB^\top (R^{-1/2})^\top$, and $V = R^{-1/2} \rB P^{*, y}$. Because $P^{*, y} \succeq 0$ and $R \succ 0$, we have
    $$
    \lambda_{min}(W) = \gamma \lambda_{min}(UV) = \gamma \lambda_{min}(VU) \geq 0.
    $$
    So $I_d + W$ is invertible. Because $R + \gamma \rB^\top P^{*, y} \rB \succ 0$ and 
    \begin{align*}
       I_d & = I_d + \gamma \rB^\top P^{*, y} \big( I_d - (I_d + W)^{-1} - W ( I + W )^{-1} \big) \rB R^{-1}
        \\
        & =  (R + \gamma \rB^\top P^{*, y} \rB). \big( R^{-1} - \gamma R^{-1} \rB^\top P^{*,y} \big( I_d + \gamma \rB R^{-1} \rB^\top P^{*, y} \big)^{-1} \rB R^{-1} \big),
    \end{align*}
    we have
    \begin{align*}
        & I_d - \gamma \rB \big( R + \gamma \rB^\top P^{*, y} \rB \big)^{-1} \rB^\top P^{*, y}
        \\
        = & I_d - \gamma \rB. \big( R^{-1} - \gamma R^{-1} \rB^\top P^{*,y} \big( I_d + \gamma \rB R^{-1} \rB^\top P^{*, y} \big)^{-1} \rB R^{-1} \big). \rB^\top P^{*,y}
        \\
        = & I_d - \gamma \rB R^{-1} \rB^\top P^{*, y} + ( \gamma \rB R^{-1} \rB^\top P^{*, y}). \big( I_d +  \gamma \rB R^{-1} \rB^\top P^{*, y} \big)^{-1}. ( \gamma \rB R^{-1} \rB^\top P^{*, y})
        \\
        = & I_d - W + W ( I_d + W)^{-1} W
        \\
        = & (I_d + W)^{-1}
        \\
        = & \big( I_d +  \gamma \rB R^{-1} \rB^\top P^{*, y} \big)^{-1}.
    \end{align*}
    Consequently, from the first equation in (DAREs)~\eqref{eq:DARE_closeloop_feedback_excat_PG}, we obtain
    \begin{align*}
         P^{*, y}  &=  Q + \gamma \rA^\top P^{*, y}  \Big(  I_d - \gamma \rB \big( R + \gamma \rB^\top P^{*, y} \rB \big)^{-1} \rB^\top P^{*, y} \Big) \rA 
         \\
         &= Q + \gamma \rA^\top P^{*,y} \big( I_d + \gamma \rB R^{-1} \rB^\top P^{*,y} \big)^{-1} \rA
         \\
         &= Q + \rA^\top P / 2.
    \end{align*}
    Moreover, we have
    \begin{align*}
        \rA - \rB R^{-1} \rB^\top P / 2 
        & = \big[ I_d - \gamma \rB R^{-1} \rB^\top P^{*, y} \big( I_d + \gamma \rB R^{-1} \rB^{\top} P^{*, y} \big)^{-1}  \big] \rA 
        \\
        & = \big(I_d + \gamma \rB R^{-1} \rB^\top P^{*, y} \big)^{-1} \rA.
    \end{align*}
    Thus, we derive the matrix Riccati equation for $P$ by
    \begin{align*}
        P = \gamma (2 P^{*,y}). \big(I_d + \gamma \rB R^{-1} \rB^\top P^{*, y} \big)^{-1} \rA =  \gamma ( 2 Q + \rA^\top P)( \rA - \rB R^{-1} \rB P / 2).
    \end{align*}
    We apply similar arguments to $(\tilde \rA, \tilde \rB, \tilde R, P^{*, z}, \bar P)$ and we derive the matrix Riccati equation for $\bar P = \gamma ( 2 \tilde{Q} + \tilde{\rA}^\top \bar P)( \tilde{\rA} - \tilde{\rB} \tilde{R}^{-1} \tilde{\rB} \bar P / 2)$. 
\end{proof}

\begin{corollary}
\label{corollary:admissibility_of_theta^*}
Under Assumption~\ref{as:finit_cost} and~\ref{as:positivity-qr}, if the two matrices $K^{*}$ and $L^{*}$ in $\RR^{\ell \times d}$ are defined by
$$
    K^{*} = R^{-1} \rB^\top P/ 2, \qquad L^{*} = \tilde R^{-1} \tilde \rB^\top \bar P / 2,
$$
as in Theorem~\ref{thm:existence_linear_control} where $(P, \bar P)$ are the matrices defined by~\eqref{eq:expression_of_solution_openloop_DARE} based on the solutions matrices $(P^{*, y}, P^{*, z})$ to (DAREs)~\eqref{eq:DARE_closeloop_feedback_excat_PG}, then
\begin{equation*}
\left\{
\begin{array}{rl}
    \sqrt{\gamma} ( \rA - \rB K^* ) &= \sqrt{\gamma} \big(I_d + \gamma \rB R^{-1} \rB^\top P^{*, y} \big)^{-1} \rA,
    \\
    \sqrt{\gamma} ( \tilde \rA - \tilde \rB L^* ) &=  \sqrt{\gamma} \big(I_d + \gamma \tilde{\rB} \tilde{R}^{-1} \tilde{\rB}^\top P^{*, z} \big)^{-1} \tilde{\rA},
\end{array}
\right.
\end{equation*}
and the control parameter $\theta^* = (K^*, L^*)$ is admissible in the sense that $\theta^* \in \Theta$.
\end{corollary}

\begin{proof}
    The definition of $(P, \bar P)$ in~\eqref{eq:expression_of_solution_openloop_DARE} implies
    \begin{align*}
        \sqrt{\gamma} ( \rA - \rB K^{* } ) & = \sqrt{\gamma} \big( \rA - \rB R^{-1} \rB^\top P / 2) 
        \\
        & = \sqrt{\gamma} \big[ I_d - \gamma \rB R^{-1} \rB^\top P^{*, y} \big( I_d + \gamma \rB R^{-1} \rB^{\top} P^{*, y} \big)^{-1}  \big] \rA
        \\
        & = \sqrt{\gamma} \big(I_d + \gamma \rB R^{-1} \rB^\top P^{*, y} \big)^{-1} \rA.
    \end{align*}
    Because $R \succ 0$ and $P^{*, y} \succeq 0$, we have 
    $$
        \lambda_{min}(\rB R^{-1} \rB^\top P^{*, y}) = \lambda_{min}( (R^{-1/2})^\top B^\top P^{*,y} B R^{-1/2} ) \geq 0.
    $$ 
    So $\lambda_{min}(I_d + \gamma B R^{-1} B^\top P^{*,y} ) \geq 1$. 
    Thus, under Assumption~\ref{as:finit_cost}, we deduce that
    $$
        \sqrt{\gamma} \| \rA - \rB K^{*} \| \leq \frac{\sqrt{\gamma}}{\lambda_{min}(I_d + \gamma B R^{-1} B^\top P^{*,y} )} \| A\| < 1.
    $$
    Similar arguments can be applied to $L^*$ and $P^{*, z}$, leading to $\sqrt{ \gamma} \| \tilde \rA - \tilde \rB L^{*} \| < 1$. 
\end{proof}

\begin{proposition}
\label{proposition:admissibility_of_u_with_feedback_control}
    Under Assumptions~\ref{as:finit_cost} and~\ref{as:positivity-qr}, the control parameter $\theta^* = (K^*, L^*) = (- \Gamma P, -\Lambda \bar P)$ with coefficients $(K^*, L^*, \Gamma, \Lambda, P, \bar P)$ defined in Theorem~\ref{thm:existence_linear_control} is admissible in the sense that $\theta^* \in \Theta$.  
    The corresponding parameterized control process $\bu^{\theta^*}$ in a feedback form of its controlled state process and conditional mean process $\big( \bX^{\bu^{\theta^*}}, \bar{\bX}^{\bu^{\theta^*}} \big)$ with policy function $\pi_{\theta^*} : (x, \bar x)  \mapsto - K^* (x - \bar x) - L^* \bar x$ is admissible in the sense that $\bu^{\theta^*} \in \cU_{ad}$.
\end{proposition}

\begin{proof}
    Corollary~\ref{corollary:admissibility_of_theta^*} states that the control parameter $\theta^* = (K^*, L^*) \in \Theta$. Because the process $\bu^{\theta^*}$ is defined in a feedback form from the policy $\pi_{\theta^*}$ in Theorem~\ref{thm:existence_linear_control}, Lemma~\ref{lemma:admissible_of_control_with_feedback_form_on_theta} implies immediately that $\bu^{\theta^*} \in \cU_{ad}$.
\end{proof}

\subsection{Duality and Pontryagin Maximum Principle}
In this subsection, we show existence and uniqueness of an optimal control with a discrete-time version of Pontryagin's maximal principle. 
For convenience, we use the notation $\check \rA = \rA - I_d$ and $\zeta = (x, \bar x, u, \bar u) \in \RR^{d} \times \RR^{d} \times \RR^{\ell} \times \RR^{\ell}$. We define the drift function $b$ by
\begin{equation}
\label{eq:drift_function_b}
    b(\zeta) = b(x, \bar x, u, \bar u) =  \check{\rA} x + \bar{\rA} \bar x + \rB u + \bar{\rB} \bar{u}
\end{equation}
and the Hamiltonian function $h$ by:
\begin{equation}
\label{eq:Hamilton_function_h}
    h(\zeta, p) = h(x, \bar x, u, \bar u, p) = b(\zeta) \cdot p + c(\zeta) - \delta x \cdot p
\end{equation}
where $p \in \RR^{d}$, $\delta = (1 - \gamma) / \gamma$, and $c(\zeta)=c(x, \bar x, u, \bar u)$ is the instantaneous cost function defined by~\eqref{eq:lq_one_step_cost}. 
Then, for an admissible control process $(u_t)_{t \geq 0}$, the state dynamics~\eqref{fo:MKV-state} of the controlled state process $(X_t)_{t \geq 0}$ can be rewritten as
\begin{equation}
\label{eq:delta_X_t_with_b}
    X_{t+1} - X_t = b(X_t, \bar X_t, u_t, \bar u_t) +  \varepsilon_{t+1}^0 + \varepsilon_{t+1} = b(\zeta_t) + \varepsilon_{t+1}^0 + \varepsilon_{t+1}.
\end{equation}
To compute the partial derivatives of $h$ with respect to its arguments $(x, \bar x, u, \bar u, p)$, we treat $\zeta$ as a $(2d + 2\ell)$ dimension vector, $\mathrm{Vect}(\zeta) = (x^\top, \bar x^\top, u^\top, \bar u^\top)^\top \in \RR^{(2d + 2 \ell)}$. For every fixed $p \in \RR^{d}$, we have 
\begin{equation}
\label{eq:gradient_h}
\nabla_{\zeta} h (\zeta, p) = 
	\begin{pmatrix}
	\partial_x h(\zeta,p)\\
	\partial_{\bar{x}} h(\zeta,p)\\
	\partial_{u} h(\zeta,p)\\
	\partial_{\bar{u}} h(\zeta,p)
	\end{pmatrix}
	=
	\begin{pmatrix}
	(\check{\rA}-\delta I_d)^\top p + 2Q (x-\bar{x}) \\
	\bar{\rA}^\top p - 2 Q (x - \bar x) + 2 \tilde Q \bar x
	\\
	\rB^\top p + 2 R (u - \bar{u})\\
	\bar{\rB}^\top p - 2 R (u - \bar u) + 2\tilde R \bar u
	\end{pmatrix},
\end{equation}
and the Hessian of $h$ w.r.t. $(u, \bar u)$ is given by 
\begin{equation}
\label{eq:hession_h_u}
   \nabla^2_{(u,\bar{u}),(u,\bar{u})} h(\zeta,p)=
	\begin{pmatrix}
	2R & -2R \\
	-2R & 2(R+\tilde{R})
	\end{pmatrix},
\end{equation}
which is positive definite because 
\begin{equation}
   \begin{pmatrix} u & \bar{u} \end{pmatrix}
    \begin{pmatrix}
	R & -R \\
	-R & R+\tilde{R}
	\end{pmatrix}
   \begin{pmatrix} u \\ \bar{u} \end{pmatrix}
=(u - \bar u) R( u - \bar u) + \bar u \tilde R \bar u > 0
\end{equation}
since, under Assumption~\ref{as:positivity-qr}, $R$ and $\tilde R$ are assumed to be (strictly) positive definite. 
We now introduce the notion of an adjoint process associated to a state process generated by an admissible control process.

\begin{definition}
\label{def:adjoint_process}
Given an admissible control process and the corresponding controlled state process $(\bu, \bX) = (u_t, X_t)_{t \geq 0}$, we say that an $\RR^d$-valued $(\cF_{t})_{t \geq 0}$-adapted process $\bp = (p_t)_{t \geq 0}$ is a corresponding adjoint process if it satisfies the backward equation:
\begin{equation}
\label{eq:adjoint}
    p_t = \EE \Big[ p_{t+1} +  \gamma \big[ (\check \rA - \delta I_d)^\top p_{t+1} + 2 Q X_{t+1} + \bar \rA^\top \bar p_{t+1} + 2 \bar Q \bar X_{t+1} \big] \big\vert \cF_t \Big],
\end{equation}
and the transversality condition:
\begin{equation}
    \label{eq:transversality}
    \EE \Bigl[ \sum_{t \geq 0} \gamma^t \| p_t \|^2 \Bigr] < \infty.
\end{equation}
\end{definition}
The transversality condition should be viewed as a \emph{terminal condition} controling the possible growth of $p_t$ when $t\to\infty$. Notice that equation~\eqref{eq:adjoint} can equivalently be written as:
\begin{equation}
   \label{eq:simpler_adjoint}
    p_t = \gamma \EE [ \rA^\top p_{t+1} + 2 Q X_{t+1} + \bar \rA^\top \bar p_{t+1} + 2 \bar Q \bar X_{t+1} \ | \  \cF_t ]. 
\end{equation}

\begin{proposition}
\label{proposition:existence_and_uniqueness_of_adjoint_process}
    Under Assumption~\ref{as:finit_cost}, for every state process generated by an admissible control process, there exists a unique adjoint process. 
\end{proposition}

\begin{proof}
Let $\bu \in \cU_{ad}$ be an admissible control process and $\bX = (X_t)_{t\geq 0}$ the corresponding controlled state process.

\textbf{(Uniqueness)} Consider two adjoint process $(p_t)_{t \geq 0}$ and $(p'_t)_{t \geq 0}$ corresponding to $(\bu, \bX)$. Taking conditional expectation $\EE[\cdot | \cF^{0}_t ]$ on both sides of equation~\eqref{eq:simpler_adjoint} and subtracting \eqref{eq:simpler_adjoint} for $\bar{\bp}$ to the corresponding equation for $\bar{\boldsymbol{p'}}$, we get 
$$
    \EE[ \| \bar p_t - \bar p_t' \| ] \leq \gamma^{1/2} \EE[ \| \bar p_{t+1} - \bar p_{t+1}' \| ]
$$ 
because $\sqrt{\gamma} \| \rA + \bar \rA \| < 1$. This implies that for every $0 \leq t \leq s$,
$$
    \gamma^{t/2} \EE [ \| \bar p_t - \bar p_t' \| ] \leq \gamma^{s/2} \EE [ \|\bar p_s -\bar p_s' \| ].
$$
From the transversality condition and Jensen's inequality for conditional expectations, we get:
$$
\lim_{s \to \infty} \gamma^{s/2} \EE \big[ \| \bar p_s - \bar p_s' \| \big]  \leq \lim_{s \to \infty} \left[(\EE[ \gamma^s \| p_s \|^2 ] )^{1/2} + (\EE[ \gamma^s \| p_s' \|^2 ] )^{1/2} \right]
 = 0.
$$
Hence $\bar p_t = \bar p_t'$, $\PP$-a.s., for every $t \geq 0$. Using similar arguments and the fact that $\gamma^{1/2} \| \rA \| < 1$,  we obtain $p_t - \bar p_t = p'_t - \bar p'_t$, $\PP$-a.s., for every $t \geq 0$. Therefore, we have the uniqueness for the adjoint process.

\textbf{(Existence)} We construct directly an adjoint process $\bp = (p_t)_{t \geq 0}$ with the following two auxiliary processes $\by = (y_t)_{t \geq 0}, \bz=(z_t)_{t \geq 0}$: for every $t \geq 0$, 
\begin{equation}
\left\{
\begin{array}{rcl}
    y_t & := & \sum_{k \geq 0}  \bigl( \gamma \rA^\top \bigr)^k \bigl( 2 \gamma Q \bigr)  \EE \bigl[ X_{t+1+k} - \bar{X}_{t+1+k}  \, | \, \cF_t \bigr]
    \\
    z_t & := & \sum_{k \geq 0} \bigl( \gamma \tilde \rA^\top \bigr)^{k} \bigl( 2 \gamma \tilde Q \bigr) \EE[ \bar{X}_{t+1 + k} \, | \, \cF_t]
    \\
    p_t & := & y_t + z_t.
\end{array}
\right.
\end{equation}
By definition, the process $\bp$ is $(\cF_t)_{t \geq 0}$-adapted and satisfies the backward adjoint equation~\eqref{eq:adjoint}.
Next, we show that the process $\bz$ (resp. $\by$) is $L^2-$discounted integrable. Under Assumption~\ref{as:finit_cost}, we introduce $\rho \in (0, 1)$ such that $\eta := \gamma / \rho < 1$ and $\xi := \rho \| \tilde \rA \| < 1$, so that $\eta \xi = \gamma \| \tilde \rA \|^2$. By Fatou's lemma and Cauchy-Schwarz inequality, for every $t \geq 0$, we have the bound: 
$$
\gamma^t \EE[ \| z_t \|^2 ] \leq \frac{4\gamma \| \tilde Q \|^2}{1 - \eta} \EE\big[ \sum_{k \geq 0} \xi^{k} (\gamma^{t+1+k} \| \bar{X}_{t+1+k} \|^2) \big].
$$
Summing over all times $t \geq 0$, and interchanging the two infinite summations involved, we get:  
$$ \displaystyle 
\EE \big[ \sum \nolimits_{t \geq 0} \gamma^t \| z_t \|^2 \big] \leq \frac{ 4\gamma \| \tilde Q \|^2}{1 - \eta} \frac{1}{1 - \xi} \EE \big[ \sum \nolimits_{t \geq 1} \gamma^t \| X_t \|^2 \big] < \infty.
$$
We then conclude that the process $\bp$ satisfies the transversality condition~\eqref{eq:transversality} and is an adjoint process corresponding to $(\bu, \bX)$.
\end{proof}

We can now use the adjoint process to compute the Gateaux derivative of $J$, and derive a necessary condition for optimality in the spirit of the classical Pontryagin's maximum principle.

\begin{proposition}
\label{prop:necessary_condition_of_PMP}
    Under Assumption~\ref{as:finit_cost}, the Gateaux derivative of $J$ at $\bu$ in the direction $\bbeta \in \cU_{ad}$ exists and is given by
    \begin{equation}
    \label{eq:Gateau_derivative_DJ_u_beta}
        DJ(\bu)(\bbeta) = \EE \Bigg[ \sum_{t \geq 0} \gamma^t \big(  p_t^\top \rB + 2 u_t^\top R + \bar{p}_t^{\top} \bar \rB + 2 \bar{u}_t^\top \bar R \big) \beta_t \Bigg],
    \end{equation}
    where $(p_t)_{t \geq 0}$ is the adjoint process of the controlled state process $\bX^\bu = (X_t^{\bu})_{t \geq 0}$.
    Moreover, if the control process $\bu$ is optimal, namely $J(\bu) \leq J(\bu')$ for any $\bu' \in \cU_{ad}$, then we have for all $t \geq 0$,
\begin{equation}
    \label{eq:necessary_PMP}
        \rB^\top p_t + 2 R u_t + \bar \rB^\top \bar{p}_t + 2 \bar R \bar{u}_t = 0, \qquad \PP-a.s. \, .
\end{equation}
\end{proposition}

\begin{proof}
    For an admissible control process $\bu \in \cU_{ad}$ and a perturbation $\bbeta \in \cU_{ad}$, the control $\bu + \lambda \bbeta = (u_t + \lambda \beta_t)_{t \geq 0}$ is still admissible for any $\lambda \in (0,1)$ by the convexity of $\cU_{ad}$. We denote by $\bX^{\bu + \lambda \bbeta} = (X_t^{\bu + \lambda \bbeta})_{t \geq 0}$ the corresponding controlled state starting from the same initial state $X_0^{\bu} = X_0^{\bu + \lambda \bbeta}$, and subject to the same noise sequences $(\varepsilon^0_{t+1})_{t \geq 0}$ and $(\varepsilon_{t+1})_{t \geq 0}$, and by $\bM^{\lambda, \bbeta} = (M_t^{\lambda, \bbeta})_{t \geq 0}$ the variation process defined as
    \begin{equation}
    \label{eq:Mt_lambda}
        M_t^{\lambda, \bbeta} := \frac{1}{\lambda} (X_t^{\bu + \lambda \bbeta} -  X_t^{\bu}). 
    \end{equation}
    The dynamics of $\bM^{\lambda, \bbeta}$ are given by:
    \begin{equation}
    \label{eq:dyn_Mt_lambda_beta}
        M_{t+1}^{\lambda, \bbeta} = \rA M_t^{\lambda, \bbeta} + \bar \rA \bar M_t^{\lambda, \bbeta} + \rB \beta_t + \bar \rB \bar \beta_t,
    \end{equation}
    where $M_0^{\lambda, \bbeta} =0$, and $\bar M_t^{\lambda, \bbeta} = \EE[ M_{t}^{\lambda, \bbeta}  | \cF^{0}_t ] $ and $\bar{\beta}_t = \EE[ \beta_t | \cF^{0}_t ] $ are the conditional expectations w.r.t. $\cF^{0}_t$. Since $\bbeta$ is $(\cF_t )_{t \geq 0}$-adapted, equation~\eqref{eq:dyn_Mt_lambda_beta} implies directly that $M_{t+1}^{\lambda, \bbeta}$ is $\cF_t$-measurable for every $t \geq 0$. We also notice that the process $\bM^{\lambda, \bbeta}$ is independent of $\lambda$, so we shall denote it by $\bM^{\bbeta} = (M_t^{\bbeta})_{t \geq 0}$ from now on.
    
   Let $(p_t)_{t \geq 0}$ be  the adjoint process of the state process $\bX^{\bu}$. From the admissibility of $\bu$ and $\bu + \lambda \bbeta$, we infer that both processes $\bX^{\bu}$ and $\bX^{\bu + \lambda \bbeta}$ are $L^2$-discounted integrable, and together with the transversality condition of the adjoint process $(p_t)_{t \geq 0}$, we deduce that
    $$
        \sum_{t \geq 0} \EE\Big[ \gamma^t \big( X_t^{\bu + \lambda \bbeta} - X_t^{\bu} \big)  \cdot p_t \Big] < \infty,
    $$ 
    and
    $$
        \sum_{t \geq 0} \EE \Big[ \gamma^t \Big(  b(X_t^{\bu+ \lambda \bbeta}, \bar{X}_t^{\bu+ \lambda \bbeta}, u_t + \lambda \beta_t, \bar u_t + \lambda \bar{\beta}_t) - b(X_t, \bar X_t, u_t, \bar u_t) \Big) \cdot p_t \Big] < \infty.
    $$
    We then compute the Gateau derivative of $J$. Let 
    $
        \zeta^{\lambda}_t = (X_t^{\bu+ \lambda \bbeta}, \bar{X}_t^{\bu+ \lambda \bbeta}, u_t + \lambda \beta_t, \bar u_t + \lambda \bar{\beta}_t)
    $ 
    and 
    $
        \zeta_t = (X_t^{\bu}, \bar X_t^{\bu}, u_t, \bar u_t)
    $, together with~\eqref{eq:Hamilton_function_h}, \eqref{eq:delta_X_t_with_b}, and $X_0^{\bu + \lambda \bbeta} = X_0^{\bu}$, we have
    \begin{equation}
    \label{eq:difference_between_cost_function_in_hamilton_func}
    \begin{split}
        & \EE \Big[ \sum_{t \geq 0} \gamma^t \big( c(\zeta_t^\lambda) - c(\zeta_t) \big) \Big] \\
         &\hskip 15pt 
         =\sum_{t \geq 0} \EE \Big[ \gamma^t \big( h(\zeta_t^\lambda, p_t) - h(\zeta_t, p_t) \big) \Big]\\ 
         &\hskip 35pt
         - \sum_{t \geq 0} \EE \Big[ \gamma^t [ ( X_{t+1}^{\bu + \lambda \bbeta} - X_t^{\bu + \lambda \bbeta}  ) - ( X_{t+1}^{\bu} - X_t^{\bu} ) ] \cdot p_t
         + \delta \gamma^t ( X_t^{\bu + \lambda \bbeta} - X_{t}^{\bu} ) \cdot p_{t} \Big] \\
        &\hskip 15pt
        =\sum_{t \geq 0} \EE \Big[ \gamma^t \big( h(\zeta_t^\lambda, p_t) - h(\zeta_t, p_t) \big) \Big] + \sum_{t \geq 0} \EE \Big[ \gamma^t \big( X_{t+1}^{\bu + \lambda \bbeta} - X_{t+1}^{\bu}) \cdot (p_{t+1} - p_t)  \big) \Big].
\end{split}
\end{equation}
    Using the formulas for the partial derivatives of the Hamilton function $\nabla_\zeta h(\zeta, p)$, and the process $(M_t^{\bbeta})_{t \geq 0}$ with $M_t^{\beta} = \lim_{\lambda \to 0} (X_t^{\bu + \lambda \bbeta} - X_t^\bu) / \lambda$, we then have 
    \begin{align}
        & DJ(\bu)(\bbeta) = \lim_{\lambda \to 0} \frac{1}{\lambda}\big[ J(\bu + \lambda \bbeta) - J(\bu) \big] 
        \nonumber \\
        = & \sum_{t \geq 0} \gamma^t \EE \Bigl[ \Big( \partial_x h(\zeta_t, p_t) \cdot M_t^{\bbeta} + \partial_{\bar x} h(\zeta_t, p_t) \cdot \bar{M}_t^{\bbeta} + \partial_u h(\zeta_t, p_t) \cdot \beta_t + \partial_{\bar u}h(\zeta_t, p_t) \cdot \bar{\beta}_t  \Big)
        \nonumber \\
        & \hskip 40pt + \Big( M_{t+1}^{\bbeta} \cdot ( p_{t+1} - p_t) \Big) \Bigr]  
        \label{eq:Gateaux}
    \end{align}
Properties of the conditional expectations imply: 
\begin{equation*}
\begin{split}
    \EE\Bigl[ \partial_{\bar x} h(\zeta_t, p_t) \cdot \bar{M}_t^{\bbeta} \Bigr] 
    & = \EE \Bigl[ \big( \bar \rA^\top p_t - 2 Q (X_t^{\bu} - \bar X_t^{\bu}) + 2 \tilde{Q} \bar X_t^{\bu} \big) \cdot \EE_{\cF^{0}_t}[M_t^{\bbeta} ] \Bigr] \\
    & = \EE\Bigl[ \EE_{\cF^{0}_t}  \bigl[ \bar \rA^\top p_t - 2 Q (X_t^{\bu} - \bar X_t^{\bu}) + 2 \tilde{Q} \bar X_t^{\bu} \bigr] \cdot \EE_{\cF^{0}_t}[M_t^{\bbeta} ] \Bigr]
        \\
        & = \EE \Bigl[ \big( \bar \rA^\top \bar{p}_t + 2 \tilde{Q} \bar{X}_t^{\bu} \big) \cdot \EE_{\cF^{0}_t}[M_t^{\bbeta} ] \Bigr] 
        \\
        & = \EE \Bigl[  \big( \bar \rA^\top \bar{p}_t + 2 \tilde{Q} \bar{X}_t^{\bu} \big) \cdot M_t^{\bbeta}  \Bigr].
\end{split}
\end{equation*}
Since $M_{t+1}^{\bbeta}$ is $\cF_t$-measurable, the backward equation~\eqref{eq:simpler_adjoint} for $(p_t)_{t\geq 0}$ implies that

\begin{equation*}
\begin{split}
        & \sum_{t \geq 0} \gamma^t\EE \Bigl[ M_{t+1}^{\bbeta}  \cdot ( p_{t+1} - p_{t} ) \Bigr]
        \\
         & \hskip 25pt
        =\sum_{t \geq 0} \gamma^t\EE \Bigl[ (M_{t+1}^{\bbeta})^\top (-\gamma) \EE \Big[  (\check \rA - \delta I_d)^\top p_{t+1} + 2 Q X_{t+1}^{\bu} + \bar \rA^\top \bar{p}_{t+1} + 2 \bar Q \bar{X}_{t+1}^{\bu} \Big| \cF_t \Big] \Bigr]
        \\
        &\hskip 25pt
        =-\sum_{t \geq 0} \gamma^{t+1} \EE \Bigl[  \big(  (\check \rA - \delta I_d)^\top p_{t+1} + 2 Q X_{t+1}^{\bu} + \bar \rA^\top \bar{p}_{t+1} + 2 \bar Q \bar{X}_{t+1}^{\bu} \big) \cdot  M_{t+1}^{\bbeta} \Bigr]
        \\
        &\hskip 25pt
        =-\sum_{t \geq 0} \gamma^{t} \EE \Bigl[  \big(  (\check \rA - \delta I_d)^\top p_{t} + 2 Q X_{t}^{\bu} + \bar \rA^\top \bar{p}_{t} + 2 \bar Q \bar{X}_{t}^{\bu} \big) \cdot  M_{t}^{\bbeta} \Bigr]
\end{split}
\end{equation*}
by shifting the time index and using the fact that $M_0^{\bbeta} = \bar{M}_0^{\bbeta} = 0$. Next, using the formula for $\partial_x h$ we get
    \begin{align*}
        & \sum_{t \geq 0} \gamma^{t} \EE \Bigl[   \partial_x h(\zeta_t, p_t) \cdot M_t^{\bbeta} + \partial_{\bar x} h(\zeta_t, p_t) \cdot \bar{M}_t^{\bbeta}  +   M_{t+1}^{\bbeta}  \cdot ( p_{t+1} - p_{t} )\Bigr] 
        \\
         = & 
         \sum_{t \geq 0} \gamma^t \EE \Bigl[  \big( (\check \rA - \delta I_d )^\top p_{t} + 2 Q (X_t^\bu - \bar X_t^\bu \big) \big) \cdot M^{\bbeta}_t + \big(\bar \rA^\top \bar p_{t} + 2 (Q + \bar Q) \bar X_t^\bu \big) \cdot M_t^{\bbeta}  
        \\
        & \hskip 60pt - \gamma \big(  (\check \rA - \delta I_d)^\top p_{t} + 2 Q X_{t}^{\bu} + \bar \rA^\top \bar{p}_{t} + 2 \bar Q \bar{X}_{t}^{\bu} \big) \cdot  M_{t}^{\bbeta} \Bigr]
        \\
     = & 0
    \end{align*}
As a result, \eqref{eq:Gateaux} reduces to
    \begin{equation*}
    \begin{split}
        DJ(\bu)(\bbeta) & = \sum_{t \geq 0} \gamma^t \EE \Bigl[ \partial_u h(\zeta_t, p_t) \cdot \beta_t + \partial_{\bar u}h(\zeta_t, p_t) \cdot \bar{\beta}_t  \Bigr] 
        \\
        & = \sum_{t \geq 0} \gamma^t \EE \Bigl[ \big(p_t^\top \rB + 2 u_t^\top R + \bar{p}_t^\top \bar \rB + 2 \bar u_t^\top \bar R \big) \beta_t \Bigr].
    \end{split}
    \end{equation*}
We conclude that if $\bu$ is optimal, we must have $DJ(\bu)({\bbeta})\geq 0$ for every $\bbeta \in \cU_{ad}$, 
implying that, for every $t \geq 0$ 
\begin{equation}
\label{eq:necessary}
         \rB^\top p_t + 2 R u_t + \bar \rB^\top \bar{p}_t + 2 \bar R \bar{u}_t = 0, \qquad \PP-a.s. \, .
\end{equation}
\end{proof}

\subsection{Identification of the optimal control}

Let us now assume that $\bu$ is an admissible optimal control and show that $u_t$ is linear in $(p_t, \bar p_t)$ where $p_t$ is the adjoint process of $\bX^{\bu}$. By taking conditional expectation $\EE[\cdot|\cF^{0}_t]$ in equation~\eqref{eq:necessary}, we get
$$
    (\rB + \bar \rB)^\top \bar p_t + 2 (R + \bar R) \bar u_t = 0
$$
from which we derive:
$$
\bar u_t=-\frac12 (R +\bar R)^{-1}(\rB+\bar \rB)^\top\bar p_t,
$$
and going back to \eqref{eq:necessary} we finally get:
\begin{equation}
    u_t = -\frac{1}{2} R^{-1} \rB^\top p_t  - \frac{1}{2} R^{-1} \big( \bar \rB - \bar R (R + \bar R)^{-1} (\rB + \bar \rB)^\top \big)  \bar p_t.
\label{eq:control_ut_expression_1}
\end{equation}
Using the notations
\eqref{eq:Gamma_Lambda} we get the formula~\eqref{eq:control_ut_expression_2} for the process $\bu = (u_t)_{t \geq 0}$.
Given the well-known fact that at the optimum of standard LQ models, the adjoint processes are affine functions of the state, one may wonder if optimal controls can also be shown to be linear in the controlled state process $\bX^\bu$ and the conditional mean process $\bar{\bX}^\bu$.  
In order to prove that it is indeed the case, we introduce the following system of forward-backward stochastic equations for two $L^2-$discounted integrable processes $(\bX, \bp) \in \cX \times \cX$:
\begin{equation}
    \label{eq:FBDSE_state_adjoint}
    \left\{
    \begin{array}{rcl}
        X_{t+1} - X_t &=& \check \rA X_t + \bar \rA \bar{X}_t + \rB\Gamma (p_t - \bar{p}_t) + (\rB + \bar \rB) \Lambda \bar{p}_t + \varepsilon_{t+1}^0 + \varepsilon_{t+1}
        \\
        p_{t+1} - p_{t} &=& - \gamma \big[ (\check \rA - \delta I_d)^\top p_{t+1} + 2Q X_{t+1} + \bar \rA^\top \bar{p}_{t+1} + 2 \bar Q \bar{X}_{t+1} \big] 
        \\
            && \hskip 180pt + Z_{t+1}^0 \varepsilon_{t+1}^0  + Z_{t+1} \varepsilon_{t+1},
    \end{array}
    \right.
\end{equation}
where the adapted processes $(Z_t)_{t \geq 0}$ and $(Z_t^0)_{t \geq 0}$ are part of the solution, and where $X_t$ satisfies the initial condition at time $t=0$, and $p_t$ satisfies the transversality condition which we view as a terminal condition at $\infty$. This justifies our terminology of \emph{forward-backward} system. Notice that the first equation in \eqref{eq:FBDSE_state_adjoint} 
is nothing but the state equation for a linear control of the form we want to restrict ourselves to. %

\noindent
The following proposition provides a solution to the forward-backward system~\eqref{eq:FBDSE_state_adjoint}.

\begin{proposition}
\label{prop:solution_to_FBSDE_with_optimal_control}
    Under Assumption~\ref{as:finit_cost} and~\ref{as:positivity-qr}, we consider a control process $\bu = (u_t)_{t \geq 0}$ in a feedback form of its controlled state process $\bX^\bu$ and its conditional mean process $\bar \bX^\bu$ with control parameter $\theta = (- \Gamma P, -\Lambda \bar P)$ and policy $\pi_\theta : (x, \bar x) \mapsto  \Gamma P (x - \bar x) + \Lambda \bar P \bar x$  such that
    \begin{equation}
        \label{eq:closed-loop-control-on-x}
        u_t = \pi_\theta(X_t^\bu, \bar X_t^\bu) =  \Gamma P (X_t^\bu - \bar{X}_t^\bu) + \Lambda \bar P \bar{X}_t^\bu,
    \end{equation}
    for every $t \geq 0$, where the coefficients $(\Gamma, \Lambda, P, \bar P)$ are defined in Theorem \ref{thm:existence_linear_control}. Let $\bp^\bu = (p_t^\bu)_{t \geq 0}$ be a process defined by:
    \begin{equation}
        \label{eq:relation_ajoint_and_state}
        p^{\bu}_t = P(X_t^\bu - \bar{X}^{\bu}_t) + \bar P \bar{X}^{\bu}_t,
    \end{equation}
    for every $t \geq 0$. Then, the pair $(\bX^\bu, \bp^\bu)$ is a solution of the forward-backward system~\eqref{eq:FBDSE_state_adjoint} with
    \begin{equation}
        \label{eq:def_Zt}
        Z_t = \gamma \big( \rA^\top P + 2 Q \big), \qquad Z_t^0 = \gamma \big[ (\rA + \bar{\rA})^\top \bar P + 2 (Q + \bar Q) \big].
    \end{equation}
    Moreover, the process $\bp^\bu$ is the unique adjoint process corresponding to $(\bu, \bX^\bu)$.
\end{proposition}

\begin{proof} 
    Proposition~\ref{proposition:admissible_of_X} and Proposition~\ref{proposition:admissibility_of_u_with_feedback_control} state that the control parameter $\theta = (-\Gamma P, -\Lambda \bar P)$ is in $\Theta$, and that $\bu \in \cU_{ad}$ and $\bX^\bu \in \cX$. By definition, $\bp^\bu$ is adapted to $(\cF_t)_{t \geq 0}$ and it is $L^2-$discounted integrable, so $\bp^\bu \in \cX$. 
    
    For any time $t \geq 0$,  we can take the conditional expectation w.r.t. $\cF^{0}_t$ and obtain that
    $u_t = \Gamma \big( p_t^\bu - \bar{p}_t^\bu \big) + \Lambda \bar{p}_t^\bu$ and $\bar{u}_t = \Lambda \bar{p}_t^\bu$. 
    From the state dynamics~\eqref{fo:MKV-state} of $\bX^{\bu}$, we deduce the forward equation in \eqref{eq:FBDSE_state_adjoint}:
    \begin{equation*}
    \begin{split}
        X_{t+1}^\bu - X_t^\bu &= \check \rA X_t^\bu + \bar \rA  \bar{X}^{\bu}_t + \rB \big( \Gamma P (X_t^\bu - \bar{X}^\bu_t) + \Lambda \bar P \bar{X}^{\bu}_t \big) + \bar \rB \Lambda \bar P \bar{X}^{\bu}_t  + \varepsilon_{t+1}^0 + \varepsilon_{t+1} 
        \\
        & = \check \rA X_t^\bu + \bar \rA  \bar{X}^{\bu}_t + \rB \Gamma ( p^\bu_t - \bar{p}^{\bu}_t)  +  (\rB + \bar \rB) \Lambda \bar p^{\bu}_t  + \varepsilon_{t+1}^0 + \varepsilon_{t+1} .
    \end{split}
    \end{equation*}
    To show that the process $\bp^\bu$ satisfies the backward equation  in \eqref{eq:FBDSE_state_adjoint}, we use the following transformation trick with the help of $\bp^\bu - \bar{\bp}^\bu = (p_t^\bu - \bar{p}_t^\bu)_{t \geq 0}$ and $\bar{\bp}^\bu = (\bar p_t^\bu)_{t \geq 0}$:
    $$
        p^\bu_{t+1} - p^\bu_t = \big( (p^\bu_{t+1} - \bar p^\bu_{t+1}) - (p^\bu_t - \bar{p}^\bu_t) \big) + ( \bar p^\bu_{t+1} - \bar{p}^\bu_t).
    $$
    By taking conditional expectations of the state dynamics~\eqref{fo:MKV-state} of $\bX^{\bu}$, we obtain a dynamical equation for $\bar{X}^\bu$: for every $t \geq 0$,
    $$
        \bar{X}^{\bu}_{t+1} = [ ( \rA + \bar \rA) + (\rB + \bar \rB) \Lambda \bar P ] \bar{X}^\bu_t + \varepsilon_{t+1}^0.
    $$
    Then, we multiply on both sides with $Z_{t+1}^0 = \gamma \big[  (\rA + \bar{\rA} )^\top \bar P + 2 (Q + \bar Q)\big]$ and apply the Riccati equation of $\bar P$ on the right-hand side of the resulting equality, we obtain 
    $$
        \gamma \big[  (\rA + \bar{\rA} )^\top \bar P + 2 (Q + \bar Q) \big] \bar{X}^{\bu}_{t+1} = \bar P \bar{X}^{\bu}_{t} + Z_{t+1}^0 \varepsilon_{t+1}^0.
    $$
    By re-arranging the terms in the above equation and using $\bar P \bar X_{t+1}^\bu = \bar p_{t+1}^\bu$, we deduce
    \begin{equation}
        \label{eq:inside_proof_diff_bar_p}
        \bar{p}^{\bu}_{t+1} - \bar{p}^{\bu}_t = - \gamma \Big[ ( \check \rA + \bar \rA - \delta I_d)^\top \bar{p}^{\bu}_{t+1} + 2 (Q + \bar Q) \bar{X}^\bu_{t+1} \Big] + Z_{t+1}^0 \varepsilon_{t+1}^0.
    \end{equation}
    Similarly, from the dynamics of $\bX^\bu$ and $\bar{\bX}^\bu$, we have 
    $$
        X_{t+1}^\bu - \bar{X}^{\bu}_{t+1} = (\rA + \rB \Gamma P) \big(X_t^\bu - \bar{X}^\bu_t \big) + \varepsilon_{t+1}.
    $$
    Then we multiply on both sides with $Z_{t+1}$ and apply the Riccati equation of $P$ on the right-hand side of the resulting equality, we deduce that
    $$
        \gamma (\rA^\top P + 2Q) (X_{t+1}^\bu - \bar{X}^{\bu}_{t+1}) = P \big(X_t^\bu - \bar{X}^\bu_t \big) + Z_{t+1} \varepsilon_{t+1},
    $$
    We re-arrange the terms and apply the condition $p^\bu_t - \bar p^\bu_t = P (X_t^\bu - \bar{X}^\bu_t)$, we have
    \begin{equation}
        \label{eq:inside_proof_diff_p_minus_bar_p}
        \begin{split}
         & (p^\bu_{t+1} - \bar p^\bu_{t+1}) - (p^\bu_t - \bar{p}^\bu_t) 
         \\
         & = - \gamma \big[ ( \check \rA - \delta I_d )^\top (p^\bu_{t+1} - \bar p^\bu_{t+1}) + 2 Q (X_{t+1}^\bu - \bar{X}^{\bu}_{t+1}) \big] + Z_{t+1}^1 \varepsilon_{t+1}.
         \end{split}
    \end{equation}
    Adding both sides of equations~\eqref{eq:inside_proof_diff_bar_p} and~\eqref{eq:inside_proof_diff_p_minus_bar_p}, we obtain immediately that $\bp^\bu$ satisfies the backward equation in~\eqref{eq:FBDSE_state_adjoint}. Thus, the pair of state-adjoint processes $(\bX^\bu, \bp^\bu) \in \cX \times \cX$ is a solution of the forward-backward system \eqref{eq:FBDSE_state_adjoint}.

    To show that the process $\bp^\bu$ is an adjoint process, we take the conditional expectations $\EE[\cdot | \cF_t]$ for every $t \geq 0$ on both sides of the backward equation of $\bp^\bu$ in \eqref{eq:FBDSE_state_adjoint}. Because the noise terms $\varepsilon_{t+1}$ and $\varepsilon_{t+1}^0$ are independent of $\cF_{t}$ and they are of means $0$, and because the terms $Z_t$ and $Z_t^0$ are deterministic constants, the process $\bp^\bu$ satisfies the backward equation~\eqref{eq:adjoint}, and $\bp^\bu \in \cX$ implies the transversality condition~\eqref{eq:transversality}, so that it is an adjoint process. Under Assumption~\ref{as:finit_cost}, we conclude that $\bp^\bu$ is the unique adjoint process corresponding to $(\bu, \bX^\bu)$.
\end{proof}

The following theorem shows that we can construct an optimal control process $\bu$ for the MFC problem~\eqref{pb:MFC_L2_admissible_control} based on a solution of the forward-backward system~\eqref{eq:FBDSE_state_adjoint}.

\begin{proposition}
\label{prop:existence_of_optimal_control_from_FBDTSE}
      We assume that Assumption~\ref{as:finit_cost} and~\ref{as:positivity-qr} hold. If $(\bX, \bp) \in \cX \times \cX$ is a solution to the forward-backward system~\eqref{eq:FBDSE_state_adjoint} with $X_0 = \varepsilon_0^0 + \varepsilon_0$, then, the control process $\bu = (u_t)_{t \geq 0}$ given by 
    \begin{equation}
    \label{eq:linearity_control_in_adjoint_process}
         u_t = \Gamma (p_t - \bar p_t) + \Lambda \bar{p}_t
    \end{equation}
    is an optimal admissible control process for the MFC problem~\eqref{pb:MFC_L2_admissible_control}. 
    Besides, the process $\bX$ is the controlled state process of $\bu$, and the process $\bp$ is the unique adjoint process corresponding to $(\bu, \bX)$.
\end{proposition}

\begin{proof}
By definition \eqref{eq:linearity_control_in_adjoint_process}, we have that $u_t = \Gamma p_t + (\Lambda - \Gamma) \bar p_t$ and $\bar u_t = \Lambda \bar p_t$ for every $t \geq 0$.
Since the process $\bp \in \cX$ is adapted to $(\cF_t)_{t \geq 0}$, $\bu$ is also adapted to $(\cF_t)_{t \geq 0}$.
Jensen's inequality for conditional expectations and Cauchy-Schwarz inequality imply that $$
\sum_{t \geq 0} \gamma^t \EE[ \| u_t \|^2 ]  \leq 2( \| \Gamma \|^2 + \| \Lambda - \Gamma \|^2) \sum_{t \geq 0} \gamma^t \EE[ \| p_t \|^2 ] <  \infty,
$$
so that $\bu \in \cU_{ad}$.
Using $\Gamma (p_t - \bar p_t) = u_t - \bar u_t$ and $\Lambda \bar p_t = \bar u_t$ in the forward equation of the forward-backward system~\eqref{eq:FBDSE_state_adjoint}, we see that the process $\bX$ satisfies the evolution dynamics~\eqref{fo:MKV-state}: 
$$
    X_{t+1} = \rA X_t + \bar \rA \bar{X}_t + \rB u_t + \bar \rB \bar u_t + \varepsilon_{t+1}^0 + \varepsilon_{t+1},
$$
for every $t \geq 0$ and starts with initial state $X_0 = \varepsilon_0^0 + \varepsilon_0$. Corollary~\ref{corollary:uniqueness_of_controlled_state_process} states that the process $\bX$ coincides with the controlled state process $\bX^\bu$ of $\bu$ in the sense that $X_t = X_t^\bu$ almost surely for every $t \geq 0$.
Meanwhile, by taking conditional expectations with respect to $\cF_{t}$ for each time $t \geq 0$ in the backward equation for $\bp$ in~\eqref{eq:FBDSE_state_adjoint}, and by applying similar arguments to those used in the last paragraph of the proof of Proposition~\ref{prop:solution_to_FBSDE_with_optimal_control}, we see that the solution process $\bp \in \cX$ is the unique adjoint process corresponding to $(\bu, \bX)$. 

With some matrix algebra on $\Lambda - \Gamma$, it is plain to see that the control process $\bu$ satisfies equation~\eqref{eq:control_ut_expression_1}, that is
$
 u_t = -\frac{1}{2} R^{-1} \rB^\top p_t  - \frac{1}{2} R^{-1} \big( \bar \rB - \bar R \tilde{R}^{-1} \tilde B^\top \big) \bar p_t
$
for all $t \geq 0$. Rearranging the terms we show that $\bu$ satisfies the necessary condition~\eqref{eq:necessary_PMP} of the Pontryagin's maximum principle with adjoint process $\bp$, that is, $\rB^\top p_t + \bar \rB^\top \bar p_t + 2 R u_t + 2 \bar R \bar u_t = 0$ for every $t \geq 0$.

To show that the process $\bu$ is an optimal control process for the MFC problem~\eqref{pb:MFC_L2_admissible_control} defined for $(\bX, \bu)$, we consider the Hamilton function $h(\zeta, p) = b(\zeta) \cdot p + c(\zeta) - \delta x \cdot p$ introduced in~\eqref{eq:Hamilton_function_h} and a perturbation direction $\bbeta \in \cU_{ad}$ for the control process. The difference between cost functions $J(\bu + \lambda \bbeta) - J(\bu)$ can then be expressed using equation~\eqref{eq:difference_between_cost_function_in_hamilton_func} like
\begin{equation}
\label{eq:difference_between_J_u_beta}
\begin{split}
     & J(\bu + \lambda \beta) - J(\bu)
    \\
    = & \sum_{t \geq 0} \EE \Big[ \gamma^t \big( h(\zeta_t^\lambda, p_t) - h(\zeta_t, p_t) \big) \Big] +  \sum_{t \geq 0} \EE \Big[ \gamma^t \big( X_{t+1}^{\bu + \lambda \bbeta} - X_{t+1}^{\bu}) \cdot (p_{t+1} - p_t)  \big) \Big]
    \\
    = &  \sum_{t \geq 0} \gamma^t \EE \big[ \nabla_{\zeta} h(\zeta_t, p_t) \cdot ( \zeta_t^\lambda - \zeta_t) \big] + \sum_{t \geq 0} \gamma^t \EE \big[ \frac{1}{2} \nabla^2_{\zeta \zeta} h(\eta_t^{\lambda}, p_t) (\zeta_t^\lambda - \zeta_t) \cdot ( \zeta_t^\lambda - \zeta_t)  \big]  \\
     & \hskip 15pt + \lambda \sum_{t \geq 0} \gamma^t \EE \big[ M_t^{\lambda, \bbeta} \cdot ( p_{t+1} - p_t ) \big] 
     \\
    = & (i) + (ii) + (iii),
\end{split}
\end{equation}
where $\zeta^{\lambda}_t = (X_t^{\bu+ \lambda \bbeta}, \bar{X}_t^{\bu+ \lambda \bbeta}, u_t + \lambda \beta_t, \bar u_t + \lambda \bar{\beta}_t)$, $\zeta_t = (X_t^{\bu}, \bar X_t^{\bu}, u_t, \bar u_t)$, and $M_t^{\lambda, \bbeta} = ( X_t^{\bu + \lambda \bbeta} - X_t^{\bu}) / \lambda$. In the second equality in the above equation we use the Taylor expansion of the Hamilton function where the terms $\nabla_\zeta h$ and $\nabla_{\zeta, \zeta}^2 h$ are the Jacobian and the Hessian matrices w.r.t. to $\zeta$, and $\eta_t^\lambda = \zeta_t + \rho_t (\zeta_t^\lambda - \zeta_t)$ for some $\rho_t \in [0, 1]$.

From the proof of Proposition~\ref{prop:necessary_condition_of_PMP} (check the computation for the Gateau derivative of $DJ(\bu)(\bbeta)$ in equation~\eqref{eq:Gateaux}), and the fact that the control process $\bu$ satisfies equation~\eqref{eq:necessary_PMP}, we have
$$
    (i) + (iii) = \lambda \sum_{t \geq 0} \gamma^t \EE \big[ (\rB^\top p_t + \bar \rB^\top \bar p_t + 2 R u_t + 2 \bar R \bar u_t) \cdot \beta_t \big] = 0.
$$
The Hessian $\nabla_{\zeta, \zeta}^2 h$ of the Hamilton function takes the form 
$$
    \nabla_{\zeta, \zeta}^2 h(\zeta, p) = \begin{pmatrix}
	2Q & -2Q & 0 & 0 \\
	-2Q & 2 Q + 2 \tilde Q & 0 & 0 \\
        0 & 0 & 2R & -2R \\
        0 & 0 & -2R & 2 R + 2 \tilde R
	\end{pmatrix}.
$$
We also notice that $\zeta^{\lambda}_t - \zeta_t = \lambda ( M_t^{\lambda, \bbeta}, \bar{M}_t^{\lambda, \bbeta}, \beta_t, \bar{\beta}_t )$. Thus,
\begin{equation*}
\begin{split}
    (ii) = & \lambda^2 \sum_{t \geq 0} \gamma^t \EE \Big[ ( M_t^{\lambda, \bbeta} - \bar{M}_t^{\lambda, \bbeta} )^\top Q  ( M_t^{\lambda, \bbeta} - \bar{M}_t^{\lambda, \bbeta} ) + (\bar{M}_t^{\lambda, \bbeta})^\top \tilde{Q} \bar{M}_t^{\lambda, \bbeta} 
    \\
    & \hskip 50pt + ( \beta_t - \bar{\beta}_t )^\top  R (\beta_t - \bar{\beta}_t ) + \bar{\beta}_t^\top \tilde{R} \bar{\beta}_t \Big].
\end{split}
\end{equation*}
Because $Q, \tilde Q \succeq 0$ and $R, \tilde R \succ 0$, we obtain that for any perturbation process $\bbeta \in \cU_{ad}$,
$$
    J( \bu + \lambda \bbeta) - J(\bu) \geq 0.
$$
Therefore, the process $\bu$ given by equation~\eqref{eq:linearity_control_in_adjoint_process} is an optimal control process to the MFC problem~\eqref{pb:MFC_L2_admissible_control}.
\end{proof}

\subsection{Proof of Theorem~\ref{thm:existence_linear_control}}

We now provide the proof of the main result of this section, namely
Theorem \ref{thm:existence_linear_control}.

\begin{proof}
\textbf{(Existence.)}
Let $\bu = (u_t)_{t \geq 0}$ be the control process specified by equation~\eqref{eq:optimal_control_linear_in_x_barx} with parameter $\theta = (- \Gamma P, - \Lambda \bar P) = (K^*, L^*)$ where $(\Gamma, \Lambda, P, \bar P, K^*, L^*)$ are defined in the statement of the theorem, and let $\bX$ be the controlled state process in $\cX$ of $\bu$. Proposition~\ref{proposition:admissibility_of_u_with_feedback_control} shows that under Assumption~\ref{as:finit_cost} and~\ref{as:positivity-qr}, $\theta$ is an admissible control parameter and $\bu \in \cU_{ad}$.

We also consider a process $\bp$ defined by equation~\eqref{eq:relation_ajoint_and_state}, that is, $p_t = P ( X_t - \bar X_t) + \bar P \bar X_t$ for every $t \geq 0$. Then, under Assumption~\ref{as:finit_cost} and~\ref{as:positivity-qr}, Proposition~\ref{prop:solution_to_FBSDE_with_optimal_control} implies that the pair of processes $(\bX, \bp)$ is a solution to the forward-backward system~\eqref{eq:FBDSE_state_adjoint}, and the process $\bp$ is the unique adjoint process associated with $(\bu, \bX)$.

Thus, from the solution process $\bp$ and the coefficients $(\Gamma, \Lambda)$, we construct another control process $\bu^* = (u_t^*)_{t \geq 0}$ defined by~\eqref{eq:linearity_control_in_adjoint_process}:
$$
    u_t^* = \Gamma ( p_t - \bar p_t) + \Lambda \bar p_t,
$$ 
for every $t \geq 0$. Under Assumption~\ref{as:finit_cost} and~\ref{as:positivity-qr}, Proposition~\ref{prop:existence_of_optimal_control_from_FBDTSE} implies that the process $\bu^*$ is admissible, and it is an optimal admissible control process for the MFC problem~\eqref{pb:MFC_L2_admissible_control} with $(\bX, \bu^*)$. Moreover, the process $\bX$ is also the controlled state process of $\bu^*$, and the process $\bp$ is the unique adjoint process corresponding to $(\bu^*, \bX)$.  From the definition of $\bp$ in~\eqref{eq:relation_ajoint_and_state}, we have $p_t - \bar p_t = P (X_t - \bar X_t)$ and $\bar p_t = \bar P \bar X_t$, so we obtain
$$
      u_t^* = \Gamma P ( X_t - \bar X_t) + \Lambda \bar P \bar X_t = u_t,
$$
for every $t \geq 0$. Hence, the control process $\bu$ considered here is an optimal control for the MFC problem~\eqref{pb:MFC_L2_admissible_control}.

\noindent \textbf{(Uniqueness.)} Suppose there exists another optimal control $\bu' \in \cU_{ad}$ for the MFC problem~\eqref{pb:MFC_L2_admissible_control}. We set $\beta_t = u'_t - u_t$ the difference between the two control processes at time $t$. Because both $\bu$ and $\bu'$ are $L^2-$discounted integrable, we have $\bbeta = (\beta_t)_{t \geq 0}$ is also $L^2-$discounted integrable by Assumption~\ref{as:finit_cost}, thus $\bbeta \in \cU_{ad}$. Because the process $\bp$ is the adjoint process corresponding to $(\bu, \bX)$, from equation~\eqref{eq:difference_between_J_u_beta}, we have 
\begin{equation*}
\begin{split}
    0 &= J(\bu') - J(\bu) 
    \\ 
    & = \sum_{t \geq 0} \gamma^t \EE \big[ (\rB^\top p_t + \bar{\rB}^\top \bar{p}_t + 2 R u_t + 2\bar{R} \bar{u}_t ) \cdot \beta_t \big] + \sum_{t \geq 0} \gamma^t \EE \big[ c( M_t, \bar M_t, \beta_t, \bar{\beta}_t ) \big]
    \\
    & = (i) + (ii),
\end{split}
\end{equation*}
where $M_t = X^{\bu'}_t - X^{\bu}_t$ follows the dynamics $M_{t+1} = \rA M_t + \bar \rA \bar M_t + \rB \beta_t + \bar \rB \bar{\beta}_t$ for every $t \geq 0$ with $M_0 = 0$. From the optimality of the control process $\bu$, the necessary condition of the Pontryagin's maximum principle in Proposition~\ref{prop:necessary_condition_of_PMP} yields that $(i) = 0$. Besides, Assumption~\ref{as:positivity-qr} implies $\EE \big[ c( M_t, \bar M_t, \beta_t, \bar{\beta}_t ) \big] \geq 0$ for every $t \geq 0$. Hence, 
\begin{equation*}
\begin{split}
    0 & = \EE\big[ c(M_t, \bar M_t, \beta_t, \bar{\beta}_t ) \big] = \EE[ (M_t - \bar M_t)^\top Q (M_t - \bar M_t) ] + \EE[ \bar{M}_t^\top \tilde{Q} \bar M_t] 
    \\
    & \hskip 120pt + \EE[ ( \beta_t - \bar{\beta}_t)^\top R (\beta_t - \bar{\beta}_t) ] + \EE[ (\bar{\beta}_t)^\top \tilde{R} \bar{\beta}_t ].
\end{split}
\end{equation*}
Because $Q, \tilde Q \succeq 0$ and $R, \tilde R \succ 0$, so we must have $\beta_t = \bar{\beta}_t = 0$ for every $t \geq 0$. This implies the uniqueness of the optimal control process.
\end{proof}

\section{\textbf{Policy gradient algorithms}}
\label{ref:PGconv}

In this section, we search for the optimal control process using PG methods. 
Theorem~\ref{thm:existence_linear_control} states that the optimal control process can be parameterized with an element $\theta^* = (K^*, L^*) \in \Theta$. Inspired by formula~\eqref{eq:optimal_control_linear_in_x_barx}, we limit our search of the optimal control process to the set of controls $\cU_{ad}^\Theta$ that can be parameterized by an element in $\Theta$. We will present a model-based and two model-free PG algorithms to approximate the best admissible control parameter $\theta^* \in \Theta$.

\subsection{Reformulation of the MFC problem}
We reformulate the MFC problem~\eqref{pb:MFC_L2_admissible_control} as an optimization problem over the set of admissible control parameters $\Theta$. We denote the discounted MF cost~\eqref{fo:MKV-discounted_cost} for a parametrized admissible control process $\bu^\theta \in \cU_{ad}^\Theta$ by
\begin{equation}
\label{eq:def-cost-C-theta}
    C(\theta) = J(\bu^\theta).
\end{equation}
Lemma~\ref{lemma:admissible_of_control_with_feedback_form_on_theta} states that for a control parameter $\theta \in \Theta$, the corresponding parameterized control process $\bu^\theta \in \cU_{ad}$, thus $\bu^\theta \in \cU_{ad}^\theta$ and the controlled state process $\bX^{\theta} \in \cX$. Proposition~\ref{proposition:reformulation_of_MF_problem} below justifies the reformulation of the search for the optimal control with parameters in $\Theta$.

\begin{proposition}
\label{proposition:reformulation_of_MF_problem}
Under Assumptions~\ref{as:finit_cost} and~\ref{as:positivity-qr}, we have
\begin{equation}
\label{eq:equivalence_C_and_J_in_MFC}
    \inf_{\theta \in \Theta} C(\theta) = \inf_{\bu \in \cU_{ad}} J(\bu).
\end{equation}
\end{proposition}

\begin{proof}
Let $\theta^* = (K^*, L^*)$ be a control parameter defined in Theorem~\ref{thm:existence_linear_control}. Proposition~\ref{proposition:admissibility_of_u_with_feedback_control} states that $\theta^* \in \Theta$. Theorem~\ref{thm:existence_linear_control} states that the control process $\bu^{\theta^*}$ minimizes the discounted MF cost $J(\bu)$ among all admissible controls in $\cU_{ad}$. Thus,
$$
      J(\bu^{\theta^*}) = C(\theta^*) \geq \inf_{\theta \in \Theta} C(\theta)
     = \inf_{\bu^\theta \in \cU_{ad}^\Theta} J(\bu^\theta) \geq \inf_{\bu \in \cU_{ad}} J(\bu) = J(\bu^{\theta^*}).
$$
This concludes the proof.
\end{proof}

From now on, we focus on the search for a minimizer of $C$ in the set $\Theta$. 
We will sometimes abuse terminology and call $\theta \in \Theta$ a ``control'' or say that the state process is ``controlled'' by $\theta$.

In the following, we cast the LQMFC for $C(\theta)$ into two standard LQ problems depending only on the parameter $K$ and $L$ respectively. For a state process $(X_t^{\theta})_{t\geq 0}$ controlled by $\theta = (K, L) \in \Theta$ following dynamics~\eqref{fo:MKV-state} and starting from initial position $X_0^\theta = \varepsilon_0^0 + \varepsilon_0$, we introduce two auxiliary processes, $\by = (y_t)_{t\geq 0}$ and $\bz = (z_t)_{t\geq 0}$ defined by: 
\begin{equation}
\label{eq:def-Y-Z-fct-X}
    y_{t} := X^{\theta}_t - \bar{X}^{\theta}_t,
    \qquad\text{and}\qquad
    z_{t} := \bar{X}^{\theta}_t,
\end{equation}
for $t \geq 0$. Their dynamics are given by:
\begin{align}
    y_{t+1} &= (\rA - \rB K) y_{t} + \varepsilon_{t+1}, && y_0 = \varepsilon_0 - \EE[\varepsilon_0];
    \label{eq:dyn_y_theta}
    \\
    z_{t+1} &= ( \tilde \rA - \tilde \rB  L) z_{t}  + \varepsilon^0_{t+1}, && z_0 = \varepsilon^0_0 + \EE[\varepsilon_0].
    \label{eq:dyn_z_theta}
\end{align}
The processes $\by$ and $\bz$ depend on the matrix parameters $K$ and $L$ respectively. Since the noise processes and their initial states involved in \eqref{eq:dyn_y_theta} and \eqref{eq:dyn_z_theta} are independent, the processes $\by$ and $\bz$ are independent. We introduce the corresponding cost functions $C_y, C_z: \RR^{\ell \times d} \to \RR$ as:
\begin{equation}
\label{eq:def_Cy_Cz}
	C_y(K) = \mathbb{E} \big[  \sum_{t \geq 0} \gamma^t  (y_t)^\top (Q + K^T R K) y_t \big] , \quad
	C_z(L) = \mathbb{E} \big[  \sum_{t \geq 0} \gamma^t   (z_t)^\top (\tilde{Q} + L^T \tilde{R} L) z_t  \big].
\end{equation}
We may sometimes use $\by^K = (y_t^K)_{t \geq 0}$ and $\bz^L = (z_t^L)_{t \geq 0}$ to emphasize the dependence on $K$ and $L$ of the processes $\by$ and $\bz$. 

Let $\Theta_K = \{ K \, | \, \gamma \| \rA - \rB K \|^2 < 1 \}$ and $\Theta_L = \{ L \, | \, \gamma \| \tilde \rA - \tilde \rB L \|^2 < 1 \}$. The following lemma shows the benefit of considering the re-parametrization \eqref{eq:def-Y-Z-fct-X} inspired by Theorem~\ref{thm:existence_linear_control}.

\begin{lemma}
\label{lem:decompose-cost-reparam}
For any $\theta = (K,L) \in \Theta$, we have $C(\theta) = C_y(K) + C_z(L),$ and 
 \begin{equation}
    \label{eq:reparametrization_of_optimality_LQMF}
     \inf_{ \theta = (K, L) \in \Theta } C(\theta)=\inf_{ K \in \Theta_K} C_y(K) +\inf_{L \in \Theta_L} C_z(L).
 \end{equation}
\end{lemma}

\begin{proof}
Any parameterized control process $\bu^\theta = (u_t^\theta)_{t\geq 0}$ with $\theta = (K, L)$, can be rewritten as $u_t^\theta = - K y_t - L z_t$ for every $t \geq 0$. Thus, the instantaneous cost at any time $t \geq 0$ satisfies $\EE \big[ c(X_t^\theta, \bar X_t^\theta, u_t^\theta, \bar u_t^\theta) \big] = \EE \big[ (y_t)^\top (Q + K^\top R K) y_t + (z_t)^\top (\tilde Q + L^\top \tilde R L) z_t \big]$, so that $C(\theta) = C_y(K) + C_z(L)$.  Because $\Theta = \{ \theta = (K, L) \, | \, K \in \Theta_K, L \in \Theta_K \} = \Theta_K \times \Theta_L$, together with the fact that $\by$ and $\bz$ are independent processes controlled by $K$ and $L$ respectively, we can then separate the optimization over $\Theta$ into two optimization problems over $\Theta_K$ and $\Theta_L$.
\end{proof}

\subsection{Additional notations}
\label{subsection:PG_additional_notations}

For an admissible control parameter $\theta = (K, L) \in \Theta$, let $\by$ and $\bz$ be the two auxiliary processes following dynamics~\eqref{eq:dyn_y_theta}--\eqref{eq:dyn_z_theta}. We define $\Sigma_K^y$ and $\Sigma_L^z$ as the infinite discounted sum of the variance/covariance matrices given by
\begin{equation}
\label{eq:variance_matrices_y_and_z}
    \Sigma_K^y := \EE \Big[ \sum \nolimits_{t \geq 0} \gamma^t y_t (y_t)^\top \Big], \qquad  \Sigma_L^z := \EE \Big[ \sum \nolimits_{t \geq 0} \gamma^t z_t (z_t)^\top \Big].
\end{equation}
By an abuse of terminology, we call the matrices $\Sigma_K^y$ and $\Sigma_L^z$ the variance matrices of the auxiliary processes $\by$ and $\bz$. To alleviate the notation, we may omit the superscript in $\Sigma_K^y$ and $\Sigma_L^z$ and use $\Sigma_K, \Sigma_L$ when the context is clear.
From the expressions of $C_y(K)$ and $C_z(L)$ in~\eqref{eq:def_Cy_Cz}, we get 
\begin{equation}
\label{eq:expression_cost_sigma}
    C_y(K) = \langle Q + K^\top R K, \, \Sigma_K \rangle_{tr}, \qquad C_z(L) = \langle \tilde Q + L^\top \tilde R L ,\, \Sigma_L  \rangle_{tr}.
\end{equation}
We denote the exact gradients of the cost $C(\theta)$ with respect to $K$ and on $L$ by $\nabla_K C(\theta)$ and $\nabla_L C(\theta)$ respectively. Lemma~\ref{lem:decompose-cost-reparam} implies immediately that 
$$
\nabla_K C(\theta) = \nabla_K C_y(K), \quad \nabla_L C(\theta) = \nabla_L C_z(L). 
$$
We recall the matrices $\Sigma_{y_0}$, $\Sigma_{z_0}$, $\Sigma^1$, $\Sigma^0$ defined in Assumption~\ref{as:non-deg}, and under the non-degenerate conditions~\eqref{eq:non_degenerate_condition} in~\ref{as:non-deg} we have
\begin{equation}
\label{eq:min_eigenvalue_of_variance_matrices}
    \lambda_y^1 := \lambda_{min} \big( \Sigma_{y_0} + \frac{\gamma}{1-\gamma} \Sigma^1 \big) > 0
    , \quad  
    \lambda_z^0 := \lambda_{min}\big( \Sigma_{z_0} + \frac{\gamma}{1-\gamma} \Sigma^0 \big) > 0.
\end{equation}
The strict positivity of the above two quantities $\lambda_y^1$ and $\lambda_z^0$ is crucial for establishing uniform upper bounds for the norms of the gradients $\nabla_K C_y(K)$ and $\nabla_L C_z(L)$ which are essential to the global convergence of the PG algorithms.

Because we consider sub-gaussian random vectors for the initial perturbations in the LQ problem, the random variables $y_0$ and $z_0$ are also sub-gaussian. We denote $C_{init, noise}$ an upper bound for the sub-gaussian norm of $y_0$, $z_0$, and the i.i.d noise terms $(\varepsilon_{t+1})_{t \geq 0}$ and $(\varepsilon_{t+1}^0)_{t \geq 0}$, namely for every $t \geq 0$,
\begin{equation}
    \max\{ \| y_0 \|_{\psi_2}, \| z_0 \|_{\psi_2}, \| \varepsilon_{t+1} \|_{\psi_2}, \| \varepsilon_{t+1}^0 \|_{\psi_2} \big\} \leq C_{init, noise}.
\end{equation}

\begin{remark}
\label{remark:assumption_on_variance}
In the absence of noise processes, if the initial perturbations $(\varepsilon_0, \varepsilon_0^0)$ are non-degenerate with bounded supports (implying Assumption~\ref{as:non-deg}), the arguments in~\cite{fazel2018global} can be applied to the processes $\by$ and $\bz$ respectively to show the convergence of the PG algorithm with gradients $(\nabla_K C(\theta), \nabla_L C(\theta) )$.
In the more general setting considered here, in particular with non-degenerate sub-gaussian distributions, the strict positivity of the eigenvalues of $\Sigma^0$ and $\Sigma^1$ allows us to derive convergence results without assumption on the initial states as enforced in \cite{fazel2018global}.
It is important to stress the fact that the presence of noise terms poses additional difficulties in the proof of the convergence of PG algorithms in model-free frameworks. The challenges arise not only from the unbounded noise process but also from the quadratic forms involving sub-gaussian random vectors $y_t$ or $z_t$ whose entries are not independent. We cannot simply rely on the classical Berstein's matrix inequality or the Hanson-Wright inequality~\cite{vershynin2018high} for independent random variables, and we need more delicate arguments to prove convergence.
\end{remark}

\subsection{Exact PG for MFC}
\label{subsection:exact_PG_for_MFC}

Assuming full knowledge of the model, we compute the optimal parameters using a standard gradient descent algorithm.
With a fixed learning rate $\eta>0$ and initial parameter $\theta^{(0)} = (K_0,L_0)$, we update the parameter from $\theta^{(k)}$ to $\theta^{(k+1)} = (K^{(k+1)}, L^{(k+1)})$ via: 
\begin{equation}    
\label{eq:one-step-update-exact-pg}
\left\{
    \begin{array}{rcl}
        K^{(k+1)} &=& K^{(k)} - \eta \nabla_K C_y(K^{(k)}), \\
        L^{(k+1)} &=& L^{(k)} - \eta \nabla_L C_z(L^{(k)}).
    \end{array}
\right.
\end{equation}
To show the convergence of $C(\theta^{(k)})$ towards the optimal cost $C(\theta^*)$ with the scheme \eqref{eq:one-step-update-exact-pg}, we adapt the gradient descent methods with Polyak-Lojasiewicz condition~\cite{karimi2016linear, polyak1963gradient, lojasiewicz1963topological}, which has also been referred as the gradient domination condition~\cite{fazel2018global,frikha2024full,yang2019provably,wang2021global}, for the non-convex functions $C_y$ and $C_z$. In the following, we first state the Polyak-Lojasiewicz condition for the auxiliary costs $C_y$ and $C_z$ in Proposition~\ref{proposition:PL_cond_Cy_Cz}, and a local smoothness condition in the direction of gradients in Proposition~\ref{proposition:local_smoothness_Cy_Cz}. The proofs of these results are deferred to  Sections~\ref{subsection:proof_proposition_PL_condition} and \ref{subsection:proof_proposition_local_smoothness} respectively.

\vskip 6pt
We recall that a differentiable function $f: D \to \RR$ 
from a  open set $D \subseteq \RR^m$ into $\RR$ is said to satisfy the \textit{$\nu$-Polyak-Lojasiewicz} ($\nu$-PL for short) condition on $D$  for a constant $\nu > 0$ if 
\begin{equation}
\label{eq:PL_condition_definition}
    \nu \big( f(x) - f(x^*) \big) \leq \| \nabla_x f(x) \|^2, \qquad \forall \, x \in D
\end{equation}
where $x^* \in D$ is a global minimal points in $D$, i.e. $f(x) \geq f(x^*)$ for all $x \in D$. 

\begin{proposition}[PL condition]
\label{proposition:PL_cond_Cy_Cz}
     The cost functions $C_y$ and $C_z$ defined in~\eqref{eq:def_Cy_Cz} satisfy the PL condition on $\Theta$: for $\theta = (K, L) \in \Theta$, we have
     \begin{equation}
        \label{eq:PL_cond_Cy_Cz_w_lambda}
        \left\{
        \begin{array}{rcl}
        \nu_{y}. \big( C_y(K) - C_y(K^*) \big) & \leq & \| \nabla_K C_y(K) \|_F^2, \\ 
        \nu_{z}. \big( C_z(L) - C_z(L^*) \big) & \leq & \| \nabla_L C_z(L) \|_F^2.
        \end{array}
        \right.
     \end{equation}
     with constants $\nu_y =  4 (\lambda_y^1)^2 \lambda_{min}(R) / \| \Sigma_{K^*} \| $ and $\nu_z =  4 (\lambda_z^0)^2 \lambda_{min}(\tilde R) / \| \Sigma_{L^*} \| $.
\end{proposition}

\begin{proposition}[Local smoothness in exact PG direction]
\label{proposition:local_smoothness_Cy_Cz}
    Consider $\theta=(K,L) \in \Theta$, $K'  = K - \eta \nabla_K C_y(K)$ and $L'  = L - \eta \nabla_L C_z(L)$ such that $\theta'=(K', L') \in \Theta$,  and let $\Delta K = K' - K$ and $\Delta L = L' - L$. Then
    \begin{equation}
    \label{eq:local_smoothness_Cy_Cz}
    \left\{
        \begin{array}{rl}
        C_y(K') - C_y(K) & \leq  \big( 1 - \lambda_{var, y} \big) \langle \nabla_K C_y(K),\, \Delta K \rangle_{tr} + \lambda_{hess, y}.  \| \Delta K \|_F^2
        \\
        C_z(L') - C_z(L) & \leq   \big( 1 - \lambda_{var, z} \big) \langle \nabla_L C_z(L),\, \Delta L \rangle_{tr} + \lambda_{hess, z}.  \| \Delta L \|_F^2
        \end{array}
    \right.
    \end{equation}
    where 
\begin{equation}
    \label{fo:lambdas}
    \begin{split}
&\lambda_{var, y} =  \lambda_{var, y}(K, K') = \| \Sigma_{K'} - \Sigma_{K} \| / \lambda_y^1\\
&\lambda_{hess, y} =  \lambda_{hess, y}(K, K') =  \| \Sigma_{K'} \| (\| R \|  + \gamma \|\rB\|^2 C_y(K) / \lambda_y^1 )        
    \end{split}
\end{equation}
and $\lambda_{var, z}, \lambda_{hess, z}$ are defined similarly with $\big(\Sigma_{L'}, \Sigma_L, \lambda_z^0, \tilde R, \tilde \rB, C_z(L) \big)$.
\end{proposition}    

For each $C_0\in\RR$, we shall use freely the  four constants $h_{var, y}(C_0)$, $h_{hess, y}(C_0)$, $h_{var, z}(C_0)$, and $h_{hess, z}(C_0)$ defined in Lemma~\ref{lemma:bound_on_coefficent_in_prop_local_smoothness} in Section~\ref{section:proof_of_exact_pg} where we conduct a detailed perturbation analysis on the control parameters. These constants are polynomials in $( C_0, 1/\lambda_y^1, 1/\lambda_z^0)$ and other model parameters.
The following Lemma~\ref{lemma:small_learning_rate} states that when the learning rate $\eta$ is small enough, the perturbed control parameter $\theta' = (K', L')$ in Proposition~\ref{proposition:local_smoothness_Cy_Cz} is admissible and the coefficients related to the local smoothness in equation~\eqref{eq:local_smoothness_Cy_Cz} are bounded. 
It is stated here in this form for the purpose of the proof of Theorem \ref{thm:exact-CV} that follows. Its contents and their proofs will be part of
Lemma~\ref{lemma:bound_on_coefficent_in_prop_local_smoothness}.

\begin{lemma}[Small learning rate]
\label{lemma:small_learning_rate}
    Consider $\theta = (K, L) \in \Theta$ with $C(\theta) \leq C_0$ for some constant $C_0 \in \RR$. If the learning rate $\eta$ is small enough and satisfies
    $$
        \eta \leq \min\{ h_{var, y}(C_0)^{-1} , h_{var, z}(C_0)^{-1} \},
    $$
    then the perturbed control parameter $\theta' = (K', L')$ in Proposition~\ref{proposition:local_smoothness_Cy_Cz} with learning rate $\eta$ is admissible in the sense that $\theta ' \in \Theta$. Moreover,
    \begin{align*}
        \lambda_{var, y}(K, K') + \eta \lambda_{hess, y}(K, K') &\leq \eta. \big( h_{var, y}(C_0) + h_{hess,y}(C_0) \big),
        \\
        \lambda_{var, z}(L, L') + \eta \lambda_{hess, z}(L, L') & \leq \eta. \big( h_{var, z}(C_0) + h_{hess, z}(C_0) \big).
    \end{align*}
\end{lemma}

We are now ready to state and prove a global convergence result based on the PL condition and a local smoothness condition on the cost functions. Consider the sequence $(\theta^{(k)})_{k \geq 0}$ generated by the PG updates~\eqref{eq:one-step-update-exact-pg}. 

\begin{theorem}[Exact PG]
\label{thm:exact-CV}
Let $C_0 = C(\theta^{(0)})$ and $\nu_{pl} = \min\{ \nu_y, \nu_z \} / 2$, and let us assume that the learning rate $\eta$ is small enough so that 
\begin{equation}
\label{eq:condition_learning_rate_exact_pg}
    \eta \leq \frac{1}{4} \min \big\{ h_{var, y}(C_0)^{-1}, h_{hess, y}(C_0)^{-1} ,h_{var, z}(C_0)^{-1} , h_{hess, z}(C_0)^{-1} \big\}.
\end{equation}
Under Assumptions~\ref{as:finit_cost},~\ref{as:positivity-qr} and~\ref{as:non-deg}, for every $k \geq 0$, we have $C(\theta^{(k+1)}) \leq C(\theta^{(k)}) \leq C_0$, and
\begin{equation}
\label{eq:linear_convergence_exact_pg}
    C(\theta^{(k+1)}) - C(\theta^*)  \leq \big( 1 - \eta \nu_{pl} \big) \big(  C(\theta^{(k)}) - C(\theta^*) \big).
\end{equation}
Moreover, for every $\varepsilon > 0$, we have 
\begin{equation}
\label{eq:epsilon_approx_of_exact_pg}
C(\theta^{(k)}) - C(\theta^*)  \leq \varepsilon  \big(  C(\theta^{(0)}) - C(\theta^*) \big), \qquad k \geq \frac{1}{\eta \nu_{pl} } \log \big( \frac{1}{\varepsilon} \big).
\end{equation}
\end{theorem}

\begin{proof}

    We assume that at iteration step $j \geq 0$, the control parameter $\theta^{(j)}$ is admissible, and $C(\theta^{(j)}) \leq C_0$. Lemma~\ref{lemma:small_learning_rate} shows that that $\theta^{(j+1)} \in \Theta$ and 
    $$
        \lambda_{var, y} + \eta \lambda_{hess, y} \leq \eta \big( h_{var, y}(C_0) + h_{hess, y}(C_0) \big) \leq \frac{1}{2}.
    $$
    Because $Q \succeq 0$ and $R \succ 0$, the expression \eqref{eq:expression_cost_sigma}  of the cost implies that 
    $$
        \lambda_y^1 \leq \| \Sigma_{K^*} \| \leq \frac{ C_y(K^*) }{ \lambda_{min}(Q) } \leq \frac{ C_0 }{ \lambda_{min}(Q)}.
    $$
    Together with the definition of $h_{hess, y}(C_0) = \frac{2 C_0}{\lambda_{min}(Q)} \big( \| R \| + \gamma \| \rB \|^2 \frac{C_0}{\lambda_y^1} \big)$ in~\eqref{eq:h_hess_y}, and the definition of $\nu_y$ in Proposition~\ref{proposition:PL_cond_Cy_Cz}, we have
    $$
        \frac{1}{2} \eta \nu_y = \frac{1}{2} \eta 4 \lambda_{min}(R) \lambda_y^1 \frac{\lambda_y^1}{\| \Sigma_{K^*} \|} \leq  2 \eta  \| R \| \frac{ C_0} { \lambda_{min}(Q) } \leq \eta h_{hess, y}(C_0) \leq \frac{1}{4} < 1.
    $$    
    Consequently, with the local smoothness condition~\eqref{eq:local_smoothness_Cy_Cz} and the $\nu_y$-PL condition~\eqref{eq:PL_cond_Cy_Cz_w_lambda}, we obtain
    \begin{align*}
        C_y( K^{ (j+1)} ) - C_y ( K^{(j)} ) & \leq - \eta \big( 1 - \lambda_{var, y} - \eta \lambda_{hess, y} \big) \| \nabla_K C_y(K^{(j)})\|_F^2
        \\
        & \leq - \frac{1}{2} \eta \nu_y \big( C_y( K^{ (j)} ) - C_y ( K^* ) \big).
    \end{align*}
    The difference in costs yields directly that $C_y(K^{(j+1)}) \leq C_y(K^{(j)})$, and
    $$
        C_y(K^{(j+1)} ) - C_y(K^*) \leq (1 - \eta \nu_{pl} ) \big( C_y( K^{(j)}) - C_y(K^*) \big).
    $$
    Applying the same arguments to $C_z$ and adding up the inequalities for $C_y$ and $C_z$, we obtain equation~\eqref{eq:linear_convergence_exact_pg} for $k=j$ and $C(\theta^{(j+1)}) \leq C(\theta^{(j)}) \leq C_0$. 
    
    Furthermore, by noticing that $ \frac{1}{\eta \nu_{pl}} \log ( \frac{1}{1 - \eta \nu_{pl}}) \geq 1$, when $k \geq \frac{1}{\eta \nu_{pl}} \log( \frac{1}{\varepsilon} )$, we have
    $$
        \log \Big( \frac{ C(\theta^{(0)}) - C(\theta^*) } {C(\theta^{(k)}) - C(\theta^*)  } \Big) \geq k \log( \frac{1}{1 - \eta \nu_{pl}} )  \geq  \log( \frac{1}{\varepsilon} ),
    $$
which concludes the proof of the theorem.
\end{proof}

\begin{remark}
\label{remark:importance_of_non_degeneracy_assumption}
It is worth highlighting that the non-degeneracy Assumption~\ref{as:non-deg} on the sources of randomness plays a crucial role in the PL condition and the smoothness condition mentioned in Propositions~\ref{proposition:PL_cond_Cy_Cz} and \ref{proposition:local_smoothness_Cy_Cz}. Notably, it allows to have convergence even when the system starts from a deterministic initial condition.
\end{remark}

\subsection{Model free PG for MFC with an MKV simulator}
\label{subsection:modelfree_MKV_simulator}

We now turn to a modification of the previous PG updates~\eqref{eq:one-step-update-exact-pg}, in which the gradient will be approximated in a model-free way. The new update scheme takes the form
\begin{equation}    
\label{eq:one-step-update-model-free-MKV}
\left\{
    \begin{array}{rcl}
        K^{(k+1)} &=& K^{(k)} - \eta \tilde{\nabla}^{T,M,\tau}_K (\theta^{(k)}), 
        \\
        L^{(k+1)} &=& L^{(k)} - \eta \tilde{\nabla}^{T,M,\tau}_L (\theta^{(k)}),
    \end{array}
\right.
\end{equation}
where $\tilde \nabla^{T, M, \tau}(\theta) = \big( \tilde{\nabla}^{T, M,\tau}_K(\theta), \tilde{\nabla}^{T, M,\tau}_L(\theta) \big)$ is to be defined below, based on samples generated by a simulator.

\vskip 6pt
So we assume that we have access to the following (stochastic) simulator, called MKV simulator $\cS^T_{MKV}$: given an admissible control parameter $\theta \in \Theta$, $\cS^T_{MKV}(\theta)$ returns a sample of the discounted MF cost (see Section~\ref{subsec:OC-MKV} and equation~\eqref{fo:MKV-discounted_cost}) for the MKV dynamics~\eqref{fo:MKV-state} controlled by the process $(u_t^\theta)_{t \geq 0}$ and truncated at a finite horizon $T$. 
In other words, it returns a realization of the truncated MF cost 
\begin{equation}
\label{eq:MKV_simulator_truncated_MF_cost}
\tilde C^T(\theta) := \sum_{t=0}^{T-1} \gamma^t c_t(\theta),
\end{equation}
where $c_t(\theta) =  c\big(X_t^{\theta}, \bar{X}_t^{\theta}, u_t^\theta, \bar{u}_t^\theta \big)$ is the instantaneous cost at time $t$ (see equation~\eqref{eq:lq_one_step_cost}).%

While the simulator uses the model to generate samples, we aim at developing a model-free algorithm to compute the optimal control parameter from the samples produced by the simulator, without using the knowledge of the model. Without full knowledge of the model, we use derivative-free techniques (see Section~\ref{sec:app-proof-modelfree-MKV-CV}) to estimate the gradient of the cost. 
We introduce some notation for the gradient estimation.
Let $\mathbb{B}_\tau\subset \mathbb{R}^{\ell \times d}$ be the ball of radius $\tau$ under the Frobenius norm centered at the origin, and $\mathbb{S}_\tau = \partial \mathbb{B}_\tau$ be its boundary.  
The uniform distributions on $\mathbb{B}_\tau$ and $\mathbb{S}_\tau$ are denoted by $\mu_{\mathbb{B}_\tau}$ and $\mu_{\mathbb{S}_\tau}$ respectively. We should think of $v = (v^{(idy)}, v^{(com)}) \in \SS_\tau \times \SS_\tau$ as a direction for a perturbation of a parameter $\theta = (K, L) \in \Theta$. For $M$ perturbation directions $\underline{v} = (v_i)_{i=1}^M = \big( (v_i^{(idy)}, v_i^{(com)}) \big)_{i=1}^M$ in which $v_i^{(idy)}$ and $v_i^{(com)}$ are independently sampled from $\mu_{\SS_\tau}$, we denote by $\theta_i = (K + v^{(idy)}_i, L + v^{(com)}_i)$ the control parameter $\theta$ perturbed by $v_i$. Then, the approximated gradient based on the MKV simulator $\cS^T_{MKV}$ and $M$ perturbation directions $\underline{v}$ on $\SS_{\tau} \times \SS_{\tau}$ is defined by 
\begin{equation}
\label{eq:definition_approx_gradient_MKV_simulator}
    \tilde{\nabla}^{T, M,\tau}_K (\theta) = \frac{\ell d}{\tau^2}\frac{1}{M} \sum_{i=1}^M \tilde{C}^T(\theta_i) v^{(idy)}_i, 
    \quad
    \tilde{\nabla}^{T, M,\tau}_L (\theta) = \frac{\ell d}{\tau^2}\frac{1}{M} \sum_{i=1}^M \tilde{C}^T(\theta_i) v^{(com)}_i.
\end{equation}
We omit the perturbation directions $\underline v$ in the notations defined in equation~\eqref{eq:definition_approx_gradient_MKV_simulator} for clarity here. 

\begin{algorithm}
	\caption{Model-free MKV-Based Gradient Estimation}
	\label{algo:MKVestim}
	\begin{algorithmic}
		\STATE {\bfseries Data:} {Parameter $\theta = (K,L)$;  truncation horizon $T$; number of perturbations $M$; perturbation radius $\tau$}
		\STATE {\bfseries Result:} {An approximated gradient $\tilde \nabla^{T, M, \tau}(\theta)$ of $\nabla C(\theta)$}
		\FOR{$i = 1, 2, \dots, M$}
		\STATE Sample $v_{i}^{(idy)}, v_{i}^{(com)}$ i.i.d. $\sim \mu_{\mathbb{S}_\tau}$\; 
		\STATE Set $\theta_i = (K_i, L_i) = (K+ v_{i}^{(idy)}, L+ v_{i}^{(com)})$\;
		\STATE Compute $\tilde{C}^T(\theta_i) = \sum_{t=0}^{T-1} \gamma^t c_t(\theta_i)$ 
		\ENDFOR
		\STATE {\bfseries Set } {$\tilde{\nabla}^{T, M, \tau}_K (\theta), \tilde{\nabla}^{T, M, \tau}_L (\theta)$ with equation~\eqref{eq:definition_approx_gradient_MKV_simulator}}
		\STATE {\bfseries Return: }{$\tilde{\nabla}^{T, M, \tau} (\theta) = \big( \tilde{\nabla}^{T, M, \tau}_K (\theta), \tilde{\nabla}^{T, M, \tau}_L (\theta) \big)$}
	\end{algorithmic}
\end{algorithm}

Algorithm~\ref{algo:MKVestim} provides a (biased) estimator of the true policy gradient $\nabla C(\theta)$. This method is in the spirit of \cite[Algorithm 1]{fazel2018global}, except that here we have two components playing different roles for the state and the conditional mean. The approximated gradient $\tilde{\nabla}^{T, M, \tau}(\theta)$ thus depends in addition on three hyper-parameters, $(T, M, \tau)$: truncation horizon, number of perturbation directions, and perturbation radius. 
In Proposition~\ref{proposition:gradient_approx_with_MKV_simulator} below, we provide conditions on $(T, M, \tau)$ so that Algorithm~\ref{algo:MKVestim} yields a sufficiently good gradient estimator to ensure convergence of the approximate PG.

\begin{proposition}[Gradient approximation with an MKV simulator]
\label{proposition:gradient_approx_with_MKV_simulator} 
Consider $\theta \in \Theta$ with $C(\theta) \leq C_0$ for some constant $C_0 \in \RR$. Let $\tilde{\varepsilon} > 0$ be a target precision and $\delta_{approx} \in (0,1)$. We assume that the parameters $(T, M, \tau)$ in Algorithm~\ref{algo:MKVestim} satisfy
    \begin{align}
        \tau^{-1} & \geq \phi_{pert, radius, MKV}(\tilde{\varepsilon}, C_0),
        \label{eq:model_free_pg_condition_pert_radius}
        \\
        T  & \geq \phi_{trunc, T, MKV}(\tilde{\varepsilon}, \tau, C_0),
        \label{eq:model_free_pg_condition_truncation_T}
        \\
        M & \geq \phi_{sample, size, MKV} \big(\tilde{\varepsilon}, \tau, T, C_0, \delta_{approx} \big),
        \label{eq:model_free_pg_condition_sample_size}
    \end{align}
    where $\phi_{pert, radius, MKV}$, $\phi_{trunc, T, MKV}$, and $\phi_{sample, size, MKV}$ are defined in Section~\ref{subsection:proof_of_proposition_approx_gradient_w_MKV_simulator}.
    Then, we have
    $$
        \PP \big( \big\| \tilde \nabla^{T, M, \tau}(\theta) - \nabla C(\theta) \big\|  > \tilde{\varepsilon} \big) \leq \delta_{approx}.
    $$
    where $\big\| \tilde \nabla^{T, M, \tau}(\theta) - \nabla C(\theta) \big\| = \big\| \tilde \nabla^{T, M, \tau}_K(\theta) - \nabla_K C(\theta) \big\| + \big\| \tilde \nabla^{T, M, \tau}_L(\theta) - \nabla_L C(\theta) \big\|$, and the randomness in $ \tilde \nabla^{T, M, \tau}(\theta)$ comes from the simulator $\cS_{MKV}^T$ and the perturbation directions $\underline v \in (\SS_\tau \times \SS_\tau)^{M}$.
\end{proposition} 
The proof of Proposition~\ref{proposition:gradient_approx_with_MKV_simulator} is deferred to Section~\ref{subsection:proof_of_proposition_approx_gradient_w_MKV_simulator} after we show the approximation of exact policy gradient with derivative-free techniques based on the uniform-smoothing on control parameters (Lemma~\ref{lemma:approx_perturbed_gradient_MKV_simulator}), the approximation of infinite-horizon cost with a truncated finite-horizon cost (Lemma~\ref{lemma:approx_truncated_policy_gradient_MKV_simulator}), and the approximation of expected truncated costs with sampled costs from a simulator (Lemma~\ref{lemma:approx_sampled_pg_with_truncated_pg}). It is worth noticing that the bounds on $(\tau, T, M)$ in~\eqref{eq:model_free_pg_condition_pert_radius},~\eqref{eq:model_free_pg_condition_truncation_T},~\eqref{eq:model_free_pg_condition_sample_size} are independent of the control parameter $\theta \in \Theta$ as long as $C(\theta) \leq C_0$.

Now, we are ready to state the convergence results in the model-free setting with an MKV simulator $\cS_{MKV}^T$. We consider a bound for the learning rate: 
\begin{align*}
    \phi_{lrate, MKV}(C_0) = & \min \Big\{ \min \{ h_{var, y}^{-1}(C_0), h_{var, y}^{-1}(C_0), h_{hess, y}^{-1}(C_0), h_{hess, y}^{-1}(C_0) \} / 4, 
    \\
     & \hspace{15pt} \min\big\{ h_{small-pert, y}^{-1}(C_0), h_{small-pert, z}^{-1}(C_0) \big\}. (2 h_{cost}(C_0)) / \nu_{pl}  \Big\}
\end{align*}
where the constants $h_{var, y}(C_0)$, $h_{var, z}(C_0)$, $h_{hess, y}(C_0)$, $h_{hess, z}(C_0)$, $h_{small-pert,y}(C_0)$, $h_{small-pert, z}(C_0)$ depending on $C_0$ are defined in Lemma~\ref{lemma:useful_bounds_I} for small perturbation of control parameters, and $h_{cost}(C_0)$ is defined in Lemma~\ref{lemma:perturbation_cost_Cy}. We consider a sequence of control parameters $(\theta^{(k)})_{k \geq 0}$ generated with the model-free PG updates~\eqref{eq:one-step-update-model-free-MKV} and the approximated gradients $\big(\tilde \nabla^{T, M, \tau} (\theta^{(k)}) \big)_{k \geq 0}$ from Algorithm~\ref{algo:MKVestim} using the MKV simulator $\cS_{MKV}^T$.

\begin{theorem}
\label{thm:modelfree-MKV-CV}
   Consider a target precision $\varepsilon \leq 1$, an initial parameter $\theta^{(0)} \in \Theta$ with cost $C(\theta^{(0)}) = C_0$, and  a learning rate satisfying
    \begin{equation}
    \label{eq:learning_rate_condition_pg_MKV}
        \eta \leq  \phi_{lrate, MKV}(C_0 + 1).
    \end{equation}
    We choose simulation parameters $\big(\tau, T, M \big)$ in Algorithm~\ref{algo:MKVestim} such that they satisfy equations~\eqref{eq:model_free_pg_condition_pert_radius}, \eqref{eq:model_free_pg_condition_truncation_T}, and \eqref{eq:model_free_pg_condition_sample_size} in Proposition~\ref{proposition:gradient_approx_with_MKV_simulator} with parameters $\tilde{\varepsilon} = \varepsilon  \nu_{pl}/ (2 h_{cost}(C_0 + 1))$ and $\delta_{approx} \in (0, 1)$.  Then, under our standing assumptions, for every iteration step $k \geq 0$, if $\theta^{(k)} \in \Theta$, then we have with high probability (at least $1 - \delta_{approx}$) that $\theta^{(k+1)} \in \Theta$, $C(\theta^{(k+1)}) \leq C_0 + 1$, and   
    \begin{equation}
    \label{eq:contraction_pg_MKV}
        C(\theta^{(k+1)}) - C(\theta^*) \leq \big( 1 -  \eta \nu_{pl}/ 2 \big). \max \Big\{ C(\theta^{(k)}) - C(\theta^*) , \, \varepsilon \Big\}.
    \end{equation}
    Moreover, when the number of iteration steps $k$ satisfies 
    $$
    k \geq \frac{2}{\eta \nu_{pl}} \log( \frac{C_0 - C(\theta^*)}{ \varepsilon} ),
    $$
    we have, with a high probability, that
    \begin{equation}
        C(\theta^{(k)}) - C(\theta^*) \leq \varepsilon.
    \end{equation}
\end{theorem}

\begin{proof}
    Suppose that we are at iteration step $k \geq 0$ of the PG algorithm with MKV simulator $S_{MKV}^T(\theta^{(k)})$. 
    
    We consider first $\theta^{(exact, k)} = (K^{(exact, k)}, L^{(exact, k)})$ generated using the exact PG updates~\eqref{eq:one-step-update-exact-pg}: $K^{(exact, k)} = K^{(k)} - \eta \nabla_K C(\theta^{(k)})$, $L^{(exact, k)} = L^{(k)} - \eta \nabla_L C(\theta^{(k)})$. Under the assumption on $\eta$ by inequality~\eqref{eq:learning_rate_condition_pg_MKV}, and the fact that the coefficient functions $h_{var,y}$, $h_{var, z}$, $h_{hess,y}$ and $h_{hess, z}$ defined in Lemma~\ref{lemma:bound_on_coefficent_in_prop_local_smoothness} are increasing function of $C_0$, the learning rate $\eta$ then becomes small enough to satisfy condition~\eqref{eq:condition_learning_rate_exact_pg} in Theorem~\ref{thm:exact-CV} for the exact-PG algorithm, which implies that
    $$
        C(\theta^{(exact, k)}) - C(\theta^{*}) \leq ( 1- \eta \nu_{pl}) \big( C(\theta^{(k)}) - C(\theta^*) \big),
    $$
    where $\nu_{pl} = \min\{ \nu_{K} , \nu_{L} \} / 2$ is a coefficient related the PL conditions. We also have that $C(\theta^{(exact, k)}) \leq C(\theta^{(k)}) \leq C_0 + 1$ . 
    
    Now, we consider a control parameter $\theta^{(k+1)} = (K^{(k+1)}, L^{(k+1)})$ generated from the one-step model-free update Scheme~\eqref{eq:one-step-update-model-free-MKV}. We apply Proposition~\ref{proposition:gradient_approx_with_MKV_simulator} to $\theta = \theta^{(k)}$ with $\tilde{\varepsilon} = \kappa \varepsilon$ where $\kappa = \nu_{pl} / \big( 2 h_{cost}(C_0 + 1) \big)$ and with some $\delta_{approx} \in (0,1)$. If $(T, M, \tau)$ satisfies equations \eqref{eq:model_free_pg_condition_truncation_T}, \eqref{eq:model_free_pg_condition_sample_size}, and \eqref{eq:model_free_pg_condition_pert_radius} for $(\tilde \varepsilon, \delta_{approx})$, we then deduce that with probability at least $1 - \delta_{approx}$, 
    $$
        \| \theta^{(k+1)} - \theta^{(exact, k)} \| = \eta  \| \tilde \nabla^{M, T, \tau}(\theta^{(k)}) - \nabla C(\theta^{(k)}) \|  \leq \eta \kappa \varepsilon.
    $$
    When the learning rate $\eta$ is small enough to satisfy~\eqref{eq:learning_rate_condition_pg_MKV}, we have 
    $$
    \eta \kappa h_{small-pert, y}(C_0 + 1) \leq 1
    \qquad\text{and}\qquad \eta \kappa h_{small-pert, z}(C_0 + 1)  \leq 1.
    $$
    From the fact that $C(\theta^{(exact, k)}) \leq C(\theta^{(k)}) \leq  C_0 + 1$, the quantity $f_{var}(K^{(exact, k)})$ (resp. $f_{var}(L^{(exact, k)})$ ) defined in \eqref{eq:f_var_K} is bounded by the constant 
    $$
    h_{small-pert,y}(C_0 + 1)\qquad\text{ (resp.}\qquad h_{small-pert, z}(C_0 + 1))
    $$
    defined in Lemma~\ref{lemma:useful_bounds_I}. 
    As a consequence, with probability at least $1 -\delta_{approx}$, we have
    \begin{align*}
        & \max \Big\{ \| K^{(k+1)} - K^{(exact, k)} \|. f_{var}(K^{(exact, k)}), \,  \| L^{(k+1)} - L^{(exact, k)} \|.f_{var}(L^{(exact, k)}) \Big\} 
        \\
        &\hskip 45pt
        \leq \| \theta^{(k+1)} - \theta^{(exact, k)} \|. \max \big\{ f_{var}(K^{(exact, k)}), f_{var}(L^{(exact, k)}) \big\} 
        \\
        &\hskip 45pt
        \le \eta \kappa \varepsilon. \max\big\{ h_{small-pert, y}(C_0 + 1), h_{small-pert, z}(C_0 + 1) \big\} 
        \\
        &\hskip 45pt
        \le \varepsilon \leq 1.
    \end{align*}
    This means that $K^{(k+1)}$ (resp. $L^{(k+1)}$) is a small perturbation of $K^{(exact, k)}$ (resp. $L^{(exact, k)}$) satisfying the condition in the stability Lemma~\ref{lemma:stability_of_small_perturbation}, and thus $K^{(k+1)} \in \Theta_{K}$ and $L^{(k+1)} \in \Theta_{L}$. Hence, with probability at least $1 - \delta_{approx}$, the control parameter $\theta^{(k+1)}$ is admissible in the sense that $\theta^{(k+1)} \in \Theta$.
    
    To prove inequality~\eqref{eq:contraction_pg_MKV}, we apply Lemma~\ref{lemma:perturbation_cost_Cy} to the cost function $C(\theta)$ with the small perturbation $\theta^{(exact, k)}$ mentioned above of $\theta^{(k+1)}$, and we obtain that
    \begin{align*}
        \Big| C \big(\theta^{(k+1)} \big) - C \big( \theta^{(exact, k)} \big) \Big| & \leq  h_{cost}(C_0 + 1). \big\| \theta^{(k+1)} - \theta^{(exact, k)} \big\| 
        \\
        & \leq h_{cost}(C_0 + 1) \eta \kappa  \varepsilon
        \\
        & = \eta \nu_{pl} \varepsilon / 2.
    \end{align*}
    Then with probability at least $1 - \delta_{approx}$,
    \begin{align*}
        C(\theta^{(k+1)}) - C(\theta^*) & =  \big( C(\theta^{(exact, k)} )- C(\theta^*) \big) +  \big( C(\theta^{(k+1)}) - C(\theta^{(exact, k)}) \big)
        \\
        & \leq \big( 1 - \eta \nu_{pl} \big). \big( C(\theta^{(k)}) - C(\theta^*) \big) + \eta \nu_{pl}/2. \varepsilon
    \end{align*}
    
    If at iteration step $k$ we have $C(\theta^{(k)}) - C(\theta^*) > \varepsilon$, we deduce that with probability at least $1 - \delta_{approx}$,
    $$
        C(\theta^{(k+1)}) - C(\theta^*) \leq  (1 - \eta \nu_{pl} / 2 ) \big( C(\theta^{(k)}) - C(\theta^*) \big).
    $$   
    In this case, we have a contraction inequality moving from $\theta^{(k)}$ to $\theta^{(k+1)}$. Because $C(\theta^*) \leq C(\theta^{(k)}) \leq C_0 + 1$, we have 
    $$
        C(\theta^{(k+1)})\leq (1 - \eta \nu_{pl}/2) C(\theta^{(k)}) + \eta \nu_{pl} / 2. C(\theta^*) \leq C(\theta^{(k)} ) \leq C_0 + 1.
    $$

    Otherwise, if at iteration step $k$, we reach the target precision $\varepsilon$ such that $C(\theta^{(k)}) - C(\theta^*) \leq \varepsilon$, then with probability at least $1 - \delta_{approx}$,
    $$
        C(\theta^{(k+1)}) - C(\theta^*) \leq ( 1 - \eta \nu_{pl} ) \varepsilon + \eta \nu_{pl} /2 \varepsilon \leq \big( 1 - \eta \nu_{pl} / 2 \big) \varepsilon \leq \varepsilon.
    $$
    In this case, the difference between $C(\theta^{(k+1)})$ and $C(\theta^*)$ is still bounded by $\varepsilon$. Because $C(\theta^*) \leq C(\theta^{(0)}) = C_0$ and $\varepsilon \leq 1$, we have
    $$
        C(\theta^{(k+1)}) \leq  C(\theta^*)  + ( 1 - \eta \nu_{pl} / 2) \varepsilon \leq C_0 + 1.
    $$
    
    Finally, we conclude the proof of the theorem by taking the number of iteration steps $k$ large enough so that $k \geq \frac{2}{\eta \nu_{pl}} \log( \frac{C_0 - C(\theta^*)}{\varepsilon} )$. We showed that $C(\theta^{(j)}) \leq C_0 + 1$ for all $j = 0, \ldots, k$ as long as $(\theta^{(j)})_{j =0, \ldots k}$ are admissible. This means that all the costs along the iteration of the PG algorithm belong to the same level-set, which justifies the use of constant coefficients $(\eta, T, M, \tau)$ in the learning process. 
    Because the admissibility of $\theta^{(j+1)}$ and the inequality ~\eqref{eq:contraction_pg_MKV} between costs depend on the admissibility of $\theta^{(j)}$ for $j = 0, \ldots, k-1$, then for a small enough coefficient $\delta_{approx}$, we have with probability at least $(1 - \delta_{approx})^{k}$ that $\theta^{(k)} \in \Theta$ and $C(\theta^{(k)}) \leq C_0 + 1$. 
    Moreover, if for some $k' < k$ we have $C(\theta^{(k')}) - C(\theta^*) \leq \varepsilon$, we must have $C( \theta^{(j)}) - C(\theta^*) \leq \varepsilon$ for all $j =k', \ldots, k$. Otherwise, if for all $j \le k$ we have $C(\theta^{(j)}) - C(\theta^*) > \varepsilon$, then 
    $$
        C(\theta^{(k)}) - C(\theta^*) \leq (1 - \eta \nu_{pl} / 2)^k \big(C_0 - C(\theta^*) \big) \leq \varepsilon.
    $$
    with high probability (at least $(1 - \delta_{approx})^{k}$ for some small coefficient $\delta_{approx}$). 
\end{proof}

\begin{remark}
    The convergence rate for model-free PG algorithm with an MKV simulator is linear, namely for $k=\mathcal{O}(\log(1/\varepsilon))$ we attain an $\varepsilon$-approximation to the optimal cost with high probability.
\end{remark}

\subsection{Model-free PG for MFC  with a population simulator}
\label{subsec:PG-popsimu}

We now turn our attention to a more realistic setting where one does not have access to an oracle simulating the MKV dynamics, but only to an oracle merely capable of simulating the evolution of $N$ agents. We then use the state sample average instead of the theoretical conditional mean, and for the social cost, the empirical average instead of the mean-field cost provided by the MKV simulator.

We rely on the following population simulator $\cS^{T,N}_{pop}$: given a control parameter $\theta$, $\cS^{T,N}_{pop}(\theta)$ returns a sample of the social cost obtained by generating realizations of the $N$-agent state trajectories controlled by $\theta$ and computing the associated cost for the population, see~\eqref{fo:N-multi_state} and~\eqref{fo:N-cost}. 
In other words, it returns a realization of 
$$
   \tilde{C}^{T,N}(\theta) = \sum_{t=0}^T \gamma^t  \frac{1}{N} \sum_{n=1}^N c \big( X^{(n), \theta}_t,  \bar{X}^{N, \theta}_t, u^{(n), \theta}_t, \bar u^{N, \theta}_t \big)
$$ 
where $c$ is the one-step cost function~\eqref{eq:lq_one_step_cost}, and $(X_t^{(n), \theta}, u_t^{(n), \theta})_{n=1}^N$ are the state-action pairs at time $t$ for the $N$ agents who adopt a same policy function so that 
$
u_t^{(n), \theta} =  -K \big( X_t^{(n), \theta} - \bar{X}_t^{N, \theta} \big) - L \bar{X}^{N, \theta}_t
$
for every $t \geq 0$ and $n=1, \ldots, N$, and where $\bar{X}_t^{N, \theta}$ and $\bar{u}_t^{N, \theta}$ are the average state and average control at time $t$.
Let $C^N(\theta) = J^N(\underline{\bU}^\theta)$ in equation~\eqref{eq:social_cost_of_population} when all agents adopt the same control parameter $\theta \in \Theta$.

Notice that this population simulator $\cS^{T,N}_{pop}$ is arguably more realistic though less powerful than the previous MKV simulator $\cS^T_{MKV}$. Indeed, the former uses only a noisy approximation of the true mean processes while the MKV simulator generates the exact means. We stress that, in this population simulator, \emph{all} agents adopt the same control, and hence in Algorithm~\ref{algo:POPestim}, they use the \emph{same} perturbed version of the control parameter. This design in the algorithm is in line with the idea that the problem corresponds to an optimization problem for a central planner or for a group of cooperative agents using the same control rule to minimize the social cost (see Section~\ref{subsec:OC-N-Agents}).

For the one-step gradient update scheme with population simulator $\cS_{pop}^{T, N}$, we replace the gradient estimator $\tilde{\nabla}^{T, M, \tau}(\theta)$ obtained from Algorithm~\ref{algo:MKVestim} with an MKV simulator by another zero-th order approximation $\tilde{\nabla}^{T, N, M,\tau}(\theta)$ of the gradient based on $\cS_{pop}^{T, N}$ generated by Algorithm~\ref{algo:POPestim}.
The term $\tilde{\nabla}^{T, N, M, \tau}(\theta) = \big( \tilde{\nabla}^{T, N, M, \tau}_K( \theta),  \tilde{\nabla}^{T, N, M, \tau}_L(\theta) \big)$ is called the \emph{sampled population policy gradient} at $\theta \in \Theta$. It is defined with $M$ perturbation directions $(v_i)_{i=1}^M = (v_i^{(idy)}, v_i^{(com)})_{i=1}^M$ on $\SS_\tau \times \SS_\tau$:
\begin{equation}
\label{eq:sampled_population_gradient}
    \tilde{\nabla}_K^{T, N, M, \tau}(\theta) =  \frac{\ell d}{\tau^2} \frac{1}{M} \sum_{i=1}^M  \tilde{C}^{T, N}(\theta_i) v_i^{(idy)}, 
    \quad
    \tilde{\nabla}_L^{T, N, M, \tau}(\theta) = \frac{\ell d}{\tau^2} \frac{1}{M} \sum_{i=1}^M \tilde{C}^{T, N}(\theta_i) v_i^{(com)}
\end{equation}
where $\theta_i = \theta + v_i$ is the perturbed control parameter with $v_i$.
The update scheme at iteration step $k$ for a parameter $\theta^{(k), pop} = (K^{(k), pop}, L^{(k), pop} )$ is defined by
\begin{equation}    
\label{eq:one-step-update-model-free-PoP}
\left\{
    \begin{array}{rcl}
        K^{(k+1), pop} &=& K^{(k), pop} - \eta \tilde{\nabla}^{T, N, M,\tau}_K (\theta^{(k), pop}), 
        \\
        L^{(k+1), pop} &=& L^{(k), pop} - \eta \tilde{\nabla}^{T, N, M,\tau}_L (\theta^{(k), pop}).
    \end{array}
\right.
\end{equation}

\begin{algorithm}
	\caption{Model-free Population-Based Gradient Estimation}
	\label{algo:POPestim}
	\begin{algorithmic}
		\STATE {\bfseries Data:} {Parameter $\theta = (K,L)$;  truncation horizon $T$; number of perturbations $M$; perturbation radius $\tau$; number of agents $N$.}
		\STATE {\bfseries Result:} {An approximation of $\nabla C(\theta)$.}
		\FOR{$i = 1, 2, \dots, M$}
		\STATE Sample $v_i^{(idy)}, v_i^{(com)}$ i.i.d. $\sim \mu_{\mathbb{S}_\tau}$\; 
		\STATE Set $\theta_i = \big( K+ v_i^{(idy)}, L+ v_i^{(com)} \big)$ \;
		\STATE \textbf{Sample $\tilde{C}^{T, N}(\theta_i) = \sum_{t=0}^{T-1} \gamma^t \frac{1}{N} \sum_{n=1}^N c( x_t^{(n), \theta_i}, \bar x_t^{N, \theta_i}, u_t^{(n), \theta_i}, \bar u_t^{N, \theta_i})$ for $N$ homogeneous agents using  $\cS^{T,N}_{pop}(\theta_i)$} \;
		\ENDFOR
		\STATE {\bfseries Set} {$\tilde{\nabla}_K^{T, N, M, \tau}(\theta)$ and $\tilde{\nabla}_L^{T, N, M, \tau}(\theta)$ with equation~\eqref{eq:sampled_population_gradient}}\;
		\STATE {\bfseries Return: }{$\tilde{\nabla}^{T, N, M, \tau}(\theta) = \Big(\tilde{\nabla}_K^{T, N, M, \tau}(\theta), \tilde{\nabla}_L^{T, N, M, \tau}(\theta)  \Big)$}
	\end{algorithmic}
\end{algorithm}

The main result of this section is Theorem~\ref{thm:modelfree-POP-CV}. It shows convergence of the above learning scheme for the optimal control parameter $\theta^* = (K^*, L^*)$ of the MFC problem with a population simulator $\cS_{pop}^{T, N}$ of the social cost of $N$ homogeneous agents. Before proving this result, we state Proposition~\ref{proposition:approx_with_social_cost_to_MF_cost} and Proposition~\ref{proposition:modelfree_population_gradient_approx} providing crucial approximation estimates of the MF cost by its $N-$agent equivalent, and the approximation of the sampled population policy gradient. Their proofs are given in Section~\ref{subsection:proof_of_approx_modelfree_pop_gradient}.

\begin{proposition}
\label{proposition:approx_with_social_cost_to_MF_cost}
Consider $\theta \in \Theta$ with $C(\theta) \leq C_0$ for some $C_0 \in \RR$. Under our standing assumptions, we have
$$
    \| C^N(\theta) - C(\theta) \| \leq \frac1N\phi_{social-cost, factor}(C_0)
$$
where $\phi_{social-cost, factor}(C_0) =  2 d C_0 C_{init, noise}^2 \big( 1/ \lambda_y^1 + 1/ \lambda_z^0 \big) / (1-\lambda)$ is a constant depending on $C_0$.
\end{proposition}

\begin{proposition}
\label{proposition:modelfree_population_gradient_approx}
Let $\theta \in \Theta$ with $C(\theta) \leq C_0$ for some $C_0 \in \RR$, choose a target precision $\tilde \varepsilon > 0$ and $\delta_{approx} \in (0,1)$, and let us assume the parameters $(T, N, M, \tau)$ in Algorithm~\ref{algo:POPestim} satisfy
\begin{align}
    \tau^{-1} &\geq \phi_{pert, radius}(\tilde \varepsilon / 4, C_0)
    \\
    N & \geq \phi_{agent, size, pop}(\tilde \varepsilon / 4, \tau, C_0)
    \\
    T & \geq \phi_{trunc, T, pop}(\tilde\varepsilon / 4 , \tau, C_0, N) + 2
    \\
    M & \geq \max \big\{ \phi_{pert, size}(\tilde\varepsilon / 4, \tau, C_0, \delta_{approx}/2), 
    \nonumber \\
    & \hspace{60pt} \phi_{sample, size, pop}(\tilde\varepsilon / 4, \tau, T, C_0, \delta_{approx} /2  ) \big\}
\end{align}
where $\phi_{pert, radius}$ defined in~\eqref{eq:phi_pert_radius}, and $\phi_{agent, size, pop}$, $\phi_{trunc, T, pop}$, $\phi_{sample, size, pop}$ in Section~\ref{subsection:proof_of_approx_modelfree_pop_gradient} are polynomials in $(d$, $\ell$, $C_0$, $1/ \lambda_{y}^1$, $1/\lambda_z^0$, $C_{init, noise}$, $1/N)$ and other model parameters. Then we have 
\begin{equation}
    \PP \big( \| \tilde \nabla^{T,N,M,\tau}(\theta) - \nabla C(\theta) \| > \tilde \varepsilon \big) \leq \delta_{approx}.
\end{equation}
\end{proposition}

Now, we show that Algorithm~\ref{algo:POPestim} leads to learning the optimal control parameter of the MFC problem with a global linear convergence rate in the social cost. To state the theorem, we consider a sequence of control parameters $(\theta^{(k), pop})_{k \geq 0}$ generated by the model-free PG update Scheme with sampled population policy gradient~\eqref{eq:one-step-update-model-free-PoP}.

\begin{theorem}
\label{thm:modelfree-POP-CV}
We consider an initial control parameter $\theta^{(0), pop} \in \Theta$ with $C(\theta^{(0), pop}) = C_0$ and a target precision $\varepsilon \leq 1$. We assume that the number of agents $N$ and the learning rate $\eta$ satisfy
\begin{align}
    \label{eq:learning_rate_modelfree_pop_simulator}
    \eta & \leq \phi_{lrate, MKV}(C_0 + 1)
    \\
    N &\geq \phi_{social-cost, factor}(C_0 + 1) / (\rho \varepsilon),
    \label{eq:N_agent_social_cost_approx_condition}
\end{align}
where $\rho = (\eta \nu_{pl}) / (16 - 4 \eta \nu_{pl}) \in (0, 1)$. And we choose simulation parameters $(T, N, M, \tau)$ in Algorithm~\ref{algo:POPestim} to satisfy conditions in Proposition~\ref{proposition:modelfree_population_gradient_approx} with $\tilde \varepsilon = ( 1 + 2 \rho) \varepsilon \nu_{pl} / \big( 2 h_{cost}(C_0 + 1) \big)$ and $\delta_{approx} \in (0,1)$. Then, under our standing assumptions, for every iteration step $k \geq 0$, if $\theta^{(k), pop} \in \Theta$, then we have with probability at least $1 - \delta_{approx}$ that $\theta^{(k+1), pop} \in \Theta$, and 
\begin{equation}
\label{eq:contraction_inequality_on_C^N_algon_iteration}
    \big| C^N( \theta^{(k+1), pop}) - C^N(\theta^*) \big| \leq \big( 1 - \eta \nu_{pl} / 4 \big). \max \Big\{ \big| C^N( \theta^{(k), pop} ) - C^N(\theta^*) \big|, \, \varepsilon \Big\}.
\end{equation}
Moreover, when the number of iteration steps $k$ satisfies 
$
k \geq \frac{4}{\eta \nu_{pl}} \log \big( \frac{ | C^N(\theta^{(0), pop}) - C^N(\theta^*) |}{ \varepsilon} \big) 
$
with high probability we have that 
$$
 \big| C^N(\theta^{(k), pop}) - C^N(\theta^*) \big| \leq \varepsilon.
$$
\end{theorem}

\begin{proof}
   The proof applies similar arguments presented in Theorem~\ref{thm:modelfree-MKV-CV} but with slightly different choices of the simulation parameters. At iteration step $k=0$, we have by assumption that $\theta^{(0), pop}\in \Theta$ and $C(\theta^{(0), pop}) = C_0 \leq C_0 + 1$. Suppose that at iteration step $k \geq 0$, we have $\theta^{(k), pop} \in \Theta$ and $C(\theta^{(k), pop}) \leq C_0 + 1$, that is, the control parameter generated by~\eqref{eq:one-step-update-model-free-PoP} at step $k$ with a population simulator $\cS_{pop}^T$ is admissible. We first observe that  Proposition~\ref{proposition:approx_with_social_cost_to_MF_cost} with a large $N$ satisfying~\eqref{eq:N_agent_social_cost_approx_condition} implies
   \begin{align*}
        &  C ( \theta^{(k), pop} ) - C (\theta^*) 
        \\
        \leq & \big|  C^N(\theta^{(k), pop}) - C^N(\theta^*) \big| + \big| C^{N}(\theta^{(k), pop}) - C(\theta^{(k), pop}) \big| + \big| C(\theta^*) - C^N(\theta^*) \big| 
        \\
        \leq & \big| C^N(\theta^{(k), pop}) - C^N(\theta^*) \big| + 2\rho \varepsilon.
   \end{align*}
   Secondly, let $\varepsilon' = ( 1 + 2\rho) \varepsilon$ and $\kappa = \frac{\nu_{pl}}{2 h_{cost}(C_0 + 1)}$. If we choose a target precision $\tilde{\varepsilon} = \varepsilon' \kappa$ and parameter $\delta_{approx} \in (0, 1)$ in Proposition~\ref{proposition:modelfree_population_gradient_approx} to approximate the true gradient using a population simulator with coefficients $(T, N, M, \tau)$ defined accordingly, and we apply the same arguments in Theorem~\ref{thm:modelfree-MKV-CV} under the condition that $C(\theta^{(k), pop}) \leq C_0 + 1$ and $\theta^{(k), pop} \in \Theta$, we have with probability at least $ 1- \delta_{approx}$ that $\theta^{(k+1), pop} \in \Theta$, $C(\theta^{(k+1), pop}) \leq C_0 + 1$, and
   \begin{align*}
        C( \theta^{(k+1), pop} ) - C (\theta^*) & \leq \big( 1 - \eta \nu_{pl} / 2 \big). \max \Big\{ C ( \theta^{(k), pop} ) - C (\theta^*) , \, (1 + 2 \rho) \varepsilon \Big\}
        \\
        & \leq  \big( 1 - \eta \nu_{pl} / 2 \big). \Big( \max \Big\{ \big| C^N(\theta^{(k), pop}) - C^N(\theta^*) \big|, \, \varepsilon \Big\} + 2 \rho \varepsilon \Big).
   \end{align*}
   Besides, from the choice of parameter $\rho$, we have
   $
        2 \rho + 2 ( 1 - \eta \nu_{pl} / 2 ) \rho  = \eta \nu_{pl}  / 4.
   $
   Therefore, by applying Proposition~\ref{proposition:approx_with_social_cost_to_MF_cost} on $\theta^{(k+1), pop} \in \Theta$ with $C(\theta^{(k+1), pop}) \leq C_0 + 1$, we have with probability at least $1 - \delta_{approx}$ that
   \begin{align*}
        & \big| C^N( \theta^{(k+1), pop} ) - C^N(\theta^*) \big| 
        \\
        \leq & 2 \rho \varepsilon + C ( \theta^{(k+1), pop} ) - C(\theta^*)
        \\
        \leq & 2 \rho \varepsilon +  \big( 1 - \eta \nu_{pl} / 2 \big). \Big( \max \Big\{ \big| C^N(\theta^{(k), pop}) - C^N(\theta^*) \big|, \, \varepsilon \Big\} + 2 \rho \varepsilon \Big)
        \\
        = & \big( 2 \rho \varepsilon + 2 ( 1 - \eta \nu_{pl} / 2 ) \rho \varepsilon \big) + \big( 1 - \eta \nu_{pl} / 2 \big) \max \Big\{ \big| C^N(\theta^{(k), pop}) - C^N(\theta^*) \big|, \, \varepsilon \Big\}
        \\
        \leq & \big( 1 - \eta \nu_{pl} / 4 \big) \max \Big\{ \big| C^N(\theta^{(k), pop}) - C^N(\theta^*) \big|, \, \varepsilon \Big\}.
   \end{align*}
    We thus show the inequality~\eqref{eq:contraction_inequality_on_C^N_algon_iteration} between the $N-$agent costs from step $k$ and $k+1$ along Algorithm~\ref{algo:POPestim} using a population simulator $\cS_{pop}^T$ and with constant parameters $(\eta, T, N, M, \tau)$. 
    
    Finally, we conclude the theorem by iterating along the PG algorithm for a long enough iteration steps $k \geq \frac{4}{ \eta \nu_{pl} } \log( \frac{  | C^N(\theta^{(0), pop}) - C^N(\theta^*) | }{\varepsilon} ) $. Because the admissibility of $\theta^{(k+1), pop}$ depends on the admissibility of $\theta^{(k), pop}$, then if there exists $k' < k$ such that $| C^N(\theta^{(k'), pop}) - C^N(\theta^*) | \leq \varepsilon$, we have $|  C^N(\theta^{(j), pop}) - C^N(\theta^*) | \leq \varepsilon$ for all $j=k', \ldots, k$ with probability at least $(1 - \delta_{approx})^{j} \geq (1 - \delta_{approx})^k$; otherwise, we have  
    $$
        | C^N(\theta^{(k), pop}) - C^N(\theta^*) | \leq (1 - \eta \nu_{pl} / 4 )^k | C^N(\theta^{(0), pop}) - C^N(\theta^*) | \leq \varepsilon 
    $$
    with high probability (at least $( 1 - \delta_{approx})^k$ for some small parameter $\delta_{approx}$).
\end{proof}

\begin{remark}
     In the approximation of the social cost with the MF cost in Proposition~\ref{proposition:approx_with_social_cost_to_MF_cost}, the coefficient $h_{social-cost, factor}$ is defined with $\lambda_y^1 > 0$ and $\lambda_z^0 > 0$ under Assumption~\ref{as:non-deg}. However, in the $N-$agent dynamics, the sample average of the states at time $t$, $\bar{X}_t^N$, would remain stochastic even if there were no common noise due to the noise term $\frac{1}{N} \sum_{n=1}^N \varepsilon^{(n)}_{t+1}$ in the dynamics equation of the sample mean process. Hence, we expect that Theorem~\ref{thm:modelfree-POP-CV} could be proved without the non-degeneracy of the noises in Assumption~\ref{as:non-deg}. 
\end{remark}

\section{\textbf{Proofs of the technical results from Section~\ref{subsection:exact_PG_for_MFC}}
\label{section:proof_of_exact_pg}}

This section provides the proofs of the technical lemmas stated without proof in Section~\ref{subsection:exact_PG_for_MFC}.
These lemmas are adapted from~\cite{fazel2018global} with necessary modifications needed to cope with the idiosyncratic noise and the common noise processes present in the state dynamics. 
To make the paper self-contained, we provide their proofs here. Because of the strong similarities between the processes $\by$ and $\bz$, we only state and prove results related to the process $\by$ and the control parameter $K$.

\subsection{Cost expression and gradient expression}
\label{subsection:cost_expression_and_gradient_expression}

\subsubsection{Cost expressions}
\label{subsection:cost_expression}

For an admissible control $\theta = (K, L) \in \Theta$, consider two $L^2$-discounted integrable processes in $\by^{K, \xi_y}, \bz^{L, \xi_z} \in \cX$ following dynamics~\eqref{eq:dyn_y_theta},~\eqref{eq:dyn_z_theta} with initial values $y_0 = y_0^{K, \xi_y} = \xi_y$ and $z_0 = z_0^{L, \xi_z}= \xi_z$. We defined the corresponding value functions $V_y$ and $V_z$ on $\RR^{\ell \times d} \times \RR^d$ by 
\begin{equation*}
V_y(K, \xi_y) = \EE_{\boldsymbol{\varepsilon}} \Big[ \sum_{t \geq 0} \gamma^t f(y_t^{K, \xi_y}, K, Q, R) \Big],
\quad
V_z(L, \xi_z) = \EE_{\boldsymbol{\varepsilon^0}} \Big[ \sum_{t \geq 0} \gamma^t f(z_t^{L, \xi_z}, L, \tilde Q, \tilde R) \Big]
\end{equation*}
where  $f(\xi, K, Q, R) := \xi^\top( Q + K^\top R K ) \xi$ and the expectations $\EE_{\boldsymbol{\varepsilon}}, \EE_{\boldsymbol{\varepsilon^0}}$ are integrations w.r.t. the idiosyncratic and the common noise processes. 

These value functions can be computed simply with the help of solutions $P^y_K$ and $P^z_L$ to the following two discrete-time  Lyapunov equations (DLEs):
\begin{equation}
\label{eq:lyapunov_eq_theta}
    \left\{
        \begin{array}{rcl}
        P^y_K & = &  Q + K^\top R K + \gamma (\rA - \rB K)^\top P^y_K (\rA - B K)
        \\
        P^z_L & = & \tilde{Q} + L^\top \tilde{R} L +\gamma ( \tilde{\rA}- \tilde \rB L)^\top P^z_L ( \tilde{\rA} - \tilde \rB L).
        \end{array}
    \right. 
\end{equation}
The existence of solutions to~\eqref{eq:lyapunov_eq_theta} is deferred to Section~\ref{sec:existence_of_DLE}. Let us also define
\begin{equation}
\label{eq:value_function_noise_constant}
\alpha^y_K = \EE \big[ \sum \nolimits_{t \geq 1} \gamma^t (\varepsilon_t)^\top P^y_K \varepsilon_t \big],
\qquad   
\alpha^z_L = \EE \big[ \sum \nolimits_{t \geq 1} \gamma^t (\varepsilon^0_t)^\top P^z_L \varepsilon^0_t \big].
\end{equation}
Then, using the dynamics of $\by^{K, \xi_y}$~\eqref{eq:dyn_y_theta} and the discrete-time Lyapunov equation~\eqref{eq:lyapunov_eq_theta} for $P_K^y$, we have 
$
    V_y(K, \xi_y) = \xi_y^\top P_K^y \xi_y + \alpha_K^y.
$
When the initial states are random such that $y_0 \sim X_0 - \bar{X}_0 = \varepsilon_0 - \EE[ \varepsilon_0]$, we have 
\begin{equation}
\label{eq:cost_expression_with_PK_y}
 C_y(K) = \EE[ V_y(K, y_0) ] = \langle P_K^y, \,  \Sigma_{y_0} + \gamma/ (1-\gamma) \Sigma^1 \rangle_{tr} = \langle Q + K^\top R K, \, \Sigma_K \rangle_{tr}
\end{equation}
Because $\lambda_y^1 = \lambda_{min}( \Sigma_{y_0} + \gamma/ (1-\gamma) \Sigma^1 )$, we derive the following bounds:
\begin{equation}
\label{eq:upper_bound_Riccati_and_Variance_matrix}
    \|P^y_K\| \leq C_y(K) / \lambda_y^1, 
    \quad  
    \|\Sigma^y_K\| \leq C_y(K) / \lambda_{min}(Q).
\end{equation}
Similar expressions for cost $C_z(L)$, value $V_z(L, \xi_z)$, and bounds for $P_L^z$ and $\Sigma_L^z$ hold true with $(P_L^z, \Sigma_{z_0}, \Sigma^0, \lambda_{z}^0, \tilde Q, \tilde R)$.

\subsubsection{Gradient expression}
\label{subsection:gradient_expression}
We express the gradient of the cost $C_y(K)$ in terms of the variance matrix $\Sigma_K$. This gradient expression is crucial as the approximation of the gradient $\nabla_K C_y(K)$ should hinge on an estimation of the variance matrix $\Sigma_K$. For completeness, we provide the proof here, even though it is analogous to~\cite[Lemma 1]{fazel2018global}.

\begin{lemma} 
    \label{lemma:policy_gradient_expression}
    For $\theta = (K, L) \in \Theta$, let $ E_K = (R + \gamma \rB^\top P^y_K \rB)K - \gamma \rB^\top P^y_K \rA$, then
    $$
        \nabla_K C(\theta) = \nabla_K C_y(K) = 2 E_K \Sigma_K, 
    $$    
\end{lemma}

\begin{proof}
    Consider a process $(y_t^{K, \xi_y})_{t \geq 0}$ starting from $y_0 = \xi_y$, the value function $V_y(K, \xi_y)$ can be rewritten as
    \begin{align*}
	V_y(K, \xi_y) =  \xi_y^\top (Q + K^\top R K) \xi_y  + \gamma \mathbb{E}_{\varepsilon_{t=1}} \big[ V_y(K, (\rA -\rB K)\xi_y + \varepsilon_{t=1} ) \big], 
	\label{eq:expression_C_y_K_ytilde}
    \end{align*}
	and since $\nabla_{\xi_y} V_y(K, \xi_y) = 2 P^y_K \xi_y$, and $\alpha^y_K$ does not depend on $\xi_y$, we have 
	\begin{align*}
    	\nabla_K V_y(K, \tilde y) =  \big(2 RK - 2 \gamma \rB^\top P^y_K (\rA-\rB K) \big)\xi_y \xi_y^\top + \gamma \EE \big[ \nabla_K  V_y (K, y_{t=1}^{k, \xi_y} ) \big].
	\end{align*}
If we replace the gradient inside the expectation by its value, and we note that $C(\theta) = C_y(K) + C_z(L)$ and $C_y(K) = \EE_{\xi_y}[V_y(K, \xi_y)]$, we obtain the result of the lemma.
\end{proof}

\subsubsection{Cost variation expression}

\begin{lemma}
\label{lemma:difference_in_C_y(K)}
    Consider $\theta'=(K', L') \in \Theta$. Let $\Delta K = K' - K$, then
    \begin{equation}
		C_y(K') - C_y(K) =  \langle \Delta K, \, 2 E_K \Sigma_{K'} + (R + \gamma \rB^\top P_K^y \rB )\Delta K  \Sigma_{K'} \rangle_{tr}.
    \end{equation}
\end{lemma}

\begin{proof}
Consider a process $(y_t^{K', y_0})_{t \geq 0}$ starting from $y_0$ and following dynamics~\eqref{eq:dyn_y_theta} with parameter $K'$, and a sequence of random costs $\big( V_y( K, y_t^{K', y_0} ) \big)_{t \geq 0}$ for parameter $K$ and initial random state $\xi_y = y_t^{K', y_0}$ for $t \geq 0$. We introduce the \emph{advantage} quantity at time $t$ given by
\begin{align*}
    A_y(t, K, y_t^{K', y_0}, K') :=& (y_t^{K', y_0})^\top (Q + K'^\top R K') y_t^{K', y_0} + \gamma V_y(K, y_{t+1}^{K', y_0}) - V_y(K, y_t^{K', y_0}) 
    \\
    =& ( y_t^{K', y_0} )^\top (K' - K)^\top \big( 2 E_K + (R + \gamma \rB^\top P_K^y \rB ) (K' - K) \big)  y_t^{K', y_0} 
\end{align*}
where the second equality is justified by $V_y(K, \xi_y) = \xi_y P_K^y \xi_y + \alpha_K^y$ for any $\xi_y \in \RR^d$. Then, the difference between the auxiliary costs $C_y$ at $K$ and $K'$ becomes
\begin{align*}
    C_y(K') - C_y(K) &= \EE_{y_0 \sim X_0 - \bar{X}_0}\big[ \sum \nolimits_{t\geq 0} \gamma^t f( y_t^{K', y_0}, K', Q, R) - V_y(K, y_0) \big]
    \\
    & = \EE_{y_0 \sim X_0 - \bar{X}_0} \big[ \sum \nolimits_{t \geq 0} \gamma^t A_y( t, K, y_t^{K', y_0}, K') \big] 
    \\
    & =  \langle \Delta K, \, 2 E_K \Sigma_{K'} + (R + \gamma \rB^\top P_K^y \rB )\Delta K  \Sigma_{K'} \rangle_{tr}
\end{align*}
which is the desired formula.
\end{proof}

\subsection{Solution of the discrete-time Lyapunov equations}
\label{sec:existence_of_DLE}

For any admissible parameters $\theta \in \Theta$, we construct solutions to the DLEs~\eqref{eq:lyapunov_eq_theta}.

\begin{proposition}
\label{prop:existence_of_lyapunov_equation}
Under Assumption~\ref{as:positivity-qr}, for any admissible control $\theta \in \Theta$, there exists a unique pair of matrices $(P^y_K, P^z_L)$ satisfying (DLEs)~\eqref{eq:lyapunov_eq_theta} and $P^y_K \succ 0$, $P^z_L \succ 0$. 
\end{proposition}

\begin{proof}
    For an admissible parameter $\theta = (K, L) \in \Theta$, the existence of a positive definite symmetric matrice $P_K^y$ (resp. $P^z_L$) as the unique solution to the corresponding Lyapunov equation in~\eqref{eq:lyapunov_eq_theta} is a direct application of~\cite{bof2018lyapunov}[Theorem 3.2] with a matrix $\sqrt{\gamma} (\rA - \rB K)$ (resp. $\sqrt{\gamma} (\tilde \rA - \tilde \rB L)$) and a positive definite matrix $Q + K^\top R K$ (resp. $\tilde Q + L^\top \tilde{R} L$). 
    Indeed, we can construct directly the solutions as: 
    \begin{equation*}
    \left\{
    \begin{array}{rcl}
       P_K^y &:=& \sum_{t=0}^\infty \gamma^t \big( (\rA - \rB K)^\top \big)^t ( Q + K^\top R K ) (\rA - \rB K)^t 
       \\
       P_L^z &:=& \sum_{t=0}^\infty \gamma^t \big( ( \tilde{\rA} - \tilde{\rB} L)^\top \big)^t ( \tilde{Q} + L^\top \tilde{R} L ) (\tilde{\rA} - \tilde{\rB} L)^t.
    \end{array}
    \right.
    \end{equation*}
    We easily deduce $P_K^y, P_L^z \succ 0$ from $Q, \tilde Q \succeq 0$ and $R, \tilde R \succ 0$.
\end{proof}

The following lemma completes the connections between the solution matrices $(P, \bar P)$ of the Riccati equations~\eqref{eq:Riccati} and the solution matrices of (DLEs)~\eqref{eq:lyapunov_eq_theta} with the specific parameters $K^*$ and $L^*$ defined in Theorem~\ref{thm:existence_linear_control}. As a result, we rediscover that the control parameter $\theta^* = (K^*, L^*)$ is admissible in the sense that $\theta^* \in \Theta$.

\begin{lemma}
We assume that Assumption~\ref{as:positivity-qr} holds, that the matrix Riccati equations~\eqref{eq:Riccati} admit solutions $(P, \bar P)$, and we define the two symmetric matrices $P^{*, y}$ and $P^{*,z}$ in $\RR^{d \times d}$ by
\begin{equation}
\label{eq:riccati_matrix_close_loop}
    P^{*,y} := ( \rA^\top P + P^\top \rA + 4Q ) / 4, \qquad P^{*,z} := ( \tilde{\rA}^\top \bar{P} + \bar{P}^\top \tilde{\rA}  + 4 \tilde Q ) / 4.
\end{equation}
Finally, we define the parameters $K^*$ and $L^*$ defined by
$$
    K^* = \frac{1}{2} R^{-1} \rB^\top P , \qquad L^* = \frac{1}{2} \tilde{R}^{-1} \tilde{\rB}^\top \bar{P}.
$$
Then, we have that $P^{*,y}$ and $P^{*,z}$ satisfy (DLEs)~\eqref{eq:lyapunov_eq_theta} for $K = K^*$ and $L = L^*$.
Moreover, the control parameter $\theta^* = (K^*, L^*) \in \Theta$.
\end{lemma}

\begin{proof}
    From the Riccati equation~\eqref{eq:Riccati} for $P$, we have
    $$
       2 R K^* =  \rB^\top P = \gamma \rB^\top ( \rA^\top P + 2 Q) (\rA - \rB K^{*})
    $$
    and then
    \begin{align*}
        \rA^\top P + P^\top \rA & = \gamma \rA^\top \big( \rA^\top P + 2 Q \big) (\rA - \rB K^*) + \gamma (\rA - \rB K^*)^\top \big(  P^\top \rA + 2 Q \big) \rA 
        \\ \
        & = \gamma (\rA - \rB K^*)^\top \big( \rA^\top P +  P^\top \rA + 4 Q \big) (\rA - \rB K^*) + 4 (K^*)^\top R K^*.
        \nonumber 
    \end{align*}
    Thus, $P^{*, y} = ( \rA^\top P + P^\top \rA + 4Q) / 4$ is a solution to the first equation in~\eqref{eq:lyapunov_eq_theta}. Moreover $P^{*, y} \succ 0$ because of $R \succ 0$. By applying similar arguments to the Riccati equation for $\bar P$ in~\eqref{eq:Riccati}, we see that $P^{*,z}$ satisfies the second equation in the (DLEs)~\eqref{eq:lyapunov_eq_theta}.

    To show $\theta^* \in \Theta$, we consider an auxiliary process with linear dynamics $\tilde{y}_{t+1} = \sqrt{\gamma} (\rA - \rB K^*) \tilde{y}_{t}$, and its associated value function $V: \RR^{d} \ni y \mapsto \sum_{t \geq 0} (\tilde{y}_t)^{\top} \big( Q + (K^*)^\top R K^* \big) \tilde{y}_t \in \RR \cup \{ \infty \}$ for the auxiliary process $(\tilde{y}_t)_{t \geq 0}$ starting from $\tilde y_0 = y$. We introduce the set $D = \{ y \in \RR^d \, | \, V(y) < \infty \} \subset \RR^{d}$ and we write the dynamics of the process $(\tilde{y}_t)_{t \geq 0}$ by $\tilde{y}_{t+1} = f(\tilde{y}_t)$ with $f : D \to \RR^{d}$. 
     Under Assumption~\ref{as:positivity-qr} with $Q \succeq 0$ and $R \succ 0$, the value function satisfies that $V (0) = 0$, $V(y) > 0$ for any $y \in D - \{ 0 \}$, and $V(f(y)) \leq V(y) < \infty$ for $y \in D$. So the system function $f$ is stable on the set $D$.
     Thus, the value function restricted to the set $D$, denoted by $V_D : D \ni y \mapsto V(y) \in \RR$, is a Lyapunov function for the linear system of $(\tilde y_t)_{t \geq 0}$ with initial status $\tilde{y}_0 = y \in D$. 
    
    Because the matrix $P^{*, y}$ from equation~\eqref{eq:riccati_matrix_close_loop} satisfies equation~\eqref{eq:lyapunov_eq_theta}, the value function $V_D$ can then be expressed into $V_D(y) = y^\top P^{*, y} y$ for any $y \in D$. As a result, if $\| y \| \to \infty$, we have $V_D(y) \rightarrow \infty$. Because $V_D(f(y)) < V_D(y) $ for any $y \in D - \{0\}$,~\cite[Theorem 1.4]{bof2018lyapunov} implies that the starting point $y=0$ is \emph{globally asymptotically stable}. As a result, we have $\gamma \| \rA - \rB K^* \|^2 < 1$ by~\cite[Theorem 3.1]{bof2018lyapunov}.
    Using similar arguments for $L^*$ and $P^{*,z}$, we have $\gamma \| \tilde \rA - \tilde \rB L^* \|^2 < 1$. We conclude that $\theta^* \in \Theta$.
\end{proof}

\subsection{Perturbation analysis}
\label{subsection:perturbation_analysis_of_exact_PG}

\subsubsection{Perturbation of the variance}
\label{subsection:perturbation_on_variance}

In this subsection, we show that when $\Delta K = K' - K$ is small enough, the norm between the corresponding variance matrices is controlled by $\|\Delta K \|$. To do so, we introduce a function $f_{var}$ defined by
\begin{equation}
    \label{eq:f_var_K}
    f_{var}(K) = 4 ( 1 + \| \rA - \rB K \| ) \| \rB \|  \frac{ C_y(K)} { \lambda_y^1 \lambda_{min}(Q) }.
\end{equation} 
We also introduce the linear operators $\cF_K$ and $\cT_K$ on symmetric matrices $\cS\cM$: for any $\Sigma \in \cS\cM$,: 
\begin{equation}
\label{eq:def-ope-F_K_y}
    \cF_K(\Sigma) = \gamma (\rA- \rB K) \Sigma (\rA-\rB K)^\top,
    \quad
    \cT_K(\Sigma) = \sum_{t \geq 0} \gamma^t (\rA - \rB K)^t \Sigma \big( (\rA - \rB K)^t \big)^\top.
\end{equation}
and we consider the operator norm for $\cF_K$ and $\cT_K$ associated with the spectral norm for symmetric matrices, denoted by $\vertiii{\cF_K}$ and $\vertiii{\cT_K}$. 
Because $\vertiii{\cF_K} = \sup_{\| \Sigma \| = 1} \| \cF_K(\Sigma) \| \leq \gamma \| \rA - \rB K \|^2 < 1$ for control parameters $K$ in $\theta = (K, L ) \in \Theta$, we then have
\begin{equation}
\label{eq:cT}
    \cT_K(\Sigma) = \sum \nolimits_{t \geq 0} (\cF_K)^t(\Sigma) = ( \bI - \cF_K)^{-1} (\Sigma).
\end{equation}

\begin{lemma} The variance matrix $\Sigma_K$ satisfies
\label{lemma:expression_Sigma_K_with_T_K}
    \begin{equation}
    \label{eq:link_vairance_y_with_T_K}
        \Sigma_K = \mathcal{T}_K \big(\Sigma_{y_0} +  \gamma / (1- \gamma) \Sigma^1 \big)
    \end{equation}
    Moreover, $\lambda_{min} (\Sigma_K) \geq \lambda_y^1$.
\end{lemma}

\begin{proof} 
Expanding the variance matrix $\Sigma_K$~\eqref{eq:variance_matrices_y_and_z} using $y_{t+1}^{K, y_0} = (\rA - \rB K) y_t^{K, y_0} + \varepsilon_{t+1}$ for $t \geq 0$, we obtain 
\begin{align*}
    \Sigma_K 
    & = \EE [ y_0 y_0^\top] + \EE \big[ \sum_{t \ge 1} \gamma^t \varepsilon_{t} (\varepsilon_{t})^\top \big] + \EE \big[ \sum_{t \ge 0} \gamma^{t+1} (\rA - \rB K) y_{t}^{K, y_0} (y_t^{K, y_0})^\top (\rA - \rB K)^\top \big]
    \\
    & 
    = \Sigma_{y_0} + \gamma / (1 - \gamma) \Sigma^1 + \cF_K ( \Sigma_K),
\end{align*}
and we conclude because $\cF_K(\Sigma_K) \succeq 0$.
\end{proof}

The following lemma provides bounds on the operator norms $\cF_k$ and $\cT_k$.

\begin{lemma} 
\label{lemma:useful-bounds-ope-T-F} 
Under  Assumption~\ref{as:non-deg}, if $\Delta K = K' - K $, we have
\begin{enumerate}[label=(\roman*)]
    \item   \label{eq:TyK-Ctheta-lambda}
    $      
        \| \Sigma_K \| / \| \Sigma_{y_0} + \gamma / (1-\gamma) \Sigma^1 \| \leq  \vertiii{\cT_K} \leq \| \Sigma_K \| /  \lambda^1_y.
    $
    
    \item \label{eq:perturbation_F_K^y}
        $
            \vertiii{\cF_{K'} - \cF_{K} } \leq \gamma \big( 2 \| \rA - \rB K \| \| \rB \| \| \Delta K \| + \| \rB \|^2 \| \Delta K \|^2 \big) .
        $
        
    \item \label{eq:perturbation_T_K^y_in_T_K^y} 
        For any $\Sigma \in \cM\cS$ and any $\eta < 1$ such that $\vertiii{ \cT_K} \vertiii{ \cF_{K'} - \cF_{K} } \leq \eta$, we have
        \begin{equation*}
            \| (\cT_{K'} - \cT_{K}) (\Sigma) \|
            \leq \frac{1}{1-\eta} \vertiii{ \cT_K} \vertiii{\cF_{K'} - \cF_{K} } \Vert \cT_{K}(\Sigma) \Vert 
            \leq \frac{\eta}{1- \eta} \Vert \cT_K (\Sigma) \Vert.
        \end{equation*}
   
    \end{enumerate}
\end{lemma}

\begin{proof} 
Inequality \ref{eq:TyK-Ctheta-lambda} is due to the definition of $\vertiii{\cT_K}$ and Lemma~\ref{lemma:expression_Sigma_K_with_T_K}. To show the second inequality, we notice that for any $U \in \cS\cM$ and $\Sigma \succ 0$, we have
\begin{align*}
    \| \cT_K(U) \| 
    &= \sup_{ \{ \xi \in \RR^d : \| \xi \| = 1 \} } \big| \xi^\top \cT(U) \xi \big|
    \\
    &= \sup_{ \{ \xi \in \RR^d : \| \xi \| = 1 \} } \Big| \sum \nolimits_{t \geq 0} \gamma^t Tr\Big( \xi^\top  (\rA - \rB K)^t  U \big( (A - B K )^t \big)^\top \xi \Big) \Big|
    \\
    &= \sup_{ \{ \xi \in \RR^d: \| \xi \| = 1 \} } \Big| \sum \nolimits_{t \geq 0} \gamma^t \big\langle (\Sigma^{1/2})^\top \big( (\rA-\rB K)^t \big)^\top \xi \xi^\top (\rA - \rB K)^t \Sigma^{1/2}, \,  \Sigma^{-1/2} U \Sigma^{-1/2}  \big\rangle_{tr} \Big|
    \\
    & \leq \sup_{\{ \xi \in \RR^d: \| \xi \| = 1 \} }  \sum \nolimits_{t \geq 0} \gamma^t Tr \Big(  (\Sigma^{1/2})^\top \big( (\rA-\rB K)^t \big)^\top \xi \xi^\top(\rA - \rB K)^t \Sigma^{1/2} \Big). \big\| \Sigma^{-1/2} U \Sigma^{-1/2} \big\|
    \\
    & = \sup_{\{ \xi \in \RR^d: \| \xi \| = 1 \} } \big|  \xi^\top \cT_K(\Sigma) \xi \big| \| \Sigma^{-1/2} U \Sigma^{-1/2} \|
    \\
    & = \| \cT_K(\Sigma) \|\|  U \Sigma^{-1} \|.
\end{align*}
where the fourth inequality above is justified by the Von Neumann's trace inequality stating that for $V \succeq 0$, 
$$
| \langle V, \Sigma^{-1/2} U \Sigma^{-1/2} \rangle_{tr} | \leq \big(\sum_{i=1}^d  \sigma_i(V) \big) \| \Sigma^{-1/2} U \Sigma^{-1/2} \| = Tr(V) \| \Sigma^{-1/2} U \Sigma^{-1/2} \|,
$$
and the last equality is due the definition of $\| \cT_K (\Sigma) \|$ and $\| \Sigma^{-1/2} U \Sigma^{-1/2} \| = \| U \Sigma^{-1} \|$.
Under Assumption~\ref{as:non-deg}, we choose $\Sigma = \Sigma_{y_0} + \gamma / (1-\gamma) \Sigma^1$ so that $\Sigma \succ 0$. Together with Lemma~\ref{lemma:expression_Sigma_K_with_T_K}, we obtain
$$
\vertiii{\cT_K} = \sup_{\| U \| = 1} \| \cT_K(U) \| \leq \| \cT_K \big( \Sigma) \|  \| \Sigma^{-1} \|   \leq \frac{ \| \Sigma_K \|} { \lambda_{min}(\Sigma) } = \frac{\| \Sigma_K \| }{\lambda_y^1}.
$$
Inequality \ref{eq:perturbation_T_K^y_in_T_K^y} is adapted from~\cite{fazel2018global}[Lemma 20] and we generalize it with coefficient $\eta < 1$.
\end{proof}

Now, we are ready to prove a perturbation result for the variance matrix. 

\begin{lemma}
\label{lemma:perturbation_variance_Cy}
Consider $\theta, \theta' \in \Theta$. If $ f_{var}(K) \| K' - K \| \leq 1$, then 
    \begin{equation}
    \label{eq:perturbation_variance}
        \| \Sigma_{K'} - \Sigma_{K} \| \leq  f_{var}(K) \| \Sigma_K \| \| K' - K \|  \leq \| \Sigma_{K} \|
    \end{equation}
\end{lemma}

\begin{proof}
    Let $\Delta K = K' - K$. Because 
$$
\|\rB\| \|\Delta K \| \leq \|\rB \| f_{var}(K)^{-1} \leq 1/4,
$$
Lemma~\ref{lemma:useful-bounds-ope-T-F}\ref{eq:TyK-Ctheta-lambda}\ref{eq:perturbation_F_K^y} and inequality~\eqref{eq:upper_bound_Riccati_and_Variance_matrix}, we have 
\begin{align*}
        \vertiii{\cT_K} \vertiii{\cF_{K'} - \cF_K} 
        &\leq \frac{ C_y(K) }{\lambda_{min}(Q) \lambda_y^1}. \gamma \big( 2 \| \rA - \rB K \|  + \frac{1}{4} \big) \| \rB\| \| \Delta K \|\\
        &\leq f_{var}(K) \|\Delta K \| \gamma /2 \leq 1/2.    
\end{align*}
    Then, by Lemma~\ref{lemma:useful-bounds-ope-T-F}~\ref{eq:perturbation_T_K^y_in_T_K^y} with $\eta = 1/2$, we obtain that
    \begin{align*}
        \| \Sigma_{K'} - \Sigma_K \| = \| ( \cT_{K'} - \cT_K) (\Sigma_{y_0} + \gamma / (1- \gamma) \Sigma^1 ) \| 
        & \leq 2 \vertiii{\cT_K} \vertiii{\cF_{K'} - \cF_K} \| \Sigma_K \|
        \\ 
        & \leq  f_{var}(K) \|\Delta K \| \| \Sigma_K \|
    \end{align*}
    which completes the proof.
\end{proof}

\subsubsection{Perturbation of the DLE solution}
\label{subsection:pertubation_on_DLE}

Consider the solutions $P_K^y$ and $P_{K'}^y$ of the first DLE equation in~\eqref{eq:lyapunov_eq_theta} for parameters $K$ and $K'$. We introduce a function $f_{riccati}$ defined by
\begin{equation}
\label{eq:f_riccati_K}
    f_{riccati}(K, \Delta K) = \frac{ C_y(K) }{\lambda_y^1} \bigg( \gamma f_{var}(K) + \frac{4\|R\|}{\lambda_{min}(Q) } \| K \| + \frac{ 2 \| R \|}{ \lambda_{min}(Q)} \| \Delta K \| \bigg)   
\end{equation}

\begin{lemma} 
    \label{lemma:perturbation_riccati_Cy}
    Consider $\theta, \theta' \in \Theta$. If $f_{var}(K) \| K' - K \| \leq 1$, then
    \begin{equation}
    \label{eq:perturbation_DLE_solution}
        \| P_{K'}^y - P_K^y \|  \leq  f_{riccati}(K, \Delta K ) \| K' - K \|.
    \end{equation}
\end{lemma}
\begin{proof} From the DLE~\eqref{eq:lyapunov_eq_theta} for $P_K^y$, we have $P_K^y = Q + K^\top R K + \cF_K(P_K^y)$, which implies that $P_K^y = \cT_K( Q + K^\top R K )$. Then by Lemma~\eqref{lemma:useful-bounds-ope-T-F}\ref{eq:perturbation_T_K^y_in_T_K^y} with $\eta = 1/2$, we have
    \begin{align*}
        \| P_{K'}^y - P_K^y \| & = \| (\cT_{K'} - \cT_K)(Q + K^\top R K )  + \cT_{K'} \big( (K')^\top R K' - K^\top R K\big) \|
        \\
        & \leq  \gamma f_{var}(K)  \| \Delta K \| \| P_K^y \| +  \vertiii{\cT_{K'}} \| R \|  ( 2 \| K \|\|\Delta K \| + \| \Delta K \|^2 ) 
    \end{align*}
    Lemma~\ref{lemma:useful-bounds-ope-T-F}~\ref{eq:TyK-Ctheta-lambda} and Lemma~\ref{lemma:perturbation_variance_Cy} imply $\vertiii{\cT_{K'}} \leq \| \Sigma_{K'} \|/ \lambda_y^1 \leq 2 \| \Sigma_K \| / \lambda_y^1$. Together with inequalities~\eqref{eq:upper_bound_Riccati_and_Variance_matrix} for $\| P_K^y \|$ and $\| \Sigma_K \|$, we conclude the lemma.
\end{proof}

\subsubsection{Perturbation of the gradient}
\label{subsection:perturbation_on_gradient}

We show that the difference between gradients $\nabla_K C_y$ computed at $K$ and $K'$ is controlled by $\| K' - K \|$. To do so, we consider the function 
\begin{equation}
    \label{eq:f_grad_K}
    f_{grad}(K) = \Big( 6 \| R \| +  \frac{ 4 \| \rB\|^2 C_y(K)}{\lambda_y^1}  + 6 f_{var}(K) \| R \| \| K \| + \frac{3}{2} ( f_{var}(K))^2 \lambda_{min}(Q) \Big)\frac{C_y(K)}{\lambda_{min}(Q)}
\end{equation}

\begin{lemma}
\label{lemma:perturbation_gradient_Cy}
    Consider $\theta, \theta' \in \Theta$. If $f_{var}(K) \| K' - K \| \leq 1$, then
    \begin{equation}
    \label{eq:perturbation_gradient}
        \| \nabla_K C_y(K') - \nabla_K C_y(K) \| \leq f_{grad}(K) \| K' - K \|.
    \end{equation}
\end{lemma}

\begin{proof}
    We use the formula $\nabla_K C(K) = 2 E_K \Sigma_K$ proven in Lemma~\ref{lemma:policy_gradient_expression}.
    $$
        \| \nabla_K C(K') - \nabla_K C(K) \| \leq 2 \| E_{K'} - E_K \|. \|\Sigma_{K'} \| + 2 \| E_{K} \|. \| \Sigma_{K'} - \Sigma_K \| 
    $$
    and we bound each term in the above right hand side. First, with inequality~\eqref{eq:upper_bound_Riccati_and_Variance_matrix} on $\| P_K^y \|$, we notice that 
    $$ 
        \| E_K \| = \| RK - \gamma \rB^\top P_K^y (\rA - \rB K) \| \leq \| R \|\| K \| + \gamma f_{var}(K) \lambda_{min}(Q) / 4.
    $$
    Next, using Lemma~\ref{lemma:perturbation_riccati_Cy} and $\| \rB \| \| \Delta K \| \leq f_{var}(K) \| \Delta K \| \leq 1$, we have 
    \begin{align*}
        \| E_{K'} - E_K \| & = \| R \Delta K - \gamma \rB^\top ( P_{K'}^y - P_K^y) ( \rA - \rB K - \rB \Delta K ) + \gamma \rB^\top P_K^y \rB \Delta K \| 
        \\
        & \leq \big( \| R \| + \gamma \| \rB \|  ( 1 + \| \rA - \rB K \| ) f_{riccati}(K, \Delta K ) + \gamma \| \rB \|^2 \| P_K^y \| \big) \| \Delta K \|
    \end{align*}
    By expanding the term $f_{riccati}(K, \Delta K)$ in~\eqref{eq:f_riccati_K}, we obtain
    \begin{align*}
        & \gamma \| \rB \| ( 1 + \| \rA - \rB K\| ) f_{riccati}(K, \Delta K)  
        \\
        &\hskip 35pt
        =\frac{\gamma^2}{4} ( f_{var}(K) )^2 \lambda_{min}(Q) + \gamma f_{var}(K) \| R \| \| K \| + \frac{\gamma \| R \|}{2} f_{var}(K) \| \Delta K \|. 
    \end{align*}
   Using Lemma~\ref{lemma:perturbation_variance_Cy} to control  $\| \Sigma_{K'} - \Sigma_K \|$, we get:
    \begin{align*}
         \| \nabla_K C(K') - \nabla_K C(K) \| \leq &  2 \Big(  \|R\|  + \gamma^2 ( f_{var}(K) )^2 \lambda_{min}(Q) / 4 + \gamma f_{var}(K) \| R \| \| K \| 
         \\
         & \hspace{10pt} + \gamma \| R \| /2 + \gamma \| \rB \|^2 C_y(K) / \lambda_y^1 \Big) \| \Delta K \|. \big( 2 \| \Sigma_K \| \big)
         \\
         &  + 2 \Big( \| R \| \| K \| + \gamma f_{var}(K) \lambda_{min}(Q) / 4 \Big)  f_{var}(K) \| \Sigma_K\|. \| \Delta K \| 
         \\
         \leq & f_{grad}(K) \| \Delta K \|.
    \end{align*}
    which is the desired result.
\end{proof}

\subsection{Admissibility of perturbed control parameters}
\label{subsection:stability_of_small_perturbation}
Here we show that when the perturbations $\| \Delta K \| = \| K' - K \|$ and $\| \Delta L \| = \| L' - L \|$ are small enough, the perturbed parameter $\theta' = (K', L')$ is still admissible, i.e. $\theta' \in \Theta$. 

\begin{lemma} 
\label{lemma:stability_of_small_perturbation}
Consider $K \in \RR^{\ell \times d}$ such that $\gamma \| \rA - \rB K \|^2 < 1$. If a matrix $K' \in \RR^{\ell \times d}$ satisfies $f_{var}(K)\| K' - K \| \leq 1$, then 
$
        \gamma \| \rA - \rB K' \|^2 < 1.
$
\end{lemma}
\begin{proof}
The set $\cB_{K} = \{ \tilde K \in \RR^{\ell \times d} \, | \, \| \tilde K - K \| \leq f_{var}(K)^{-1} \}$ is convex. 
     For any matrix $\tilde K \in \RR^{\ell \times d}$, let $\rA - \rB \tilde K = U M U^{-1}$ be the Schur decomposition of square matrix $\rA -\rB \tilde K \in \RR^{d \times d}$, where $U$ is a unitary matrix and $M$ is an upper triangular matrix. The diagonal elements of $M$, denoted by $m_{ii} = \lambda_{i}$ for $i = 1, \ldots, d$, are the eigenvalues $(\lambda_1, \ldots, \lambda_d)$ of the matrix $\rA - \rB \tilde K$.
     For every $t \geq 1$ and $i=1, \ldots, d$, the $i-$th diagonal term of an upper triangular matrix $M^t$ satisfies $(M^t)_{ii} = (M_{ii})^t = \lambda_{i}^t$, Consequently
     $$
        \| ( \rA - \rB \tilde K )^t \|^2_F = \| U M^t U^{-1} \|_F^2 = \| M^t \|_F^2 \geq \sup_{i=1,\ldots, d} | (M^t)_{ii} |^2 = \sup_{i=1,\ldots d} (\lambda_{i}^t)^2 = \| \rA - \rB K \|^{2t}.
     $$
     Hence, for any matrix $\tilde{K} \in \RR^{\ell \times d}$ satisfying $\gamma \| \rA - \rB \tilde K \|^2 < 1$, Lemma~\ref{lemma:expression_Sigma_K_with_T_K} implies that
    \begin{align*}
        Tr(\Sigma_{\tilde K} ) & = \sum_{t \geq 0} Tr\Big( \cF_K^t \big(\Sigma_{y_0} + \frac{\gamma}{1-\gamma} \Sigma^1 \big) \Big)
        \\
        & = \sum_{t \geq 0} \gamma^t Tr \Big( (\rA - \rB \tilde K)^t \big(\Sigma_{y_0} + \frac{\gamma}{1-\gamma} \Sigma^1 \big). \big( (\rA - \rB \tilde K )^\top \big)^t \Big)
        \\
        & \geq \sum_{t \geq 0} \gamma^t \lambda_{min}\big(\Sigma_{y_0} + \frac{\gamma}{1-\gamma} \Sigma^1 \big) Tr \Big( \big( (\rA - \rB \tilde K )^\top \big)^t. (\rA - \rB \tilde K)^t \Big)
        \\
        & \geq \lambda_y^1 \sum_{t \geq 0} \gamma^t \| ( \rA - \rB \tilde K )^t \|_F^2
        \\
        & \geq \lambda_y^1 \sum_{t \geq 0} \gamma^t \| \rA - \rB \tilde K \|^{2t}
        \\
        & = \frac{\lambda_y^1}{ 1 - \gamma \| \rA - \rB \tilde K \|^2 }.
    \end{align*} 
     Moreover, if $\tilde K \in \cB_K$, Lemma~\ref{lemma:perturbation_variance_Cy} implies that $\| \Sigma_{\tilde K} - \Sigma_K \| \leq \| \Sigma_K \|$. Together with~\eqref{eq:upper_bound_Riccati_and_Variance_matrix} we have $Tr ( \Sigma_{\tilde K } ) \leq d \| \Sigma_{\tilde K} \| \leq 2 d \| \Sigma_K \| \leq 2 d C_y(K) / \lambda_{min}(Q)$. Consequently, if $\tilde K \in \cB_K$ and $\gamma \| \rA - \rB \tilde K \|^2 < 1$, we also have that 
    $$
        \gamma \| \rA - \rB \tilde K \|^2 \leq  1 - \frac{ \lambda_y^1 }{ Tr( \Sigma_{\tilde K} ) } \leq 1 - \frac{\lambda_y^1 \lambda_{min}(Q) }{ 2 d C_y(K) } = 1 - \zeta
    $$
    where $\zeta = ( \lambda_y^1 \lambda_{min}(Q) ) / ( 2 d C_y(K) )$ is independent of $\tilde K$. 
    
    In the following, we show that $\tilde K \in \cB_K$ implies $\gamma \| \rA - \rB \tilde K \|^2 < 1$. Otherwise, suppose there exists $K' \in \cB_K$ with $\gamma \| \rA - \rB K' \|^2 \geq 1$. Consider $f_K(\lambda) = 1 - \gamma \| \rA - \rB \big( (1- \lambda) K + \lambda K' \big) \|^2$. $f_K$ is continuous because the norm $\| . \|$ is continuous. Besides, we have $f_K(0) = \gamma \|\rA - \rB K \|^2 \leq 1 - \zeta < 1 - \zeta /2 $ and $f_K(1) = \gamma \| \rA - \rB K' \|^2 > 1 - \zeta / 2$. Thus, by continuity of $f_K$, there exists $\lambda^{\zeta} \in [0,1]$ such that $f_K(\lambda^{\zeta}) = 1 - \zeta / 2 < 1$. By convexity of $\cB_K$, $K^{\lambda^\zeta} = (1 - \lambda^{\zeta}) K + \lambda^{\zeta} K' \in \cB_K$, so that 
    $$ 1 - \zeta / 2 = f_K(\lambda^\zeta) =  1 - \gamma \| \rA - \rB K^{\lambda^\zeta} \|^2 \leq 1 - \zeta, $$
    which leads to a contradiction.
\end{proof}

\subsection{Proof of Proposition~\ref{proposition:PL_cond_Cy_Cz}}
\label{subsection:proof_proposition_PL_condition}

\begin{proof} As in \cite[Lemma 11]{fazel2018global}, we first notice that for real values $a, b, c \in \RR$ with $c \neq 0$, $2ab + a^2 c = (a + b/c)^2 c - b^2 / c \geq - b^2 / c$, and this inequality still holds for matrices with appropriate dimensions. For $t \geq 0$, let $a = (K' - K) y_t^{K', y_0}$, $b=E_K y_t^{K', y_0} $, $c=R + \gamma \rB^\top P_K^y \rB$. The \emph{advantage} quantity $A_y(t, K, y_t^{K', y_0}, K')$ in Lemma~\ref{lemma:difference_in_C_y(K)} satisfies
$$
    A_y(t,K, y_t^{K', y_0}, K')  \geq - ( E_K y_t^{K', y_0})^\top (R + \rB^\top P_K^y \rB)^{-1} E_K y_t^{K', y_0}.
$$
Then, using Lemma~\ref{lemma:policy_gradient_expression}, we deduce that
\begin{align*}
    C_y(K) - C_y(K^*) &= - \EE \Big[ \sum_{t \geq 0} \gamma^t A_y(t, K, y_t^{K^*, y_0}, K^*) \Big] 
    \\
    &\leq \| (R + \gamma \rB^\top P_K^y \rB)^{-1} \| \| \Sigma_{K^*} \|.  \| E_K \|_F^2
    \\
    &\leq  \frac{\Vert \Sigma_{K^*} \Vert}{4 \lambda_{min}(R) \lambda_{min}(\Sigma_K)^2} \| \nabla_K C_y(K) \|_F^2.
\end{align*}
Because $\lambda_y^1 \leq \lambda_{min}( \Sigma_K)$ (see Lemma~\ref{lemma:expression_Sigma_K_with_T_K}),  we conclude that $C_y$ satisfies the $\nu_y-$PL condition for $\nu_y = 4 \lambda_{min}(R). (\lambda_y^1)^2  / \| \Sigma_{K^*} \| $.
\end{proof}

Notice that we can also derive a lower bound for $C_(y)(K) - C_y(K^*)$. Indeed, if we consider $K' = K - (R + \rB^\top P_K^y \rB)^{-1} E_K$ so that the \emph{advantage} quantities in Lemma~\ref{lemma:difference_in_C_y(K)} attain their lower bounds, then with $R \succ 0$ and $P_K^y \succ 0$ (Assumption~\ref{as:positivity-qr} and Proposition~\ref{prop:existence_of_lyapunov_equation}), we have 
\begin{align}
    C_y(K) - C_y(K^*) \geq C_y(K) - C_y(K') &= \langle (R + \gamma \rB^\top P_K^y \rB)^{-1} E_K, \,  E_K \Sigma_{K'}  \rangle_{tr}
    \nonumber \\
    & \geq  \lambda_y^1 \| E_K \|_F^2  / \| R + \gamma \rB^\top P_K^y \rB \|.
    \label{eq:lower_bound_for_Cy_minus_Copt}
\end{align}

\subsection{Proof of Proposition~\ref{proposition:local_smoothness_Cy_Cz}}
\label{subsection:proof_proposition_local_smoothness}

\begin{proof} Consider $\Delta K = K' - K = - \eta \nabla_K C_y(K)$ with $\eta > 0$. By Lemma~\ref{lemma:policy_gradient_expression} on the gradient $\nabla_K C_y(K)$, we have 
\begin{align*}
   \langle \Delta K , 2 E_K \Sigma_{K'} \rangle_{tr}  & =  \langle\Delta K,  \nabla_K C_y(K) + 2 E_K (\Sigma_{K'} - \Sigma_K) \rangle_{tr}
    \\
    & = \langle \Delta K,  \nabla_K C_y(K) \rangle_{tr} -  \langle 2 \eta E_K \Sigma_K, 2E_K (\Sigma_{K'} - \Sigma_K) \rangle_{tr}
     \\
     & \leq \langle \Delta K,  \nabla_K C_y(K) \rangle_{tr} + 4\eta \big| \langle (E_K)^\top E_K \Sigma_K, \Sigma_{K'} - \Sigma_K \rangle_{tr} \big|
     \\
     & \leq \langle \Delta K,  \nabla_K C_y(K) \rangle_{tr} + 4\eta \langle (E_K)^\top E_K \Sigma_K, \Sigma_K \rangle_{tr} \frac{\| \Sigma_{K'} - \Sigma_{K} \|}{\lambda_{min}(\Sigma_K)}
     \\
     & \leq \big( 1 - \lambda_{var, y}(K, K') \big)  \langle \Delta K, \nabla_K C_y(K) \rangle_{tr}.
\end{align*}
where $\lambda_{var, y}(K, K') = \| \Sigma_{K'} - \Sigma_K \| / \lambda_y^1$. Then, together with Lemma~\ref{lemma:difference_in_C_y(K)}, we obtain the local smoothness inequality~\eqref{eq:local_smoothness_Cy_Cz} for $C_y$ in Proposition~\ref{proposition:local_smoothness_Cy_Cz}.
\end{proof}

\subsection{Perturbation analysis bounds}
\label{subsection:bound_on_local_smoothness_coefficient}

The following bounds on the coefficients $\lambda_{var, y}$ and $\lambda_{hess,y}$ introduced in Proposition~\ref{proposition:local_smoothness_Cy_Cz} are needed in the convergence analysis of the exact PG descent algorithm~\ref{thm:exact-CV}.

\begin{lemma}
\label{lemma:useful_bounds_I}
    Consider $\theta \in \Theta$ and $C(\theta) \leq C_0$  for some constant $C_0 \in \RR$. We define the following constants
    \begin{align}
        & h_{small-pert, y}(C_0) = 4 \|\rB \| \big( 1 + 1 / \sqrt{\gamma} \big) \frac{C_0}{\lambda_y^1 \lambda_{min}(Q)}
        \label{eq:bound_f_var_K}
        \\
        & h_{grad,y}(C_0)  = 2 \frac{C_0}{\lambda_{min}(Q)} \sqrt{\Big( \| R \| + \gamma \| \rB \|^2 \frac{C_0}{ \lambda_y^1} \Big) \frac{C_0}{\lambda_y^1}}
        \label{eq:bound_grad_Cy_K}
        \\
        & h_{K}(C_0)  = \frac{1}{\lambda_{min}(R)} \Big(   \sqrt{\Big( \| R \| + \gamma \| \rB \|^2 \frac{C_0}{ \lambda_y^1} \Big) \frac{C_0}{\lambda_y^1}} + \gamma \|\rB \| \| \rA \| \frac{C_0}{\lambda_y^1} \Big)
        \label{eq:bound_K}
    \end{align}
    Then we have
    \begin{equation}
    \label{eq:inequality_useful_bounds_I}
        f_{var}(K) \leq h_{small-pert,y}(C_0), \quad \| \nabla_K C_y(K) \| \leq h_{grad,y}(C_0), \quad \| K \| \leq h_{K}(C_0).
    \end{equation}
    Similar functions $h_{small-pert, z}(C_0)$, $h_{grad, z}(C_0)$ and $h_{L}(C_0)$ can be defined by replacing $(\rB, Q, R, \lambda_y^1)$ with $(\tilde \rB, \tilde Q, \tilde R, \lambda_z^0)$ in~\eqref{eq:bound_f_var_K},~\eqref{eq:bound_grad_Cy_K}, and ~\eqref{eq:bound_K} and bounds analog to \eqref{eq:inequality_useful_bounds_I} hold true.
\end{lemma}
\begin{proof}
    By definition of $\theta \in \Theta$, we have $\| \rA - \rB K \| \leq 1 / \sqrt{\gamma}$, then equation~\eqref{eq:f_var_K} for $f_{var}(K)$ implies the first inequality in~\eqref{eq:inequality_useful_bounds_I}. From equation~\eqref{eq:lower_bound_for_Cy_minus_Copt}, we have $\| E_K \|^2 \leq \| E_K \|_F^2 \leq ( \|R\| + \gamma \|\rB\|^2 \| P_K^y \|) C_y(K) / \lambda_y^1$. Together with the bounds on $\| P_K^y \|$ and $\| \Sigma_K \|$ in~\eqref{eq:upper_bound_Riccati_and_Variance_matrix} and the gradient expression from Lemma~\ref{lemma:policy_gradient_expression}, we derive the second inequality in~\eqref{eq:inequality_useful_bounds_I}. Moreover, we notice that
    $$
        \| K \| = \| (R + \gamma \rB^\top P_K^y \rB)^{-1} ( R + \gamma \rB^\top P_K^y \rB) K \| \leq \frac{1}{\lambda_{min}(R)} ( \| E_K \| + \gamma \| \rB^\top P_K^y \rA \| ),
    $$
    from which we obtain the third inequality in~\eqref{eq:inequality_useful_bounds_I}.
\end{proof}

\begin{lemma}
\label{lemma:bound_on_coefficent_in_prop_local_smoothness}
    Consider $\theta = (K, L) \in \Theta$ such that $C_y(K) \leq C(\theta) \leq C_0$ with a constant $C_0 \in \RR$. Consider the following two constants in $\RR$ defined by
    \begin{align}
       h_{var, y}(C_0) & = 8 \frac{ \big( 1 + 1 / \sqrt{\gamma} \big) \|\rB \| ( C_0 )^3 }{ (\lambda_y^1)^2 (\lambda_{min}(Q) )^3} \sqrt{\Big(\| R \| + \gamma \| \rB \|^2 \frac{C_0}{\lambda_y^1} \Big) \frac{C_0}{\lambda_y^1}} 
        \label{eq:h_var_y}
        \\
        h_{hess, y}(C_0) & = 2 \frac{C_0}{\lambda_{min}(Q)}\Big( \| R \| + \gamma \| \rB \|^2 \frac{C_0}{\lambda_y^1} \Big).
        \label{eq:h_hess_y}
    \end{align}
    Then
    $$
        f_{var}(K) \| \nabla_K C_y(K) \|\frac{ \| \Sigma_K \| }{\lambda_y^1} \leq h_{var, y}(C_0).
    $$
    Moreover, if $K' = K - \eta \nabla_K C_y(K)$ for some $\eta > 0$, if $\eta \leq h_{var, y}^{-1}(C_0)$, then 
    \begin{align*}
        & \| K' - K \| f_{var}(K) \leq 1
        \\
        & \gamma \| A - B K' \|^2 < 1
        \\
        & \lambda_{var, y}(K, K') := \| \Sigma_{K'} - \Sigma_K \| / \lambda_y^1  \leq \eta h_{var, y}(C_0)
        \\
        & \lambda_{hess, y}(K, K') := \| \Sigma_{K'} \|. \big( \|R \| + \gamma \|B \|^2 C_y(K) / \lambda_y^1 \big) \leq h_{hess, y}(C_0).
    \end{align*}
\end{lemma}
\begin{proof}
    From inequality~\eqref{eq:bound_f_var_K} and~\eqref{eq:bound_grad_Cy_K} and $ 1 \leq \frac{ \| \Sigma_K \| }{ \lambda_y^1} \leq \frac{ C_0 }{ (\lambda_y^1 \lambda_{min}(Q)) }$, we have 
    $$
        f_{var}(K)\| \nabla_K C_y(K) \| \frac{ \| \Sigma_K \| }{\lambda_y^1}  \leq h_{small-pert,y}(C_0) h_{grad}(C_0) \frac{C_0}{\lambda_y^1 \lambda_{min}(Q)} = h_{var, y}(C_0).
    $$
    For such a  $K'$, we have 
    $$
        \| K' - K \| f_{var}(K) = \eta f_{var}(K) \| \nabla_K C_y(K) \| \leq \eta h_{var, y}(C_0) \leq 1.
    $$
    Then, Lemma~\ref{lemma:stability_of_small_perturbation} implies that $\gamma \| A - B K' \|^2 < 1$, and Lemma~\ref{lemma:perturbation_variance_Cy} implies
    $$ 
        \| \Sigma_{K'} - \Sigma_K \| \leq f_{var}(K) \| \Sigma_K \| \| \Delta K \| \leq \| \Sigma_K \|.
    $$
    As a result, 
    $$
        \lambda_{var, y}(K, K') \leq  \eta  f_{var}(K) \| \nabla_K C_y(K) \| \frac{ \| \Sigma_K \|}{\lambda_y^1} \leq \eta h_{var, y}(C_0)
    $$
    and
    $$
        \lambda_{hess, y}(K, K') \leq \big( \| \Sigma_K \| + \| \Sigma_{K'} - \Sigma_K \| \big). \big( \|R \| + \gamma \|B \|^2 C_y(K) / \lambda_y^1 \big) \leq h_{hess, y}(C_0),
    $$
    which completes the proof.
\end{proof}

\subsection{Perturbations of the cost function}
\label{section:perturbation_on_cost_function}

\begin{lemma}
\label{lemma:perturbation_cost_Cy}
    Consider $\theta \in \Theta$ and $C(\theta) \leq C_0$ for some $C_0 \in \RR$. Let
    \begin{align}
        \label{eq:h_raccati}
        h_{riccati, y}(C_0) & = \frac{C_0}{\lambda_y^1} \Big( \gamma h_{small-pert,y}(C_0) + \frac{4 \| R \|}{\lambda_{min}(Q)} h_{K}(C_0) + \frac{2 \| R \| }{\| \rB \| \lambda_{min}(Q) } \Big)
        \\
        h_{cost, y}(C_0) & = \frac{2d\; h_{riccati, y}(C_0)\; C_{init, noise}^2}{1- \gamma} 
        \label{eq:h_cost_y}
        \\
        h_{f\_grad\_K}(C_0) & = \frac{C_0}{\lambda_{min}(Q)}\big( 6 \| R \| ( 1+ h_{small-pert, y}(C_0) h_K(C_0) )
        \nonumber \\
        & \hspace{60pt} + 4 \| \rB \|^2 C_0 / \lambda_y^1 + 3h_{small-pert,y}^2(C_0) \lambda_{min}(Q) / 2 \big)
        \label{eq:h_f_grad_K}
    \end{align}
    If $f_{var}(K) \| K' - K \| \leq 1$, then
    \begin{align*}
         | C_y(K') - C_y(K) | &\leq h_{cost, y}(C_0) \| K' - K \|, 
         \\
         \| \nabla_K C_y(K') - \nabla_K C_y(K) \| & \leq h_{f\_grad\_K}(C_0) \| K' - K \|.
    \end{align*}
    Similarly, one can define constants $h_{riccati, z}(C_0)$, $h_{cost, z}(C_0)$, $h_{f\_grad\_L}(C_0)$ with $(\lambda_z^0, \tilde \rB, \tilde Q, \tilde R)$ and with $h_{small-pert, z}(C_0), h_L(C_0)$ from Lemma~\ref{lemma:useful_bounds_I}, and obtain similar upper bounds. 
    Let $h_{cost}(C_0) = h_{cost, y}(C_0) + h_{cost, z}(C_0)$ and $\| \theta' - \theta \| = \| K' - K \| + \| L' - L \|$, then 
    \begin{equation}
    \label{eq:perturbation_cost}
        | C(\theta') - C(\theta) | \leq h_{cost}(C_0) \| \theta' - \theta \|
    \end{equation}
\end{lemma}

\begin{proof}
    From the cost expression~\eqref{eq:cost_expression_with_PK_y} $C_y(K) = \langle P_K^y, \Sigma_{y_0} + \frac{\gamma}{1-\gamma} \Sigma^1 \rangle_{tr}$, and from Lemma~\ref{lemma:perturbation_riccati_Cy} on the perturbation between $\| P_{K'}^y - P_K^y \|$, we have
      \begin{align*}
         | C_y(K') - C_y(K) | &=  \Big| \langle P_{K'}^y - P_{K}^y,  \Sigma_{y_0} + \gamma \Sigma^1 / (1-\gamma) \rangle_{tr} \Big|
        \\
        &\leq h_{riccati, y}(C_0) \| \Delta K \| ( \| y_0 \|^2_{L^2(\RR^d)} + \gamma \| \varepsilon_{t=1} \|^2_{L^2(\RR^d)} / ( 1- \gamma) )
    \end{align*}
    Then, by the fact that $\| \xi \|_{L^2(\RR^d)} \leq \sqrt{2d} \| \xi \|_{\psi_2} $ for a sub-gaussian random vector $\xi \in \RR^d$, and $\| y_0 \|_{\psi_2}, \| \varepsilon_{t} \|_{\psi_2} \leq C_{init, noise}$, we obtain inequality $| C_y(K') - C_y(K) | \leq h_{cost, y}(C_0) \| \Delta K \|$. By Lemma~\ref{lemma:perturbation_gradient_Cy} and the definition of $f_{grad}(K)$ in~\eqref{eq:f_grad_K}, we obtain inequality $\| \nabla_K C_y(K') - \nabla_K C_y(K) \| \leq h_{f\_grad\_K}(C_0) \| \Delta K \|$. 
\end{proof}

\section{\textbf{Proofs of technical results from Section~\ref{subsection:modelfree_MKV_simulator}}}
\label{sec:app-proof-modelfree-MKV-CV}

\subsection{Intuition for zeroth order optimization}

Derivative-free optimization (see e.g.~\cite{MR2487816,MR3627456}) attempts to find an optimizer of a function $f: x \mapsto f(x)$ using only pointwise values of $f$, without having access to its gradient. The basic idea is to express the gradient of $f$ at a point $x$ by using values of $f$ in a small region around point $x$. This can be achieved by introducing an approximation of $f$ using Gaussian smoothing.  Formally, if we define the perturbed function for some $\sigma >0$ as
$
f_{\sigma^2}(x) = \mathbb{E}[f(x+ \varepsilon)],
$ where
$\varepsilon \sim \mathcal{N}(0,\sigma^2 I)$,
then its gradient can be written as
$
\frac{d}{dx}f_{\sigma^2}(x) = \frac{1}{\sigma^2} \mathbb{E}[f(x+\varepsilon) \varepsilon].
$
We can replace the Gaussian smoothing term $\varepsilon$ by a uniform random variable with bounded norms. 

Consider $U = (U^{(idy)}, U^{(com)})$ with $U^{(idy)},U^{(com)}$ being two independent random matrices in $\RR^{\ell \times d}$  uniformly distributed on the ball with radius $\tau$: $\mathbb{B}_{\tau} = \{ W \in \RR^{\ell \times d} : \| W \|_F \leq \tau \}$.
The unit ball under Frobenius norm for matrices in $\RR^{\ell \times d}$ can be identified as a unit ball for vectors in $\RR^{\ell d}$. Let $\mathbb{S}_\tau = \{ W \in \RR^{\ell \times d}: \| W \|_F = \tau \}$ be the boundary of $\mathbb{B}_\tau$, viewed as a sphere of dimension $\ell d -1$. The ratio between the sphere area of a ball with radius $\tau$ in $\RR^{\ell d}$  to its volume equals to $(\ell d) / \tau$~\cite{MR2298287}. The spectral norm of $W \in \mathbb{S}_\tau$ satisfies $\| W \| \leq \| W \|_F = \tau$.

For any $\theta = (K,L)$ and $ \tau > 0$, we introduce the following smoothed version $C_\tau$ of $C$ defined by
\begin{equation}
\label{eq:formula_perturbation_definition_C}
    C_\tau(\theta) = \mathbb{E}_U[C(\theta + U)] = \mathbb{E}_{(U^{(idy)},U^{(com)})}[C_y(K+U^{(idy)}) + C_z(L+U^{(com)})],
\end{equation}
where $\mathbb{E}_U$ means an expectation over $U$, and the costs $C_y(K+U^{(idy)})$ and $C_z(L + U^{(com)})$ are evaluated with  auxiliary processes $\by^{K + U^{(idy)}}$ and $\bz^{L + U^{(com)}}$ as defined in equations~\eqref{eq:def_Cy_Cz}.
Let $V^{(idy)} \sim \tau U^{(idy)} / \| U^{(idy)} \|_F$ and $V^{(com)} \sim \tau U^{(com)} / \| U^{(com)} \|_F$ be two random matrices independently and uniformly distributed on the sphere $\SS_\tau$, then $V = (V^{(idy)}, V^{(com)}) \in \SS_\tau \times \SS_\tau$. From~\cite[Lemma 2.1 ]{MR2298287} or~\cite[Lemma 26]{fazel2018global}, we have the following expressions for the gradients of $C_\tau(\theta)$:
\begin{align}
    \label{eq:gradCK-zero-order}
    &\nabla_K C_\tau(\theta) = \frac{\ell d}{\tau^2} \mathbb{E}_{V^{(idy)}} \big[ C_y(K + V^{(idy)}) V^{(idy)} \big] = \frac{\ell d}{\tau^2}\mathbb{E}_V \big[C(\theta + V) V^{(idy)} \big],
    \\
    \label{eq:gradCL-zero-order}
    &\nabla_L C_\tau(\theta) = \frac{\ell d}{\tau^2} \mathbb{E}_{V^{(com)}} \big[ C_z(L + V^{(com)}) V^{(com)} \big] = \frac{\ell d}{\tau^2} \mathbb{E}_V \big[C(\theta+V) V^{(com)} \big].
\end{align}
The last equalities in~\eqref{eq:gradCK-zero-order} and~\eqref{eq:gradCL-zero-order} are justified by Lemma~\ref{lem:decompose-cost-reparam} and by the fact that $V^{(idy)}$ and $V^{(com)}$ are independent mean-zero random variables.

To show that the gradient estimator from Algorithm~\ref{algo:MKVestim} provides a good approximation of the exact policy gradient, we introduce two more gradient estimators before we dive into the perturbation analysis.  Let $M$ be the number of perturbation directions, $\tau > 0$ be the perturbation radius, and $T$ be the truncation horizon for the infinite discounted MF cost. The following four forms of gradient will be used in the next subsections:
\begin{enumerate}[label=(\roman*)]
    \item \textbf{[Policy gradient]} $\nabla C(\theta) = \big(\nabla_K C(\theta), \nabla_L C(\theta) \big)$ represents the exact gradient of the cost $C(\theta)$ with respect to $\theta = (K, L)$.
    
    \item \textbf{[Perturbed policy gradient]}	$\hat{\nabla}^{M,\tau}(\theta, \underline v)$ defined in~\eqref{eq:def_perturbed_PG_MKV} below is an approximation of the exact PG based on the zeroth order approximation  with $M$ perturbation directions $\underline v = (v_i)_{i=1}^M$ where $v_i \sim \mu_{\mathbb{S}_\tau} \otimes \mu_{\mathbb{S}_\tau}$.

    \item \textbf{[Truncated policy gradient]} $\hat{\nabla}^{T, M, \tau}(\theta, \underline v)$ defined in~\eqref{eq:def_truncated_gradient_MKV} below is an approximation of the above perturbed policy gradient using a modified cost obtained by truncation at a finite horizon $T$.
	
    \item \textbf{[Sampled policy gradient]} $\tilde{\nabla}^{T, M,\tau} \big(\theta, \underline v \big)$ defined in equation~\eqref{eq:definition_approx_gradient_MKV_simulator} is the output of Algorithm~\ref{algo:MKVestim}.
\end{enumerate}

\subsection{Approximation by a smoothed cost}
\label{subsection:approximation_with_smoothed_cost}
In this subsection, we approximate the exact policy gradient $\nabla C(\theta)$ by the perturbed policy gradient $\hat{\nabla}^{M, \tau}(\theta, \underline v)$ defined in~\eqref{eq:def_perturbed_PG_MKV} below. 
To start, we consider the gradient of a smooth cost $C_\tau(\theta)$ defined in~\eqref{eq:formula_perturbation_definition_C} as $\nabla C_\tau(\theta) = \big( \nabla_K C_\tau(\theta) , \nabla_L C_\tau(\theta) \big)$. Equations~\eqref{eq:gradCK-zero-order} and~\eqref{eq:gradCL-zero-order} provide expression for $\nabla C_\tau(\theta)$ without computing the expectations. Now, we approximate these expectations $\EE_V$ therein using $M$ independent samples $(v_i^{(idy)})_{i=1}^M$ from $V^{(idy)}$ and $M$ independent samples $(v_i^{(com)})_{i=1}^M$ from $V^{(com)}$ respectively. Let $\underline v = (v_i)_{i=1}^M$ with $v_i = (v_{i}^{(idy)}, v_{i}^{(com)}) \in \RR^{\ell \times d} \times \RR^{\ell \times d}$ and $v_i \sim \mu_{\mathbb{S}_\tau} \otimes \mu_{\mathbb{S}_\tau}$. We define the ``perturbed policy gradient'' in $K$ and in $L$ by
\begin{equation}
\label{eq:def_perturbed_PG_MKV}
 \hat{\nabla}_K^{M, \tau}(\theta, \underline v) = \frac{\ell d}{\tau^2} \frac{1}{M} \sum_{i=1}^M C(\theta + v_i) v_i^{(idy)}, 
 \quad 
 \hat{\nabla}_L^{M, \tau}(\theta, \underline v) = \frac{\ell d}{\tau^2} \frac{1}{M} \sum_{i=1}^M C(\theta + v_i) v_i^{(com)}
\end{equation}
Let $\hat{\nabla}^{M, \tau}(\theta, \underline v) = \big( \hat{\nabla}_K^{M, \tau}(\theta, \underline v),  \hat{\nabla}_L^{M, \tau}(\theta, \underline v) \big)$. 
In~\eqref{eq:def_perturbed_PG_MKV}, we consider $\tau$ small enough so that by the stability Lemma~\ref{lemma:stability_of_small_perturbation}, we have $\theta + v_i \in \Theta$ for all $i=1, \ldots, M$, and thus $\hat{\nabla}^{M, \tau}(\theta, \underline v)$ is well-defined. It is clear that the randomness of the perturbed policy gradient comes from the random directions $\underline v \in (\mathbb{S}_{\tau} \times \mathbb{S}_\tau )^{M}$. 

\begin{lemma}
\label{lemma:approx_perturbed_policy_gradient_K}
    Consider $\theta \in \Theta$ with $C(\theta) \leq C_0$ for some constant $C_0 \in \RR$. For any $\varepsilon > 0$ and $\delta_{pert, grad} \in (0,1)$, if 
    \begin{align}
        \tau & \leq \min\{ C_0 h_{cost}^{-1}(C_0), h_{small-pert, y}^{-1}(C_0), \varepsilon h_{f\_grad\_K}(C_0)^{-1} / 2 \}
        \label{eq:condition_tau_K}
    \end{align}
    where $h_{cost}(C_0), h_{small-pert, y}(C_0), h_{f\_grad\_K}(C_0)$ are defined in Lemma~\ref{lemma:useful_bounds_I} and Lemma~\ref{lemma:perturbation_cost_Cy}, then, we have
    $$
        \PP_{\underline v} \left( \| \hat\nabla^{M, \tau}_K(\theta, \underline v) - \nabla_K C(\theta) \| > \varepsilon \right) \leq (\ell + d) \exp \Big( -\frac{M \varepsilon^2 \tau^2 } { 64 \ell^2 d^2 C_0^2 + 16 \ell d C_0 \varepsilon / 3 } \Big).
    $$
\end{lemma}

\begin{proof} 
    For $i=1, \ldots, M$, let $\theta_i = \theta + v_i$ and 
    $$
    W_i = \frac{\ell d}{\tau^2 M} \big( C(\theta_i) \tilde{v}_i^{(idy)} - \EE_{v_i}[ C(\theta + v_i) \tilde{v}_i^{(idy)} ] \big)
    $$ 
    where $\tilde{v_i}^{(idy)}$ is the dilation matrix of $v_i^{(idy)}$ and thus $\| \tilde v_i^{(idy)} \| = \| v_i^{(idy)} \| \leq \| v_i^{(idy)} \|_F = \tau$.
    For $\tau$ small enough so that $\tau \leq C_0 h_{cost}(C_0)^{-1}$ where $h_{cost}(C_0)$ was defined in Lemma~\ref{lemma:perturbation_cost_Cy}, we have 
    $$
    | C(\theta_i) - C(\theta) | \leq h_{costs, y}(C_0) \| v_i^{(idy)} \| + h_{cost, z}(C_0) \| v_i^{(com)} \| \leq  h_{cost}(C_0) \tau \leq C_0
    $$ 
    and so $\| C(\theta_i) \tilde{v}_i^{(idy)} - \EE_{v_i}[C(\theta_i) \tilde{v}_i^{(idy)}] \| \leq 4C_0 \tau$. Thus, for $i=1, \ldots, M$,
    $$
        \| W_i \| \leq \frac{4 \ell d C_0 }{M \tau} =: \lambda_{W},
        \qquad
        \Big \| \sum_{i=1}^M \EE[ W_i^2 ] \Big\| \leq M \sup_{i=1,\ldots, M} \EE[ \| W_i \|^2 ] \leq \frac{8 \ell^2 d^2 C_0^2}{M \tau^2} = :\sigma^2_W.
    $$
    Using  Bernstein's matrix inequality with bounded random matrices~\cite{MR2946459}[Theorem 6.1], and expressions~\eqref{eq:gradCK-zero-order} for $\nabla_K C_\tau(\theta)$, we obtain that for any $\tilde \varepsilon > 0$,
    \begin{align*}
        \PP \big( \| \hat{\nabla}_K^{M, \tau}(\theta, \underline v) - \nabla_K C_\tau(\theta) \| > \tilde{\varepsilon} \big) & = 
        \PP \big( \big\| \sum_{i=1}^M W_i \big\|  > \tilde{\varepsilon} \big) \leq (\ell + d) \exp \Big(\frac{- \tilde{\varepsilon}^2 / 2}{ \sigma^2_W + \tilde{\varepsilon} \lambda_W / 3}  \Big).
    \end{align*}
    Moreover, when $\| (K + U^{(idy)}) - K \| = \| U^{(idy)} \| \leq \| U^{(idy)} \|_F = \tau \leq h_{small-pert,y}^{-1}(C_0)$, Lemma~\ref{lemma:perturbation_gradient_Cy} on the perturbation of gradients implies that
    \begin{align*}
        \| \nabla_K C_\tau(\theta) - \nabla_K C(\theta) \| & = \| \EE_{U}\big[ \nabla_K C_y( K+ U^{(idy)}) - \nabla_K C_y(K) \big] \|
        \\
        & \leq \EE_U[ f_{grad}(K) \| U \| ] \leq h_{f\_grad\_K}(C_0) \tau, 
    \end{align*}
    where $h_{f\_grad\_K}(C_0)$ is defined in~\eqref{eq:h_f_grad_K}.
    Hence, with high probability, we have
    \begin{align*}
        \| \hat\nabla^{M, \tau}_K( \theta, \underline v) - \nabla_K C(\theta) \| & \leq 
         \|  \hat\nabla^{M, \tau}_K(\theta, \underline v) - \nabla_K C_\tau(\theta) \| + \| \nabla_K C_\tau(\theta) - \nabla_K C(\theta) \|
         \\
         & \leq  \tilde{\varepsilon} + h_{f\_grad\_K}(C_0) \tau.
    \end{align*}
    By choosing $\tilde \varepsilon = \varepsilon /2$, whenever $\tau$ satisfies \eqref{eq:condition_tau_K}, we have $ \tilde{\varepsilon} + h_{f\_grad\_K}(C_0) \tau \leq \varepsilon$, the desired result follows.
\end{proof}

Prompted by the conditions under which Lemma~\ref{lemma:approx_perturbed_policy_gradient_K} holds, we define the following two functions for the perturbation radius and the number of perturbation directions:
\begin{align}
    & \phi_{pert, radius}(\varepsilon, C_0) = \max \big\{ h_{cost}(C_0) / C_0,  h_{small-pert, y}(C_0), h_{small-pert, z}(C_0),
    \nonumber \\
    & \hspace{120pt}  4 h_{f\_grad\_K}(C_0) / \varepsilon, 4 h_{f\_grad\_L}(C_0) / \varepsilon \big\},
    \label{eq:phi_pert_radius}
    \\
    & \phi_{pert, size}(\varepsilon, \tau, C_0, \delta_{pert} ) =   \Big( \frac{16 \ell d  C_0 }{ \tau \varepsilon } + 1 \Big)^2 \log \Big( \frac{ 2(\ell + d) }{\delta_{pert}} \Big)
    \label{eq:phi_pert_size}
\end{align}
where $h_{cost}(C_0)$, $h_{small-pert, y}(C_0)$, $h_{small-pert, z}(C_0)$, $h_{f\_grad\_K}(C_0)$, $h_{f\_grad\_L}(C_0)$ are defined in Lemma~\ref{lemma:useful_bounds_I} and Lemma~\ref{lemma:perturbation_cost_Cy}.

\begin{lemma}
\label{lemma:approx_perturbed_gradient_MKV_simulator}
Consider $\theta \in \Theta$ with $C(\theta) \leq C_0$ for some constant $C_0 \in \RR$. For a target precision $\varepsilon > 0$, if 
\begin{equation}
\label{eq:condition_tau_M_approx_perturbed_gradient}
    \tau^{-1} \geq \phi_{pert, radius}(\varepsilon, C_0)
    \quad \text{and} \quad
    M \geq \phi_{pert, size}(\varepsilon, \tau, C_0, \delta_{pert}),
\end{equation}
then
$$
    \PP_{\underline v}\big( \| \hat \nabla^{M, \tau}(\theta, \underline v) - \nabla C(\theta) \| > \varepsilon \big) \leq \delta_{pert},
$$
where $\| \hat \nabla^{M, \tau}(\theta, \underline v) - \nabla C(\theta) \|  = \| \hat \nabla^{M, \tau}_K(\theta, \underline v) - \nabla_K C(\theta) \| + \| \hat \nabla^{M, \tau}_L(\theta, \underline v) - \nabla_L C(\theta) \|$.
\end{lemma}

\begin{proof} We apply Lemma~\ref{lemma:approx_perturbed_policy_gradient_K} with target precision $\varepsilon / 2$, then for $\tau$ small enough and $M$ large enough satisfying condition~\eqref{eq:condition_tau_M_approx_perturbed_gradient}, we have
$$
    \PP_{\underline v}\big( \| \hat \nabla^{M, \tau}_K(\theta, \underline v) - \nabla_K C(\theta) \| > \varepsilon / 2 \big) \leq \delta_{pert}.
$$
Together with the similar concentration inequality for $\| \hat \nabla^{M, \tau}_L(\theta, \underline v) - \nabla_L C(\theta) \|$, we conclude the proof.
\end{proof}

\subsection{Approximation by a truncated cost} 
\label{subsection:approximation_C_with_C_T}
For $T \geq 0$, we define the truncated variances and the truncated costs by:
\begin{align}
    \label{eq:truncation_Sigma_y_Sigma_z_T}
    \Sigma_{K}^{y,T}  &= \mathbb{E}\big[\sum_{t=0}^{T-1} \gamma^t (y_t^{K} (y_t^K)^\top \big],
    \quad 
    \Sigma_{L}^{z,T}  = \mathbb{E} \big[ \sum_{t=0}^{T-1} \gamma^t (z_t^L) (z_t^L)^\top \big],
    \\
    \label{eq:truncation_Cy_Cz_T} 
    C^T_y(K) &= \mathbb{E} \big[  \sum_{t=0}^{T-1}  \gamma^t f(y_t^{K}, K, Q, R) \big] ,
    \quad
    C^T_z(L) = \mathbb{E} \big[  \sum_{t=0}^{T-1}  \gamma^t f(z_t^{L}, L, \tilde{Q}, \tilde{R} ) \big]
\end{align}
where $f(\xi, \varphi, q, r) = \xi^{\top} (q + \varphi^\top r \varphi)\xi$. We define $C^T(\theta) = \EE[ \sum_{t=0}^{T-1} \gamma^t c_t(\theta) ]$ where $c_t(\theta) = c(x_t^\theta, \bar{x}_t^\theta, u_t^\theta, \bar u_t^\theta)$ is the instantaneous cost at time $t$. From the definition of the auxiliary process $(\by^K, \bz^L)$ from~\eqref{eq:dyn_y_theta} and~\eqref{eq:dyn_z_theta}, we have $C^T(\theta) = C_y^T(K) + C_z^T(L)$.

For $\theta = (K, L) \in \Theta$ with $C(\theta) \leq C_0$, from the proof of Lemma~\ref{lemma:stability_of_small_perturbation} on the stability of small perturbations of control parameters, we have 
$$
    \gamma \| \rA - \rB K \|^2 \leq  1 -  \frac{ \lambda_y^1 \lambda_{min}(Q) } {2 d C_0 }, \quad \gamma \| \tilde \rA - \tilde \rB L \|^2 \leq 1 - \frac{ \lambda_z^0 \lambda_{min}(\tilde Q) }{ 2 d C_0}.
$$
We define the constant $\gamma_{pert}$ and the function $\phi_{truncation}(\varepsilon)$ of $\varepsilon$ for the lower bound of the truncation time horizon:
\begin{align}
    & \gamma_{pert} = \max \big\{ \gamma, 1 - \lambda_y^1 \lambda_{min}(Q) / (2 d C_0),\, 1 - \lambda_z^0 \lambda_{min}(\tilde Q) / (2 d C_0) \big\},
    \label{eq:gamma_pert}
    \\
    & \phi_{truncation}(\varepsilon) = \Big( \frac{\log( 2d C_{init, noise}^2) - \log(\varepsilon (1 - \gamma_{pert})^2 ) + 1 }{\log(1/\gamma_{pert})} \Big)^2.
    \label{eq:phi_truncation}
\end{align}
The truncated cost $C^T(\theta)$ and the MF cost $C(\theta)$ are related by the following..

\begin{lemma}
\label{lemma:approx_truncated_cost}
    Consider $\theta \in \Theta$ with $C(\theta) \leq C_0$. For any $\varepsilon > 0$ and $T \geq 2$, if $T \geq \phi_{truncation}(\varepsilon)$, we have
    \begin{equation}
        C(\theta) - C^T(\theta)  \leq \varepsilon h_{T}(C_0)
    \end{equation}
    where $h_T(C_0) = d(\| Q \| + \| \tilde Q\| + \|R \| h_K(C_0)^2 + \| \tilde R \| h_L(C_0)^2 )$ with $h_K(C_0)$, $h_L(C_0)$ defined in Lemma~\ref{lemma:useful_bounds_I}.
\end{lemma}

\begin{proof} 
    We first show that if $T \geq \max\{ 2, \phi_{truncation}(\varepsilon)\}$, then $ \big\| \Sigma_K^y - \Sigma_K^{y,T} \big\| \leq \varepsilon$. We notice that  
    $
        \mathbb{E}\left[ \gamma^t y^K_t (y^K_t)^\top \right] = (\cF_K^y)^t \left( \Sigma_{y_0} \right) + \sum_{s=0}^{t-1} (\cF_K^y)^s \left( \gamma^{t-s} \Sigma^1 \right).	
    $
    So, by Lemma~\ref{lemma:expression_Sigma_K_with_T_K} and the definition of $\cF_K$ and $\cT_K$ in~\eqref{eq:def-ope-F_K_y}, we have
    \begin{align*}
       \Sigma_K^y - \Sigma_K^{y, T} = & \Sigma_K^y - \mathbb{E} \big[ \sum_{t=0}^{T-1} \gamma^t y_t^K (y_t^K)^\top \big]  
        \\
        = & \cT_K \big(\Sigma_{y_0} + \frac{\gamma}{1-\gamma}\Sigma^1 \big) -  \Big[ \sum_{t=0}^{T-1} (\cF_K)^t (\Sigma_{y_0})  + \sum_{t=1}^{T-1} \Big(\sum_{s=0}^{t-1} (\cF_K)^s (\gamma^{t-s} \Sigma^1) \Big) \Big]
         \\
        = & \frac{\gamma}{1-\gamma}\cT_K( \Sigma^1 ) +  \sum_{t=T}^\infty (\cF_K)^t(\Sigma_{y_0})  - \sum_{s=1}^{T-1}  \gamma^s \Big( \cT_K(\Sigma^1) - \sum_{t=T-s}^\infty (\cF_K)^t(\Sigma^1) \Big).
         \\
         = & \frac{\gamma^T}{1- \gamma} \cT_K \left( \Sigma^1 \right) + (\cF_K)^T \left( \cT_K(\Sigma_{y_0})  \right) 
         + \sum_{s=1}^{T-1} \gamma^s \left(\cF_K \right)^{T-s} \left( \cT_K (\Sigma^1) \right).
    \end{align*}
Using $\gamma \leq \gamma_{pert}$, $\vertiii{\cF_K} = \gamma \| \rA - \rB K \|^2 \leq \gamma_{pert}$, $\vertiii{ \cT_K } \leq (1- \gamma_{pert})^{-1}$, $\| \Sigma_{y_0} \| \leq \EE[\| y_0\|^2] \leq  \EE[  2 d \| y_0 \|_{\psi_2})^2 ] \leq 2d C_{init, noise}^2$ and $\| \Sigma^1 \| \leq 2d C_{init, noise}^2$ we obtain
    \begin{align*}
        \big\| \Sigma_K^y - \Sigma_K^{y,T} \big\|
        & \leq  \Big( \frac{\gamma_{pert}^T}{1- \gamma_{perc}} + \gamma_{perc}^T + (T-1) \gamma_{perc}^T \Big)  \frac{2d C_{init, noise}^2}{1- \gamma_{perc}}  
        \\
        & \leq  2d \gamma_{perc}^T (T +1) C_{init, noise}^2 / (1- \gamma_{perc})^2 .
    \end{align*} 
    For $T \geq 2$ and $T \geq \phi_{truncation}(\varepsilon, \gamma_{\theta})$, we have
    $$
          T \log( 1/ \gamma_{perc}) - \log(T + 1) \geq \sqrt{T} \log(1/\gamma_{perc}) -1 \geq \log \Big( \frac{2d C_{init, noise}^2}{\varepsilon ( 1- \gamma_{perc})^2 } \Big),
    $$
    which implies directly $\| \Sigma_K^y - \Sigma_K^{y, T} \| \leq \varepsilon$. 
    Similar arguments imply $\| \Sigma_L^z - \Sigma_{L}^{z, T} \| \leq \varepsilon$.
    From equation~\eqref{eq:cost_expression_with_PK_y} and the definition of $C_y^T(K)$, $C_z^T(L)$ in~\eqref{eq:truncation_Cy_Cz_T}, we have 
    \begin{align*}
        | C(\theta) - C^T(\theta) | &= \Big| \langle Q + K^\top R K,\, \Sigma_K^y - \Sigma_K^{y, T} \rangle_{tr} + \langle \tilde Q + L^\top \tilde R L,\, \Sigma_L^z - \Sigma_L^{z, T} \rangle_{tr} \Big|
        \\
        & \leq \varepsilon h_T(C_0)
    \end{align*}
which is the desired estimate.
\end{proof}

We define the truncated policy gradients in $K$ and in $L$ at time $T$ with perturbation size, perturbation radius, and perturbation directions $(M, \tau, \underline v)$ by
\begin{equation}
\label{eq:def_truncated_gradient_MKV}
    \hat{\nabla}_K^{T, M, \tau}(\theta, \underline v) = \frac{\ell d}{M \tau^2} \sum_{i=1}^M C^T(\theta + v_i) v_i^{(idy)}, 
    \quad
    \hat{\nabla}_L^{T, M, \tau}(\theta, \underline v) = \frac{\ell d}{M \tau^2} \sum_{i=1}^M C^T(\theta + v_i) v_i^{(com)},  
\end{equation}
and we set $\hat{\nabla}^{T, M, \tau}(\theta, \underline v) = \big( \hat{\nabla}_K^{T, M, \tau}(\theta, \underline v), \,  \hat{\nabla}_L^{T, M, \tau}(\theta, \underline v)  \big)$. Here, we assume that the perturbation radius $\tau$ is small enough so that $\theta + v_i \in \Theta$ for all $i=1, \ldots, M$, and thus the truncated policy gradient $\hat{\nabla}^{T, M, \tau}(\theta, \underline v)$ is well-defined.

To control the approximation of the perturbed policy gradient $\hat{\nabla}^{M, \tau}$ by a truncated policy gradient $\hat{\nabla}^{T, M, \tau}$,  based on Lemma~\ref{lemma:approx_truncated_cost}, we define
\begin{align}
    h_T(C_0) &= d(\| Q \| + \| \tilde Q\| + \|R \| h_K(C_0)^2 + \| \tilde R \| h_L(C_0)^2 ),
    \label{eq:h_T}
    \\
    \phi_{trunc, T}(\varepsilon, \tau, C_0) & = \phi_{truncation}\big( \frac{\tau  \varepsilon }{2 \ell d. h_T(C_0)} \big)
    \nonumber 
    \\
    & = \Big( \frac{ \log( 4 \ell d^2 C_{init, noise}^2 h_T(C_0) ) - \log( \varepsilon \tau (1 - \gamma_{pert})^2 ) + 1 }{\log(1/\gamma_{pert})} \Big)^2,
    \label{eq:phi_trunc_T}
\end{align}
where $h_{K}(C_0)$ and $h_L(C_0)$ are defined in Lemma~\ref{lemma:useful_bounds_I}, and $\gamma_{pert}$ and $\phi_{truncation}$ are defined in equations~\eqref{eq:gamma_pert} and~\eqref{eq:phi_truncation}.

\begin{lemma}
\label{lemma:approx_truncated_policy_gradient_MKV_simulator}
If $\theta \in \theta$ with $C(\theta) \leq C_0$ for some constant $C_0 \in \RR$, for any given $\varepsilon > 0$, if the truncation time horizon satisfies $T \geq \max \{2, \phi_{trunc, T}(\varepsilon, \tau, C_0) \}$, then we have,
\begin{equation}
        \big\| \hat{\nabla}^{T, M, \tau}(\theta, \underline v) -  \hat{\nabla}^{M, \tau}(\theta, \underline v) \big\| \leq \varepsilon
\end{equation}
where $ \big\| \hat{\nabla}^{T, M, \tau}(\theta, \underline v) -  \hat{\nabla}^{M, \tau}(\theta, \underline v)  \big\| =  \big\| \hat{\nabla}_K^{T, M, \tau}(\theta, \underline v) -  \hat{\nabla}_K^{M, \tau}(\theta, \underline v)  \big\| +  \big\| \hat{\nabla}_L^{T, M, \tau}(\theta, \underline v) -  \hat{\nabla}_L^{M, \tau}(\theta, \underline v)  \big\|$.
\end{lemma}

\begin{proof}
    Let $\theta_i = \theta + v_i$ for $i=1,\ldots, M$. From the definitions~\eqref{eq:def_perturbed_PG_MKV} and~\eqref{eq:def_truncated_gradient_MKV}, Lemma~\ref{lemma:approx_truncated_cost} implies that
    $$
        \| \hat{\nabla}_K^{T, M, \tau}(\theta, \underline v) -  \hat{\nabla}_K^{M, \tau}(\theta, \underline v)  \big\| 
        \leq  \frac{\ell d}{M \tau^2} \sum_{i=1}^M | C(\theta_i) - C^T(\theta_i) |. \| v_i^{(idy)} \| \leq  \frac{\varepsilon}{2}.
    $$
    We then conclude with similar inequalities for the gradient in $L$.
\end{proof}

\subsection{Costs with MKV simulator}
 \label{subsection:approx_sampled_costs_with_MKV_simulator}
In this subsection, we will first prove a matrix concentration inequality for a quadratic form $(y_t^K)^\top ( Q + K^\top R K) y_t^K$ in an auxiliary process $\by^{K}$ at time $t$, being a sub-gaussian random vector with dependent coordinate entries. The quadratic forms for $\by^K$ as well as for $\bz^L$ are closely related to the sample costs simulated from an MKV simulator $\cS_{MKV}^T$. Based on this result, we will quantify the approximation error between the sampled policy gradients as outputs of Algorithm~\ref{algo:MKVestim} and the truncated police gradient defined in equation~\eqref{eq:def_truncated_gradient_MKV}.
To start, we define the constant
\begin{equation}
\label{eq:h_psi_2, y}
    h_{\psi_2, y}(\tau, C_0) = \frac{2 \tau d^2 C_{init, noise}^2  \Big( \| Q \| + \|R\| ( \tau + h_K(C_0) )^2 \Big) }{(1 - \sqrt{\gamma})^2}.
\end{equation}
A similar constant $h_{\psi_2, z}(\tau, C_0)$ is defined with $(\tilde Q, \tilde R, h_L(C_0))$ in the above equation.
As before, we shall find it convenient to use the function $f$ given by $f(\xi, \varphi, q, r) \mapsto \xi^\top ( q + \varphi^\top r \varphi ) \xi$. 

The following result will be needed in the proof of the main probabilistic estimate of this section in Lemma~\ref{lemma:concentration_ineq_for_yt_K}. Its proof is similar to \cite[Proposition 2.4]{zajkowski2020bounds}.

\begin{lemma}
\label{lemma:subexponential_quadratic_form_bound}
If $\xi \in \RR^{d}$ is a sub-gaussian random vector  and $U \succeq 0$ a positive semi-definite matrix  in $\RR^{d \times d}$, then $\xi^\top U \xi$ is a sub-exponential random variable, and
\begin{equation}
    \| \xi^\top U \xi \|_{\psi_1} \leq Tr(U) \| \xi \|_{\psi_2}^2.
\end{equation}
\end{lemma}
\begin{proof}
   Let $U = V^\top \Sigma V$ be the singular value decomposition of $U$ with  $\Sigma$ being diagonal, i.e. $\Sigma = \diago(\lambda_{1}, \lambda_{2}, \ldots, \lambda_{d})$ where $\{\lambda_i;\;i=1,\ldots,d\}$ are the eigenvalues of $U$. Let $\zeta = V \xi = [\zeta_1, \ldots, \zeta_d]^\top$, then the $d$ entries $\zeta_{i=1,\ldots d}$ are sub-gaussian random variables, and $\zeta$ is also a sub-gaussian random vector such that $\| \zeta \|_{\psi_2} = \sup_{ s \in \mathbb{S}^{d-1}}\| \langle V \xi, s\rangle\|_{\psi_2} = \| \xi \|_{\psi_2}$. Consequently, we have $\xi^\top U \xi = \zeta^\top \Sigma \zeta =  \sum_{i=1}^d \lambda_i \zeta_i^2 $, and~\cite{vershynin2018high}[Lemma 2.7.7] implies that $\xi^\top U \xi$ is sub-exponential. Moreover,
    \begin{align*}
        \| \xi^\top U \xi \|_{\psi_1} = \Big\| \sum_{i=1}^d \lambda_i \zeta_i^2 \Big\|_{\psi_1} \leq \sum_{i=1}^d\lambda_i  \| \zeta_i^2 \|_{\psi_1} = \sum_{i=1}^d \lambda_i  \| \zeta_i \|_{\psi_2}^2 \leq \big( \sum_{i=1}^d \lambda_i \big) \| \zeta \|_{\psi_2}^2,
    \end{align*}
    where the second equality is justified by the definitions of $\| \cdot \|_{\psi_1}$ and $\| \cdot \|_{\psi_2}$ norms (see~\cite{vershynin2018high}[Lemma 2.7.6]), and the last inequality is due to $\| \zeta_i \|_{\psi_2} = \| \langle \zeta, e_i \rangle \|_{\psi_2} \leq \| \zeta \|_{\psi_2}$ with $e_i \in \mathbb{S}^{d-1}$ being the $i$-th vector of the canonical basis in dimension $d$. We conclude the proof of the lemma noticing that $Tr(U) = \sum_{i=1}^d \lambda_{i}$.
\end{proof}

\begin{lemma}
\label{lemma:concentration_ineq_for_yt_K}
    Consider $\theta = (K, L) \in \Theta$ with $C(\theta) \leq C_0$ for some constant $C_0 \in \RR$, let $\by^K = (y_t^K)_{t \geq 0}$ be the auxiliary process with dynamics~\eqref{eq:dyn_y_theta}, and for any $t \geq 0$, consider the quadratic form $\zeta_t^{K} = f(y_t^K, K, Q, R) = (y_t^K)^\top (Q + K^\top R K ) y_t^K$. Then, the random vector $y_t^K$ is sub-gaussian in $\RR^d$, and the real-valued random variable $\zeta_t^K$ is sub-exponential. Their norms are bounded by
    \begin{equation}
    \label{eq:bound_on_y_t_K_and_quadratic_form}
        \| y_t^K \|_{\psi_2} \leq \gamma^{-t/2} d C_{init, noise} / (1 - \sqrt{\gamma}),
        \qquad
        \| \zeta_t^K \|_{\psi_1} \leq Tr( Q + K^\top R K) \| y_t^{K} \|_{|\psi_2}^2.
    \end{equation}
    Moreover, for $M$ perturbed control parameters $K_i = K + v_i^{(idy)}$ with perturbation directions $(v_i^{(idy)})_{i=1, \ldots, M}$ on the sphere $\SS_{\tau}$ independently sampled from $\mu_{\SS_\tau}$, we have, for any $\varepsilon > 0$, 
    \begin{equation}
    \label{eq:concentration_ineq_for_quadratic_yt_K}
        \PP \Big( \Big\| \frac{\gamma^{t}}{M} \sum_{i=1}^M ( \zeta_t^{K_i} -  \EE[ \zeta_t^{K_i} ] ) v_i^{(idy)} \Big\| \geq \varepsilon \Big) \leq (\ell + d)  \exp \Big( \frac{-M \varepsilon^2 }{ 4 h_{\psi_2, K}^2(\tau, C_0) + \varepsilon h_{\psi_2, K}(\tau, C_0) } \Big) 
    \end{equation}
    where the expectation $\EE$ and the probability $\PP$ are with respect to the randomness in $(y_t^{K_i})_{i=1,\ldots, M}$, but not the randomness in the perturbation directions $(v^{(idy)}_i)_{i=1, \ldots, M}$.
\end{lemma}

\begin{proof}
We first show the two bounds on $y_t^K$ and $\zeta_t^K$ in equation~\eqref{eq:bound_on_y_t_K_and_quadratic_form}.
From the dynamics~\eqref{eq:dyn_y_theta}, we have $y_t^K = (\rA - \rB K)^{t} y_0 + \sum_{s=1}^{t} ( \rA - \rB K)^{t-s} \varepsilon_s$ where $y_0 = \varepsilon_0 - \EE[ \varepsilon_0 ]$. Because $(\varepsilon_s)_{s \geq 0}$ are all sub-gaussian random variables, then $y_t^K$ is also sub-gaussian. Now. Let $\xi_0 = y_0$, $\xi_{s+1} = \varepsilon_{s+1}$ and $H_s = (\rA - \rB K)^{t-s}$ for $s=0, \ldots, t$. Let $\xi_{s,j} \in \RR$ and $H_{s, j} \in \RR^d$ be the $j$-th entry of $\xi_s$ and the $j-$th column of $H_s$, for all $j=1, \ldots, d$. Then $\| \xi_{s, j} \|_{\psi_2} \leq \| \xi_s \|_{\psi_2} \leq C_{init, noise} $, $\| H_{s, j} \| \leq \| H_s \|$, and
$
y_t^K = \sum_{s=0}^t H_s \xi_s.
$
We deduce that
\begin{align*}
    \gamma^t \| y_t^{K}\|_{\psi_2}^2 
    \leq \gamma^t \Big( \sum_{s=0}^{t} \sum_{j=1}^d \| H_{s, j} \xi_{s, j} \|_{\psi_2} \Big)^2 
    &  \leq \gamma^t \Big( \sum_{s=0}^t  \sum_{j=1}^d \| H_{s, j} \|_2. \| \xi_{s, j} \|_{\psi_2} \Big)^2
    \\
    &  \leq   d^2 C_{init, noise}^2 \Big( \sum_{s=0}^t \gamma^{s/2} ( \gamma^{(t-s)/2} \| H_s \| ) \Big)^2  
    \\
    &  \leq d^2 C_{init, noise}^2 / (1 - \sqrt{\gamma})^2.
\end{align*}
Besides, from Lemma~\ref{lemma:subexponential_quadratic_form_bound}, we know that the quadratic form $\zeta_t = f(y_t^K, K, Q, R)$ is a sub-exponential random variable and $\| \zeta_t^K \|_{\psi_1} \leq Tr(Q + K^\top R K) \| y_t^K \|_{\psi_2}$. 

Now, we will show the concentration inequality~\eqref{eq:concentration_ineq_for_quadratic_yt_K} with small perturbation directions $(v_i^{(idy)})_{i=1}^M$.
Consider a small enough perturbation radius $\tau$, so that the perturbed control parameter $K_i = K + v_i^{(idy)}$ satisfies $\gamma \| \rA - \rB K_i \|^2 < 1$ for all $i = 1, \ldots, M$. For each $K_i$, we define an auxiliary process $\by^{K_i} = (y_t^{K_i})_{t \geq 0}$ following dynamics~\eqref{eq:dyn_y_theta}, but with an independent copy $\boldsymbol{\varepsilon^{i}}$ of the noise process $\boldsymbol{\varepsilon}$ and an independent copy of the initial status $y_0^{K_i} \sim y_0$. By consequence, for each time $t \geq 0$, $(\by^{K_1}_t, \ldots, \by^{K_M}_t)$ are $M$ i.i.d. random vectors in $\RR^{d}$. 
For a chosen time $t \geq 0$ and $i=1, \ldots, M$, consider the quadratic form of $y_t^{K_i}$ given by $\zeta_t^{K_i} = f(y_t^{K_i}, K_i, Q, R)$ and its unbiased term $W_t^i = \zeta_t^{K_i} - \EE[ \zeta_t^{K_i}]$. 

Let $\tilde v_i^{(idy)} \in \RR^{ (\ell + d) \times (\ell + d) }$ be the dilatation matrix of $v_i^{(idy)} \in \RR^{\ell \times d}$. Such a matrix is symmetric and has the same norm, i.e. $\|\tilde v_i^{(idy)}\|=\|v_i^{(idy)}\|$. See for example \cite{MR2946459}. For sub-exponential random variables, we have 
$
     \| \zeta_t^{K_i} \|_{L^p}^p \leq 2 p!  \| \zeta_t^{K_i} \|_{\psi_1}^p 
$
so that
$$
 \EE\big[ ( W_t^i \tilde{v}_i^{(idy)})^p \big] = \tau^p \big\| W_t^i \big\|_{L^p}^p \big( \tilde{v}_i^{(idy)} / \tau \big)^p \preceq \tau^p 2^p. ( 2 p! \| \zeta_t^{K_i} \|_{\psi_1}^p ). \mathbf{I}_{(\ell + d)}
$$
where $\mathbf{I}_{(\ell + d)}$ is an identity matrix on $\RR^{(\ell + d) \times (\ell + d)}$. Let's define two time-dependent coefficients:
$$
    r_t = \gamma^{-t} h_{\psi_2}(\tau, C_0) , 
    \quad
    \sigma_t = 2  \gamma^{-t} \sqrt{M} h_{\psi_2}(\tau, C_0).
$$
We observe that for all $i=1, \ldots, M$, $2 \tau \| \zeta_t^{K_i} \|_{\psi_1} \leq \gamma^{-t} h_{\psi_2}(\tau, C_0)$, so $\EE\big[ ( W_t^i \tilde{v}_i^{(idy)})^p \big] \preceq ( p! / 2 ) r_t^{p-2} \big( 4 \tau \| \zeta_t^{K_i} \|_{\psi_1} \mathbf{I}_{(\ell + d)} \big)^2$, and $\| \sum_{i=1}^M \big( 4 \tau \| \zeta_t^{K_i} \|_{\psi_1} \mathbf{I}_{(\ell + d)} \big)^2 \| \leq \sigma_t^2$. Thus, by the Matrix Bernstein inequality for sub-exponential case~\cite{MR2946459}[Theorem 6.2], we have 
\begin{align*}
     \PP \Big( \Big\| \frac{\gamma^{t}}{M} \sum_{i=1}^M ( \zeta_t^{K_i} -  \EE[ \zeta_t^{K_i} ] ) v_i^{(idy)} \Big\| \geq \varepsilon \Big) & = \PP \Big( \big\| \sum_{i=1}^M W_t^i \tilde{v}_i^{(idy)} \big\| \geq \frac{\varepsilon M}{\gamma^t} \Big) 
     \\
     & \leq (\ell + d) \exp \Big( \frac{- (\varepsilon M \gamma^{-t} )^2 }{ \sigma_t^2 + r_t (\varepsilon M \gamma^{-t}) }\Big). 
\end{align*}
Replacing the coefficient $r_t$ and $\sigma_t$ with $h_{\psi_2, y}(\tau, C_0)$, we obtain the concentration inequality~\eqref{eq:concentration_ineq_for_quadratic_yt_K}.
\end{proof}

\begin{remark}
It should be noted that the concentration inequality~\eqref{eq:concentration_ineq_for_quadratic_yt_K} is related to a quadratic form of sub-gaussian random vector with dependent coordinate entries, so we cannot apply the Hanson-Wright concentration inequality~\cite{vershynin2018high} or the classical Bernstein's inequality with bounded variables as considered in~\cite{fazel2018global}. Lemma~\ref{lemma:concentration_ineq_for_yt_K} provided a tail bound with the help of a discount coefficient $\gamma$.
\end{remark}

Now, we quantify the quality of the approximation of the sampled policy gradients from Algorithm~\ref{algo:MKVestim} by the truncated policy gradients defined in~\eqref{eq:def_truncated_gradient_MKV}. We recall that the sampled policy gradients at $\theta = (K, L)$ with $M$ perturbation directions $\underline v = ( (v_i^{(idy)}, v_i^{(com)} ) )_{i=1}^M$ on $\SS_{\tau} \times \SS_{\tau}$ with radius $\tau > 0$ and from an MKV simulator $\cS_{MKV}^{T}$ are given by 
\begin{equation}
\label{eq:def_sampeld_policy_gradient_MKV_restate}
    \tilde{\nabla}^{T, M,\tau}_K (\theta, \underline v) = \frac{\ell d}{\tau^2}\frac{1}{M} \sum_{i=1}^M \tilde{C}^T(\theta_i) v^{(idy)}_i, 
    \quad
    \tilde{\nabla}^{T, M,\tau}_L (\theta, \underline v) = \frac{\ell d}{\tau^2}\frac{1}{M} \sum_{i=1}^M \tilde{C}^T(\theta_i) v^{(com)}_i, 
\end{equation}
where $\theta_i = (K + v_i^{(idy)}, L + v_i^{(com)})$ for $i=1, \ldots, M$. We also define the following constants to bound the perturbation size $M$ in the sampled costs:
\begin{align}
     \phi_{sample,size}(\varepsilon, \tau, T, C_0, \delta_{sample}) = \log \Big( \frac{4T (\ell + d)}{\delta_{sample}} \Big)\Big( \frac{8 \ell d T}{\varepsilon \tau^2 }   h_{\psi_2}(\tau, C_0)  + \frac{1}{4} \Big)^2
    \label{eq:phi_sample_size}
\end{align}
where $h_{\psi_2}(\tau, C_0) = \max \big\{ h_{\psi_2, y}(\tau, C_0),\, h_{\psi_2, z}(\tau, C_0) \big\}$ with $h_{\psi_2, y}(\tau, C_0)$ and $h_{\psi_2, y}(\tau, C_0)$ defined in~\eqref{eq:h_psi_2, y}. 

\begin{lemma}
\label{lemma:approx_sampled_pg_with_truncated_pg}
If for $\theta \in \Theta$ with $C(\theta) \leq C_0$ for some $C_0 \in \RR$, the perturbation size $M$ satisfies 
$$
    M \geq \phi_{sample, size}(\varepsilon, \tau, T, C_0, \delta_{sample}),
$$
then for any sequence of $2M$ perturbation directions $\underline v$ independently and uniformly sampled from $\SS_\tau$, we have, for any $\varepsilon > 0$,
\begin{equation}
\label{eq:approx_sampled_pg_with_truncated_pg}
    \PP \big( \| \tilde \nabla^{T, M, \tau}(\theta, \underline v) - \hat{\nabla}^{T, M, \tau}(\theta, \underline v) \| \geq \varepsilon \big) \leq \delta_{sample}.
\end{equation}
\end{lemma}

\begin{proof}
    Let $\theta_i = \theta + v_i = ( K + v_i^{(idy)}, L + v_i^{(com)}) = (K_i, L_i)$ with perturbation direction $v_i = (v_i^{(idy)}, v_i^{(com)})$, for $i=1, \ldots, M$. As before, we use the auxiliary processes $(\by^{K_i}, \bz^{L_i})$ and for every $t\geq 0$ the quadratic forms $\zeta_t^{K_i} = f(y_t^{K_i}, K_i, Q, R)$ and $\tilde \zeta_t^{L_i} = f(z_t^{L_i}, L_i, \tilde Q, \tilde R)$. Then 
    \begin{align*}
        & \|  \tilde \nabla^{T, M, \tau}_K(\theta, \underline v) - \hat{\nabla}^{T, M, \tau}_K(\theta, \underline v)  \|
        \\
        &\hskip 15pt
        =\Big\| \frac{\ell d}{M \tau^2} \sum_{i=1}^M \sum_{t=0}^{T-1} \gamma^t \big[ \big( \zeta_t^{K_i} - \EE[\zeta_t^{K_i}] \big) + \big( \tilde{\zeta}_t^{L_i} - \EE[ \tilde{\zeta}_t^{L_i} ] \big) \big] v_i^{(idy)} \Big\|
        \\
        &\hskip 15pt
        \leq  \frac{\ell d}{\tau^2} \sum_{t=0}^{T-1} \Big( \Big\| \frac{\gamma^t}{M} \sum_{i=1}^M  \big( \zeta_t^{K_i} - \EE[\zeta_t^{K_i}] \big) v_i^{(idy)} \Big\| + \Big\| \frac{\gamma^t}{M} \sum_{i=1}{M}  \big( \tilde{\zeta}_t^{L_i} - \EE[ \tilde{\zeta}_t^{L_i} ] \big) v_i^{(idy)} \Big\| \Big)
        \\
        &\hskip 15pt
        \leq \frac{\ell d}{\tau^2}.2T \tilde{\varepsilon}
    \end{align*}
    with probability at least $1 - 2T (\ell + d) \exp\big( \frac{- M \tilde{\varepsilon}^2 }{4 h_{\psi_2}^2(\tau, C_0) + \tilde{\varepsilon} h_{\psi_2}(\tau, C_0) } \big)$. A similar inequality holds for $\|  \tilde \nabla^{T, M, \tau}_L (\theta, \underline v) - \hat{\nabla}^{T, M, \tau}_L(\theta, \underline v)  \| $. Choosing $\tilde{\varepsilon}$ and $M$ so that
    $$
        \Big( \frac{4 \ell d T}{\tau^2} \Big).\tilde{\varepsilon} = \varepsilon,
        \quad
        4 T (\ell + d) \exp \Big( \frac{-M \tilde{\varepsilon}^2 }{4 h_{\psi_2}^2(\tau, C_0)  + \tilde{\varepsilon} h_{\psi_2}(\tau, C_0) } \Big) \leq \delta_{sample},
    $$
    or letting $M \geq \phi_{sample, size}(\varepsilon, \tau, T, C_0, \delta_{sample})$, we conclude that with probability at least $1 - \delta_{sample}$, $\| \tilde \nabla^{T, M, \tau}(\theta, \underline v) - \hat{\nabla}^{T, M, \tau}(\theta, \underline v) \| \leq \varepsilon$. 
\end{proof}

\subsection{Proof of Proposition~\ref{proposition:gradient_approx_with_MKV_simulator}}
\label{subsection:proof_of_proposition_approx_gradient_w_MKV_simulator}
In this subsection, we show that when the truncation horizon $T$ is large enough, the perturbation radius $\tau$ is small enough, and the perturbation size $M$ is large enough in Algorithm~\ref{algo:MKVestim}, the output of the approximated policy gradient $\tilde \nabla^{T, M, \tau}$ from the algorithm is close to the exact policy gradient $\nabla C$ with high probability.  We shall use the constants:
\begin{align}
    \phi_{pert, radius, MKV}( \varepsilon, C_0) &= \phi_{pert, radius}( \varepsilon / 3, C_0) 
    \label{eq:phi_MKV_pert_radius}
    \\
    \phi_{trunc, T, MKV}(\varepsilon, \tau,  C_0) &= \phi_{trunc, T}( \varepsilon / 3, \tau, C_0) + 2
    \label{eq:phi_MKV_trunc_T}
    \\
    \phi_{sample, size, MKV}(\varepsilon, \tau, T, C_0, \delta_{approx}) &= \phi_{pert, size} \big(\varepsilon / 3, \tau, C_0, \delta_{approx} / 2 \big) 
    \nonumber \\ 
    & \hspace{10pt} + \phi_{sample, size} \big(\varepsilon / 3, \tau, T, C_0, \delta_{approx} / 2\big) 
    \label{eq:phi_MKV_sample_size}
\end{align}
where $\phi_{pert, radius}$~\eqref{eq:phi_pert_radius}, $\phi_{pert, size}$~\eqref{eq:phi_pert_size}, $\phi_{trunc, T}$~\eqref{eq:phi_trunc_T}, $\phi_{sample, size}$~\eqref{eq:phi_sample_size} are polynomial in $(d,$ $\ell,$ $C_0,$ $(\lambda_y^1)^{-1},$ $(\lambda_z^0)^{-1},$ $C_{init, noise},$ $\varepsilon^{-1},$ $\tau^{-1})$ and model parameters defined in Sections~\ref{subsection:approximation_with_smoothed_cost},~\ref{subsection:approximation_C_with_C_T}, and~\ref{subsection:approx_sampled_costs_with_MKV_simulator}.\\

\begin{proof}
From Lemma~\ref{lemma:approx_perturbed_gradient_MKV_simulator}, Lemma~\ref{lemma:approx_truncated_policy_gradient_MKV_simulator}, and Lemma~\ref{lemma:approx_sampled_pg_with_truncated_pg}, we have 
\begin{align*}
    & \PP \big( \big\| \tilde \nabla^{T, M, \tau}(\theta, \underline v) - \nabla C(\theta) \big\|  > \varepsilon \big)
    \\
   &\hskip 15pt
   \le \PP \big( \big\| \hat{\nabla}^{M, \tau}(\theta, \underline v) - \nabla C(\theta) \big\| + \big\| \hat{\nabla}^{T, M, \tau}(\theta, \underline v) - \hat{\nabla}^{M, \tau}(\theta, \underline v) \big\| 
    \\
    & \hspace{25pt} + \big\| \tilde{\nabla}^{T,M, \tau}(\theta, \underline v) - \hat{\nabla}^{T, M, \tau}(\theta, \underline v) \big\| > \varepsilon \big)  
    \\
    &\hskip 15pt
    \le \PP \big(  \big\| \hat{\nabla}^{M, \tau}(\theta, \underline v) - \nabla C(\theta) \big\| > \varepsilon / 3 \big) + \PP\big( \big\| \hat{\nabla}^{T, M, \tau}(\theta, \underline v) - \hat{\nabla}^{M, \tau}(\theta, \underline v) \big\|  > \varepsilon / 3 \big)
    \\
    & \hspace{25pt} + \PP \big( \big\| \tilde{\nabla}^{T,M, \tau}(\theta, \underline v) - \hat{\nabla}^{T, M, \tau}(\theta, \underline v) \big\| > \varepsilon / 3 \big) 
    \\
    &\hskip 15pt
    \le\delta_{approx} / 2 + 0 + \delta_{approx} / 2
    \\
    &\hskip 15pt
    \le\delta_{approx},
\end{align*}
which is what needed to be proved.
\end{proof}

\section{\textbf{Proofs of technical results from Section~\ref{subsec:PG-popsimu}}}
\label{sec:app-proof-modelfree-POP-CV}

\subsection{Auxiliary costs for the $N$-agent problem}

In order to show that the model-free population-based gradient estimator $\tilde{\nabla}^{T,N,M, \tau}$ from Algorithm~\ref{algo:POPestim} approximates well the exact policy gradient $\nabla C$, like in Section~\ref{sec:app-proof-modelfree-MKV-CV} for the MKV simulator, we introduce the \emph{perturbed population policy gradient} $\hat{\nabla}^{N, M, \tau}$, an the \emph{truncated population policy gradient} $\hat{\nabla}^{T, N, M, \tau}$ by replacing the perturbed MF cost $C(\theta + v_i)$ in~\eqref{eq:def_perturbed_PG_MKV} and the truncated perturbed cost $C^T(\theta + v_i)$ in~\eqref{eq:def_truncated_gradient_MKV} by the corresponding $N-$agent social costs adapted from equation~\ref{eq:social_cost_of_population}. 

Let $[N] = \{1, \ldots, N\}$. For $\theta = (K, L) \in \Theta$, we introduce auxiliary processes $(\by^{K, n, N})_{n \in [N]}$ and $\bz^{L, N}$ for the $N$-agents by $ y_t^{K, n, N} = X_t^{(n), \theta} - \bar X_t^{N, \theta}$ and $z_t^{L, N} = \bar X_t^{N, \theta} = \frac{1}{N} \sum_{n=1}^N X_t^{(n), \theta}$ for all $t \geq 0$ and $n \in [N]$, where $(X_t^{(n), \theta})_{t\geq 0}$ is the state process for agent $n$ following dynamics~\eqref{fo:N-multi_state} with control at time $t$ given by $u_t^{(n), \theta} = -  K (X_t^{(n), \theta} - \bar{X}_t^{N, \theta}) - L \bar{X}_t^{N, \theta}$. The dynamics of $(\by^{K, n, N})_{n \in [N]}$ and $\bz^{L, N}$ are as follows
\begin{align}
    & y_{t+1}^{K, n, N} = ( \rA - \rB K ) y_t^{K, n, N} + \tilde{\varepsilon}_{t+1}^{(n)}, \quad && y_0^{n, N} = \varepsilon_0^{(n)} - \sum \nolimits_{n \in [N]} \varepsilon_0^{(n)} / N,
    \label{eq:def_dyn_y_K_n_N}
    \\
    & z_{t+1}^{L, N}  = ( \tilde \rA - \tilde \rB L ) z_{t}^{L, N} + \varepsilon_{t+1}^0 + \overline{{\varepsilon}_{t+1}^{N}}, \quad && z_0^{N} = \varepsilon_0^0 +  \sum \nolimits_{n\in [N]} \varepsilon_0^{1, (n)} / N,
    \label{eq:def_z_L_N}
\end{align}
where $\overline{{\varepsilon}_{t+1}^{N}} = \frac{1}{N} \sum_{n=1}^N \varepsilon_{t+1}^{(n)}$ 
and $\tilde{\varepsilon}_{t+1}^{(n)} = \varepsilon_{t+1}^{(n)} - \overline{{\varepsilon}_{t+1}^{N}}$ for $ t \geq 0$. We define the variance matrices $\Sigma_K^{y, n, N}$ and $\Sigma_L^{z, N}$ by replacing $y_t$ and $z_t$ in~\eqref{eq:variance_matrices_y_and_z} with $y_t^{K, n, N}$ and $z_t^{L, N}$, and define $\Sigma_K^{y, N} = \frac{1}{N} \sum_{n=1}^N \Sigma_K^{y, n, N}$. Similar modifications are applied to the initial perturbation variances $\Sigma_{y_0}^N = \frac{1}{N} \sum_{n=1}^N \EE[ y_0^{n, N} (y_0^{n,N})^\top ] $ and $\Sigma_{z_0}^N = \EE[ z_0^N (z_0^N)^\top ]$, and for the noise variances 
$
\Sigma^{1, N} = \frac{1}{N} \sum_{n=1}^N \EE\big[ \tilde{\varepsilon}_{t=1}^{(n)}. \big( \tilde{\varepsilon}_{t=1}^{(n)} \big)^\top \big]$ and $\Sigma^{0, N} = \EE\big[ \big(  \varepsilon_{t=1}^0 + \overline{{\varepsilon}_{t=1}^{N}} \big). \big(  \varepsilon_{t=1}^0 + \overline{{\varepsilon}_{t=1}^{N}} \big)^\top \big] $.  For the auxiliary costs in the $N-$agent problem, we adjust equation~\eqref{eq:def_Cy_Cz} with $(\by^{K, n, N})_{n\in [N]}$ and $\bz^{L, N}$ so that  
\begin{align*}
    C_y^{N}(K) & = \frac{1}{N} \sum \nolimits_{n \in [N]} \EE \big[ \sum \nolimits_{t\geq 0} \gamma^t (y_t^{K, n, N})^\top (Q + K^\top R K ) y_t^{K, n, N} \big]\\
    & = \langle Q + K^\top R K, \Sigma_{K}^{y, N} \rangle_{tr}
\end{align*}
and
\begin{align*}
    C_z^{N}(L) &= \EE \big[ \sum \nolimits_{t\geq 0} \gamma^t (z_t^{L, N})^\top ( \tilde{Q} + L^\top \tilde{R} L ) z_t^{L, N} \big]\\
    & = \langle \tilde{Q} + L^\top \tilde{R} L, \Sigma_{L}^{z, N} \rangle_{tr}.
\end{align*}
For $\theta = (K, L) \in \Theta$, $C^N(\theta) = C_y^N(K) + C_z^N(L)$. 
Because all agents adapt the same control parameter $\theta$, the $N-$agent auxiliary costs $C_y^N(K)$ and $C_z^{N}(L)$ can then be expressed using the matrices $P_K^y$ and $P_L^z$ solving the DLEs~\eqref{eq:lyapunov_eq_theta}:
\begin{equation}
\label{eq:N_agent_Cy_Cz_expression_PK_PL}
    C_y^{N}(K) = \langle P_K^y, \, \Sigma_{y_0}^N + \frac{\gamma}{1 - \gamma} \Sigma^{1, N} \rangle_{tr}, 
    \quad
    C_z^N(L) = \langle P_L^z, \, \Sigma_{z_0}^N + \frac{\gamma}{1 - \gamma} \Sigma^{0, N} \rangle_{tr}.
\end{equation}
The truncated versions at horizon $T$ for the variances and the costs $\Sigma_K^{y, N, T}$, $\Sigma_{L}^{z, N, T}$, $C^{N, T}_y(K)$, $C^{N, T}_z(L)$, $C^{N, T}(\theta)$, and the sampled truncated costs with population simulator $\tilde C_y^{N, T}(K), \tilde C_z^{N, T}(L), \tilde C^{N, T}(\theta)$ are defined accordingly. 

By definition of the matrices $\Sigma^{1,N}$, $\Sigma^{0, N}$, $\Sigma_{y_0}^{N}$ and $\Sigma_{z_0}^N$, we have
\begin{align}
    \Sigma^{1, N} &= (1 - 1 / N ) \Sigma^1, & \Sigma^{0, N} &= \Sigma^0 + \Sigma^1 / N,
    \label{eq:variance_noise_N_agent}
    \\
    \Sigma_{y_0}^{N} & = ( 1 - 1 / N)\Sigma_{y_0}, & \Sigma_{z_0}^{N} &=  \Sigma_{z_0} + \Sigma_{y_0} / N.
    \label{eq:variance_init_N_agent}
\end{align}
Moreover, because the sub-gaussian norms of $y_0$, $z_0$ and the noise processes are bounded by $C_{init,noise}$, we have for all $t \geq 0$ and $n \in [N]$,
\begin{equation}
\label{eq:bound_phi2_y0_z0_N_agent}
    \max \Big\{ \| y_0^{n, N} \|_{\psi_2}, \| z_0^{N} \|_{\psi_2}, 
     \| \tilde{\varepsilon}_{t+1}^{(n)} \|_{\psi_2}, \big\| \varepsilon_{t+1}^0 + \overline{{\varepsilon}_{t+1}^{N}} \big\|_{\psi_2}  \Big\}  \leq 2 C_{init, noise}.
\end{equation}
We define the minimum eigenvalues of the initial and noise variance matrices by:
\begin{equation}
    \lambda_y^{1,N} = \lambda_{min}\big( \lambda_{y_0}^N + \gamma \Sigma^{1, N } / (1 - \gamma) \big), \quad \lambda_{z_0}^N = \lambda_{min}\big( \lambda_{z_0}^N + \gamma \Sigma^{0, N} / ( 1- \gamma) \big).
\end{equation}

\subsection{Proof of Proposition~\ref{proposition:approx_with_social_cost_to_MF_cost} }
\label{subsection:proof_of_approx_modelfree_pop_gradient}

\begin{proof}
From the expressions~\eqref{eq:N_agent_Cy_Cz_expression_PK_PL} for the auxiliary costs $C_y^N$ and $C_z^N$ for the $N$-agent problem, and the expressions~\eqref{eq:cost_expression_with_PK_y} for $C_y$ and $C_z$ for the MFC problem, we have 
\begin{align*}
    | C^N(\theta) - C(\theta) | & = \Big| \frac{1}{N} \langle P_L^z - P_K^y, \, \frac{1}{N} \big( \Sigma_{y_0} + \frac{\gamma}{1- \gamma} \Sigma^1 \big) \rangle_{tr} \Big| 
    \\
    & \leq \frac{1}{N} ( \|P_K^y\| + \|P_L^z \|). \frac{ \max \big\{ Tr(\Sigma_{y_o}), Tr(\Sigma^1) \big\} }{1-\gamma} .
\end{align*}
Using $Tr(\Sigma_{y_0}) = \EE[ \| y_0 \|^2 ] \leq 2d C_{init,noise}^2$ and $Tr(\Sigma^1) \leq 2 d C_{init, noise}^2$, and the bounds~\eqref{eq:upper_bound_Riccati_and_Variance_matrix}, we conclude by setting 
$
    h_{social-cost, factor}(C_0) = \frac{2 d C_{init,noise}^2 C_0}{(1 - \gamma)} \Big( \frac{1}{\lambda_y^1} + \frac{1}{\lambda_z^0} \Big).
$
\end{proof}

Using the same arguments as in Section~\ref{sec:app-proof-modelfree-MKV-CV}, we now show that the output of Algorithm~\ref{algo:POPestim} is close enough to the exact policy gradient $\nabla C$. In the following, we only focus on the differences in the corresponding proof arguments when we use a social cost instead of an MF cost. We shall use the following constants.
\begin{align}
    \phi_{agent, size, pop}(\varepsilon, \tau, C_0) 
    &=  \big( \ell d / (\varepsilon \tau) \big). h_{social-cost, factor}(C_0) 
    \nonumber \\
    & = \frac{2 \ell d^2 C_{init, noise}^2 C_0}{\tau \varepsilon ( 1- \gamma)}  \Big( \frac{1}{\lambda_y^1} + \frac{1}{\lambda_z^0} \Big)
    \label{eq:phi_agent_size_N}
\end{align}

\begin{align}
\phi_{sample,size, pop}(\varepsilon, \tau, T, C_0, \delta_{sample}) &= \phi_{sample, size}(\varepsilon / 4, \tau, T, C_0, \delta_{sample}) 
    \nonumber \\
    & =  \log \Big( \frac{4T (\ell + d)}{\delta_{sample}} \Big).\Big(  \frac{32 \ell d T}{\varepsilon \tau^2}  (  h_{\psi_2}(\tau, C_0) )  + \frac{1}{4} \Big)^2
    \label{eq:phi_sample_size_N-agent}
\end{align}
\begin{align}
    \phi_{trunc, T, pop}(\varepsilon, \tau, C_0, N)& = \phi_{trunc, T}(\varepsilon / (1 + 1/N), \tau, C_0) 
    \nonumber \\
    &\hskip -25pt = \Big( \frac{1}{\log(1 / \gamma_{pert})} \Big)^2 \Big[ \log \Big( \frac{4 \ell d C_{init, noise}^2 h_T(C_0)}{\varepsilon \tau ( 1- \gamma_{pert})^2} (1 + 1 / N) \Big) + 1 \Big]^2
    \label{eq:phi_trunc_T_N_agent}
\end{align}
where $h_{social-cost, factor}(C_0)$ is defined in Proposition~\ref{proposition:approx_with_social_cost_to_MF_cost} , and $\phi_{sample, size}$ and $\phi_{trunc, T}$ are defined in equations~\eqref{eq:phi_sample_size} and~\eqref{eq:phi_trunc_T}.

\begin{lemma}
\label{lemma:approx_gradient_pop_simulator}
If $\theta \in \Theta$ is such that $C(\theta) \leq C_0$ for some $C_0 \in \RR$, $\varepsilon > 0$ and $\delta_{approx} \in (0,1)$, and if 
    \begin{align*}
        \tau^{-1} &\geq \phi_{pert, radius}(\varepsilon, C_0)
        \\
        N &\geq \phi_{agent, size}^N(\varepsilon, \tau, C_0)
        \\
        T & \geq \phi_{trunc, T}^N(\varepsilon, \tau, C_0, N)
        \\
        M & \geq \max \big\{ \phi_{pert,size}(\varepsilon, \tau, C_0, \delta_{approx}),  \phi_{sample, size}^N( \varepsilon, \tau, T, C_0, \delta_{approx})  \big\},
    \end{align*}
    then with $M$ perturbation directions $\underline v = (v_i^{(idy)}, v_i^{(com)})_{i=1}^M$ on $\SS_\tau \times \SS_\tau$, we have
    \begin{align}
         & \PP \big( \| \tilde{\nabla}^{T, N, M, \tau}(\theta, \underline v) - \hat{\nabla}^{T, N, M, \tau}(\theta, \underline v) \| > \varepsilon \big) \leq \delta_{approx},
         \label{eq:ineq_approx_sample_gradient_pop}
        \\
        & \| \hat{\nabla}^{T, N, M, \tau}(\theta, \underline v) - \hat{\nabla}^{N, M, \tau}(\theta, \underline v) \|  \leq \varepsilon,
        \label{eq:ineq_approx_truncated_gradient_pop}
        \\
        & \| \hat{\nabla}^{N, M, \tau}(\theta, \underline v) - \hat{\nabla}^{M, \tau} (\theta, \underline v) \|  \leq \varepsilon.
        \label{eq:ineq_approx_perturbed_gradient_pop}
    \end{align}
\end{lemma}

\begin{proof}
We first show the inequality~\eqref{eq:ineq_approx_sample_gradient_pop}. For $\theta = (K, L) \in \Theta$ and $K_i = K + v_i^{(idy)}$, from the bounds in~\eqref{eq:bound_phi2_y0_z0_N_agent} and Lemma~\ref{lemma:concentration_ineq_for_yt_K}, we deduce that 
$$
\gamma^t \| y_t^{K_i, n, N} \|_{\psi_2}^2 \leq \frac{4 d^2 C_{init, noise}^2}{( 1- \sqrt{\gamma})^2}.
$$
So by setting $\zeta_t^{K_i, n, N} = (y_t^{K_i, n, N})^\top ( Q + K_i^\top R K_i) y_t^{K_i, n, N}$, using the same arguments as in Lemma~\ref{lemma:concentration_ineq_for_yt_K}, we get 
\begin{align*}
    & \PP \Big( \big\| \frac{\gamma^t}{M} \sum_{i=1}^M \big( \frac{1}{N} \sum_{n=1}^N ( \zeta_t^{K_i, n, N} - \EE[\zeta_t^{K_i, n, N} ] \big) v_i^{(idy)}  \big\| \geq \varepsilon \Big) 
    \\
    &\hskip 35pt\leq  (\ell + d ) \exp \Big( (- M \varepsilon^2 ) / (64 h_{\psi_2, K}^2(\tau, C_0) + 8 \varepsilon h_{\psi_2, K}(\tau, C_0) ) \Big).
\end{align*}
Thus, like in Lemma~\ref{lemma:approx_sampled_pg_with_truncated_pg}, we obtain~\eqref{eq:ineq_approx_sample_gradient_pop} if $M \geq \phi_{sample, size}^N( \varepsilon, \tau, T, C_0, \delta_{approx})$.

Now we show the inequality~\eqref{eq:ineq_approx_truncated_gradient_pop} for the truncated gradient with social cost. We first derive a bound similar to the one in Lemma~\ref{lemma:approx_truncated_cost}, namely: 
\begin{align*}
    \big\| \Sigma_K^{y, N} - \Sigma_K^{y,T, N} \big\|
    & \leq  \Big( \frac{\gamma_{pert}^T}{1- \gamma_{perc}} + \gamma_{perc}^T + (T-1) \gamma_{perc}^T \Big)  \frac{\max\big\{  \| \Sigma_{y_0}^N \|, \| \Sigma^{1, N} \| \big\} }{1- \gamma_{perc}}   
    \\
    & \leq ( 1 + 1 / N)  2d \gamma_{perc}^T (T +1) C_{init, noise}^2  / (1- \gamma_{perc})^2,
\end{align*}
and then, like in Lemma~\ref{lemma:approx_truncated_policy_gradient_MKV_simulator}, we obtain inequality~\eqref{eq:ineq_approx_truncated_gradient_pop}.

Finally, to show inequality~\eqref{eq:ineq_approx_perturbed_gradient_pop} for the perturbed policy gradient with social cost, we notice that 
$$
    \| \hat{\nabla}^{N, M, \tau}(\theta, \underline v) - \hat{\nabla}^{M, \tau} (\theta, \underline v) \|
    \leq \frac{\ell d}{M \tau} \sum_{i=1}^M | C^{N}(\theta_i) - C(\theta_i) |. \| v_i^{(idy)} \|, 
$$
and we conclude using Proposition~\ref{proposition:approx_with_social_cost_to_MF_cost}.
\end{proof}

\subsection{Proof of Proposition~\ref{proposition:modelfree_population_gradient_approx}}
\begin{proof}
This is an immediate consequence of Lemma~\ref{lemma:approx_gradient_pop_simulator}, Lemma~\ref{lemma:approx_perturbed_gradient_MKV_simulator} and: 
    \begin{align*}
        & \big\| \tilde \nabla^{T, N, M, \tau}(\theta, \underline v) - \nabla C(\theta) \big\| 
        \\
        &\hskip 25pt
        \le \| \tilde{\nabla}^{N, T, M, \tau}(\theta, \underline v) - \hat{\nabla}^{N, T, M, \tau}(\theta, \underline v) \|
        + 
        \big\| \hat{\nabla}^{N, T, M, \tau}(\theta, \underline v) - \hat{\nabla}^{N, M, \tau}(\theta, \underline v) \big\|
        \\
        &\hskip 45pt
        + \big\| \hat{\nabla}^{N, M, \tau}(\theta, \underline v) - \hat{\nabla}^{M, \tau}(\theta, \underline v) \big\|
        + \big\| \hat{\nabla}^{M, \tau}(\theta, \underline v) - \nabla C(\theta) \|.
\end{align*}
\end{proof}

\section{\textbf{Numerical Results}}
\label{sec:numerics}

Although the main thrust of our work is theoretical, we demonstrate the accuracy of our convergence results by numerical tests. Also, in order to illustrate the \emph{robustness} of our computational algorithms, we depart from the strict exchangeability assumption of the mean field models and consider  cases heterogeneous agents by perturbing the coefficients of the mean field cost: to be specific, for a set of $N$ agents, we shall use for agent $n$, a one-step cost function $c^{(n)}(x, \bar x, u, \bar u)$  which deviates from equation~\eqref{eq:lq_one_step_cost} by using a coefficient $Q^{(n)} \succ 0$ such that 
$$
    \| Q^{(n)} - Q \| \leq \tilde h
$$
for $n = 1, \ldots, N$. We call the factor $\tilde h$ is called the \emph{degree of heterogeneity} of the $N-$agent problem.
In this section, we consider a numerical example with $1$-dimensional state and action spaces because it is easier to visualize and allows us to focus on the main difficulty of MFC problems, namely the MF interactions. 

\subsection{Details on the numerical results}

\begin{table}[!h]
\begin{center}
    \begin{adjustbox}{width=0.5\columnwidth,center}
        \begin{tabular}{cccccccccc}
            
            \multicolumn{10}{c}{Model parameters}\\
            
            \hline \hline
            
            $a$ & $\overline{a}$ & $b$ & $\overline{b}$ & $q$ & $\overline{q}$ & $r$ & $\overline{r}$ & $\gamma$ & $\tilde{h}$ \\
            
            \hline 
            
            0.5 & 0.5 & 0.5 & 0.5 & 0.5 & 0.5 & 0.5 & 0.5 & 0.9 & 0.1 \\
            
            \hline \hline
            
            \\
            
            \multicolumn{10}{c}{Initial distribution and noise processes}\\
            
            \hline \hline
            
            & \multicolumn{2}{c}{$\varepsilon_0^0$} &  \multicolumn{2}{c}{$\varepsilon_0$} &  \multicolumn{2}{c}{$\varepsilon^0_t$} &  \multicolumn{2}{c}{$\varepsilon_t$} & \\
            
            \hline
            
            & \multicolumn{2}{c}{$\mathcal{U}([-1, 1])$} & \multicolumn{2}{c}{$\mathcal{U}([-1, 1])$} & \multicolumn{2}{c}{$\mathcal{N}(0, 0.01)$} & \multicolumn{2}{c}{$\mathcal{N}(0, 0.01)$} &\\
            
            \hline \hline
            
            \\
            
            \multicolumn{10}{c}{Learning parameters}\\
            
            \hline \hline
            
            & $T$ & $M$ & $\tau$  & $\eta$ & $\beta_1$ & $\beta_2$ & $K_0$ & $L_0$ & \\
            
            \hline
            
            & 50 & 10000 & 0.1 & 0.01 & 0.9 & 0.999 & 0.0 & 0.0 \\
            
            \hline \hline
            
        \end{tabular}
        
    \end{adjustbox}
\end{center}
\caption{Simulation parameters for the numerical computations.} 
\label{tab:simulation_parameters}    
\end{table}

Recall that the MFC dynamics and the one-step cost function are given respectively by~\eqref{eq:lq_sys_function_F} and~\eqref{eq:lq_one_step_cost}, with parameters denoted, in the one-dimensional case, by $a, \bar a, b, \bar b$ for the dynamics, and $q, \bar q, r, \bar r$ for the cost. The discount factor is still denoted by $\gamma$, and the perturbed cost coefficients for $N$ heterogeneous agents are denoted by $( q^{(n)} )_{n=1,\ldots, N}$. The parameters of the model-free PG method are $(T, M, \tau)$. 

The numerical results presented below are based on the numerical values specified in Table~\ref{tab:simulation_parameters}.  
The randomness in the state process is given in the second part of Table~\ref{tab:simulation_parameters}, where $\mathcal{U}([a,b])$ stands for the uniform distribution on interval $[a,b]$ and $\mathcal{N}(\mu, \sigma^2)$ is the one-dimensional Gaussian distribution with mean $\mu$ and standard deviation $\sigma$. 
Note that in the simulations, the social cost for the $N$-agent problem takes into account the heterogeneity from the agents by randomly sampling the model coefficient $q^{(n)}$ for agent $n$ using a uniform distribution $\cU( [ q - \tilde h, q + \tilde h])$.

We used the Adam optimizer (which provided slightly smoother results than a constant learning rate $\eta$ as in our theoretical analysis) initialized with initial learning rate $\eta$ and exponential decay rates given by $\beta_1$ and $\beta_2$. The number of perturbation directions is $M$, and their radius is $\tau$. These values were chosen based on the theoretical bounds we found, and further tuned after a few tests using the exact PG method.  

The computations have been run on one 2.4 GHz Skylake processor. For the parameters described here, the model-free PG with MKV simulator took roughly $10$ hour for $5000$ iterations. For the same number of iterations, the model-free PG with $N-$agent population simulator took roughly $48$ hours for $N=10$. The displayed curves correspond to the mean over $10$ different realizations. The shaded region, if any, corresponds to the mean $\pm$ standard deviation.

\subsection{Simulation with model-free PG algorithms}

\begin{figure}[!hbtp]
	\centering
		\begin{subfigure}{0.4\columnwidth}
			\centering
			\includegraphics[width=\columnwidth]{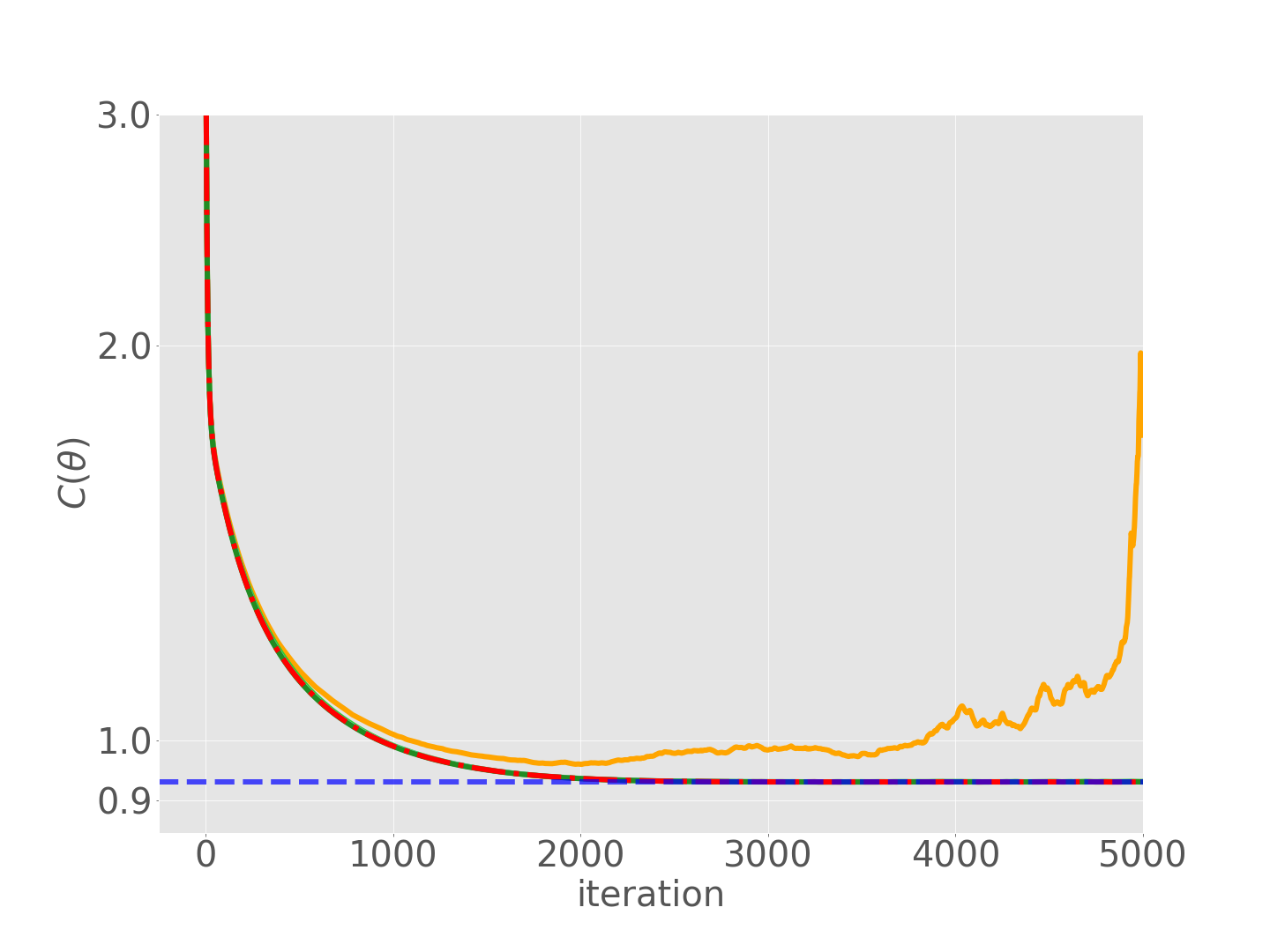}
			\caption{$C(\theta)$}
			\label{fig:1d-adam-cost-cv-MF}
		\end{subfigure}%
		\hspace{-0.3cm}%
		\begin{subfigure}{0.4\columnwidth}
			\centering
			\includegraphics[width=\columnwidth]{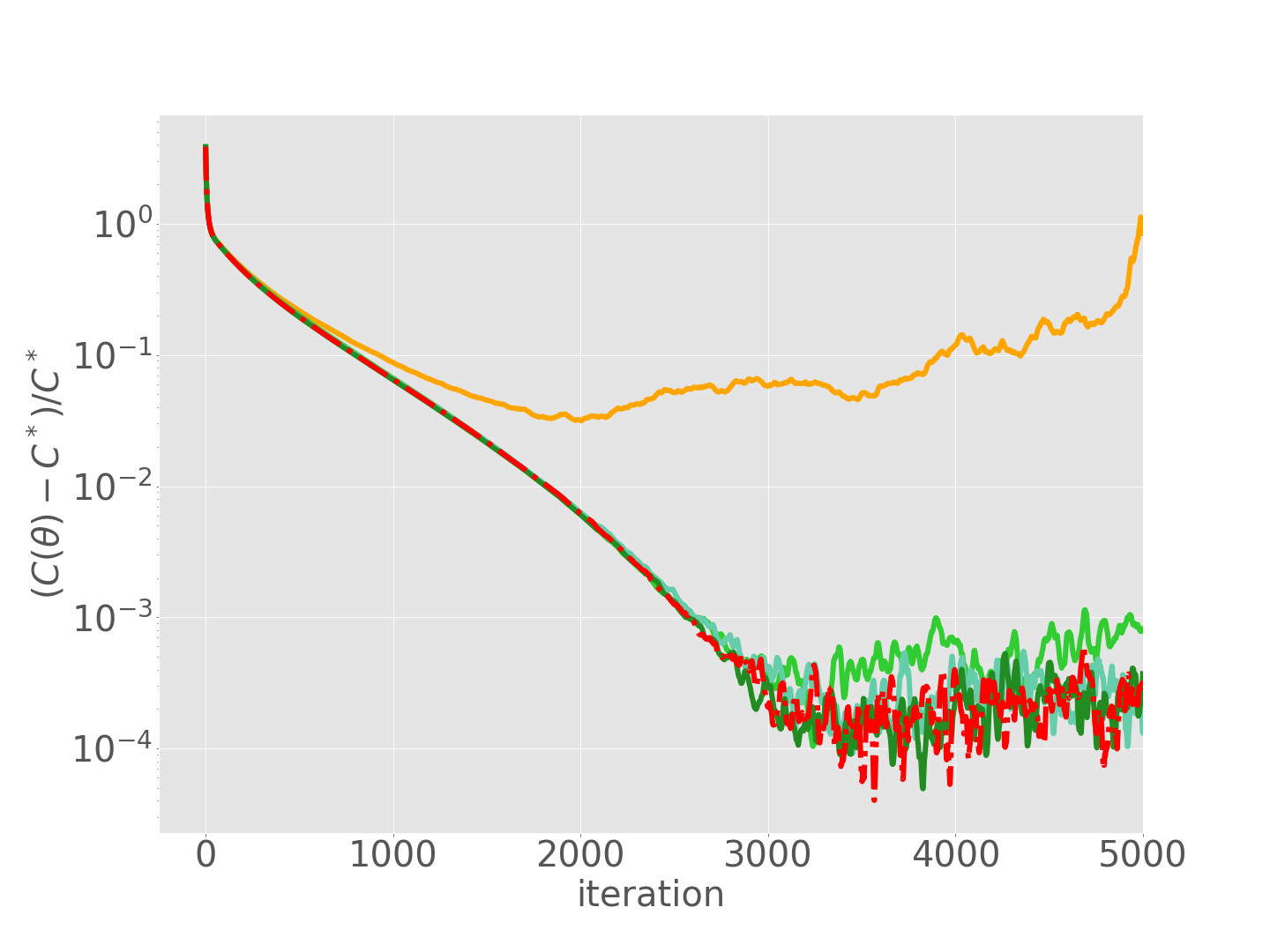} 
			\caption{$ \frac{ C(\theta) - C^* }{ C^* }$}%
			\label{fig:1d-adam-cost-difference-MF}
		\end{subfigure}%
		\hspace{-0.3cm}%
		\begin{subfigure}{.4\columnwidth}
			\centering
			\includegraphics[width=\columnwidth]{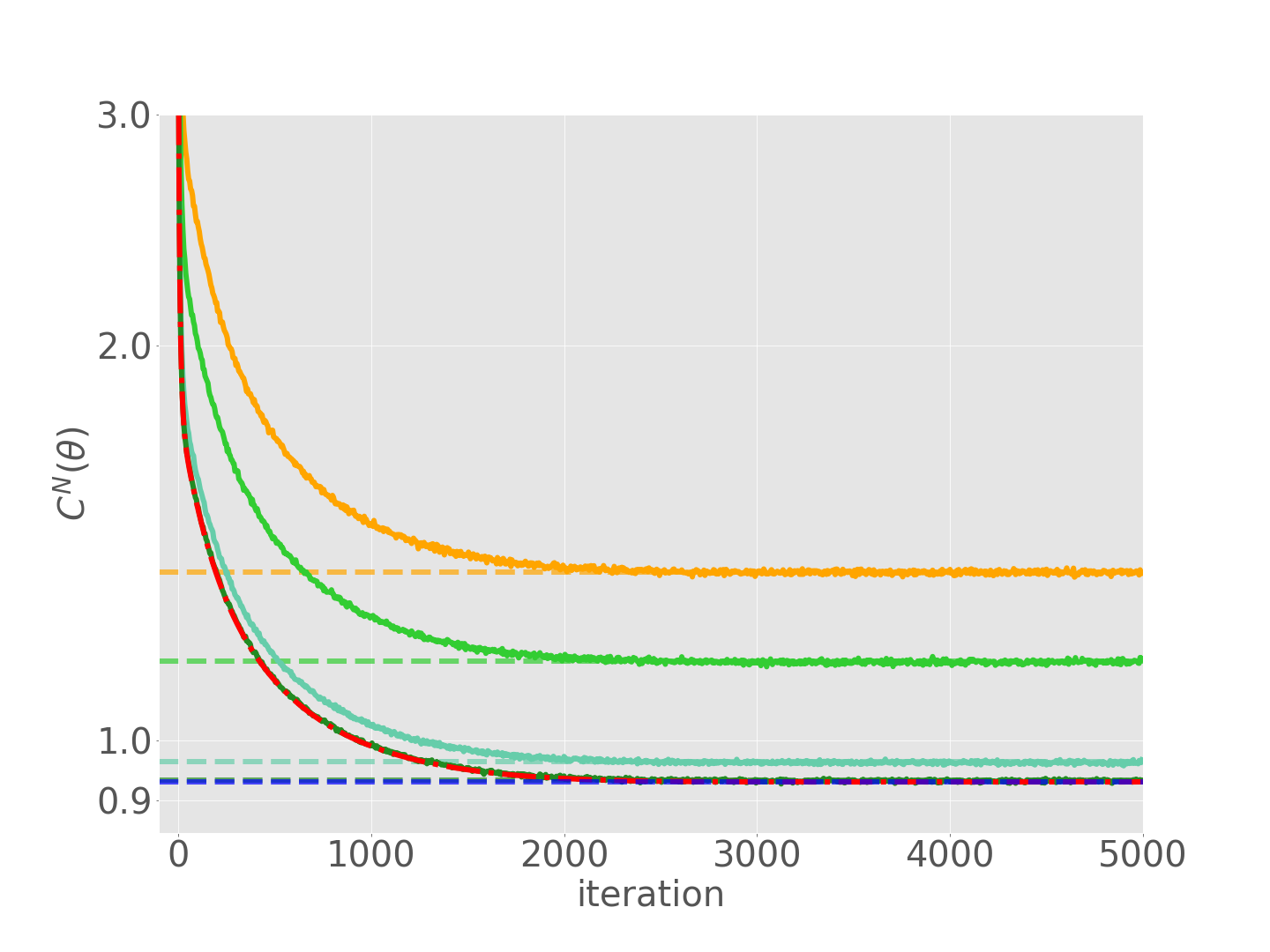}
			\caption{$C^N(\theta)$}
			\label{fig:1d-adam-cost-cv-N}
		\end{subfigure}%
		\hspace{-0.3cm}%
		\begin{subfigure}{.4\columnwidth}
			\centering %
			\includegraphics[width=\columnwidth]{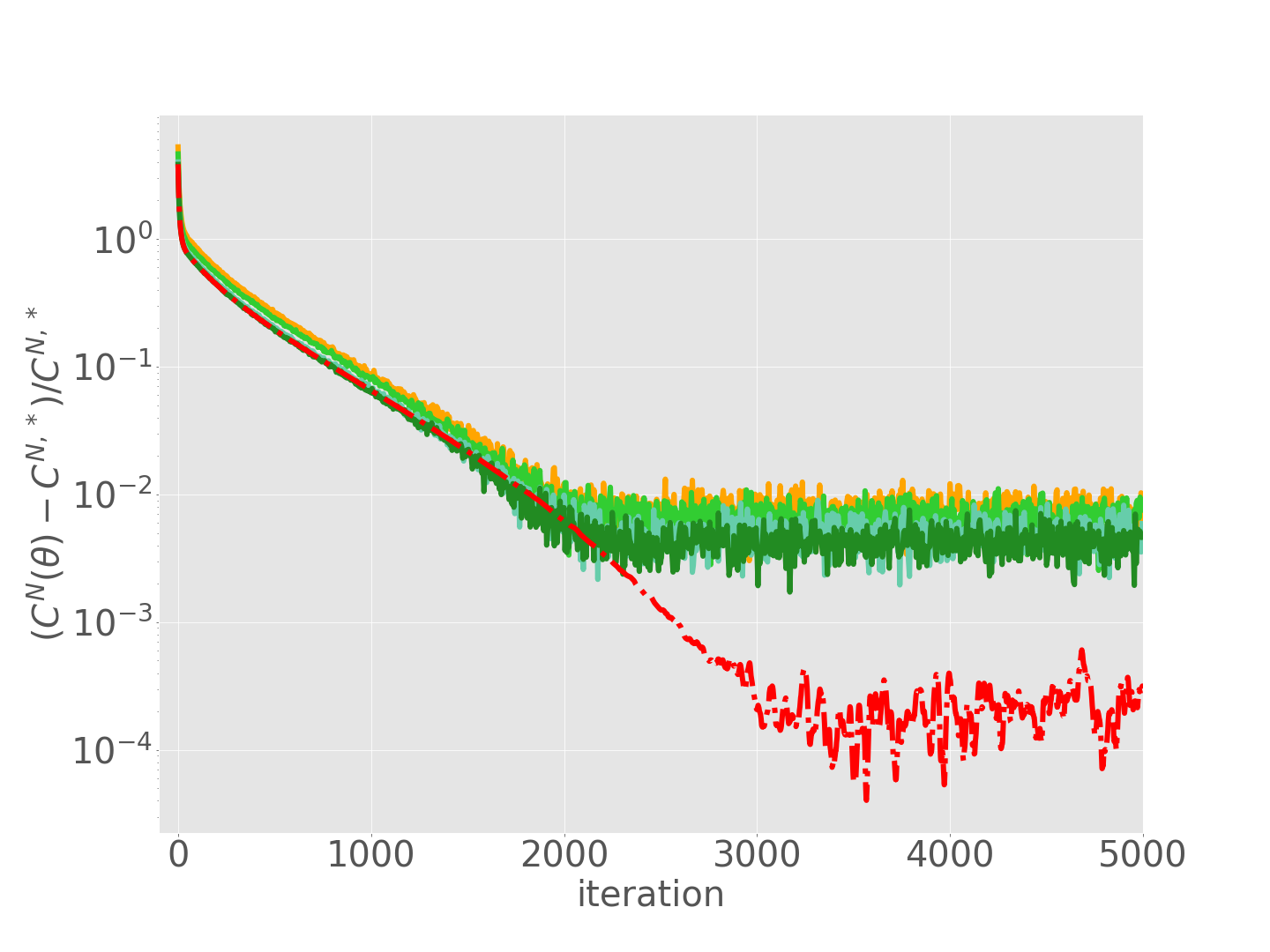}
			\caption{$\frac{C^{N}(\theta) - C^{*,N}}{C^{*,N}}$} %
			\label{fig:1d-adam-cost-difference-N}
		\end{subfigure}%
	\begin{subfigure}{0.15\columnwidth}
		\includegraphics[width=\textwidth]{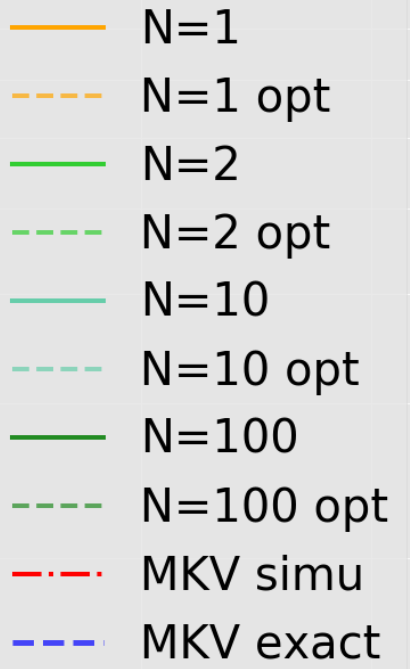}
	\end{subfigure}
	\caption{Comparison of costs achieved by the controls learned with MKV and population simulator using model-free PG. 
	(a)~:~MF costs; 
	(b)~:~relative error in MF cost w.r.t. MF optimum; 
	(c)~:~$N$-agent costs and their respective optimal value;
	(d)~:~relative error in $N$-agent costs w.r.t. $N$-agent optimal cost.
	}
    \label{fig:1d-adam-cost-MF-N}
\end{figure}

We investigate the performance of control parameters learned using an MKV simulator and a population simulator, on both the MF cost and the associated $N-$agent social cost. The results are shown in Figure~\ref{fig:1d-adam-cost-MF-N}.

\begin{figure}[!hbtp]
	\begin{subfigure}{.4\columnwidth}
		\centering
		\includegraphics[width=1.05\columnwidth]{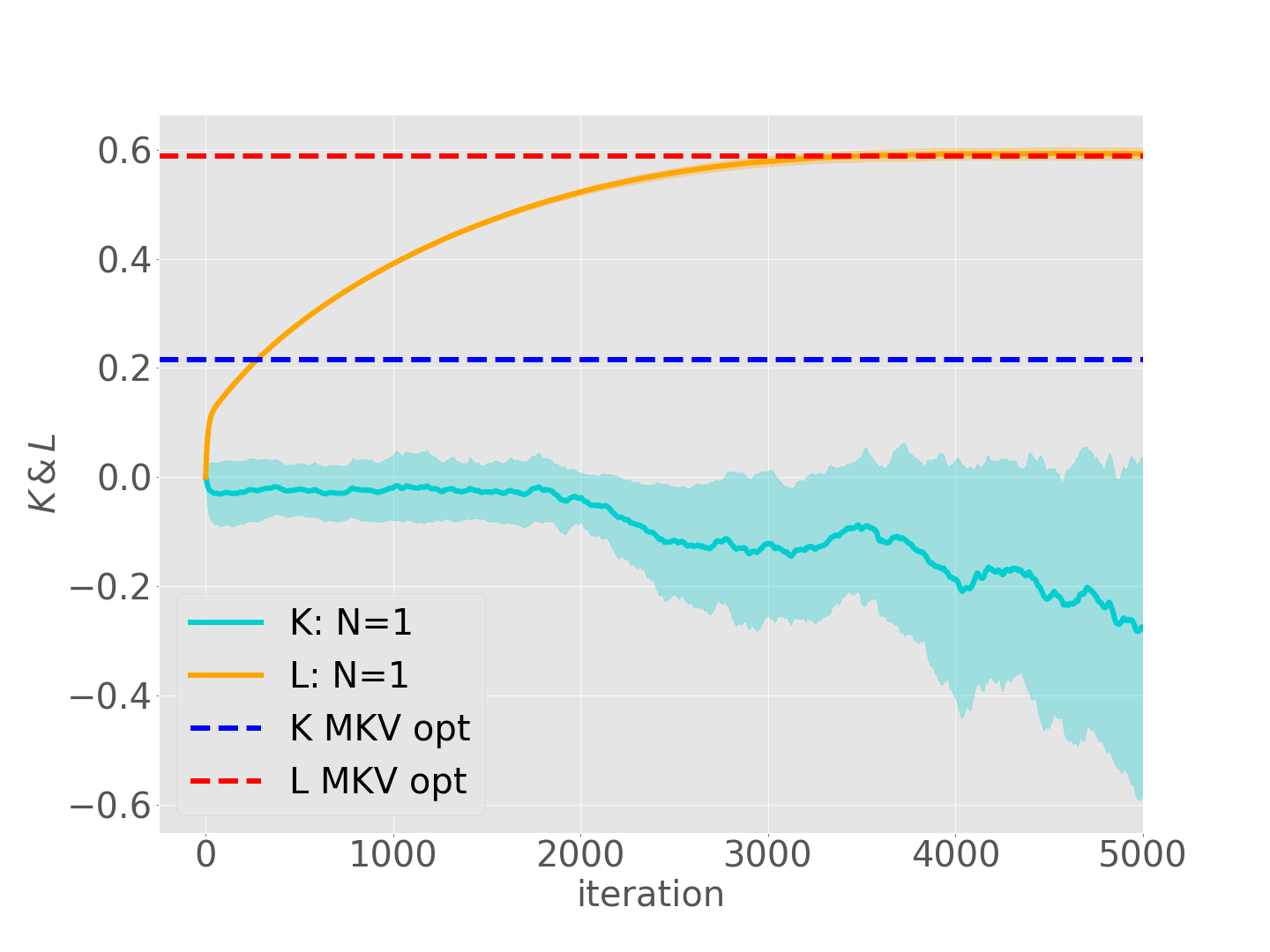}
		\caption{\,}
		\label{fig:1d-adam-controls-N1}
	\end{subfigure}%
	\begin{subfigure}{.4\columnwidth}
		\centering %
		\includegraphics[width=1.05\columnwidth]{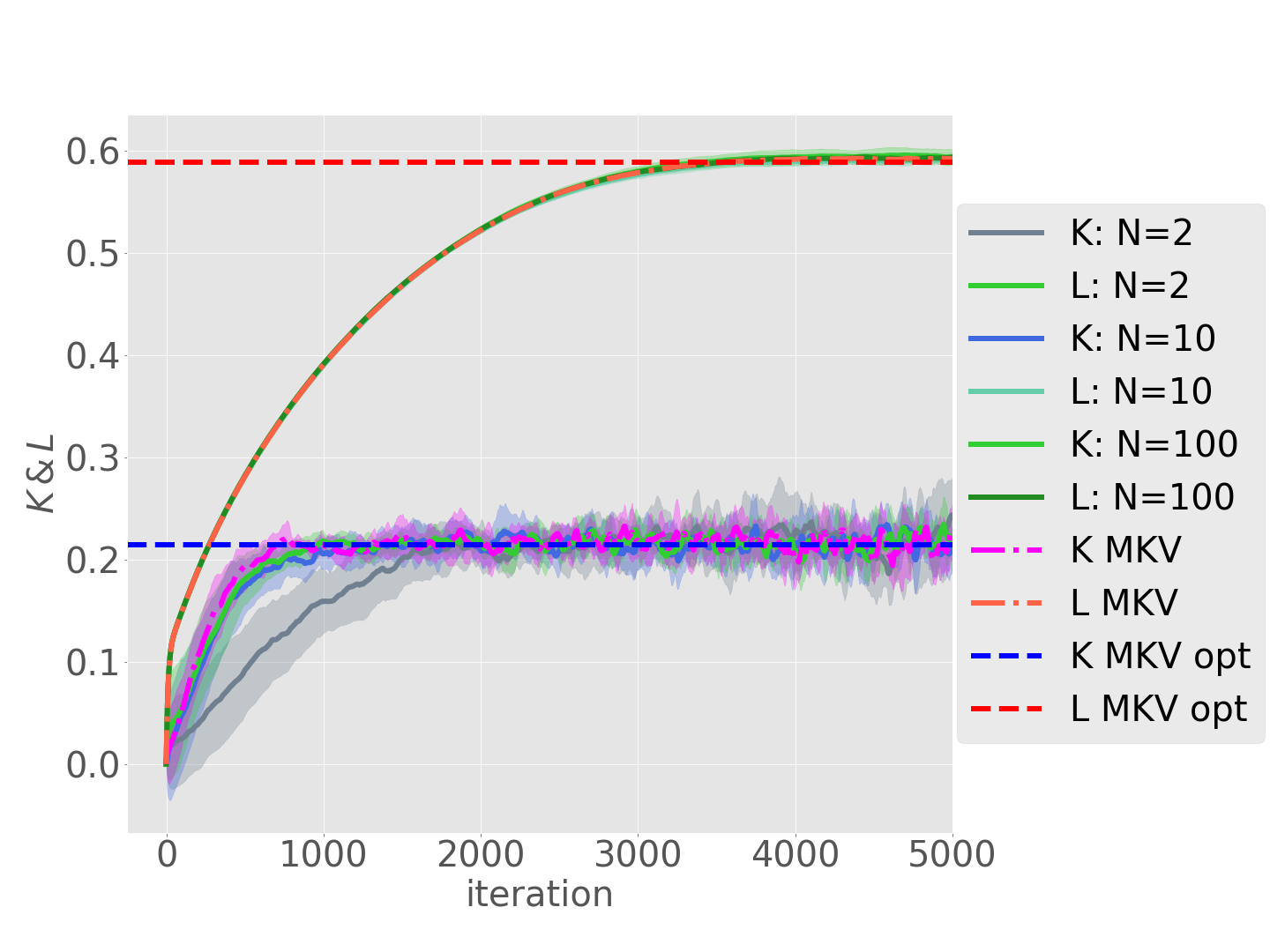}
		\caption{\,}
		\label{fig:1d-adam-controls-Nlarge}
	\end{subfigure}
	\caption{Convergence of the control parameters. (a)~$1$-agent simulator. (b)~MKV and population simulators.}
	\label{fig:global}
\end{figure}

In Figure~\eqref{fig:1d-adam-cost-cv-MF} and~\eqref{fig:1d-adam-cost-difference-MF}, one first sees that the MF costs evaluated with the learned parameters using either an MKV simulator or a population simulator approximate the optimal MF cost (blue line), provided that the number of agents is strictly larger than $1$. 
Indeed, as shown in Figure~\eqref{fig:1d-adam-controls-N1}, a single agent is able to learn the second component $L^*$ of the optimal control parameter $\theta^* = (K^*, L^*)$ but not the first component. This observation can be explained by the fact that $X^{(1)}_t - \bar X^N_t = 0$ for $t \geq 0$ when $N=1$. So, the agent does not have any estimate of the non-common component of the state (see~\eqref{eq:def-Y-Z-fct-X}) and is thus unable to learn the associated control parameter $K$. Figure~\eqref{fig:1d-adam-cost-difference-MF} displays the relative error of the MF cost, computed by $(C(\theta^{(k)}) - C^*) / C^*$ at iteration step $k$ where $C^* = C(\theta^*)$ is the optimal MF cost.
In addition to the fact that a single agent cannot learn $K^*$, we note in Figure~\eqref{fig:1d-adam-controls-Nlarge} that for $N=2$ the convergence is almost as good as with $N \ge 10$ as with MKV simulator, but it is slower for $N=2$. This certainly comes from the propagation of chaos, and we expect some dependence on the dimension of the state in accordance with our theoretical bounds.

Figures~\eqref{fig:1d-adam-cost-cv-N} and~\eqref{fig:1d-adam-cost-difference-N} display the values of the $N-$agent social cost and the error with the optimal social cost as functions of the number of iterations, for a population size $N=1,2, 10, 100$. 
In Figure~\eqref{fig:1d-adam-cost-cv-N}, the dashed line for each $N$ corresponds to the value of $C^{*,N} = J^N(\Phi^{*,N})$, which is the optimal $N-$agent social cost (i.e., the social cost evaluated at the optimal $N-$agent control $\Phi^{*,N}$ introduced in Section~\ref{subsec:OC-N-Agents}). Figure~\eqref{fig:1d-adam-cost-difference-N} shows the relative error, defined as $(C^N(\theta^{(k)}) - C^{*,N}) / C^{*,N}$ at iteration $k$. Here, we can see that the $N$ agents manage to learn a control that is approximately optimal for \emph{their} own social cost (even when $N=1$).

\subsection{Approximate optimality of the MFC for the $N-$agent problem}
\label{sec:approx-opt-MKV}

Starting from the optimal control parameter $\theta^* = (K^*, L^*)$ for the MFC problem, we consider a control parameter $\Phi^{*, N}_{MKV}$ for the control vector $\underline{U}^{N, MKV, *}_t \in \RR^{\ell N}$ and the state vector $\underline{X}_t \in \RR^{dN}$ in the dynamics~\eqref{fo:N-vector_state} when all the agents adopt the same control parameter $\theta^*$, namely $ \underline{U}^{N, MKV, *}_t =  \Phi^{*,N}_{MKV} \underline{X}_t$ for $t \geq 0$ with
$$
    \Phi^{*,N}_{MKV} = -I_N \otimes K^*   -  \frac{1}{N} \bone_N \otimes (L^* - K^*).
$$
We provide numerical evidence showing that even though the optimal control $\underline{U}^{*, N}_t = \Phi^{*, N} \underline{X}_t$ of the $N$ agents problem differs from the control vector $\underline{U}^{N, MKV, *}_t$, the latter provides an approximately optimal social cost, where the quality of the approximation depends on $N$ and on the degree of heterogeneity $\tilde h$ (recall that the variations $q^{(n)} - q$ are of size at most $\tilde h$).

\begin{figure}[h!]
	\centering
	\begin{subfigure}{.5\textwidth}
		\centering
		\includegraphics[width=1.0\textwidth]{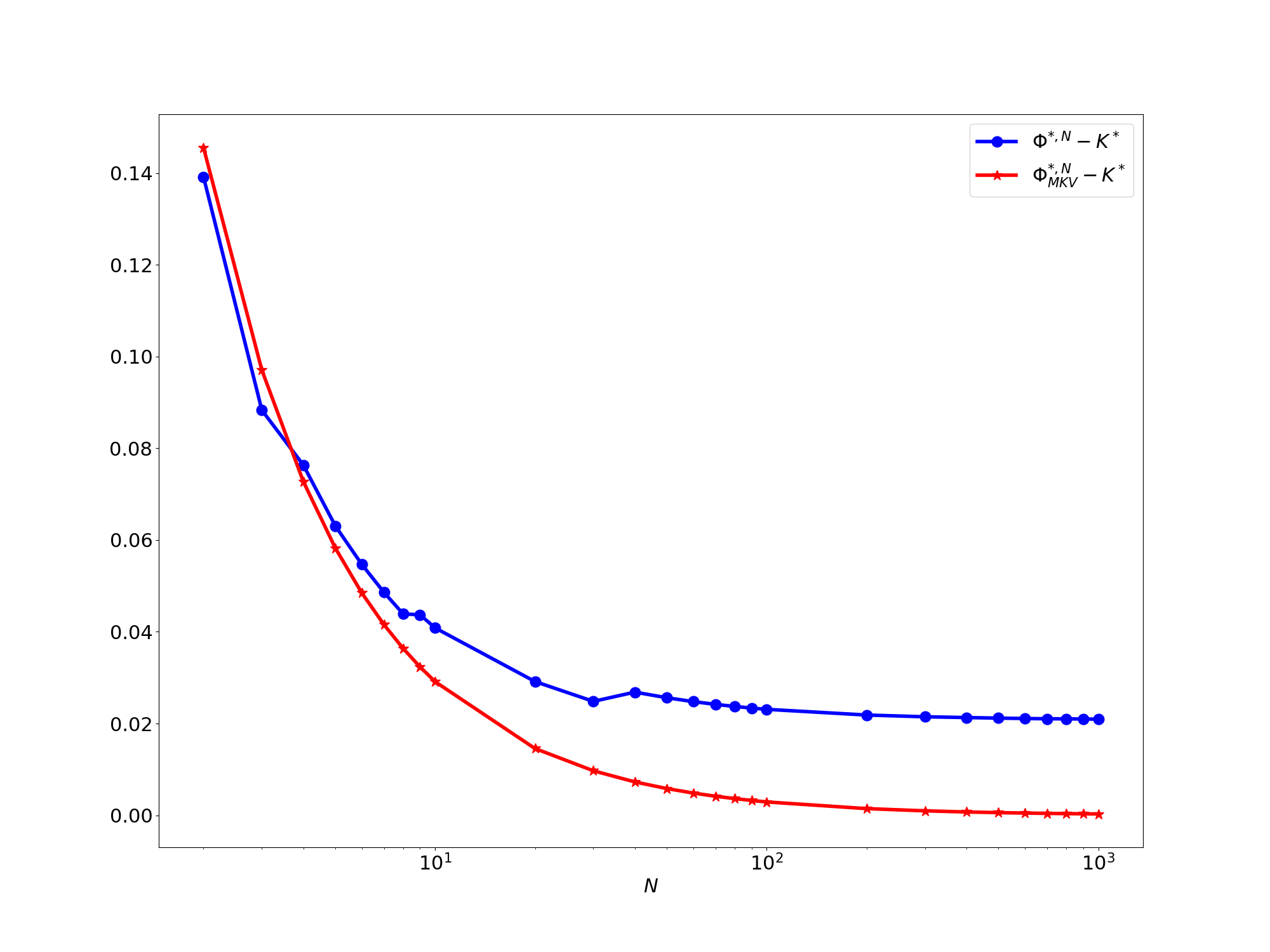}
		\caption{}
		\label{fig:cmp-N-MKV-opt-diagterm}
	\end{subfigure}%
	\begin{subfigure}{.5\textwidth}
		\centering
		\includegraphics[width=1.0\textwidth]{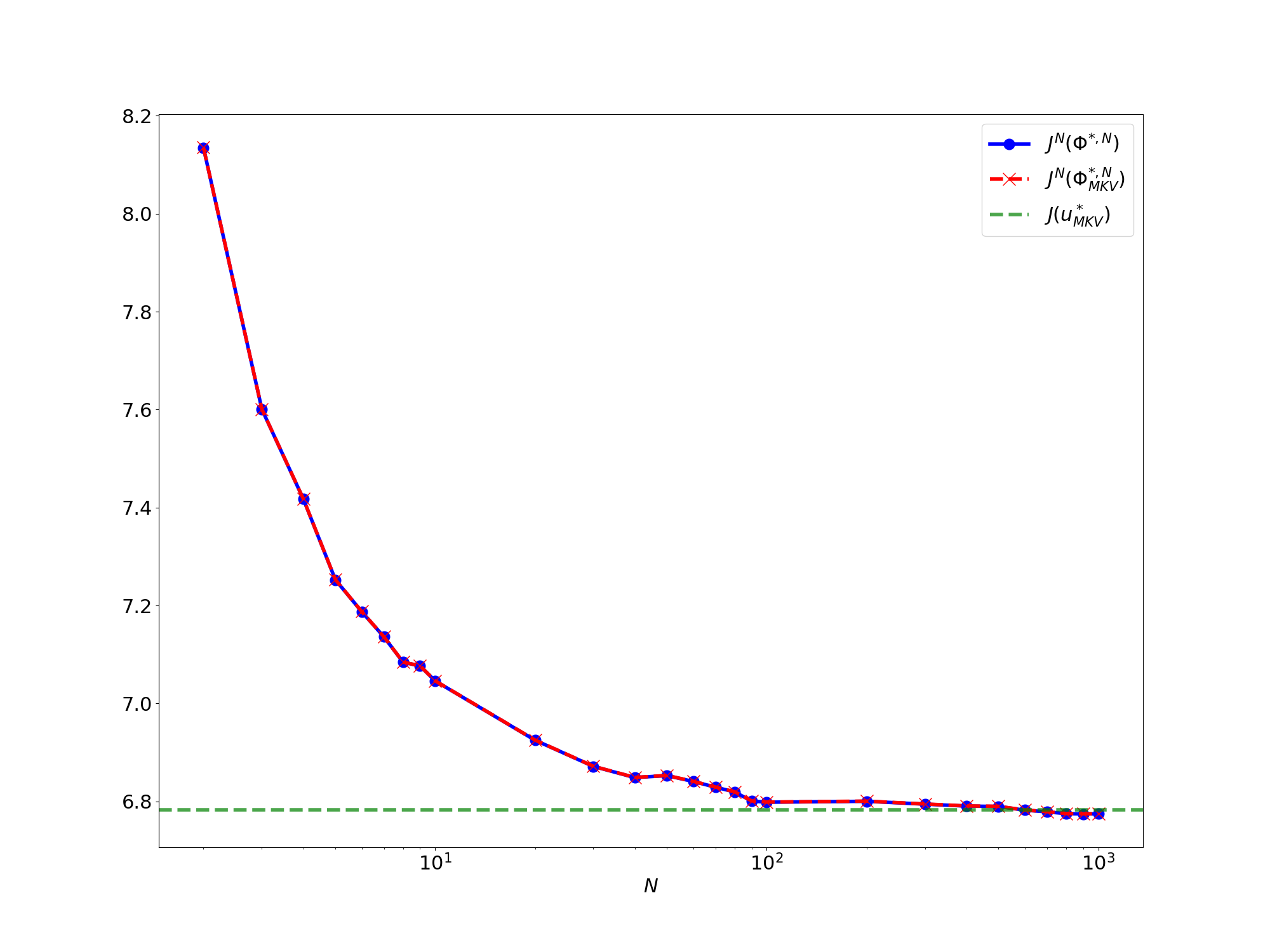}
		\caption{}
		\label{fig:cmp-N-MKV-opt-cost}
	\end{subfigure}
	\caption{Comparison of the $N$ agent optimal control and the MKV optimal control. (a): Maximum difference between diagonal terms; (b) Social cost for each control.}
\end{figure}

When the coefficients of the model are known, the matrices $\Phi^{*,N}$ (for the $N-$agent problem) and $(K^*,L^*)$ (for the MFC problem)  can be computed by solving corresponding Riccati equations. In the one dimensional case with $d = \ell = 1$, we compare the gap between the diagonal coefficients of $\Phi^{*,N}$ and $K^* \in \RR$.
Figure~\eqref{fig:cmp-N-MKV-opt-diagterm} shows (in blue) the graph of the function 
$$
N \mapsto \max_{n = 1,\dots, N} \big|(- \Phi^{*,N})_{n,n} - K^*  \big| .
$$
While this quantity decreases with $N$, it does not converge to $0$. This is due to the fact that the optimal control with heterogeneous agents does not converge to the optimal MKV control. The figure was produced using one sample of coefficients $(q^{(n)})_{n=1}^N$ for the matrix $\Phi^{*,N}$. 
For the sake of comparison, we also show (see the red curve in Figure~\eqref{fig:cmp-N-MKV-opt-diagterm}) that the absolute value of diagonal coefficients of $\Phi^{*,N}_{MKV}$, defined by $K^* + \frac{1}{N}(L^* - K^*)$, converge to $K^*$.

On the other hand, instead of comparing the controls, we can compare the associated costs.
Indeed, once the matrices $\Phi^{*,N}$ and $\Phi^{*,N}_{MKV}$ are computed, if we use them as feedback functions, the corresponding $N-$player social cost can be readily computed using an analytical formula on the feedback coefficients. Figure~\eqref{fig:cmp-N-MKV-opt-cost} shows the difference between the two social costs as $N$ increases. One can see that whether the agents use the real optimal control (blue line) or the one coming from the optimal MKV control (red line), the values of the social cost are almost the same as they seem to converge to the value of the theoretical optimal MKV social cost (green dashed line). However, a small discrepancy remains due to the heterogeneity of the population.

\begin{figure}[h!]
	\centering
	\begin{subfigure}{0.4\textwidth}
		\centering
		\includegraphics[width=1.0\textwidth]{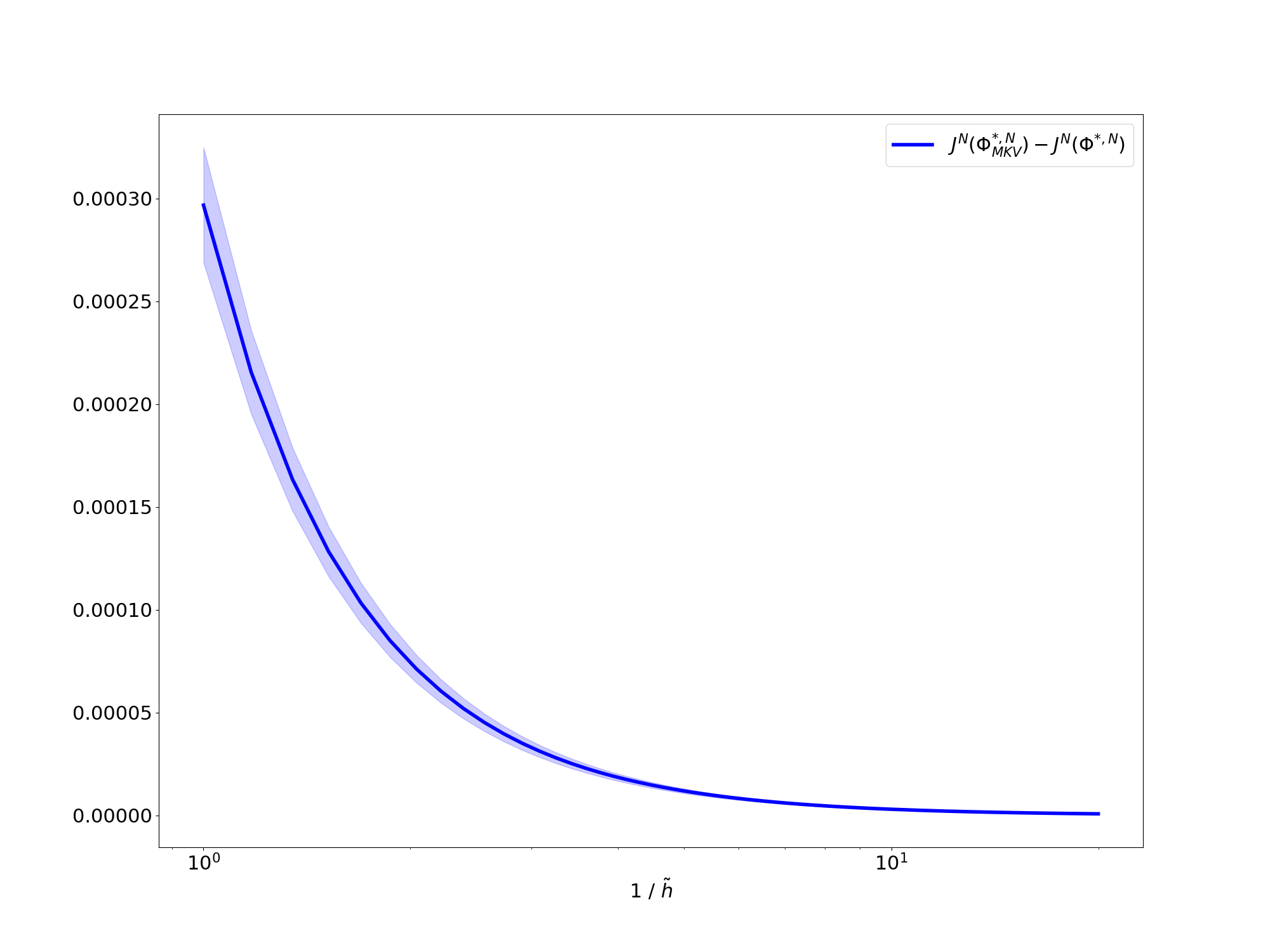}
	\end{subfigure}%
	\caption{Influence of the heterogeneity degree on the social cost.}
	\label{fig:influence-hetero-radius}
\end{figure}

\begin{figure*}[!htp]
	\centering
	\begin{minipage}[h]{0.92\textwidth}
		\begin{subfigure}{0.26\textwidth}
			\centering
			\includegraphics[width=\textwidth]{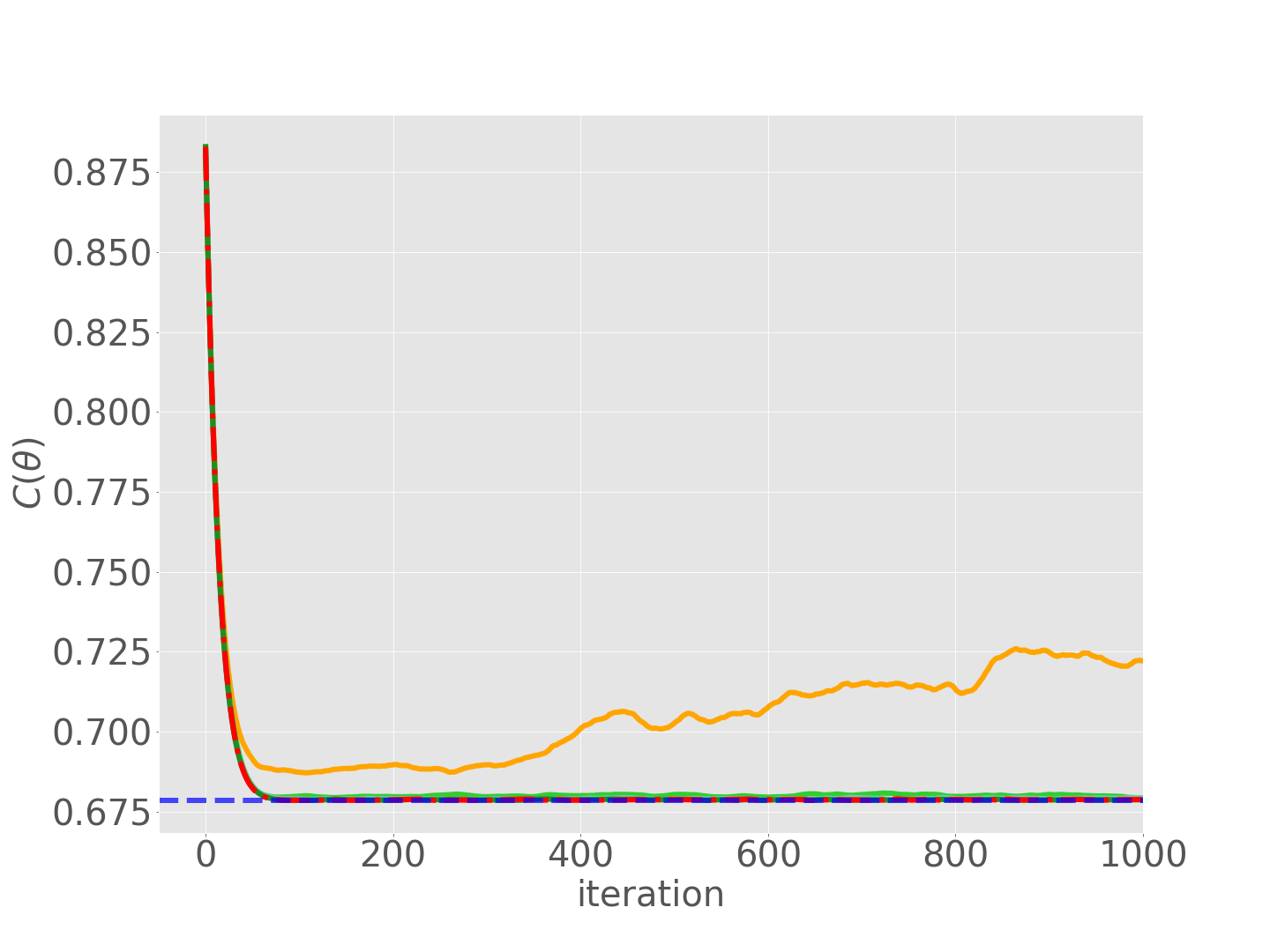}
			\caption{$C(\theta)$}
		\end{subfigure}%
		\hspace{-0.3cm}%
		\begin{subfigure}{.26\textwidth}
			\centering
			\includegraphics[width=\textwidth]{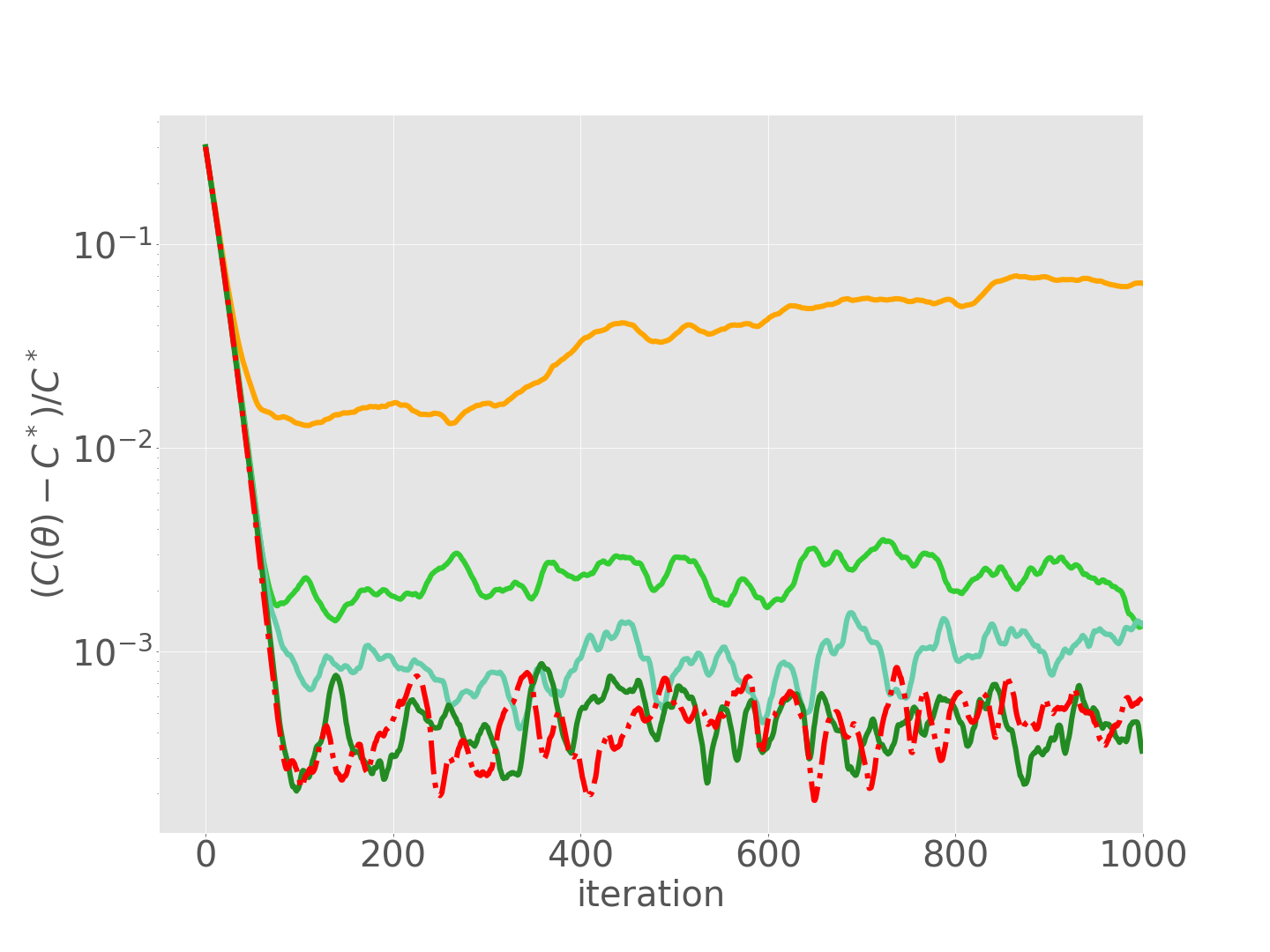} 
			\caption{$ \frac{ C(\theta) - C^* }{ C^* }$}%
		\end{subfigure}%
		\hspace{-0.3cm}%
		\begin{subfigure}{.26\textwidth}
			\centering
			\includegraphics[width=\textwidth]{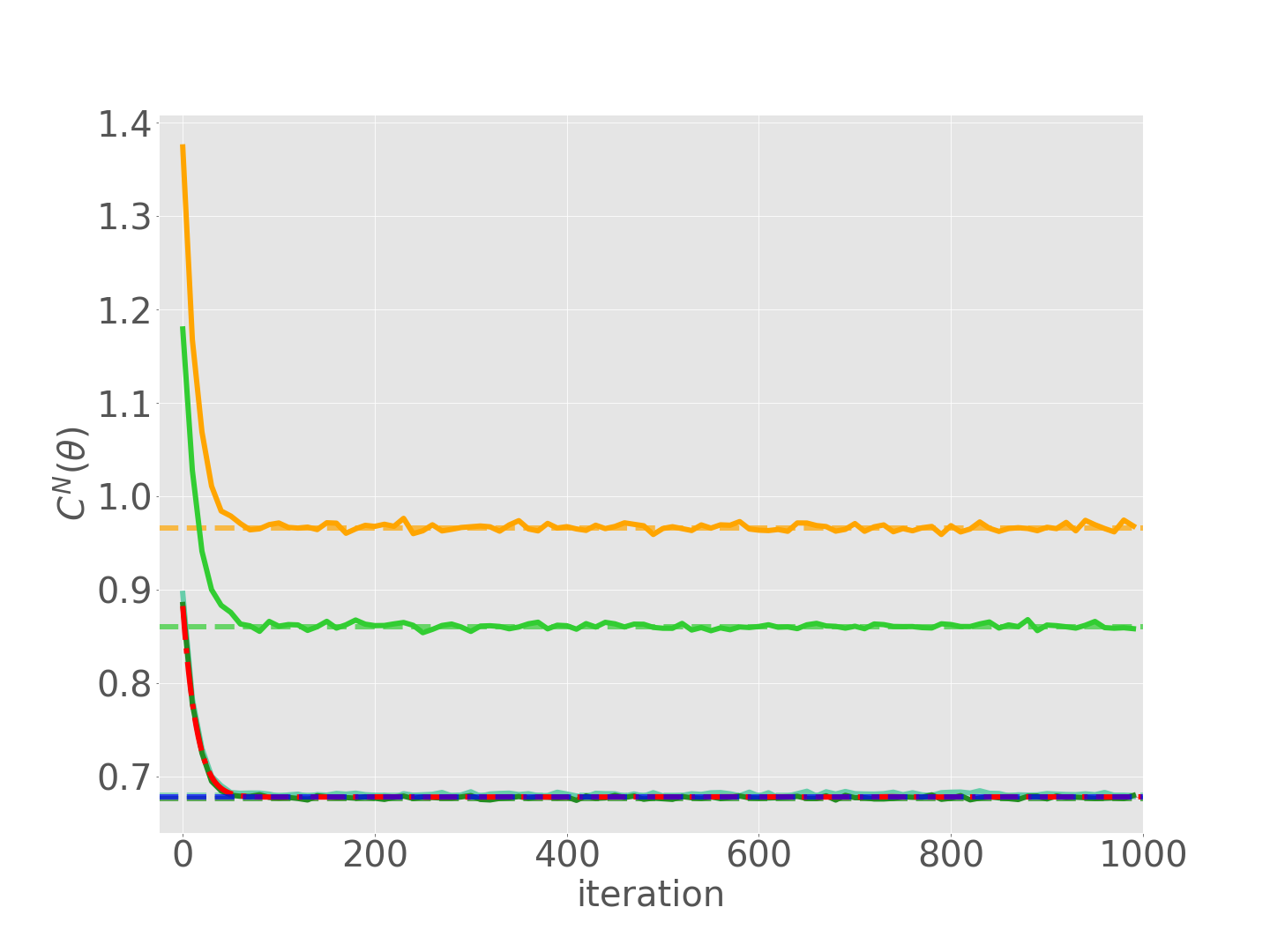}
			\caption{$C^N(\theta)$}
		\end{subfigure}%
		\hspace{-0.3cm}%
		\begin{subfigure}{.26\textwidth}
			\centering\captionsetup{width=.8\linewidth}
			\includegraphics[width=\textwidth]{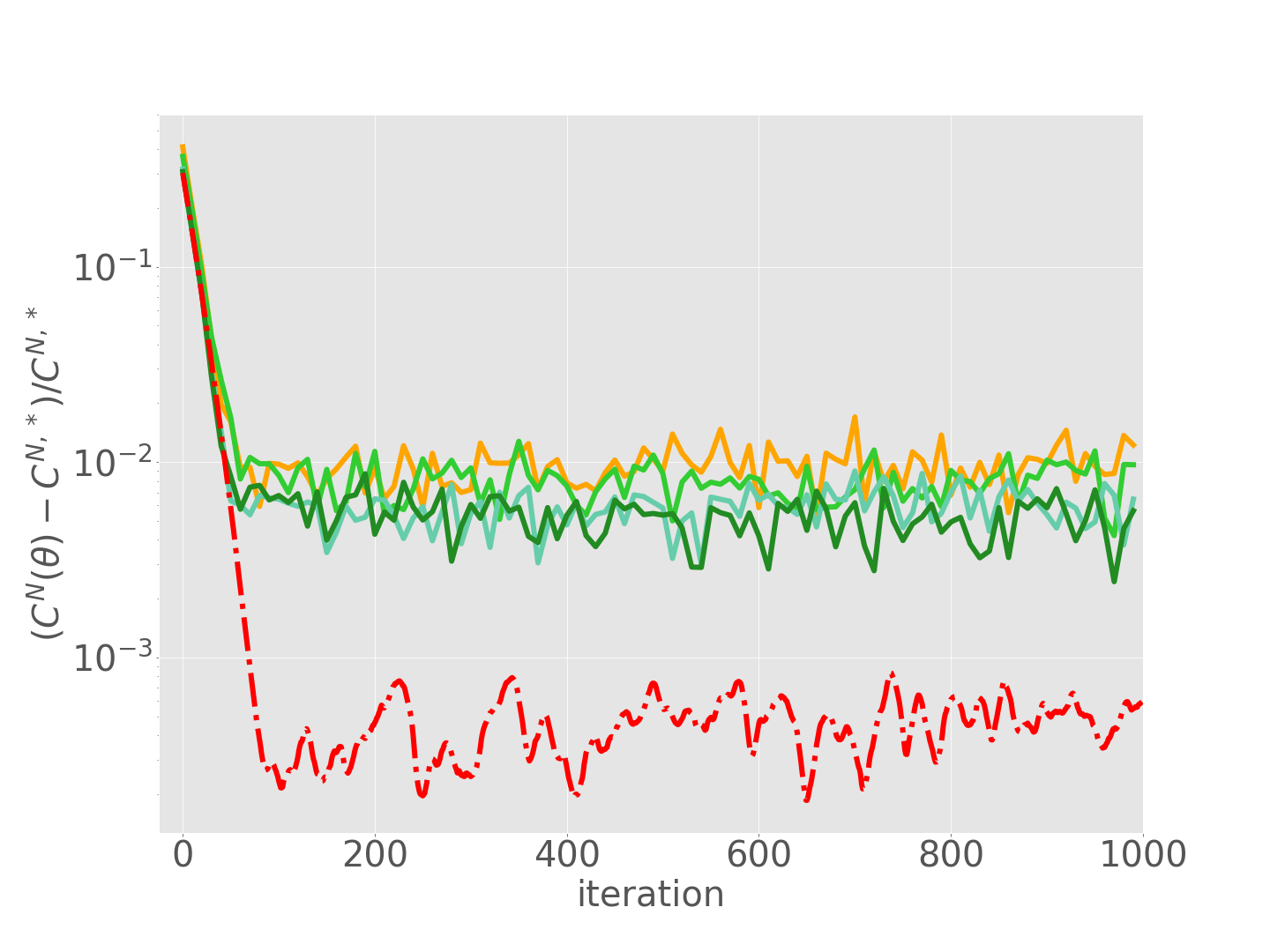}
			\caption{$\frac{C^{N}(\theta) - C^{*,N}}{C^{*,N}}$} %
		\end{subfigure}%
	\end{minipage}%
	\hspace{-0.3cm}%
	\begin{minipage}[t]{0.08\textwidth}
		\centering
		\begin{subfigure}{\textwidth}
			\includegraphics[width=\textwidth]{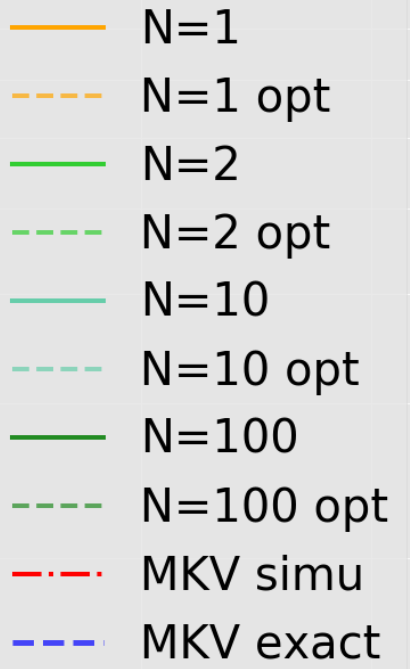}
		\end{subfigure}
	\end{minipage}
	\caption{ Discount factor $\gamma = 0.5$, heterogeneous perturbation size $\tilde{h}=0.4$. }
	\label{fig:gamma_0.5_hetero_0.4}
\end{figure*}

\begin{figure*}[!htp]
	\centering
	\begin{minipage}[h]{0.92\textwidth}
		\begin{subfigure}{0.26\textwidth}
			\centering
			\includegraphics[width=\textwidth]{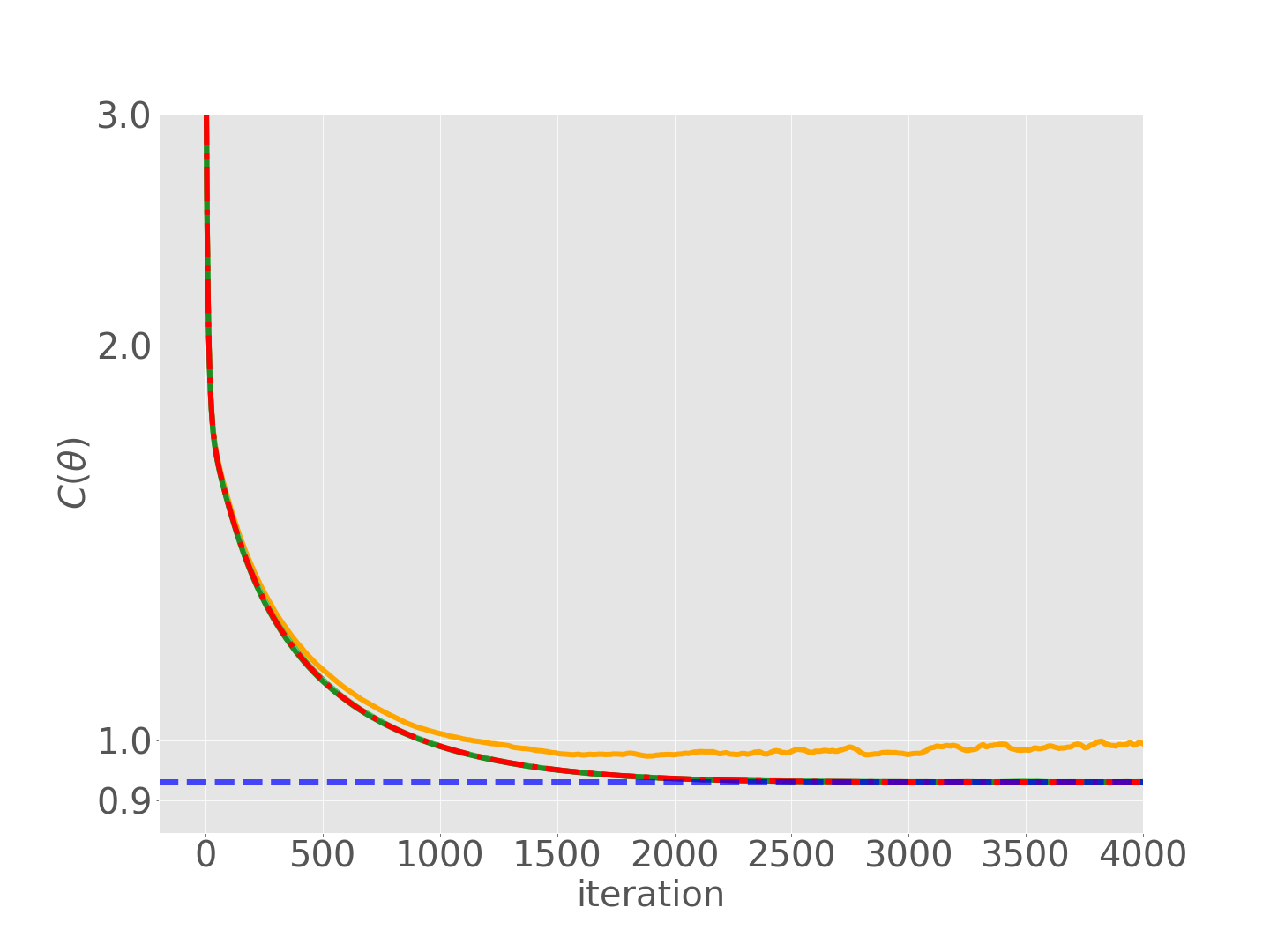}
			\caption{$C(\theta)$}
		\end{subfigure}%
		\hspace{-0.3cm}%
		\begin{subfigure}{.26\textwidth}
			\centering
			\includegraphics[width=\textwidth]{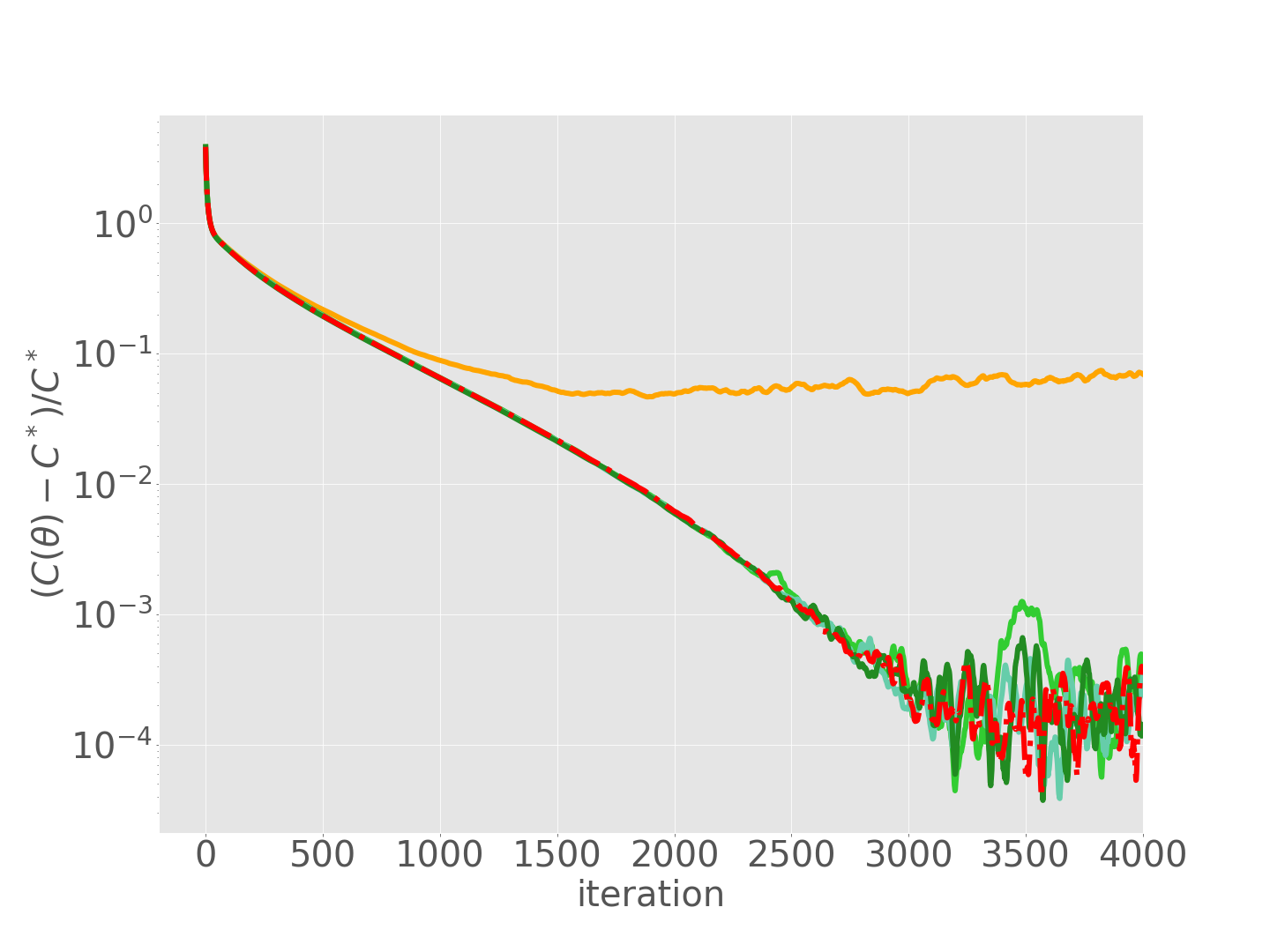} 
			\caption{$ \frac{ C(\theta) - C^* }{ C^* }$}%
		\end{subfigure}%
		\hspace{-0.3cm}%
		\begin{subfigure}{.26\textwidth}
			\centering
			\includegraphics[width=\textwidth]{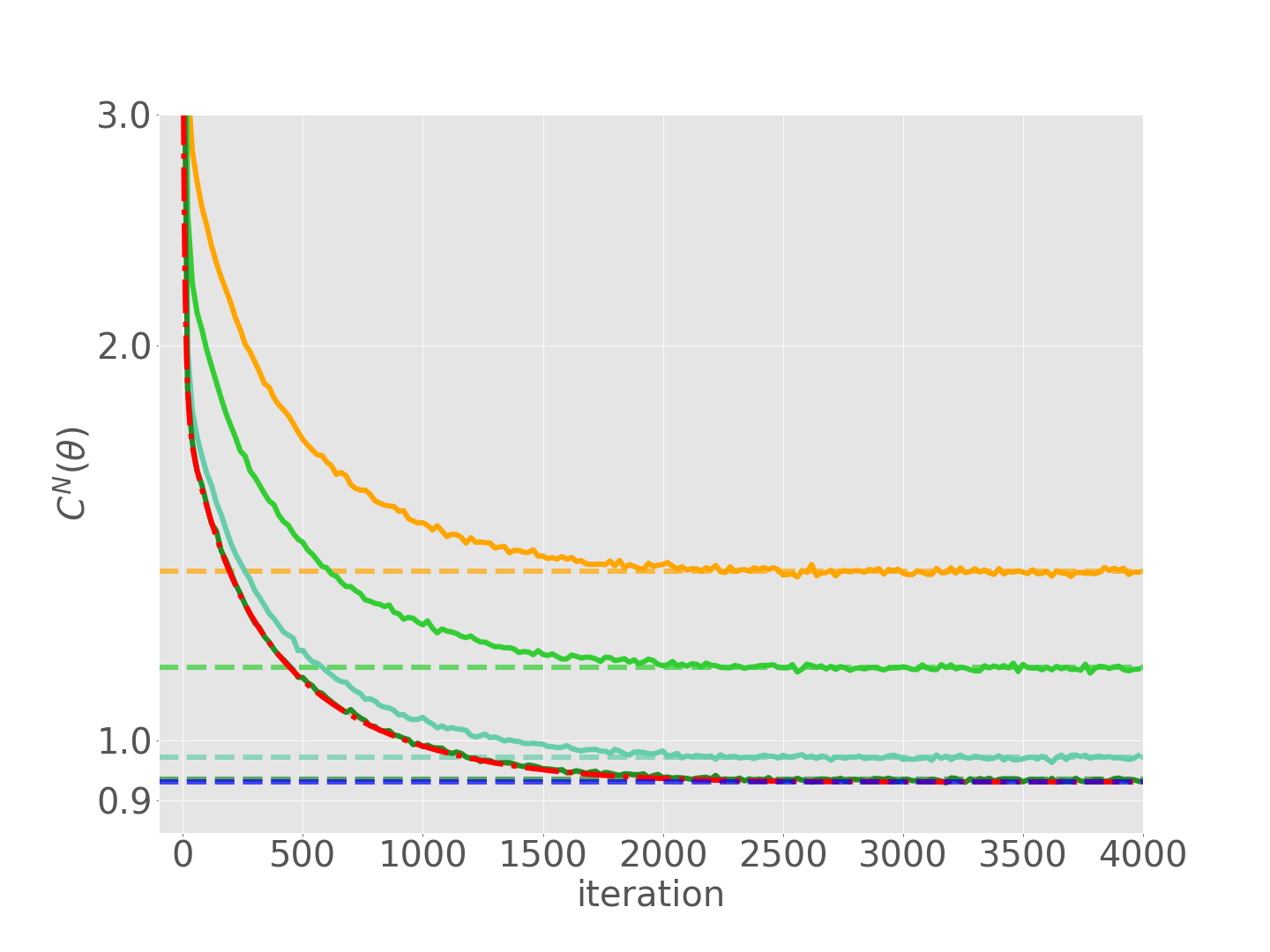}
			\caption{$C^N(\theta)$}
		\end{subfigure}%
		\hspace{-0.3cm}%
		\begin{subfigure}{.26\textwidth}
			\centering\captionsetup{width=.8\linewidth}
			\includegraphics[width=\textwidth]{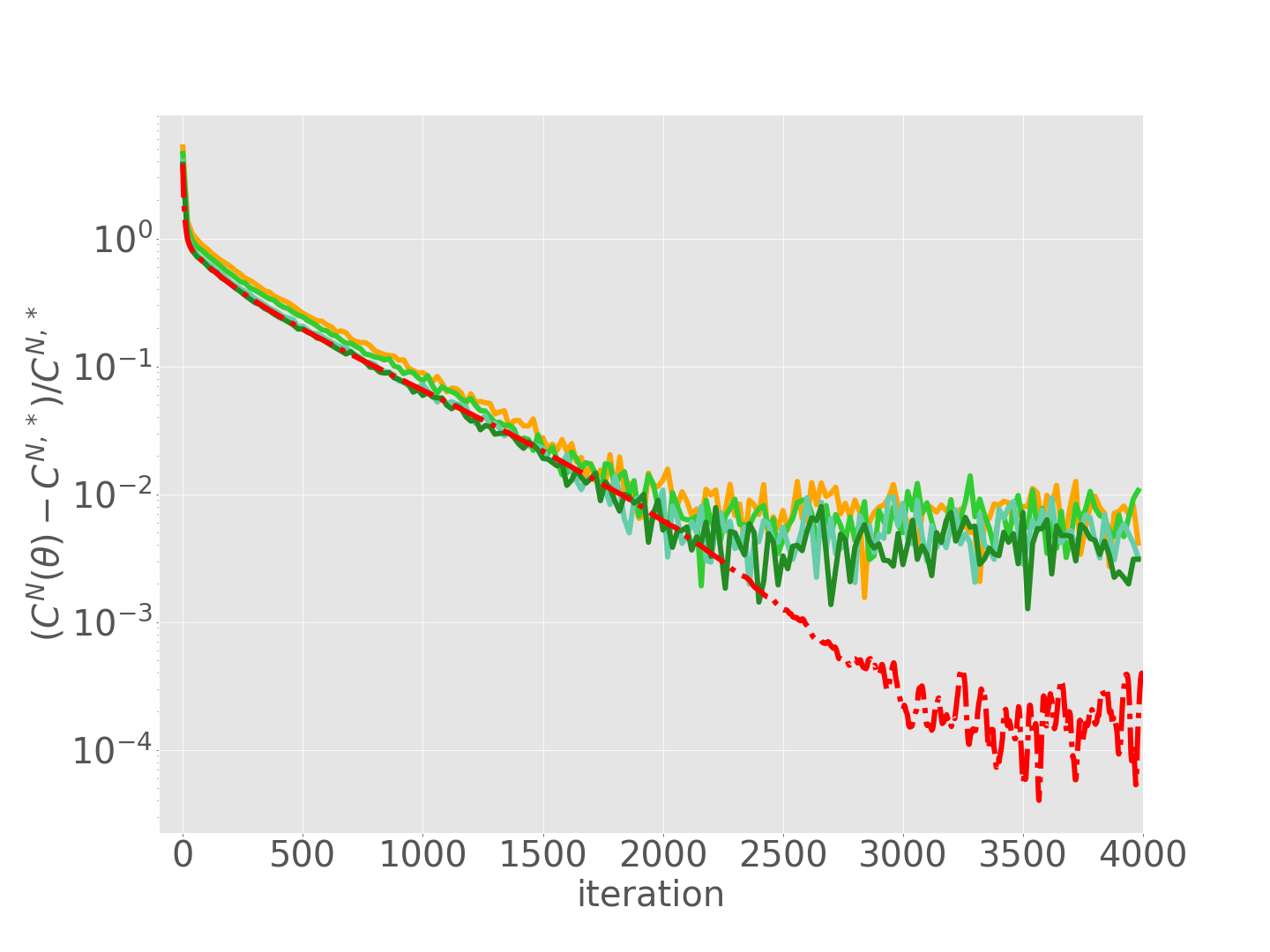}
			\caption{$\frac{C^{N}(\theta) - C^{*,N}}{C^{*,N}}$} %
			\label{fig:gamma_0.9_homo_relcost}
		\end{subfigure}%
	\end{minipage}%
	\hspace{-0.3cm}%
	\begin{minipage}[t]{0.08\textwidth}
		\centering
		\begin{subfigure}{\textwidth}
			\includegraphics[width=\textwidth]{{FIGURES/FIGURES_APPENDIX/gamma_0.5_hetero_0.4/legend_cost}.png}
		\end{subfigure}
	\end{minipage}
	\caption{ Discount factor $\gamma = 0.9$, homogeneous case $\tilde{h}=0$ }
	\label{fig:gamma_0.9_homo}
\end{figure*}

Figure~\eqref{fig:influence-hetero-radius} illustrates the influence of the heterogeneity of the agents. For fixed $N=100$, as the degree of heterogeneity $\tilde h$ vanishes, the difference in social costs evaluated with $\Phi^{*,N}$ and $\Phi^{*,N}_{MKV}$ decreases. In this figure, the curve is obtained by averaging over $5$ random realizations of the model coefficients $(q^{(n)})_{n=1}^N$, and the shaded region corresponds to the mean $\pm$ one standard deviation. 
We also observe the effect of heterogeneity on the PG convergence algorithm when $\tilde h$ decreases. We see clearly a concentration towards the MFC problem in terms of the relative error when $\tilde{h}$ decreases from $0.4$ to $0.1$ in Figure~\eqref{fig:gamma_0.5_hetero_0.4}(b) and~\ref{fig:gamma_0.5_heteo_0.1}(b), and when $\tilde{h}$ decrease from $0.1$ to $=0$ in Figures~\ref{fig:gamma_0.9_heteo_0.1}(d) and~\ref{fig:gamma_0.9_homo}(d).

\subsection{Effect of the discount coefficient $\gamma$}

\begin{figure*}[!htp]
	\centering
	\begin{minipage}[h]{0.92\textwidth}
		\begin{subfigure}{0.26\textwidth}
			\centering
			\includegraphics[width=\textwidth]{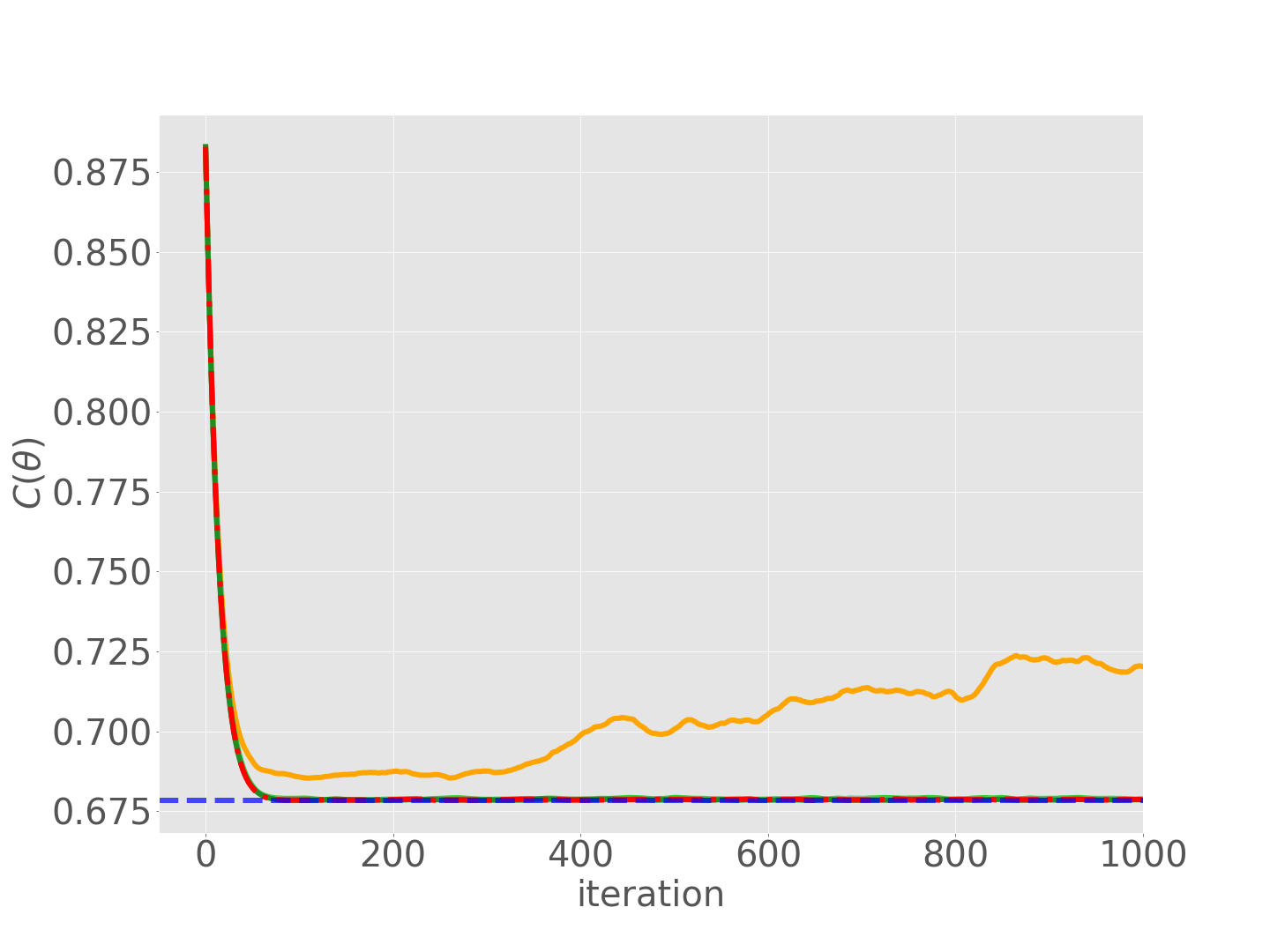}
			\caption{$C(\theta)$}
		\end{subfigure}%
		\hspace{-0.3cm}%
		\begin{subfigure}{.26\textwidth}
			\centering
			\includegraphics[width=\textwidth]{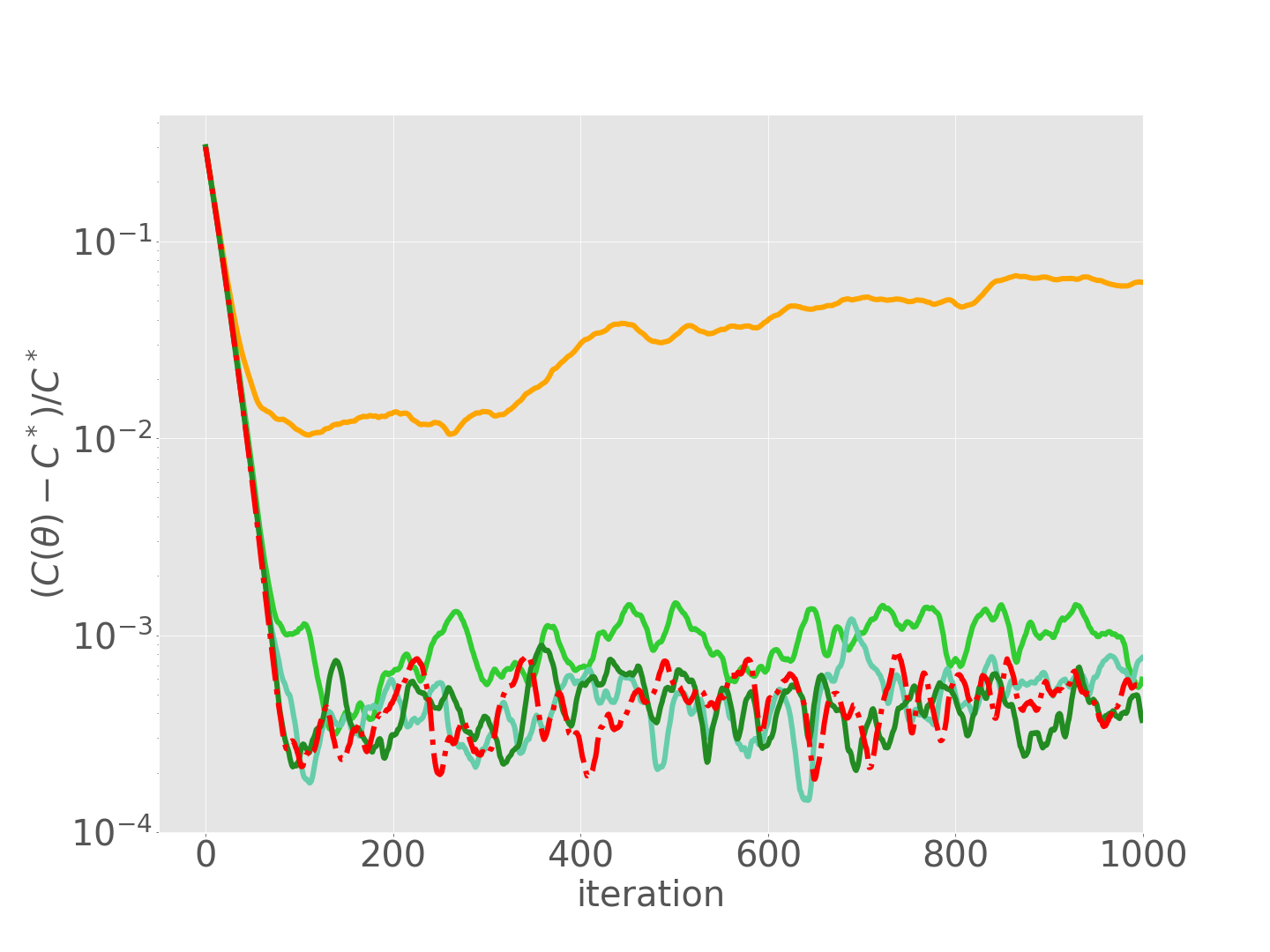} 
			\caption{$ \frac{ C(\theta) - C^* }{ C^* }$}%
		\end{subfigure}%
		\hspace{-0.3cm}%
		\begin{subfigure}{.26\textwidth}
			\centering
			\includegraphics[width=\textwidth]{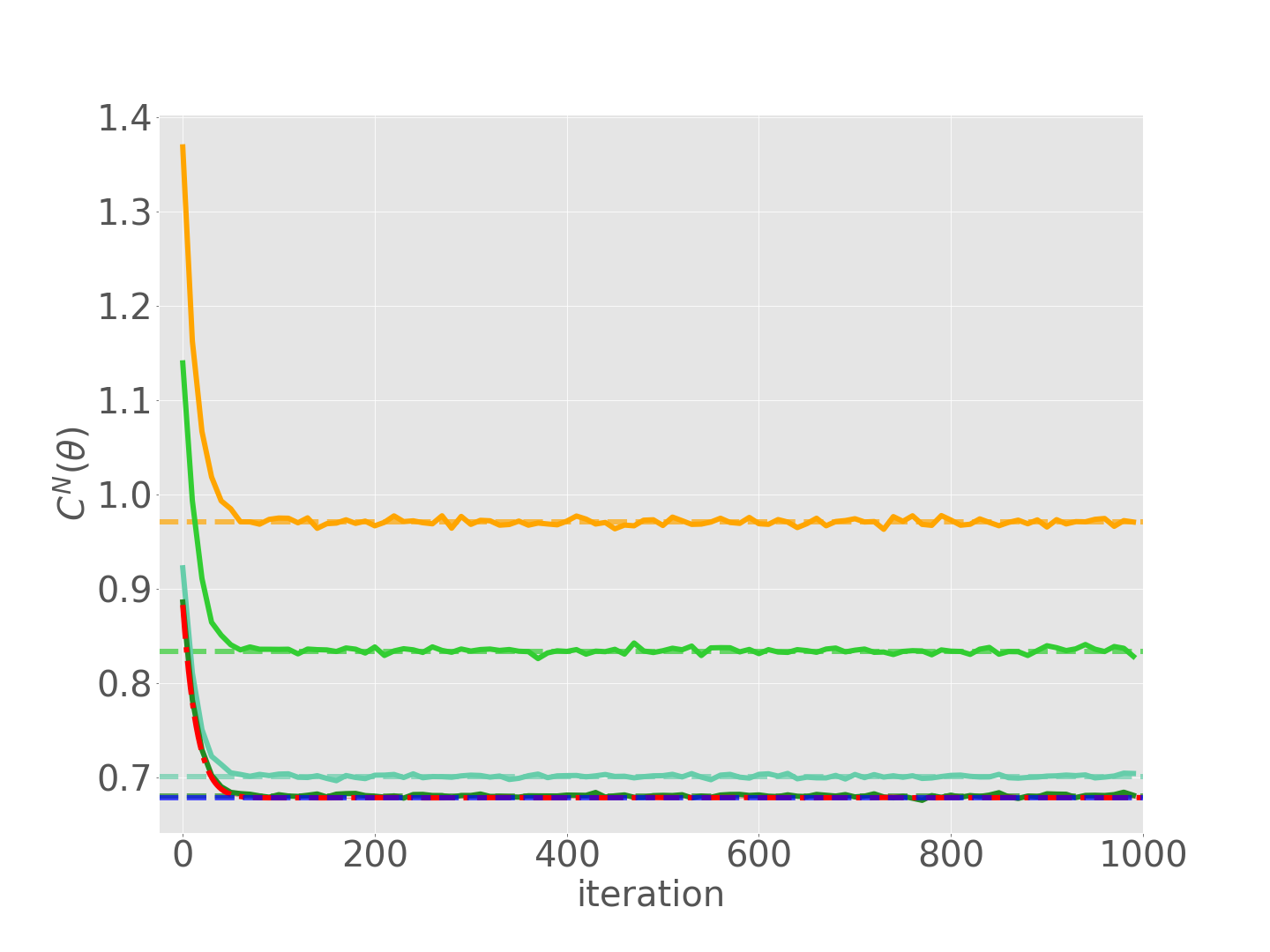}
			\caption{$C^N(\theta)$}
		\end{subfigure}%
		\hspace{-0.3cm}%
		\begin{subfigure}{.26\textwidth}
			\centering\captionsetup{width=.8\linewidth}
			\includegraphics[width=\textwidth]{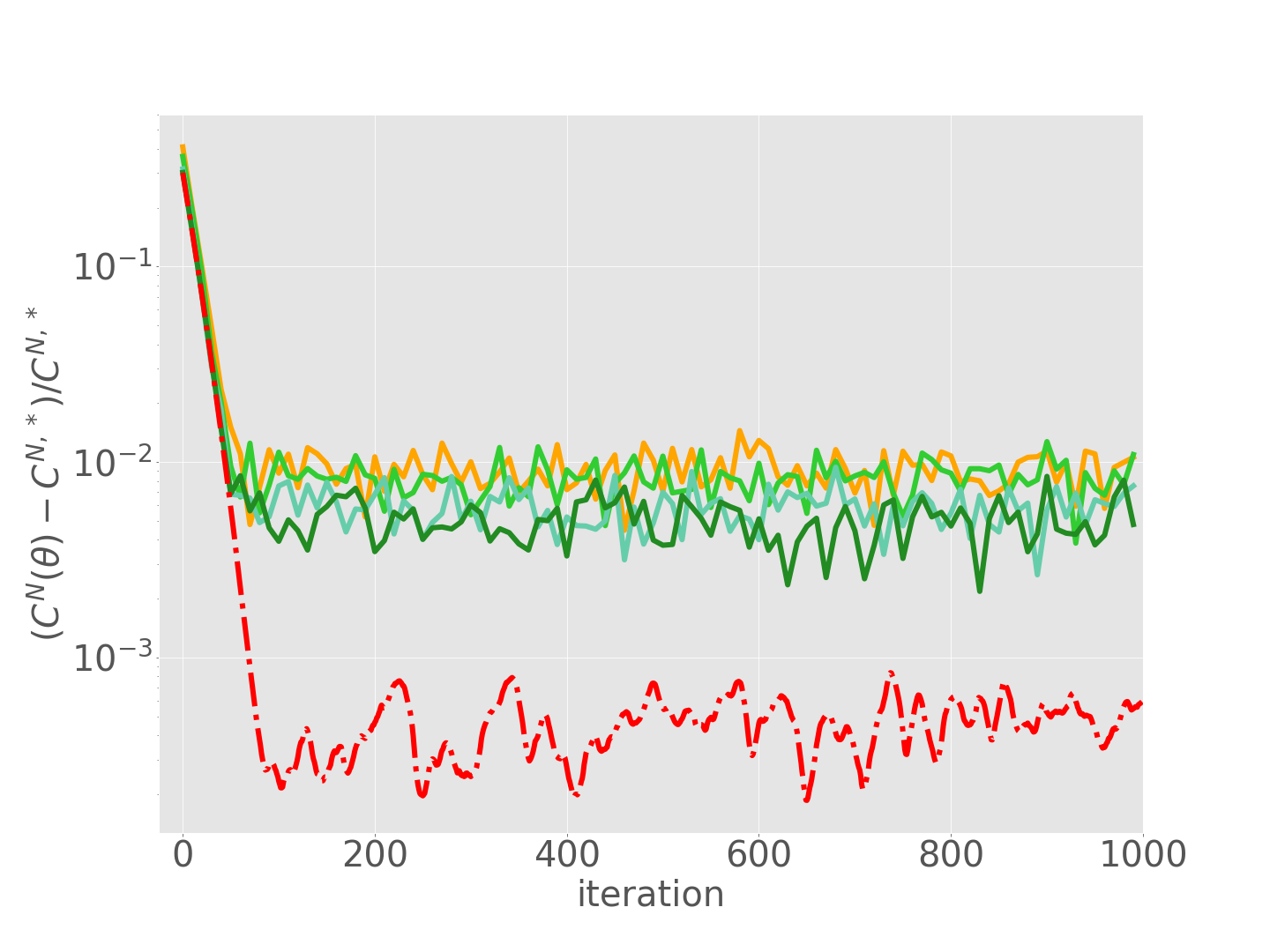}
			\caption{$\frac{C^{N}(\theta) - C^{*,N}}{C^{*,N}}$} %
		\end{subfigure}%
	\end{minipage}%
	\hspace{-0.3cm}%
	\begin{minipage}[t]{0.08\textwidth}
		\centering
		\begin{subfigure}{\textwidth}
			\includegraphics[width=\textwidth]{{FIGURES/FIGURES_APPENDIX/gamma_0.5_hetero_0.4/legend_cost}.png}
		\end{subfigure}
	\end{minipage}
	\caption{ Discount factor $\gamma = 0.5$, heterogeneous perturbation size $\tilde{h}=0.1$. }
	\label{fig:gamma_0.5_heteo_0.1}
\end{figure*}

\begin{figure*}[!htp]
	\centering
	\begin{minipage}[h]{0.92\textwidth}
		\begin{subfigure}{0.26\textwidth}
			\centering
			\includegraphics[width=\textwidth]{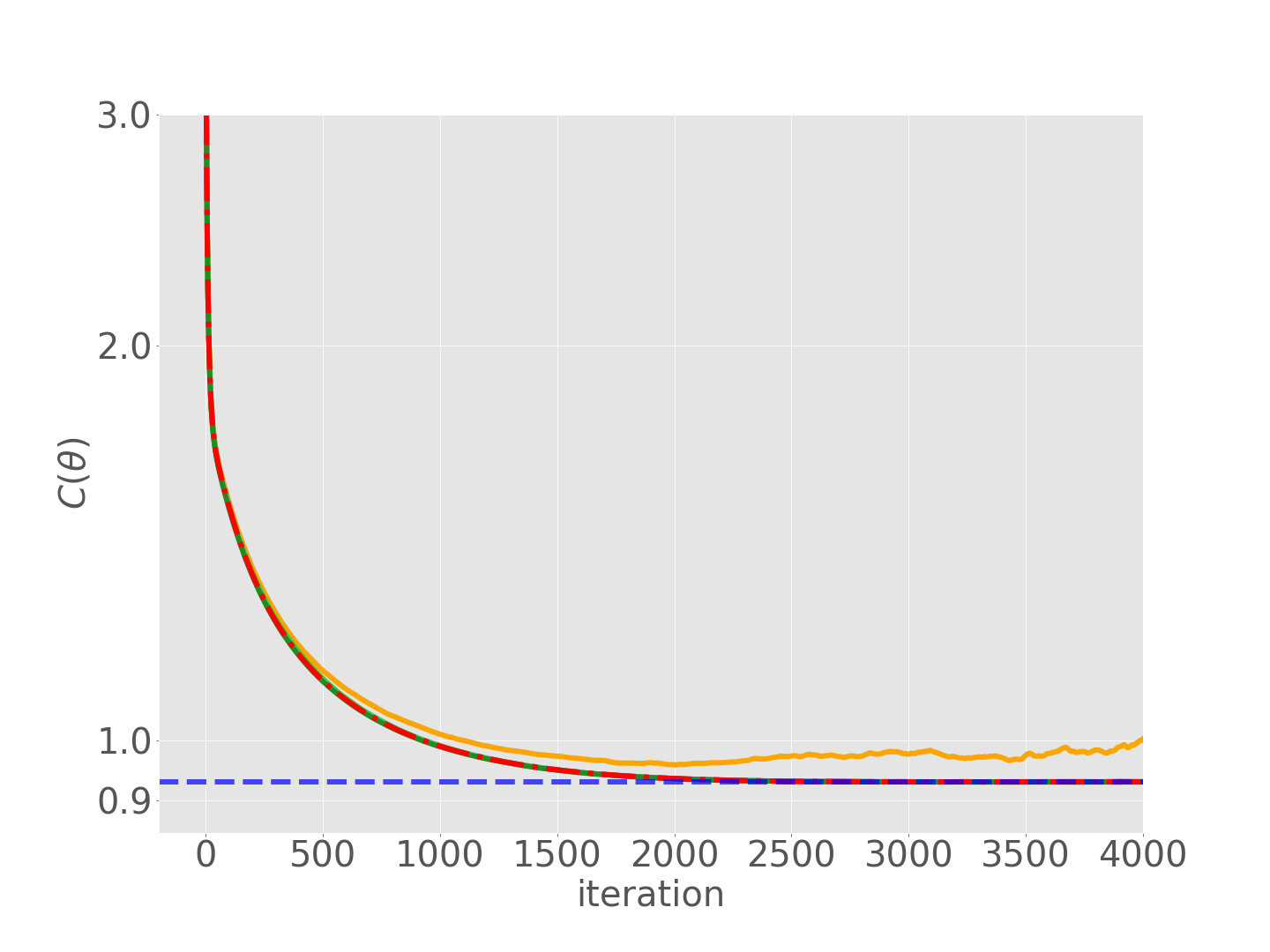}
			\caption{$C(\theta)$}
		\end{subfigure}%
		\hspace{-0.3cm}%
		\begin{subfigure}{.26\textwidth}
			\centering
			\includegraphics[width=\textwidth]{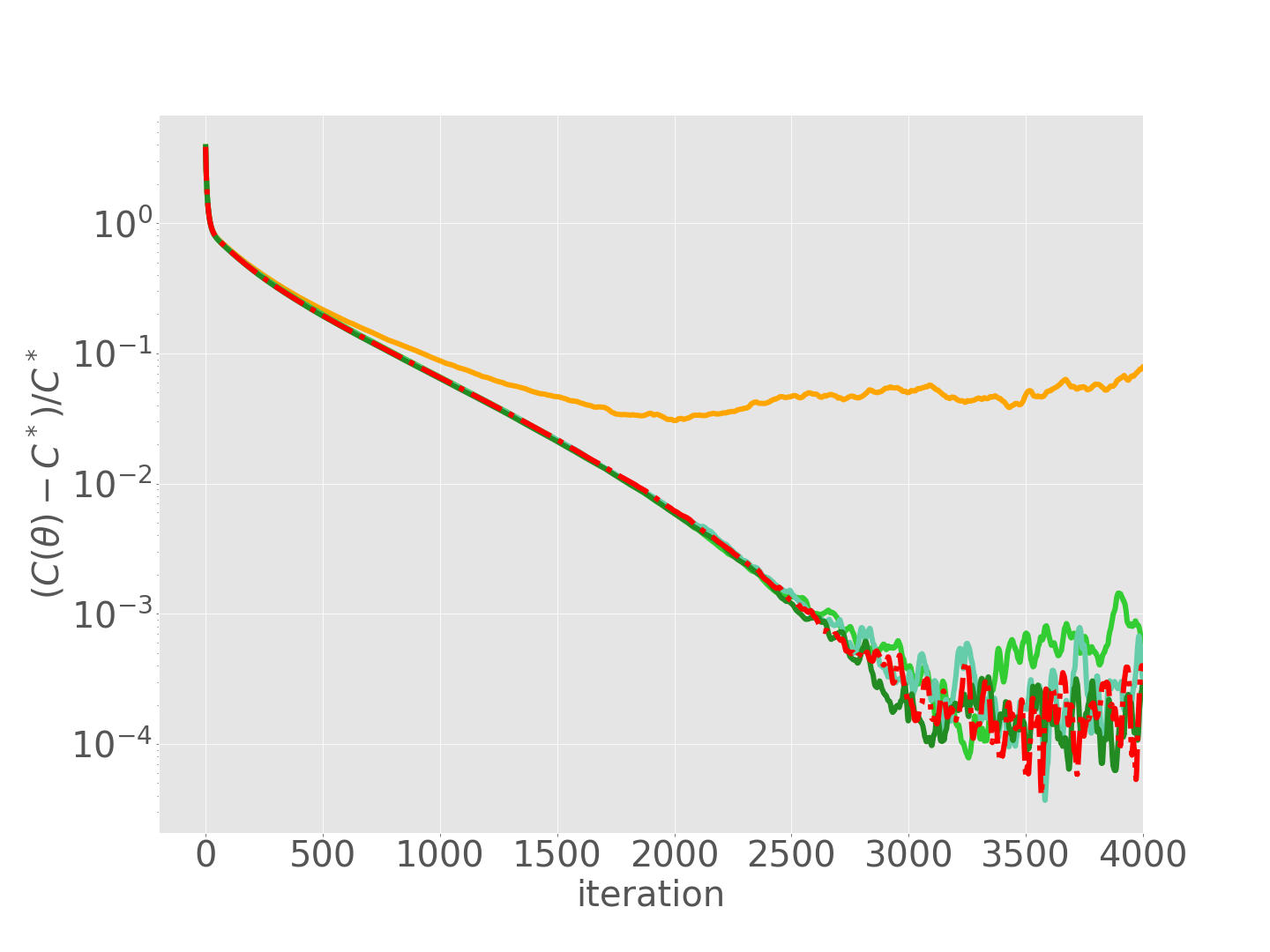} 
			\caption{$ \frac{ C(\theta) - C^* }{ C^* }$}%
		\end{subfigure}%
		\hspace{-0.3cm}%
		\begin{subfigure}{.26\textwidth}
			\centering
			\includegraphics[width=\textwidth]{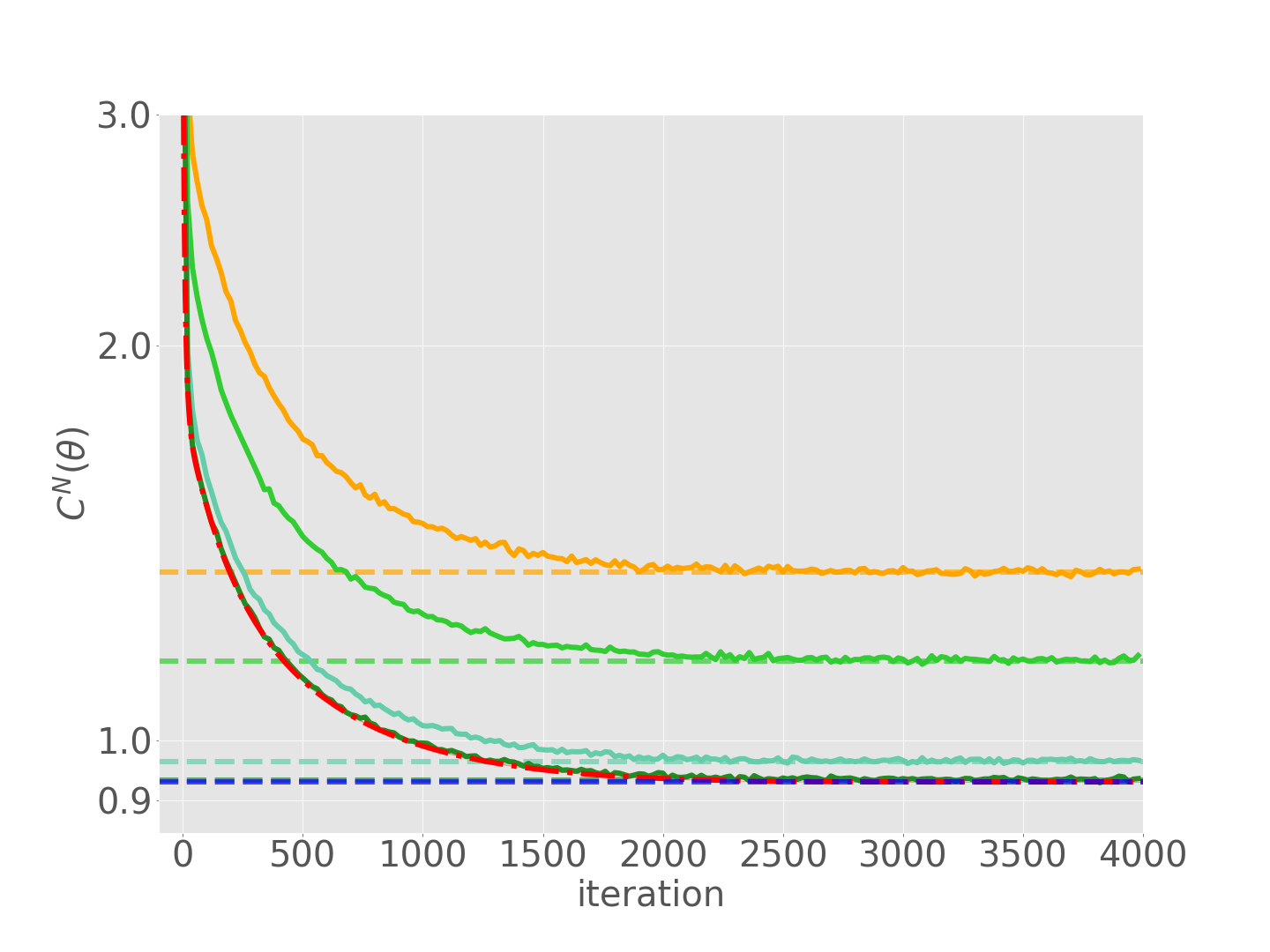}
			\caption{$C^N(\theta)$}
		\end{subfigure}%
		\hspace{-0.3cm}%
		\begin{subfigure}{.26\textwidth}
			\centering\captionsetup{width=.8\linewidth}
			\includegraphics[width=\textwidth]{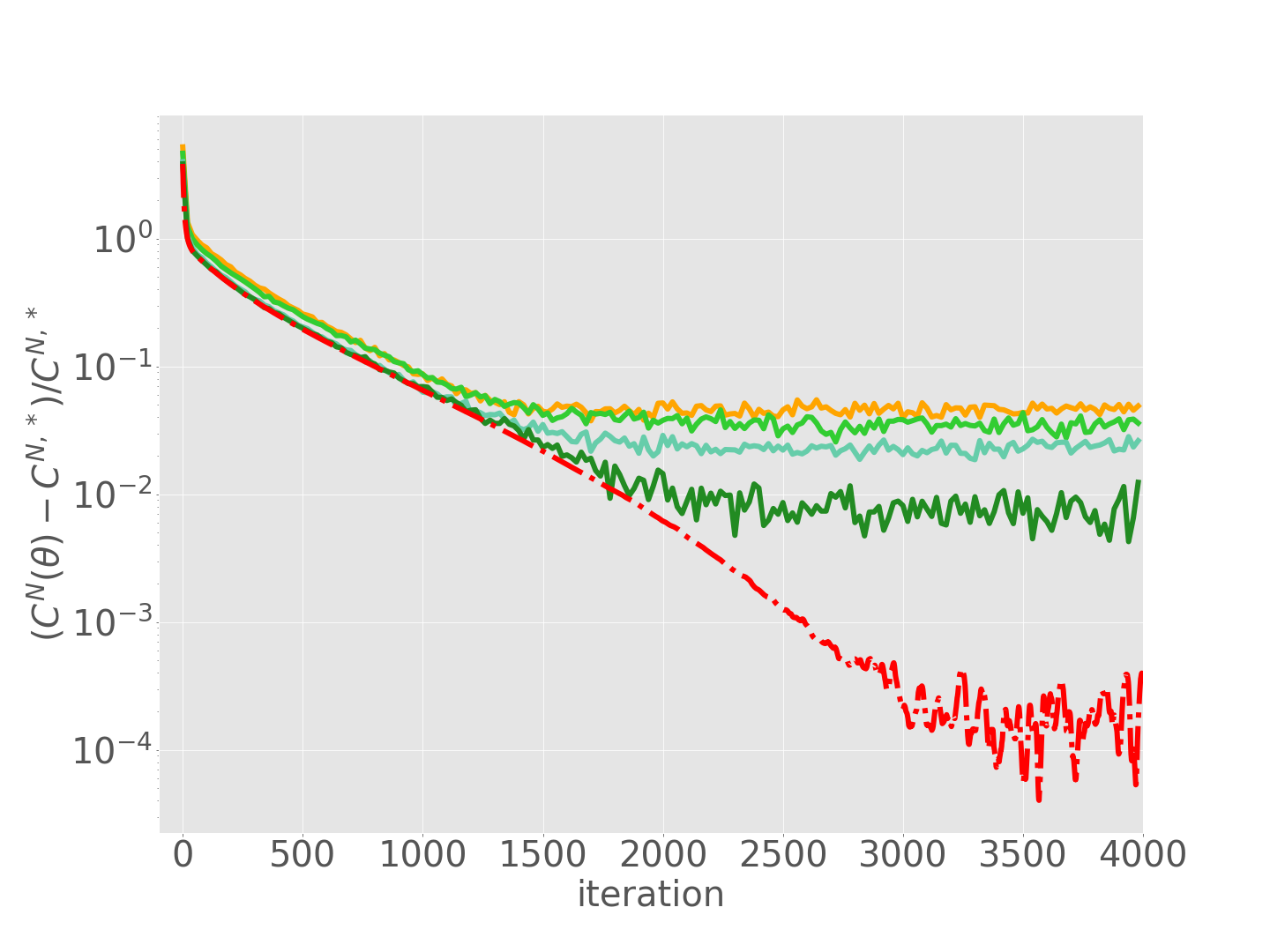}
			\caption{$\frac{C^{N}(\theta) - C^{*,N}}{C^{*,N}}$} %
			\label{fig:gamma_0.9_heteo_0.1_relcost}
		\end{subfigure}%
	\end{minipage}%
	\hspace{-0.3cm}%
	\begin{minipage}[t]{0.08\textwidth}
		\centering
		\begin{subfigure}{\textwidth}
			\includegraphics[width=\textwidth]{{FIGURES/FIGURES_APPENDIX/gamma_0.5_hetero_0.4/legend_cost}.png}
		\end{subfigure}
	\end{minipage}
	\caption{Discount factor $\gamma = 0.9$, heterogeneous perturbation size $\tilde{h}=0.1$. }
	\label{fig:gamma_0.9_heteo_0.1}
\end{figure*}

The discount factor also plays an important role in the learning procedure. The higher it is,
the harder the problem becomes. Indeed, as $\gamma$ gets close to $1$, the future times are less and less discounted
and one needs to find an optimal behavior over a longer time horizon, which is more challenging.
For the numerical results presented previously, we took $\gamma= 0.9$, which is quite high and in particular, it is higher than e.g. the value considered in \cite{guo2019learning} (where the authors used the value $0.8$). 
With a larger $\gamma$, the model-free algorithm needs a larger number of perturbation directions $M$ for the estimation of gradients, and it also needs a longer truncated time horizon $T$ as well as more iteration steps to converge. Figure~\eqref{fig:gamma_0.5_heteo_0.1} and Figure~\eqref{fig:gamma_0.9_heteo_0.1} provide a simple illustration of the effect of $\gamma$ along the iteration of the algorithm. For $\gamma = 0.5$, it takes approximately $100$ iterations to achieve $10^{-3}$ in terms of the relative error with respect to the optimal MF cost, whereas for $\gamma=0.9$, the number of iterations goes up to $3000$.

\bibliographystyle{apalike}
\bibliography{MF_RL.bib}

\end{document}